\documentclass[5p,times]{elsarticle}

\usepackage{graphicx} % Required for inserting images
\usepackage{hyperref}
\usepackage{subcaption}
\usepackage{amsmath}
\usepackage{mathtools,amsfonts}
\usepackage{framed}
\usepackage{color}
\usepackage[export]{adjustbox}
\usepackage{multirow}
\usepackage[table,xcdraw]{xcolor}
\usepackage{tabularx}
\usepackage[shortlabels]{enumitem}
\usepackage{ragged2e}
\usepackage[format=plain, labelfont={bf}]{caption}
\usepackage{tikz}
\usepackage{pgfplots}
\pgfplotsset{compat=1.18} 
\usepackage{booktabs}
\usepackage{lineno}
\usepackage{rotating}

\usepackage[normalem]{ulem}

\newcolumntype{+}{!{\vrule width 2pt}}

\newlength\savedwidth

\newcommand\thickhline{\noalign{\global\savedwidth\arrayrulewidth\global\arrayrulewidth 2pt}%
\hline
\noalign{\global\arrayrulewidth\savedwidth}}
\begin{document}
\begin{frontmatter}

\title{\textcolor{black}{\emph{In Silico} Study for Optimizing Intensity and Focality Electrode Configurations for Directional DBS Under Uncertainty Using Metaheuristic L1L1 Method}}
%PREVIOUS: Balancing Intensity and Focality in Directional DBS Under Uncertainty: A Simulation Study of Electrode Optimization via Metaheuristic L1L1 Method

\author[inst1]{Fernando {Galaz Prieto}\corref{cor1}}
\ead{fernando.galazprieto@tuni.fi}
\cortext[cor1]{Corresponding author at: Sähkötalo building, Korkeakoulunkatu 3, Tampere, 33720, FI}
\author[inst1]{Antti Lassila}
\author[inst1,inst2]{Maryam Samavaki}
\author[inst1]{Sampsa Pursiainen}
\affiliation[inst1]{Computing Sciences Unit, Faculty of Information Technology and Communication Sciences, Tampere University, Tampere, Finland}
\affiliation[inst2]{Department of Physics and Mathematics, Faculty of Science, Forestry and Technology, University of Eastern Finland, Joensuu, Finland}

\begin{abstract}
\noindent \textbf{Background and Objective:} As Deep Brain Stimulation (DBS) advances toward directional leads and optimization-based current steering, selecting electrode contact configurations becomes complex. This study formulates configuration selection as an inverse mapping between target activation and electrode currents using metaheuristic L1-norm regularized L1-norm fitting (L1L1). L1L1 incorporates lead-field uncertainty arising from electrode placement, tissue properties, and forward modeling assumptions.

\noindent \textbf{Methods:} The framework introduces lead-field perturbations and restricts the controllable domain through a sensitivity-based feasibility criterion within a finite element formulation derived using the Complete Electrode Model. Current distributions were optimized for 8- and 40-contact leads. Performance was evaluated using focused current density, nuisance current density, and their ratio under safety and sparsity constraints.

\noindent \textbf{Results:} L1L1 was evaluated using noiseless and noisy lead fields, with noise selected to reflect attenuation within the volume of tissue activated. The method produced sparse, spatially selective stimulation patterns across perturbation levels. Hyperparameter optimization yielded bipolar or multipolar configurations. Compared with the Reciprocity Principle (RP), which produced strictly bipolar configurations, and Tikhonov-regularized least squares (TLS), which produced more distributed solutions, L1L1 enabled controlled transitions between sparse and multipolar patterns. It concentrated stimulation within the target while limiting unintended current spread, particularly under noisy conditions.

\noindent \textbf{Conclusions:} L1L1 can assist specialists in optimizing DBS configurations. By incorporating uncertainty directly into optimization, it provides robust and interpretable current steering across lead configurations while accounting for forward-model variability.

\end{abstract}


\end{frontmatter}

\section{Introduction}
\label{sec:Intro}

% Introduction
Deep brain stimulation (DBS) is a neuromodulatory therapy in which electrical stimulation is delivered to targeted sub-cortical brain regions through implanted electrode leads to modulate pathological neural activity \cite{volkmann2002introduction, Gardner_2013_history_DBS}. The clinical procedure involves stereotactic implantation of a multi-contact electrode lead into a selected deep brain target, followed by the delivery of controlled electrical stimulation through specific contacts along the lead. Since its clinical introduction, DBS has become an established treatment for several neurological and movement disorders. In Parkinson’s disease, DBS has been widely adopted to alleviate motor symptoms \cite{Limousin_1998_Parkinson, picillo2016programming, krack2019deep, rodriguez2005bilateral, steffen2020bipolar}. The therapy has also been successfully applied in the treatment of essential tremor \cite{benabid1996chronic, Rehncrona2003, Wong2020, Kremer2021}, refractory epilepsy \cite{LEHTIMAKI_2016_DBS_targetting, Jarvenpaa_2020_Improving_ANTDBS, fasano2021experience, Salanova_2021_SANTE10years}, and dystonia \cite{johnson2013neuromodulation}.

%% VTA
\textcolor{black}{
The electro-physiological effects of the stimulation is commonly evaluated through the concept of the \textit{volume of tissue activated} (VTA) \cite{butson_2007_patient, Duffley2019, patrick2024modeling} which represents a spatial region of neural elements influenced by an electric field. The electric field is used to evaluate stimulation patterns which are obtained from a forward model that describes how the current-injected pattern, delivered through a set of active electrode contacts, propagates through the surrounding brain tissue. This forward model representation, often solved via Finite Element Model (FEM), is determined by a set of multiple assumptions including tissue conductivity distributions, behavior of the electrode–tissue interface, modeling assumptions used to represent neural activation, and complex electrode geometries \cite{Alonso2018, McCann2019, FrankemolleGilbert2021, zhang2020steering}. While there is no standardized version for computing VTA explicitly, most pipelines discretize the space around the implanted electrode into a regular grid and calculate which nodes are active according to an independent excitability metric.  Some of these excitability metrics used are axon modelling, activation functions, or through electric-field norm modelling \cite{Duffley2019}. Several DBS studies has applied to some degree an interpretation of VTA modelling in their pipelines, for instance, in the works of \cite{aastrom2014relationship} with Lead-DBS \cite{horn2019lead}, \cite{madler2012activated} with DBSproc \cite{Lauro_2016_DBSProc}, \cite{chaturvedi_2013_ann} with StimVision \cite{noecker2021stimvision}, and \cite{Mikroulis_2025_DBS_MRI} with StimFit \cite{roediger2022stimfit}.
}

%%% Uncertainty
\textcolor{black}{
Uncertainty in DBS modeling, however, reflects the limited precision with which the physical system represented in the computational model can be specified \cite{schmidt2016uncertainty, athawale2019statistical}. In our framework, we take into account this matter in the form of deviations of spatial position of the implanted electrodes which, for instance, this may occur during surgical implantation or when reconstructing electrode trajectories \cite{willsie2015computational, Lasica2025LocalizationPrecision}. When stimulation configurations are evaluated through mathematical optimization methods, these form of variability appear as deviations in the electric fields predicted by the forward model \cite{pena_2017_particle, cubo2019dbs_optimization, Anderson_2021, Mikroulis_2025_DBS_MRI}. In our work, we model uncertainty through a stochastic perturbation of the forward representation linking electrode currents to spatial electric-field distributions. The perturbations are modeled using a \textit{Gaussian noise} assumption, which provides a neutral representation of aggregate modeling variability while allowing our optimization framework to evaluate stimulation configurations under controlled deviations of the forward model.
}
\textcolor{black}{Our optimization framework consists of retrieving optimal electrode configurations by solving an inverse problem that maps desired target activation to stimulation amplitudes using the metaheuristic \textit{L1-norm regularized L1-norm fitting} (L1L1) method~\cite{galazprieto_2022_L1vsL2} for 8- and 40-contact DBS lead models corresponding to segmented and high-density directional lead designs reported in clinical and experimental studies, respectively. Segmented 8-contact directional leads have been widely investigated in clinical DBS~\cite{schnitzler2022directional, steigerwald2019directional}, whereas high-density 40-contact directional leads have been evaluated for intraoperative stimulation~\cite{contarino2014directional, Anderson_2018, masuda2022surgical}.} The L1L1 formulation \textcolor{black}{minimizes the mismatch between the induced and prescribed target fields while} penalizes both the control vector and the data misfit using L1 norms, which promotes sparse stimulation patterns and \textcolor{black}{provides robustness to modeling variability and noise compared with least-squares approaches}. For comparison, configurations based on the \textit{Reciprocity Principle} (RP) are also considered, where reciprocity theory is used to infer stimulation patterns that maximize current density at a target location \cite{FERNANDEZCORAZZA2020116403} \textcolor{black}{but do not explicitly control off-target stimulation}. In addition, \textit{Tikhonov Regularized Least Squares} (TLS) solutions are included as a classical regularization approach that balances data fidelity and solution smoothness through L2-norm penalization \cite{dmochowski2011optimized} \textcolor{black}{and typically yields more spatially distributed stimulation patterns}. In this context, the terms ``L1'' and ``L2'' refer to the discrete 1- and 2-norms of finite-dimensional vectors. The capitalized form \textit{L1L1} denotes the algorithm name and does not refer to the Lebesgue functional spaces $L^{1}$ or $L^{2}$.

\section{Materials and Methods}
\label{sec:Materials}
The quasi-static forward problem was modeled using the Complete Electrode Model (CEM) formulation~\cite{pursiainen2017forward} and solved with the Finite Element Method (FEM). \textcolor{black}{The resulting electric-field distributions were evaluated over the volumetric mesh using interpolation procedures consistent with}~\cite{pursiainen2016, miinalainen2019}. All modeling, meshing, and optimization procedures were performed using the Zeffiro Interface (ZI) version~6.1.0~\cite{pursiainen2025zeffiro}, \textcolor{black}{an open-source MATLAB-based framework for constructing unstructured, boundary-fitted, multi-compartment tetrahedral FE meshes from MRI-derived segmentations}.

\subsection{Realistic Head Model}
\label{sec:head_model}
We generated a tetrahedral-based mesh with a nominal element size of 3.0~mm from an open access T1 weighted MRI data set from a neurologically healthy adult \cite{piastra_maria_carla_2020_3888381}. We performed cortical segmentations using FreeSurfer \cite{fischl2012freesurfer}, supported by functions from SPM12 package \cite{ASHBURNER20091163}. Deep nuclei labels were extracted according to the Aseg atlas (Fig.~\ref{fig:brain_probe}) resulting in a realistic head model with 18 anatomically distinct tissue compartments. \textcolor{black}{A cylindrical compartment with a 1.27~mm outer diameter was incorporated to represent the DBS lead geometry and interpolated using in-house scripts, with the electrode–tissue interface represented through the CEM formulation, where interfacial effects are incorporated via a lumped impedance parameter at the boundary.}

%A double-layered cylindrical compartment with a 1.27~mm outer diameter, and a 0.5~mm thick encapsulation layer \cite{anderson_2019_anodic} was further added and interpolated using in-house scripts. The cylinder functions as a proxy for setting the electrodes with their respective geometric properties. 

As \textcolor{black}{the target sub-compartment is} the \textit{anterior nucleus of the thalamus} (ANT), consistent with ANT-DBS applications in epilepsy \cite{LEHTIMAKI_2016_DBS_targetting, Jarvenpaa_2020_Improving_ANTDBS}, the thalamus compartment underwent three levels of nested \textcolor{black}{refinement, where each tetrahedron was subdivided at each step into four congruent sub-elements} \cite{galazprieto_2023_mesh}. This procedure yielded element diameters ranging from 0.1725 to 0.375~mm \textcolor{black}{within the refined region} and \textcolor{black}{resulted in} a final mesh comprising approximately 4.4 million nodes and 25.4 million tetrahedral elements (Fig.~\ref{fig:thalamus_three_nested}).

\textcolor{black}{Electrical conductivities were assigned as follows:} 0.14 S/m for white matter, 0.33 S/m for gray matter, 0.0064 S/m for skull, 1.79 S/m for cerebrospinal fluid (CSF), 0.33 S/m for subcortical structures, and 0.33 S/m for skin, \textcolor{black}{following} \cite{dannhauer2011modeling}. Subcortical nuclei were \textcolor{black}{assigned the conductivity of gray matter}, while ventricular regions were \textcolor{black}{assigned the conductivity of CSF}, in accordance with \cite{rezaei2021reconstructing} and its referenced sources.

\subsection{Lead Field Matrix}
\begin{figure}[!ht]
\begin{scriptsize}
    \centering
        \begin{subfigure}[t]{4.0cm}
            \centering
            \includegraphics[width = 3.5cm]{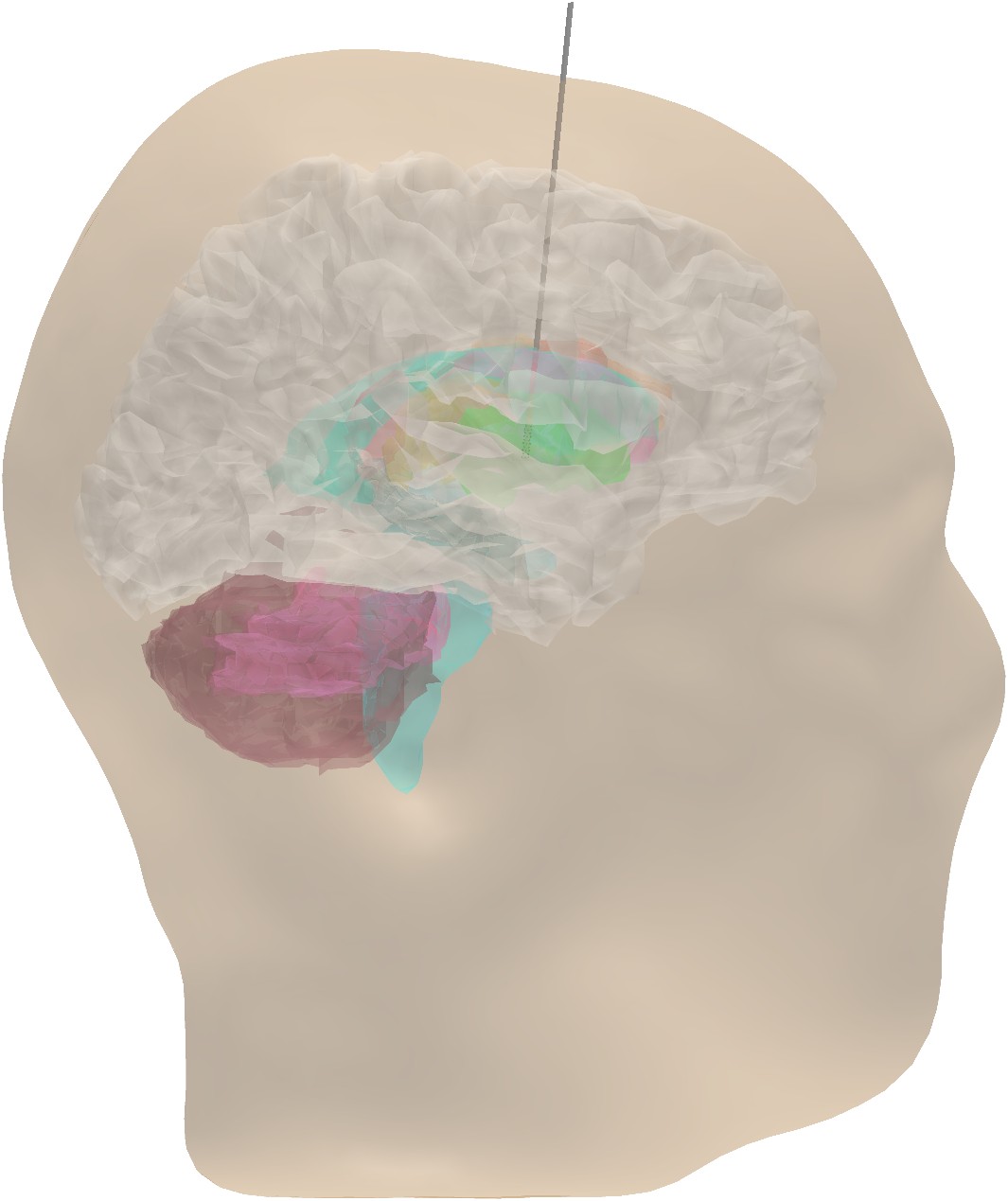}
            \caption{Segmented head model domain depicting main subcortical structures.}
            \label{fig:brain_probe}
        \end{subfigure}
        \hskip0.50cm
         \begin{subfigure}[t]{4.0cm}
            \centering
            \includegraphics[width = 3.50cm, height = 2.50cm]{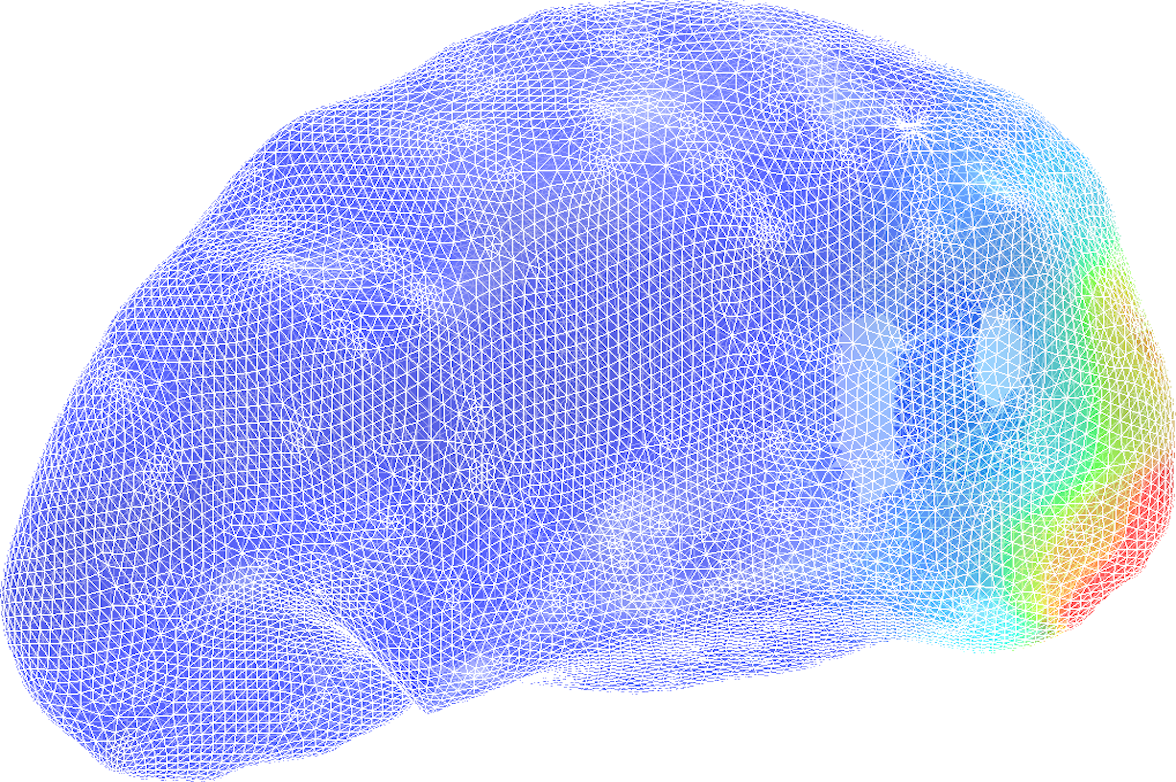} \\
            \caption{Refined unstructured boundary-fitted thalamus sub-compartment with activity (smoothed).}
            \label{fig:thalamus_three_nested}
        \end{subfigure}
        \vskip0.5cm
        \begin{subfigure}[t]{8.0cm}
            \centering
            \includegraphics[width = 3.5cm, height = 2.85cm, angle = -6]{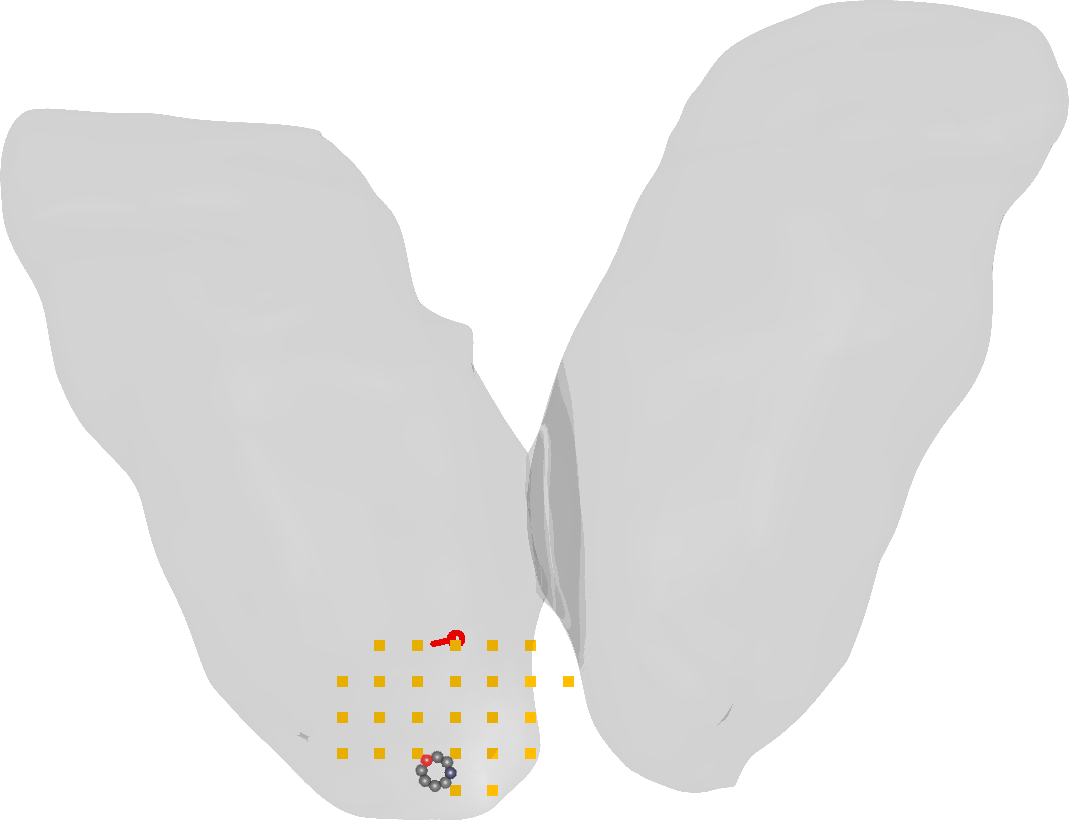}
            \hskip0.25cm
            \includegraphics[width = 3.75cm, height = 2.65cm, angle = 0]{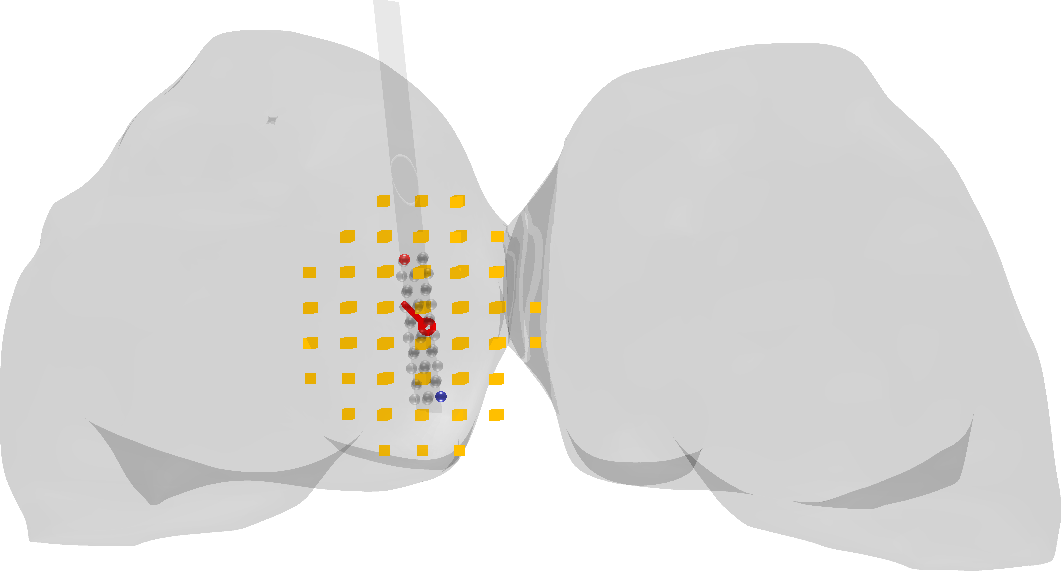} \\
            \vskip0.25cm
            \includegraphics[width = 4.25cm, height = 3.50cm, angle = -3]{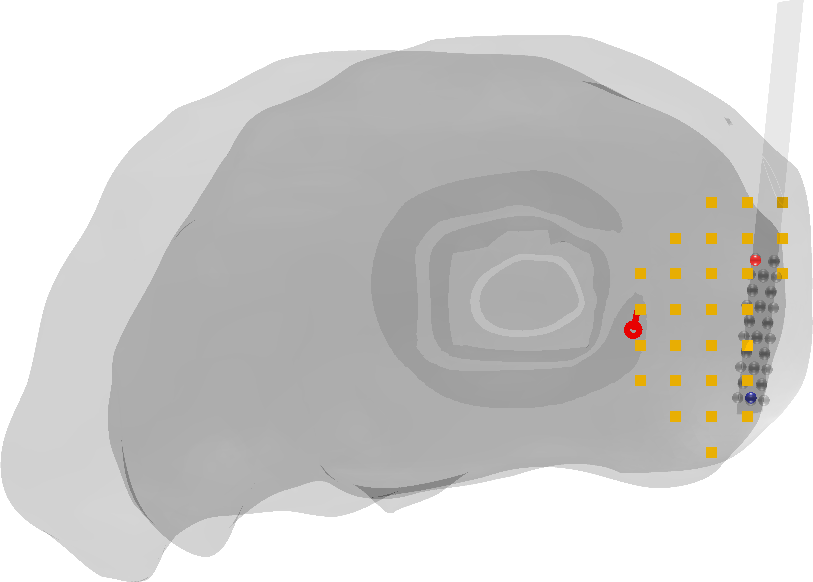}
            \caption{Superior, coronal, and sagittal view of target domain. The interthalamic adhesion can be easily noticed.}
            \label{fig:thalamus_views}
        \end{subfigure}
        \vskip0.35cm
\caption{
{\bf (\subref{fig:brain_probe})}: A head model including skull, forebrain, midbrain, hindbrain structures as well as a directional deep brain stimulation (DBS) lead positioned at the right-hemispheric anterior nucleus of thalamus (ANT). {\bf (\subref{fig:thalamus_three_nested})}: Refined unstructured boundary-fitted thalamus sub-compartment with the distribution of tetrahedral mesh edges and an interpolated current intensity distribution shown on the surface. {\bf (\subref{fig:thalamus_views})}: Placement of a 40-contact dDBS lead with a bipolar electrode configuration (anodal as red, and cathodal as blue). The dipolar target current (red) is depicted by the red pin, indicating the target position and orientation; a 6.0 millimeter (mm) fitting region is displayed as an arrangement of source points (yellow squares) for High-Resolution Lead Field (HR-LF) matrix.
}
\label{fig:reconstruction_head}
\end{scriptsize}
\end{figure}

To establish a baseline representation for spatial attenuation of the lead field and its correspondence to feasible VTA coverage, we derived an initial low-resolution lead field (LR-LF) comprising $244 \times 3$ dipolar DOFs \textcolor{black}{uniformly distributed} throughout the thalamus. The spatial extent for detailed current steering was then guided by reported activation radii around the electrodes, where typical stimulation produces tissue activation on the order of approximately 2.0--4.0~mm \cite{butson_2007_patient, miocinovic_2006_subthalamic}, and may extend to roughly 3.8--4.0~mm under clinically relevant parameters \cite{madler2012activated}. \textcolor{black}{Based on these observations,} we then \textcolor{black}{defined} a spherical region of 6.0~mm radius centered on the lead to represent the domain in which fine-grained modulation is physiologically meaningful. A high-resolution lead field (HR-LF) was subsequently generated by increasing the dipole density to $3620 \times 3$ DOFs across the same anatomical region \textcolor{black}{and restricting the set to} those sources falling within the 6.0~mm domain (Fig.~\ref{fig:thalamus_views}). \textcolor{black}{Further details of the lead-field formulation are provided in} \ref{App:sec:Math_Lyx}.

%% Location of the probe

        % \begin{subfigure}[t]{4.0cm}
        %     \centering
        %     \includegraphics[width = 0.5cm, height = 3.5cm, angle = -1]{Figure_brain_reconstructions/AB_probe.png}
        %     \hskip0.75cm
        %     \includegraphics[width = 0.5cm, height = 3.5cm, angle = -2]{Figure_brain_reconstructions/MED_probe.png}
        %     \caption{8- and 40-electrode contacts, respectively.}
        %     \label{fig:dDBS_probes}
        % \end{subfigure}

\subsection{Optimization Framework}
\label{sec:Optimization}
The applied optimization framework aims to return a \emph{focal} anodal–cathodal electrode configuration that best reproduces the volumetric current distribution associated with a prescribed target \emph{current dipole model}~\citep{demunck1988, mosher2002multiple, medani_2012_ModellingDipoles}. \textcolor{black}{The formulation promotes alignment of the induced current field with the target dipole while suppressing components associated with off-target activation.} Active electrodes are selected iteratively by adjusting the current pattern $\mathbf{y}$ under a set of predefined hyperparameters and an expected number of $K$~active contacts, depending on the chosen algorithm (e.g., Reciprocity Principle, L1L1 Method, or Tikhonov-regularized least squares). \textcolor{black}{The resulting configurations correspond to the candidate stimulation pattern evaluated according to their ability to reproduce the prescribed target field.} 

To restrict the number of admissible solutions, a metaheuristic \textcolor{black}{design is employed} to explore the feasible design space and identify current patterns that best satisfy the required \textcolor{black}{computational} performance characteristics~\cite{galazprieto_2024_lattice}. We impose a lower bound on the focused current density $\Gamma$ ($\Gamma \ge \Gamma_0 = 0.80~\mathrm{mA}$) \textcolor{black}{to enforce a minimum level of target engagement}, as well as an upper bound on the nuisance current density $\Xi$ \textcolor{black}{to limit unintended activation in surrounding regions}. The \emph{candidate solution} is thus defined as the current pattern that maximizes the \textit{field ratio} between target- and nuisance-region current densities while satisfying these constraints.

As safety constraints, the total (absolute) injected current is constrained to $\mu = 4.0~\mathrm{mA}$ \eqref{eq:L1L1Method_2}, with per-contact limits $\gamma = \pm 2.0~\mathrm{mA}$ \eqref{eq:L1L1Method_3}. \textcolor{black}{Charge balance is enforced through the net-current constraint} \eqref{eq:L1L1Method_4}. The stimulation is treated as time-invariant, \textcolor{black}{focusing on the spatial distribution of the induced current density}, rather than \textcolor{black}{temporal waveform characteristics or frequency-dependent neural responses}.

\textcolor{black}{As reference magnitude,} $\|{\bf x}_1\|_2 = 3.85$~A/m\textsuperscript{2} \textcolor{black}{is used as the \textit{activation threshold} for the current dipole moment density which represents a rough approximate threshold for the excitation of nerve fibers}. \textcolor{black}{This value is derived by dividing} the average \textcolor{black}{current} dipole moment density of pyramidal neurons (0.77~nAm/mm\textsuperscript{2}) \textcolor{black}{by an approximate dendritic axis length of} 0.2~mm \textcolor{black}{for patients with epileptogenic responses} \cite{Oishi20021390, MURAKAMI201549}. In this study, it is considered as a reference to normalize the present optimization framework.

\subsubsection{Electrode Configurations}
\label{sec:ele_conf}
Two clinically-inspired electrode configurations are considered, corresponding to $K = 8$ and $K = 40$ contact arrangements. The 8-contact configuration consists of four axially cylindrical bands (1.5~mm length, 0.5~mm spacing), with the two central bands subdivided into three segments, yielding a ``1–3–3–1'' directional layout. The 40-contact configuration comprises eight longitudinal rows of five ellipsoidal contacts (0.8~mm width, 0.66~mm height), spaced 1.5~mm along the lead axis and arranged with a 0.75~mm helical offset between adjacent rows.

% \begin{figure}[!ht]
%     \centering
%         \includegraphics[width = 0.75cm, height = 4.5cm, angle = 87]{Figure_brain_reconstructions/AB_probe.png}
%         \hskip1.25cm
%         \includegraphics[width = 0.75cm, height = 4.5cm, angle = 86]{Figure_brain_reconstructions/MED_probe.png}
%         \caption{8- and 40-electrode contacts, respectively.}
%         \label{fig:dDBS_probes}
% \end{figure}

\textcolor{black}{While physical electrode designs follow cylindrical and segmented geometries, in the present framework the electrode contact interfaces are modeled as equipotential spherical surfaces with a lumped contact impedance \cite{butson2006sources,Lempka_2009}, consistent with the CEM formulation in which electrode--tissue interactions are represented through a boundary impedance rather than an explicit geometric layer \cite{pursiainen2012complete}. Consequently, the distinction between the two lead models arises from the number and spatial arrangement of the contacts rather than from the detailed directional current distribution of individual segmented electrodes. The present comparison therefore focuses on the effect of contact layout on the optimization process instead of reproducing the exact steering characteristics of physical directional DBS leads.} 

A uniform impedance value of 2.0~k$\Omega$ was adopted for the present numerical model because it yields the greatest lead field amplitudes within the feasible impedance interval from 50~$\Omega$ to 2.0~k$\Omega$ suggested in \cite{butson2006sources}. \textcolor{black}{In our numerical comparison, increasing the impedance across this interval changed the median lead field strength by only approximately 14--19\%, whereas increasing the contact count from 8 to 40 produced an approximately 40\% increase. These results indicate that, within the present current-driven framework, the optimization is considerably more sensitive to the spatial contact layout than to feasible variations in contact impedance. A uniform impedance value of 2.0~k$\Omega$ was therefore adopted to provide a consistent modeling basis while preserving the relative differences between probe configurations.}

\subsubsection{Bipolar Configurations}
\label{sec:bipolar_conf}

\begin{figure*}[ht]
    \begin{scriptsize}
    \centering
    %%%%%%%%%%%%%%%%%%%%%%%%%%%%%%%%%%%%%%%%
    \begin{minipage}{18.0cm}
    \centering
    {\bf 8-contact Electrode Configuration} \\
        \begin{tabular}{cc|cc|cc|cc}
        %\hline
        \begin{minipage}[b]{1.0cm}
            \RaggedLeft
            \includegraphics[width=0.90cm, height = 4.5cm]{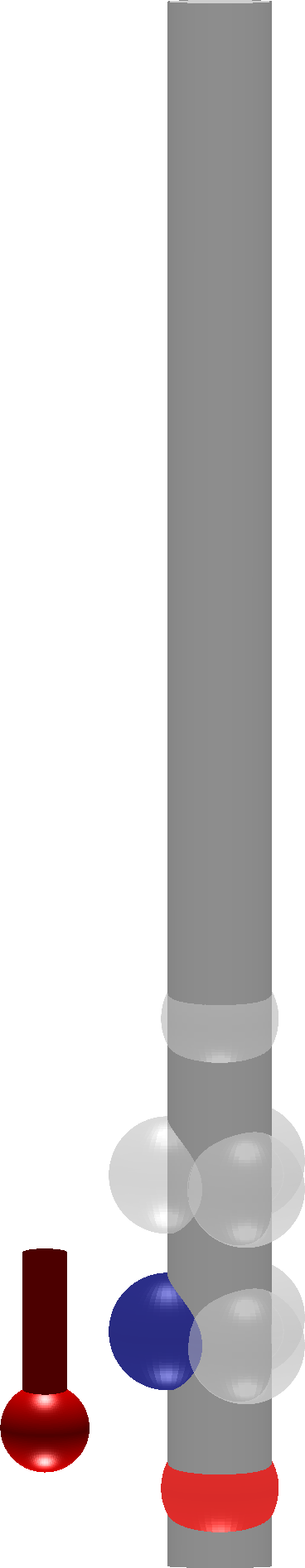}     
        \end{minipage} &
        \begin{minipage}[b]{2.85cm}
            \centering
            \vskip0.2cm
            RP \\
            \includegraphics[width=2.85cm, height=2.00cm]{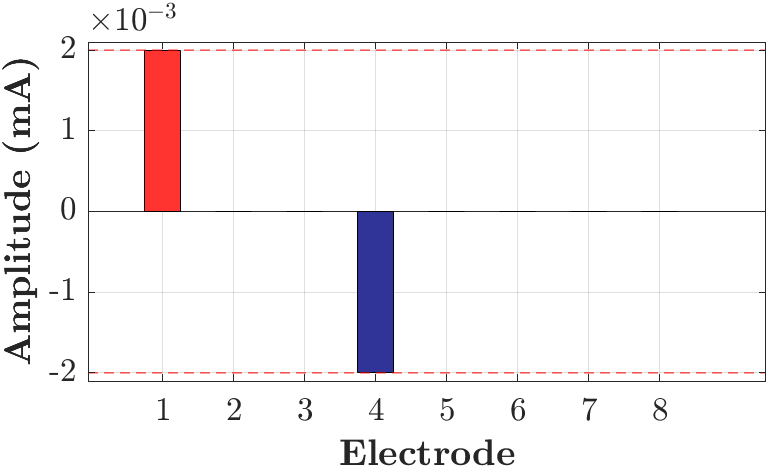} \\ \vskip0.2cm
            \includegraphics[width=1.50cm, height=2.00cm]{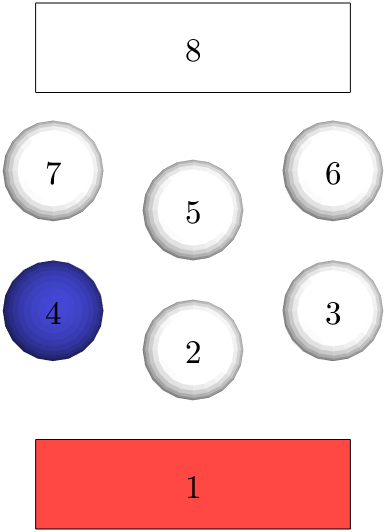} \vskip0.2cm
        \end{minipage} &
        \includegraphics[width=0.50cm, height = 4.5cm]{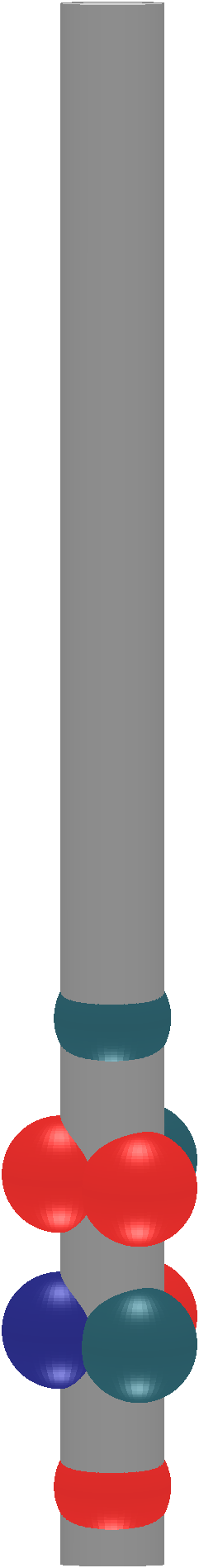} &
        \begin{minipage}[b]{2.85cm}
            \centering
            \vskip0.2cm
            L1L1(A) \\
            \includegraphics[width=2.85cm, height=2.00cm]{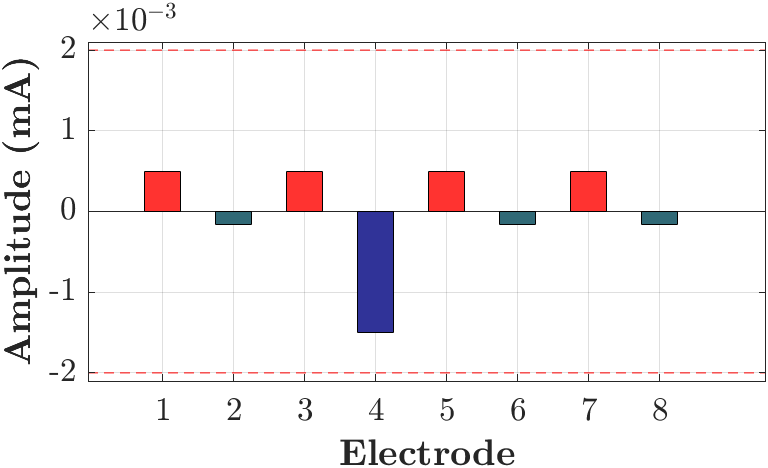} \\ \vskip0.2cm
            \includegraphics[width=1.50cm, height=2.00cm]{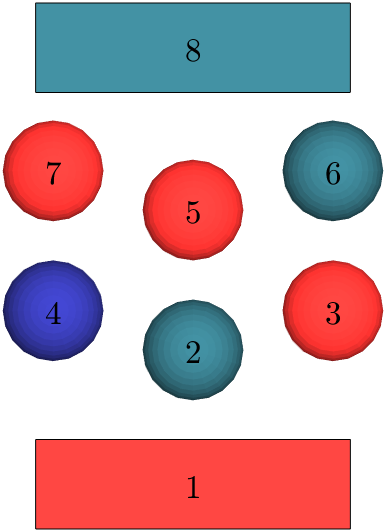} \vskip0.2cm 
        \end{minipage} &
        \includegraphics[width=0.50cm, height = 4.5cm]{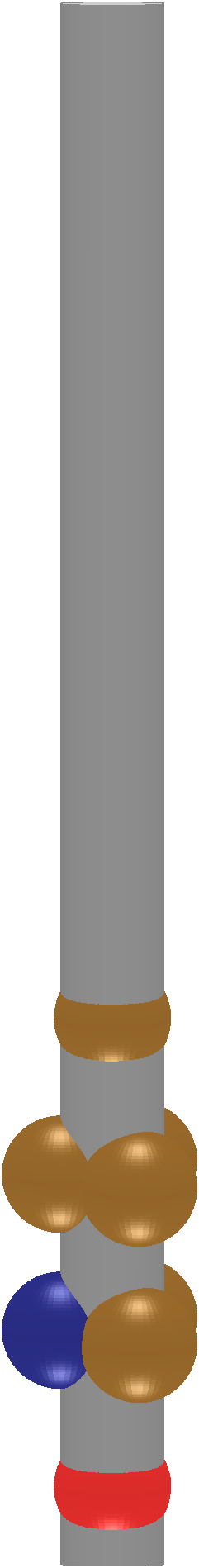} &
        \begin{minipage}[b]{2.85cm}
            \centering
            \vskip0.2cm
            L1L1(B) \\
            \includegraphics[width=2.85cm, height=2.00cm]{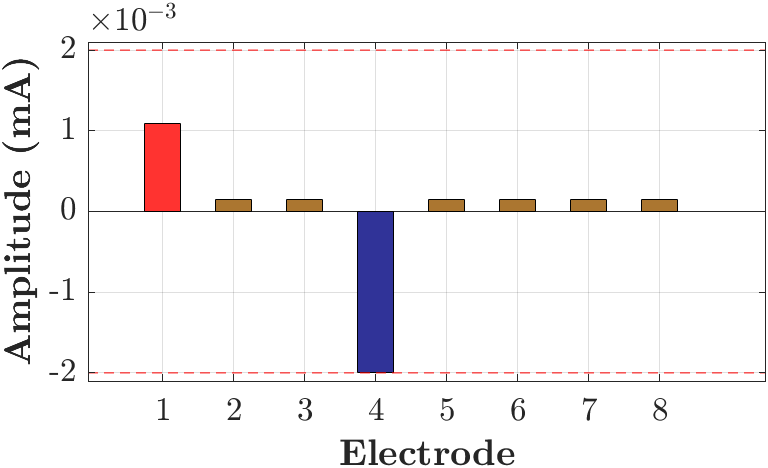} \\ \vskip0.2cm
            \includegraphics[width=1.50cm, height=2.00cm]{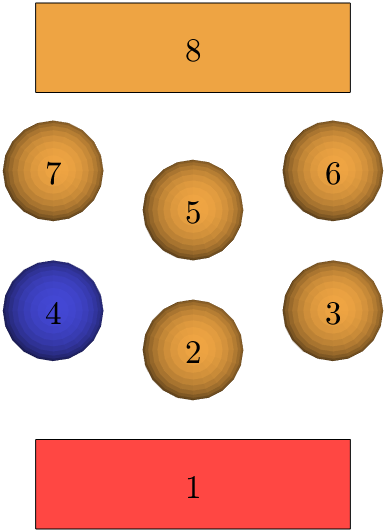} \vskip0.2cm
        \end{minipage} &
        \includegraphics[width=0.50cm, height = 4.5cm]{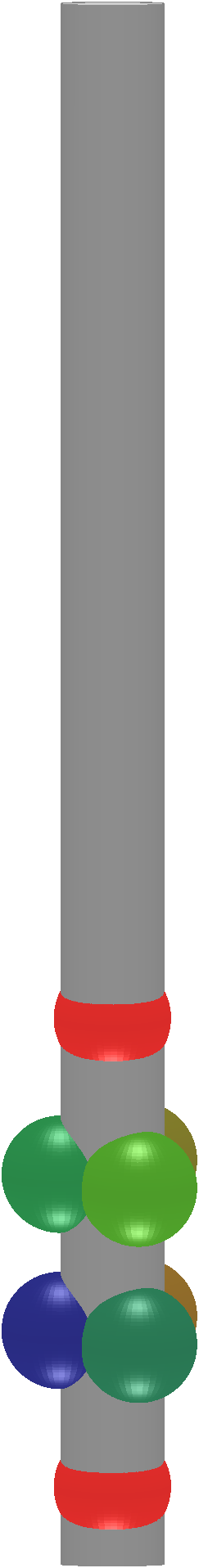} &
        \begin{minipage}[b]{2.85cm}
            \centering
            \vskip0.2cm
            TLS \\
            \includegraphics[width=2.85cm, height=2.00cm]{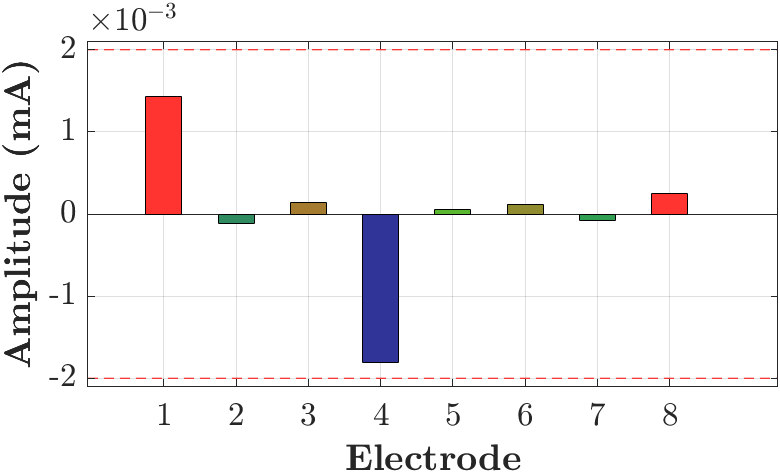} \\ \vskip0.2cm
            \includegraphics[width=1.50cm, height=2.00cm]{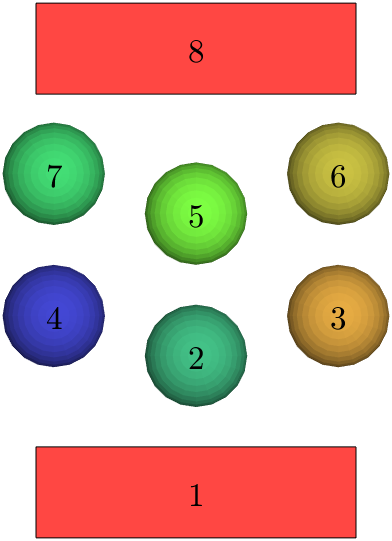} \vskip0.2cm
        \end{minipage} \\
        %\hline
        \end{tabular}
    \end{minipage}
    %%%%%%%%%%%%%%%%%%%%%%%%%%%%%%%%%%%%%%%%
    \vskip0.5cm
    %%%%%%%%%%%%%%%%%%%%%%%%%%%%%%%%%%%%%%
    \begin{minipage}{18.0cm}
    \centering
    {\bf 40-contact Electrode Configuration} \\
        \begin{tabular}{cc|cc|cc|cc}
            %\hline
        \begin{minipage}[b]{1.0cm}
            \RaggedLeft
            \includegraphics[width=1.0cm, height = 4.5cm]{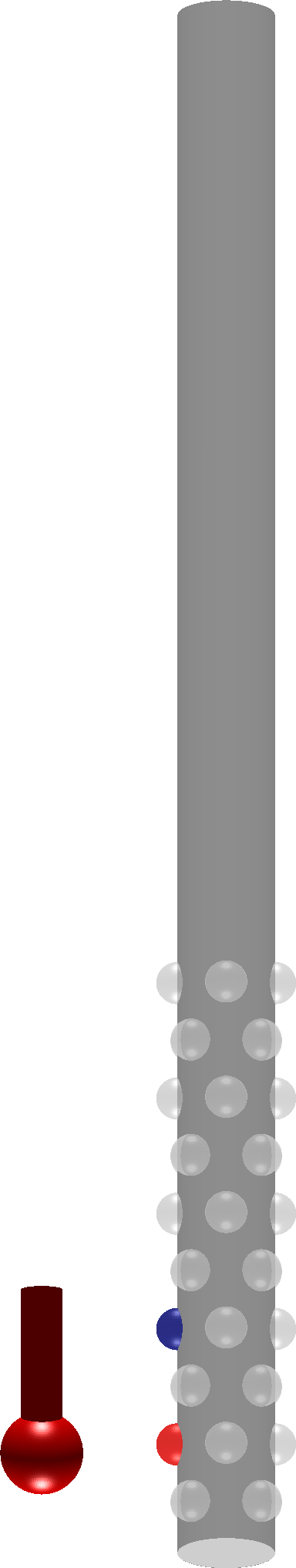}     
        \end{minipage} &
            \begin{minipage}[b]{2.85cm}
                \centering
                \vskip0.2cm
                RP \\
                \includegraphics[width=2.85cm, height=2.00cm]{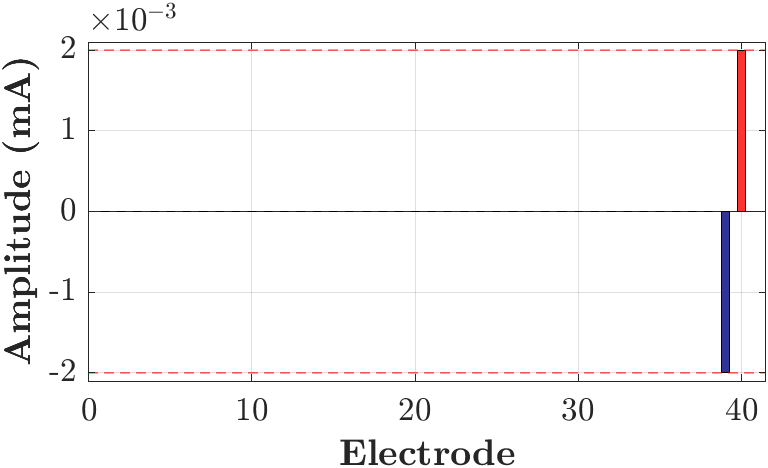} \\ \vskip0.2cm
                \includegraphics[width=2.85cm, height=2.00cm]{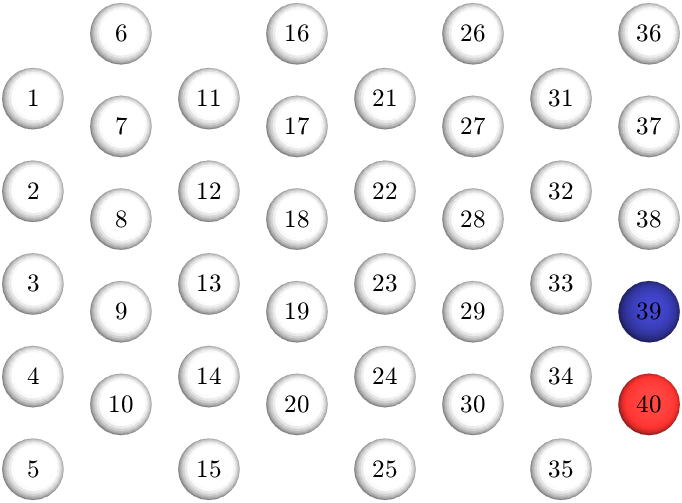} \vskip0.2cm
            \end{minipage} &
            \includegraphics[width=0.50cm, height = 4.5cm]{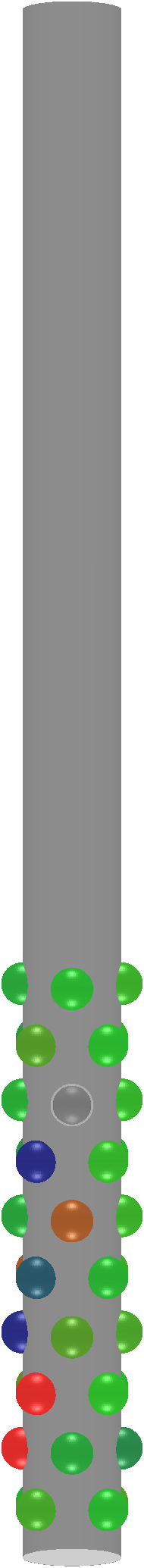} &
            \begin{minipage}[b]{2.85cm}
                \centering
                \vskip0.2cm
                L1L1(A) \\
                \includegraphics[width=2.85cm, height=2.00cm]{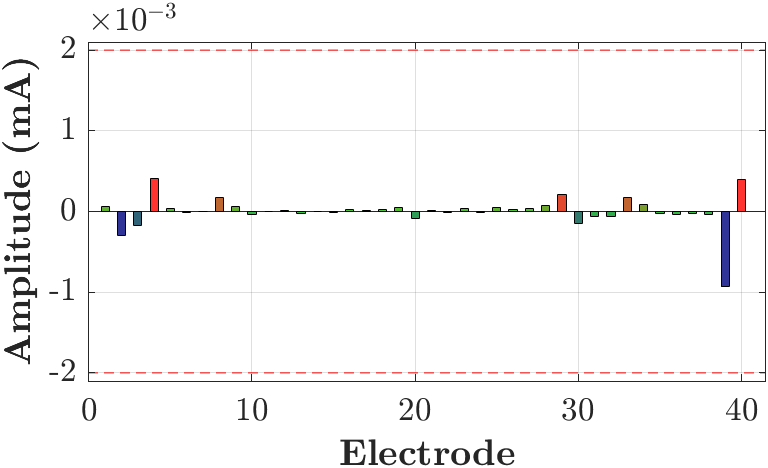} \\ \vskip0.2cm
                \includegraphics[width=2.85cm, height=2.00cm]{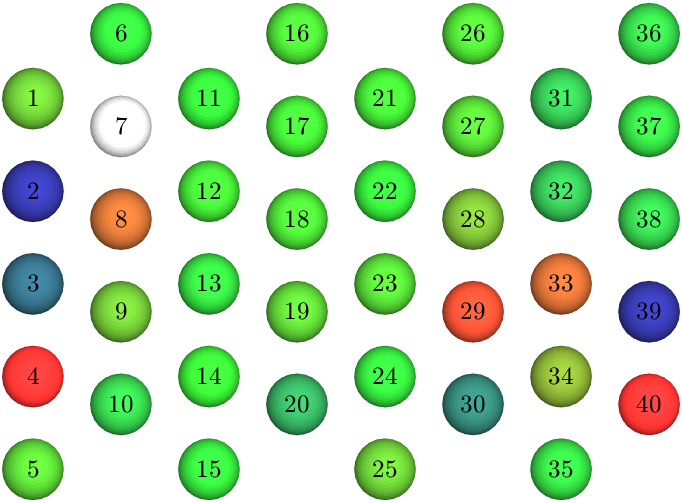} \vskip0.2cm 
            \end{minipage} &
            \includegraphics[width=0.50cm, height = 4.5cm]{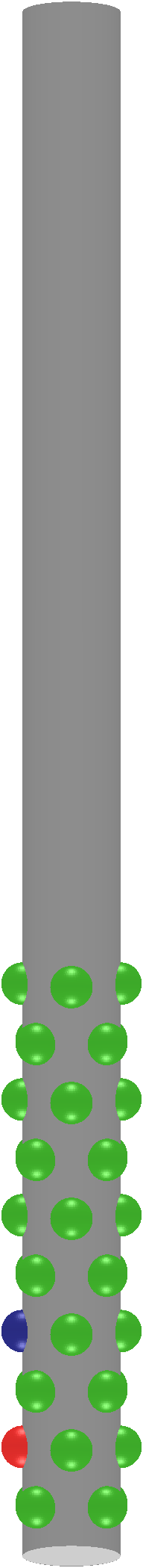} &
            \begin{minipage}[b]{2.85cm}
                \centering
                \vskip0.2cm
                L1L1(B) \\
                \includegraphics[width=2.85cm, height=2.00cm]{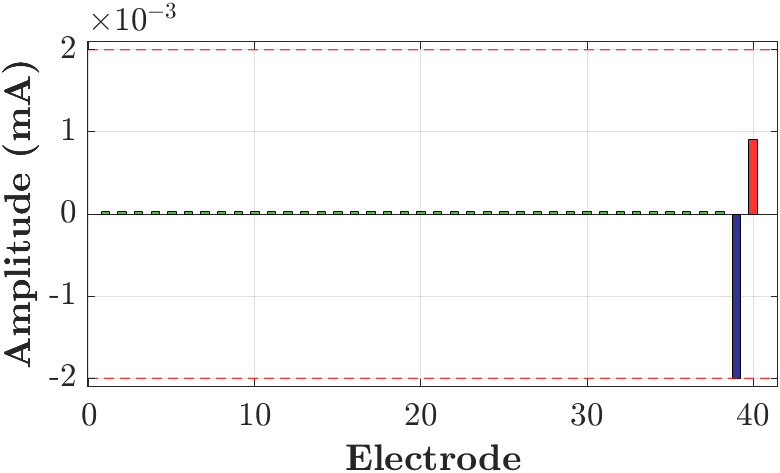} \\ \vskip0.2cm
                \includegraphics[width=2.85cm, height=2.00cm]{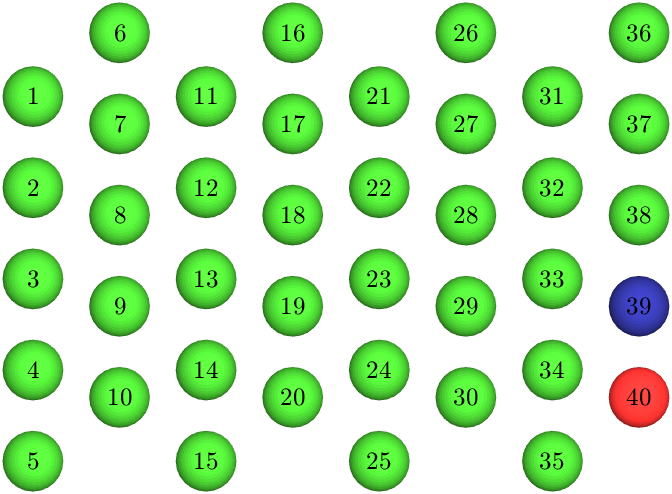} \vskip0.2cm
            \end{minipage} &
            \includegraphics[width=0.50cm, height = 4.5cm]{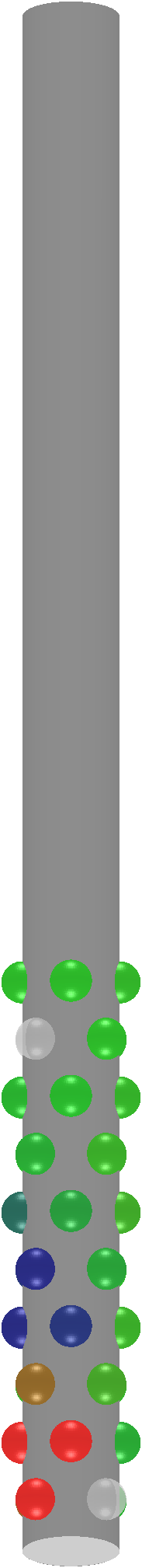} &
            \begin{minipage}[b]{2.85cm}
                \centering
                \vskip0.2cm
                TLS \\
                \includegraphics[width=2.85cm, height=2.00cm]{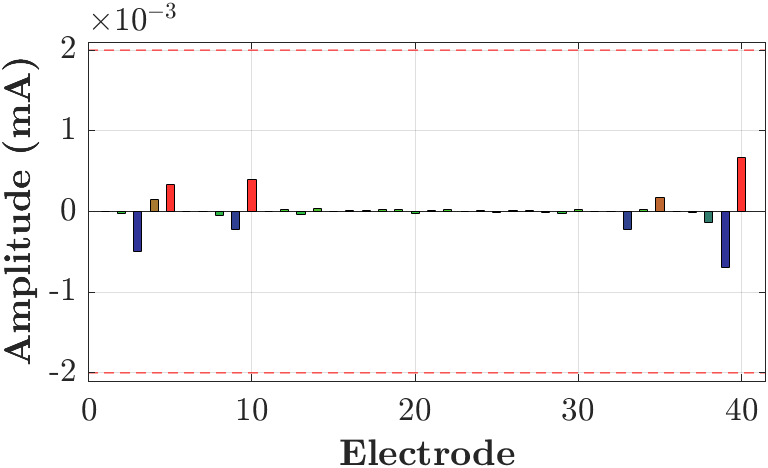} \\ \vskip0.2cm
                \includegraphics[width=2.85cm, height=2.00cm]{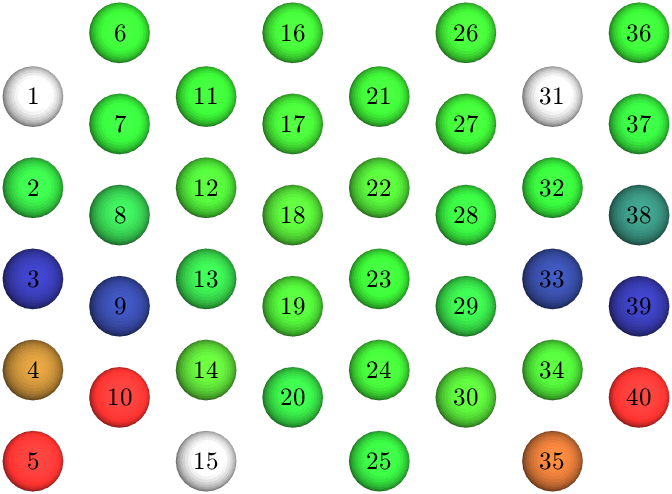} \vskip0.2cm
            \end{minipage} \\
            %\hline
        \end{tabular}
    \end{minipage}
    %%%%%%%%%%%%%%%%%%%%%%%%%%%%%%%%%%%%%%
    \vskip0.15cm
    %%%%%%%%%%%%%%%%%%%%%%%%%%%%%%%%%%%%%%
    \begin{minipage}{4.0cm}
        \RaggedLeft
        \includegraphics[height = 0.65cm, width = 3.0cm]{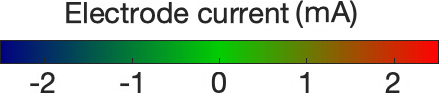}
    \end{minipage}
    %%%%%%%%%%%%%%%%%%%%%%%%%%%%%%%%%%%%%%%
    \caption{
    Schematic interpretation of the current injection patterns obtained using the Reciprocity Principle (RP), metaheuristic L1-norm regularized L1-norm fitting (L1L1) method under two different nuisance current thresholds (L1L1(A) with $\varepsilon = {0 : -160}$ and L1L1(B) with $\varepsilon = {0 : -10}$), and the Tikhonov-Regularized Least Square (TLS) method. \textcolor{black}{The target dipolar source, defined by its position and orientation (here, parallel to the lead), is indicated by the red arrow.} The bar plots display the resulting electrode current amplitudes, while the contact maps visualize their spatial distribution over the electrodes; cathodal (blue), anodal (red), and inactive (gray/green). The RP method yields a bipolar configuration maximizing focused current density $\Gamma$, while the remaining methods search for multipolar configurations maximizing focality $\Theta = \Gamma / \Xi$ under an \textit{a priori} constraint $\varepsilon$ on the nuisance current density $\Xi$, determined by the peak signal-to-noise ratio (PSNR) $\delta^{-1}$ and the lead field attenuation. The L1L1(B) result reflects more strongly constrained $\Xi$ (smaller dynamic range), activating fewer contacts whereas L1L1(A) and TLS allow higher $\Theta$ values and broader excitation patterns. \textcolor{black}{The current on the electrode contacts is bounded by $\pm 2.0$ mA. Electrode–tissue interfaces are modeled using the Complete Electrode Model (CEM), where each contact is treated as an equipotential surface with a uniform lumped impedance (2.0~k$\Omega$) across all contacts; contact-specific impedance variations due to surface area or segmentation are not explicitly modeled.}
    }
    \label{fig:probes_examples}
    \end{scriptsize}
\end{figure*}

Bipolar patterns are obtained using the reciprocity principle (RP)~\cite{FERNANDEZCORAZZA2020116403}. The principle states that the maximal focused current density $\Gamma_{\mathrm{max}}$ for a linear quasi-static system is achieved by selecting a single anode–cathode pair. In practice, the target dipolar field is multiplied by the transpose of the lead field matrix, and the contacts corresponding to the strongest positive and negative entries are chosen as anode and cathode, respectively, with all remaining contacts set to zero (Fig.~\ref{fig:probes_examples}). \textcolor{black}{The resulting stimulation pattern aligns with the direction of the target field at the location of interest.} RP maximizes the delivered current density at the target \textcolor{black}{through concentration of current between a single electrode pair}, but \textcolor{black}{it does not explicitly control the spatial distribution of the induced field outside the target region}, and thus may yield suboptimal field selectivity.

\subsubsection{Multipolar Configurations}
\label{sec:multipolar_conf}

Multipolar configurations are derived using the \emph{L1-norm regularized L1-norm fitting} (L1L1) method~\cite{galazprieto_2022_L1vsL2} and \emph{Tikhonov-regularized least squares} (TLS)~\cite{dmochowski2011optimized}. Both methods assume a predetermined set of candidate active contacts and estimate the corresponding stimulation amplitudes under polarity, charge, and safety constraints (detailed in~\ref{App:sec:L1L1}). The L1L1's formulation combines an $L_1$ penalty on the stimulation vector with an $L_1$-norm fitting term applied to \textcolor{black}{the focused-field mismatch $\mathbf{L}_1\mathbf{y}-\mathbf{x}_1$, representing the deviation between the induced electric field in the target region and the prescribed target field template (see Eq.~\ref{eq:L1L1Method_1}), promoting sparse and interpretable stimulation supports while improving robustness to modeling uncertainty in the lead-field representation}. TLS, in contrast, \textcolor{black}{employs an $L_2$ regularization term that stabilizes the solution in the presence of correlated columns and numerical ill-conditioning, typically producing smoother and more distributed amplitude patterns}. Because the full combinatorial search space of possible contact subsets can grow exponentially with the variation in the hyperparameters (Table~\ref{tab:variables_hyperparameters}), both formulations employ a metaheuristic support-selection strategy to efficiently explore only a tractable subset of candidate configurations while maintaining practical computational cost~\cite{galazprieto_2024_lattice}.

%%%%

\begin{figure}[ht!]
\begin{subfigure}[t]{8.5cm} 
    \begin{minipage}{3.5cm}
        \centering
        \begin{tikzpicture}[scale=0.9]
            \draw[thick] (-1,0) to[out=90,in=180] (0,1.5) to[out=0,in=90] (2.5,0) to[out=-90,in=0] (0,-1.5) to[out=180,in=-90] cycle;
            \node[scale = 0.75] at (0, 1.1) {\textbf{Thalamus}};
            \draw[gray, line width=1.5mm] (1.1,0) -- (0.85,2); 
            \node[gray, scale = 0.75] at (1.7, 1.7) {\textbf{DBS probe}};
            \draw[dashed, blue, thick] (1.1, 0) circle (0.75);
            \draw[solid, <->, blue] (1.1, 0) -- (1.475, 0.6495) ;
            \node[blue, scale = 0.75] at (1.63, 0.85) {\textbf{$d_\delta$}};
            \node[blue, scale = 0.75] at (2.0, 0.5) {\textbf{$\mathcal{V}_\delta$}};
            \def\targetX{1.7}
            \def\targetY{0.3}
            \draw[->, red, thick] (1.6, 0) -- (\targetX, \targetY);
            \draw[solid, <->, red] (1.1,0) -- (1.6,0);
            \node[red, scale = 0.75] at (1.4, 0.15) {\textbf{$d_T$}};
            \draw[<->, violet] (0.35, 0) -- (0.6, 0);
            \node[violet, scale = 0.75] at (0.49, 0.20) {\textbf{$d_\varepsilon$}};
            \node[red, scale = 0.75] at (1.1, -0.2) {\textbf{Target}};
            \draw[dashed, purple, thick] (1.1,0) circle (0.5);
        \end{tikzpicture}
    \end{minipage}
    \hskip0.75cm
    \begin{minipage}{4.0cm}
        \centering
        \begin{tikzpicture}[scale=0.80]
            \draw[thick,->] (1.1, 0.1) -- (5.8, 0.1) ;
            %\node[below, scale=0.5] at (5.5, 0) {\textbf{Distance}};
            \node[below, scale = 0.75] at (3.5, -0.5) {\textbf{Relative Distance from dDBS}};
            \draw[thick,->] (1.1, 0.1) -- (1.1, 4.5) node[above right, brown, scale=0.5] {\textbf{Field Strength}};
            
            %\draw[gray, line width=1.5mm] (1.1,0.1) -- (1.1,3.5);
            \draw[thick, brown, domain = 1.25:5.9, samples=100] plot (\x, {0.25+1.5/((\x-0.65)*(\x-0.65))});
            
            %\node[gray, scale=0.5] at (2,3) {\textbf{DBS probe}};
            
            \draw[thick, dashed, red] (2,0.1) -- (2,1.1);
            \filldraw[red] (2,1.1) circle (0.07);
            \node[red, below, scale = 0.75] at (2,0) {\textbf{Target} \( d_T \)};
            
            \draw[thick, dashed, blue] (4,0.1) -- (4,0.4);
            \filldraw[blue] (4,0.4) circle (0.07);
            \node[blue, below, scale = 0.75] at (4,0) {\( d_\delta \) \textbf{range}};
            
            \draw[<->, solid,  blue] (4,0.4) -- (4,4.5);
            \node[blue, right, scale = 0.75] at (3.7,2.25) {\textbf{\( {\delta} \)} };
            
            \draw[<->, solid,  violet] (4,0.4) -- (4,1.1);
            \draw[dashed, thick, blue] (2,0.4) -- (4,0.4);
            \draw[<->,solid, purple] (2,1.1) -- (4,1.1);
            \node[purple, above, scale = 0.75] at (3,1.1) {\textbf{\( d_{\varepsilon} \)} };
            \node[purple, right, scale = 0.75] at (3.7,0.75) {\textbf{\( {\varepsilon} \)} };
        \end{tikzpicture}
    \end{minipage}
    \caption{ {\bf Left:} The set $\mathcal{V}_\delta$ corresponding to PSNR $\delta^{-1}$ is shown by the blue dashed circle with radius $d_\delta$. The target is at the distance $d_T$ from the probe. Between the target and the boundary of $\mathcal{V}_\delta$ is a range $d_\varepsilon$ in which the lead field decays by the factor $\varepsilon$, which approximates the maximum achievable dynamic range $\varepsilon$ between the {\em focused current density} at the target position and in its surroundings. {\bf Right:} One-dimensional interpretation of the lead field decay.}
    \label{fig:eps_effect_graphs}
\end{subfigure}
\vskip0.5cm
\begin{subfigure}[t]{8.5cm}
    \begin{minipage}{2.0cm}
        \centering
        \includegraphics[width = 1.9cm, trim = 70 0 130 190, clip]{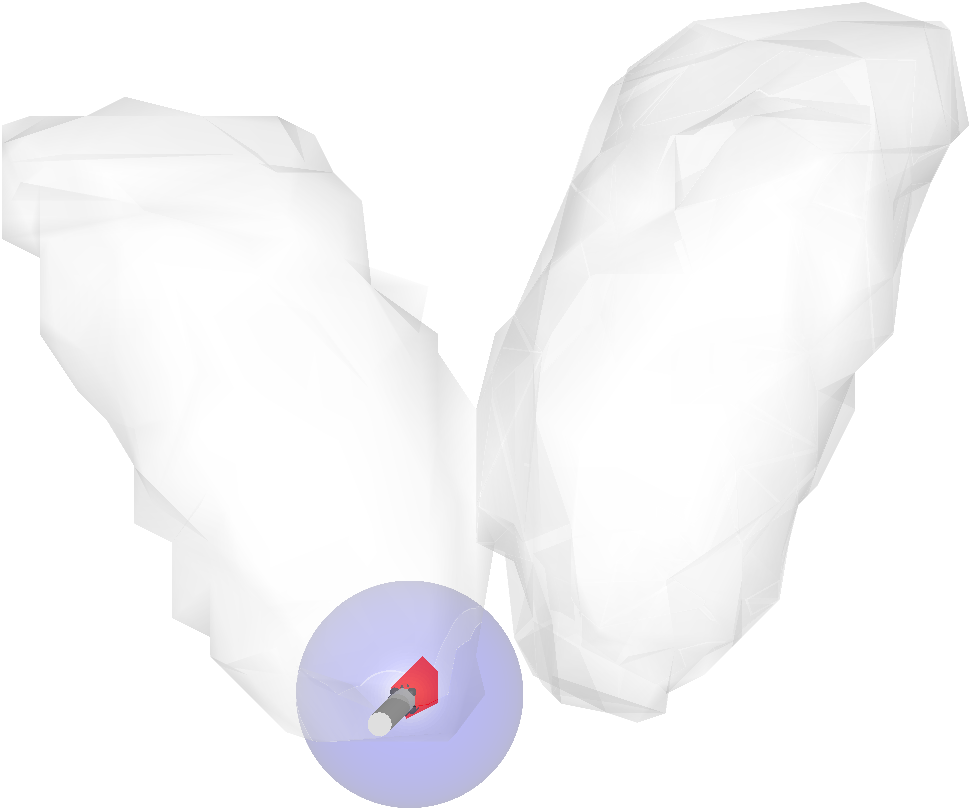} \\
        $\mathcal{V}_{-10 \text{ dB}}$
        \end{minipage}
        \begin{minipage}{2.0cm}
        \centering
        \includegraphics[width = 1.9cm, trim = 70 0 130 190, clip]{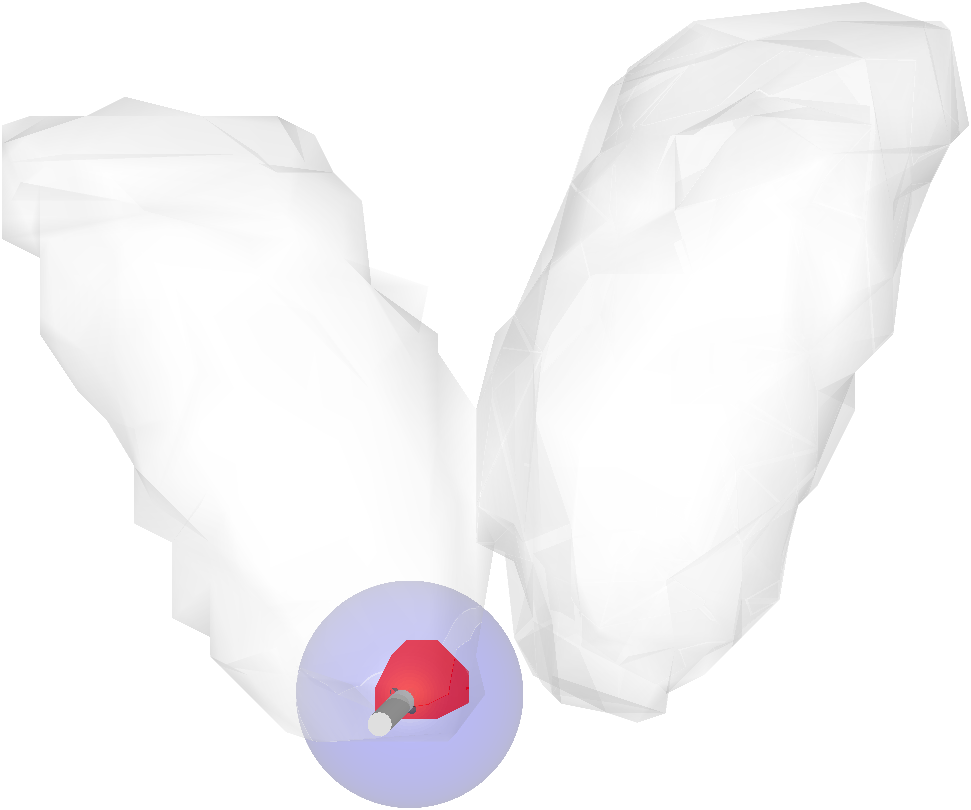} \\ 
        $\mathcal{V}_{-20 \text{ dB}}$
        \end{minipage}
        \begin{minipage}{2.0cm}
        \centering
        \includegraphics[width = 1.9cm, trim = 70 0 130 190, clip]{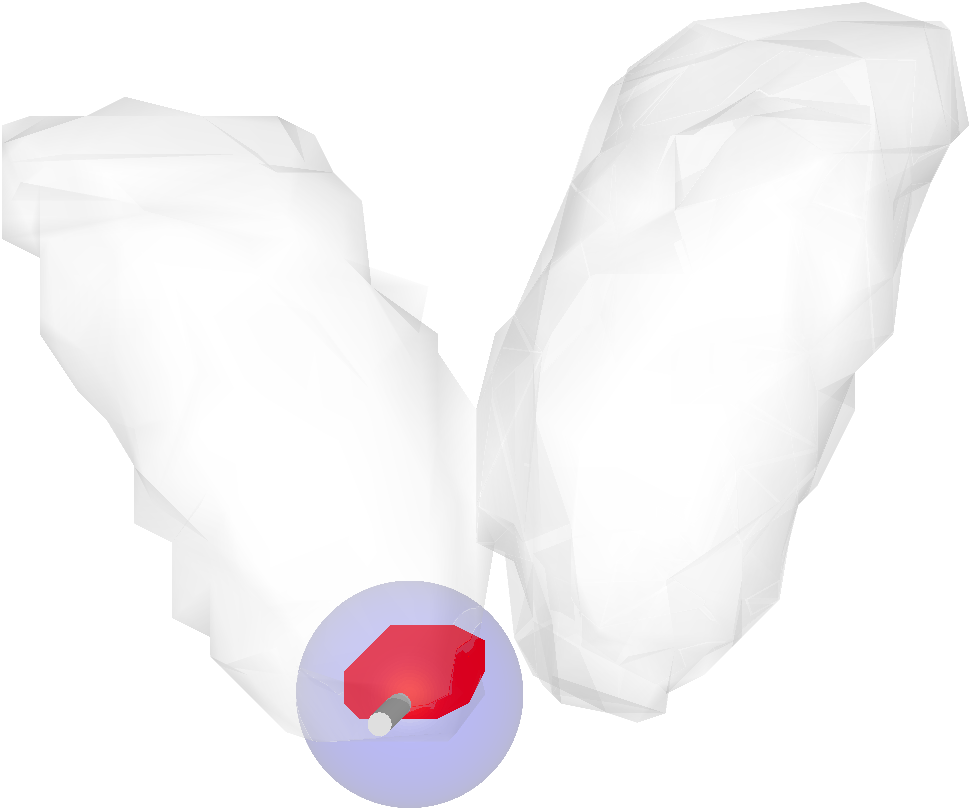} \\ 
        $\mathcal{V}_{-30 \text{ dB}}$
        \end{minipage}
        \begin{minipage}{2.0cm}
        \centering
        \includegraphics[width = 1.9cm, trim = 70 0 130 190, clip]{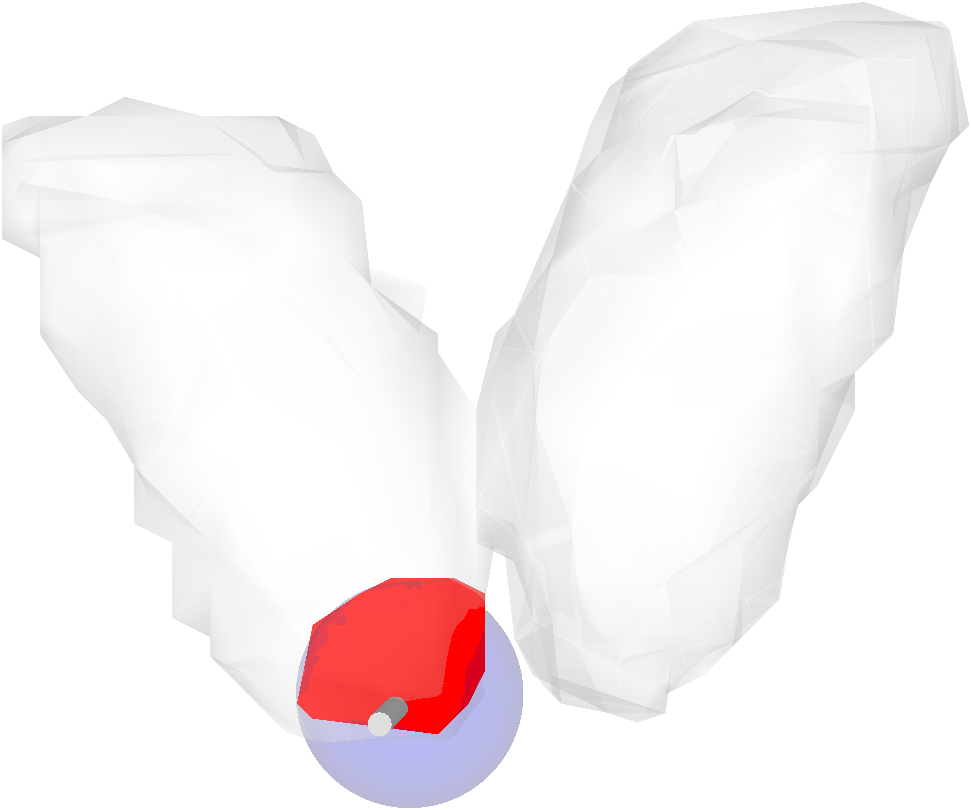} \\ 
        $\mathcal{V}_{-40 \text{ dB}}$
    \end{minipage}
    \caption{The set $\mathcal{V}_\delta$ illustrated for the volumetric model w.r.t PSNR $\delta^{-1}$ levels from 10 to 40 dB.}
    \label{fig:eps_effect_VTA}
\end{subfigure}
\caption{{\bf (\subref{fig:eps_effect_graphs})}: A schematic of 
how  the peak signal-to-noise ratio (PSNR) $\delta^{-1}$ relates to the set  $\mathcal{V}_\delta$ and to the approximate maximum dynamic range $\varepsilon^{-1}$ of the volumetric stimulation current field.  
{\bf (\subref{fig:eps_effect_VTA})}: The set $\mathcal{V}_\delta$ shown for PSNR 10--40 dB, providing a rough approximation of the maximum feasible coverage of volume of tissue activated (VTA) for an activation threshold set by $\delta$.}
\label{fig:epsilon_effect_schematic}
\end{figure}

%%%%
\subsection{Uncertainty Modeling}
\label{sec:uncertain}
\textcolor{black}{
Uncertainty in DBS optimization framework arises from physiological and technical factors that affect the accuracy of the lead field matrix $\mathbf{L}$. These factors include electrode placement imprecision, patient-specific anatomical variation, and inhomogeneities in tissue conductivity~\cite{schmidt2016uncertainty, athawale2019statistical}. These effects introduce deviations in the LF entries, which can be modeled as additive perturbations. The perturbation magnitude reflects aggregate forward‑model variability arising from sources such as electrode localization error, tissue conductivity uncertainty, and modeling discrepancies summarized in Table~\ref{tab:uncertainty}; the variance parameter $\sigma^2$ therefore represents the expected magnitude of these combined perturbations. Under this approximation, uncertainty is represented using a \textit{Gaussian noise model} applied to the LF rows. Specifically, this is modeled as
\[
\tilde{\mathbf{l}}_i^{(\mathrm{obs})}
=
\tilde{\mathbf{l}}_i + \mathbf{n}_i,
\qquad
\mathbf{n}_i \sim \mathcal{N}(\mathbf{0}, \sigma^2 \mathbf{I}),
\]
where $\mathbf{n}_i$ represents uncertainty-induced perturbations arising from forward model variability. Although several uncertainty sources (e.g., electrode localization errors or tissue conductivity variations) may introduce spatially correlated perturbations in the forward model, Gaussian noise model provides a tractable approximation of aggregate modeling variability and is commonly used for robustness analysis in inverse and forward modeling problems.}

\subsubsection{\textcolor{black}{Sensitivity and Row Reduction}}
\textcolor{black}{
Each row corresponds to a spatial degree of freedom (DOF) and describes how the electrode contacts influence the current density at that location. However, only variations in the lead field across contacts contribute to differential stimulation effects. Therefore, each row is first reduced by subtracting its mean value, yielding the zero-averaged LF row
\[
\tilde{\mathbf{l}}_i =
\left(
L_{i,1} - \overline{L}_i,\,
L_{i,2} - \overline{L}_i,\,
\dots,\,
L_{i,K} - \overline{L}_i
\right),
\]
where
\[
\overline{L}_i = \frac{1}{K}\sum_{j=1}^{K} L_{i,j}.
\]
The transformation removes uniform offsets and isolates the effective stimulation sensitivity associated with differential electrode activation. The magnitude of the reduced row is quantified using its Euclidean norm, $\|\tilde{\mathbf{l}}_i\|_2$, which provides a scalar measure of the controllability of the DOF under electrode stimulation.}

\subsubsection{\textcolor{black}{PSNR-Based Feasibility Constraint}}
\textcolor{black}{
Under the additive perturbation model, the detectability of a DOF depends on the ratio between its signal magnitude and the noise level. This relationship is characterized using peak signal-to-noise ratio (PSNR), which defines the minimum signal magnitude required for reliable discrimination from the perturbations. The feasible set of DOFs (or rather, the \textit{feasible coverage range for VTA}), determined by this detectability constraint, is therefore defined as
\[
\mathcal{V}_\delta =
\left\{
i \,\middle|\,
\|\tilde{\mathbf{l}}_i\|_2
\ge
\delta
\max_k \|\tilde{\mathbf{l}}_k\|_2
\right\},
\]
where $\delta \in (0,1]$ is the \textit{activation threshold parameter} and defines the minimum signal‑to‑uncertainty ratio required for a spatial DOF to remain distinguishable under forward‑model perturbations. Smaller values of $\delta$ permit inclusion of weaker DOFs but increase sensitivity to uncertainty, whereas larger values restrict the feasible set to DOFs whose stimulation sensitivity remains reliably distinguishable above the noise floor (Fig.~\ref{fig:eps_effect_graphs})}.

\subsubsection{\textcolor{black}{Dynamic Range Limitation}}
\textcolor{black}{
The restriction of the feasible set to DOFs with sufficiently large magnitude introduces an inherent limitation on the achievable stimulation contrast. Because the optimization operates only within $\mathcal{V}_\delta$, the relative stimulation strength between the target location and surrounding regions is governed by the structure of the reduced lead field itself. In particular, if the target corresponds to index $\ell$, the achievable dynamic range satisfies
\[
\varepsilon
\ge
\delta
\|\tilde{\mathbf{l}}_\ell\|_2^{-1}
\max_k \|\tilde{\mathbf{l}}_k\|_2.
\]
This bound follows directly from the feasibility condition defining the set $\mathcal{V}_\delta$, which restricts the optimization to spatial DOFs whose reduced lead‑field magnitude exceeds the detectability threshold. Consequently, locations with weaker lead‑field magnitude cannot be independently controlled below this detectability level. This expresses a fundamental property of the forward model: the attainable spatial selectivity is governed by the relative scaling of the reduced lead‑field rows rather than by the optimization procedure alone. As a result, the achievable stimulation contrast is determined by the sensitivity distribution encoded in the lead‑field matrix (Fig.~\ref{fig:eps_effect_VTA}).}

\begin{table}[ht]
\centering
\scriptsize
\textcolor{black}{
\caption{\textcolor{black}{\textbf{Magnitude of modeling variability in DBS activation prediction.}
Compilation of quantitative differences reported in DBS modeling studies affecting electric-field (EF) distributions or neural activation predictions across modeling variants, electrode designs, tissue properties, electrode localization, and clinical validation comparisons. Each row corresponds to a controlled perturbation investigated in the cited study. Reported relative differences are converted to decibels using $E_{\mathrm{dB}} = 20\log_{10}(\mathrm{err})$, where $\mathrm{err}$ is the fractional relative difference. Across studies, variability in activation metrics, electrode geometry, conductivity modeling, and localization uncertainty consistently falls within approximately $-30$ to $-10$~dB, with more extreme cases extending toward $-6$~dB. Interpreting these values statistically, this range corresponds to relative deviations of roughly $3\%$ to $30\%$, which defines a practical regime of modeling uncertainty observed in DBS applications. Accordingly, the selected range of $\varepsilon = [-10, -40]$~dB in the optimization framework spans moderate-to-conservative uncertainty levels, covering both typical experimental variability and lower-noise scenarios while avoiding unrealistically high perturbations.}}
\label{tab:uncertainty}
\begin{tabular}{llrc}
\toprule
\textbf{Study} & \textbf{Modeling factor} & \textbf{rel. diff (\%)} & \textbf{dB} \\
\hline
\multicolumn{4}{l}{\textit{Neural activation modeling}} \\
\hline
\cite{FrankemolleGilbert2021} &
individual axon thresholds &
24--47 &
-12.40 to -6.56 \\
— &
population activation metrics &
7--8 &
-23.10 to -21.94 \\
— &
activation volume comparison &
9--12 &
-20.92 to -18.42 \\
\hline
\multicolumn{4}{l}{\textit{Electrode design / geometry}} \\
\hline
\cite{butson_2007_patient} &
electrode geometry influence on VTA &
15--20 &
-16.48 to -13.98 \\
— &
stimulation configuration sensitivity &
10 &
-20.00 \\
\hline
\multicolumn{4}{l}{\textit{Tissue conductivity modeling}} \\
\hline
\cite{zhang2020steering} &
tissue conductivity &
3.8--5.6 &
-28.40 to -25.04 \\
\hline
\multicolumn{4}{l}{\textit{Model vs clinical measurement}} \\
\hline
\cite{Alonso2018} &
Dice overlap (0.92) &
8 &
-21.94 \\
— &
Dice overlap (0.85) &
15 &
-16.48 \\
\hline
\multicolumn{4}{l}{\textit{Electrode localization uncertainty}} \\
\hline
\cite{Lasica2025LocalizationPrecision} &
post-op. electrode localization error &
10--30 &
-20.00 to -10.46 \\
\cite{athawale2019statistical} &
electrode positional uncertainty &
10--40 &
-20.00 to -7.96 \\
\bottomrule
\end{tabular}
}
\end{table}
 
\begin{table*}[!ht]
\centering
\scriptsize
\caption{\textcolor{black}{Summary of the variables used in the forward model, uncertainty formulation, and optimization framework, together with their corresponding roles and parameter ranges. The reciprocity principle (RP) provides a direct bipolar solution obtained from the transpose of the lead-field matrix, without hyperparameter tuning. The L1L1 method performs metaheuristic optimization using an $L_1$-norm regularized $L_1$-norm fitting formulation, where $\alpha$ controls sparsity of the stimulation pattern and $\varepsilon$ defines a threshold on nuisance-field activation. The TLS approach employs an $L_2$-norm regularized least-squares formulation, where $\alpha$ controls the overall regularization strength and $\beta$ weights the suppression of off-target fields. All parameter ranges are expressed in decibel (dB) scale and follow the same logarithmic parametrization used in~\cite{galazprieto_2022_L1vsL2}. For L1L1, the $\alpha$ regularization parameter is explored over $\alpha \in [-100,-30]$~dB, whereas the $\varepsilon$ threshold weight parameter is explored over in two alternative ranges: L1L1(A) explores a wide range $\varepsilon \in [-160,0]$~dB to allow weak nuisance components, while L1L1(B) restricts the range to $\varepsilon \in [-10,0]$~dB, enforcing stronger suppression of low-magnitude off-target fields. For TLS, $\alpha \in [-200,-110]$~dB and $\beta \in [-50,40]$~dB.}}
\label{tab:variables_hyperparameters}
\scriptsize
\textcolor{black}{
\begin{tabular}{llll}
\toprule
\textbf{Var.} & \textbf{Role} & \textbf{Clinical meaning} & \textbf{Range / Notes} \\
\midrule
\multicolumn{4}{c}{\textit{Forward model}} \\
\midrule
$\mathbf{L}$ & Lead-field matrix & Maps electrode currents to the model & -- \\
$\mathbf{L}_1$ & Target rows of $\mathbf{L}$ & Maps currents to target region & -- \\
$\mathbf{L}_2$ & Nuisance rows of $\mathbf{L}$ & Maps currents to off-target regions & -- \\
$\mathbf{y}$ & Current pattern & Number of electrode contacts available & Maximize/Minimize \\
$\mathbf{x}$ & Target field & Discretized coordinate vector & Dipolar DOF \\
$\|\mathbf{x}_1\|_2$ & Reference magnitude & Maximum current dipole moment density & 3.85 A/m$^2$ \cite{Oishi20021390, MURAKAMI201549} \\
$K$ & Electrodes contacts index & Number of electrode contacts available & 8 or 40 contacts \\
\midrule
\multicolumn{4}{c}{\textit{Uncertainty}} \\
\midrule
$\tilde{\mathbf{l}}_i^{(\mathrm{obs})}$ & Perturbed LF row & Observed sensitivity & -- \\
$\mathbf{n}_i \sim \mathcal{N}(0,\sigma^2 \mathbf{I})$ & Additive perturbation & Models forward uncertainty & $\sigma^2$ set from LF variability \\
$\sigma^2$ & Noise variance & Magnitude of uncertainty & Problem-dependent \\
$\tilde{\mathbf{l}}_i$ & Reduced LF row & Differential sensitivity at DOF $i$ & -- \\
$\delta$ & Activation threshold parameter & Minimum signal-to-uncertainty ratio & $(0,1]$ \\
$\mathcal{V}_\delta$ & Feasible DOF set & Feasible coverage range for VTA & Defined by $\delta$ \\
\midrule
\multicolumn{4}{c}{\textit{Electric safety constraints}} \\
\midrule
$\|\mathbf{y}\|_1 \le \mu$ & Total current bound & Limits overall injected current & $\mu = 4.0$ mA \\
$|y_k| \le \gamma$ & Per-contact bound & Limits current per electrode contact & $\gamma = 2.0$ mA \\
$\sum_{k=1}^{K} y_k = 0$ & Charge conservation & Enforces net-zero current injection & -- \\
\midrule
\multicolumn{4}{c}{\textit{Optimization (L1L1)}} \\
\midrule
$\alpha_{(L1L1)}$ & Regularization parameter & Controls sparsity of $\mathbf{y}$ & $[-100,-30]$ dB \\
$\varepsilon$ & Threshold parameter & Controls nuisance-field activation & $[-160,0]$ dB \\
$\zeta=\|\mathbf{L}\|_1$ & Scaling factor & Normalizes optimization terms & Derived from LF \\
$\nu=\|\mathbf{x}_1\|_\infty$ & Target scaling & Defines reference magnitude & Derived from target \\
$\Psi_\varepsilon(\cdot)$ & Thresholding operator & Selects nuisance components above $\varepsilon$ & -- \\
\midrule
\multicolumn{4}{c}{\textit{Optimization (TLS)}} \\
\midrule
$\alpha_{(TLS)}$ & Regularization parameter & Controls solution magnitude & $[-200,-110]$ dB \\
$\beta$ & Nuisance weighting & Controls off-target suppression & $[-50,40]$ dB \\
$\varsigma=\|\mathbf{L}\|_2$ & Scaling factor & Normalizes optimization terms & Derived from LF \\
$\nu=\|\mathbf{x}_1\|_\infty$ & Target scaling & Defines reference magnitude & Derived from target \\
\midrule
\multicolumn{4}{c}{\textit{Decision variables}} \\
\midrule
$\Gamma$ & Focused current density & Target alignment measure (Intensity) & Maximized \\
$\Gamma_0$ & Lower bound for $\Gamma$ & Minimum level of target engagement & $\Gamma \ge 0.8$ mA \\
$\Xi$ & Nuisance current density & Off-target field measure & Minimized \\
$\Theta$ & Field ratio ($\Theta = \Gamma / \Xi$) & Selectivity metric (Focality) & Maximized \\
\bottomrule
\end{tabular}
}
\end{table*}

\subsection{Hardware}
Dell Precision 5820 Workstation equipped with 256~GB of RAM, a 10-core Intel i9–10900X CPU, and an NVidia Quadro RTX 4000 GPU. The GPU were used to accelerate both forward and inverse calculations.

\subsection{Numerical Experiments}
\label{sec:num_exp}
Our numerical experiments evaluate the proposed optimization methods with respect to: (i) 8- and 40-contact DBS lead configurations, (ii) noiseless and noisy forward mappings, (iii) dipole orientation relative to the lead array (parallel or perpendicular) \cite{Janson_2020_Activation}, and (iv) algorithm efficiency. The L1L1 method is tested under two alternative {\em a priori} assumptions for the feasible $\varepsilon$ threshold range: a case with [-160, 0]~dB (refereed to as L1L1(A)) and a case with [-10, 0]~dB (refereed to as L1L1(B)). The first case neglects uncertainty in the field ratio $\Theta$, while the second accounts for variability within a 10~dB dynamic range. TLS uses regularization parameter $\gamma$ and L2-norm weighting parameter $\beta$ intervals as in \cite{galazprieto_2022_L1vsL2} (Table \ref{tab:variables_hyperparameters}). Detailed formulation of both the L1L1 and TLS methods is provided in \ref{App:sec:L1L1} and \ref{App:sec:TLS}, respectively.

Feasible solutions are computed for each target DOF of the LR-LF under PSNR 40.0~dB corruption. While anatomical variations can introduce spatially correlated errors, this additive noise provides a tractable baseline to assess robustness to model perturbations. Although many sources of uncertainty (e.g., electrode localization error or conductivity variation) introduce spatially correlated perturbations in the lead field, the additive Gaussian noise provides a tractable approximation of aggregate forward‑model variability and is commonly used for robustness analysis.

\begin{figure}[t]
    \centering
    \includegraphics[width = 6.0cm, height = 5.5cm]{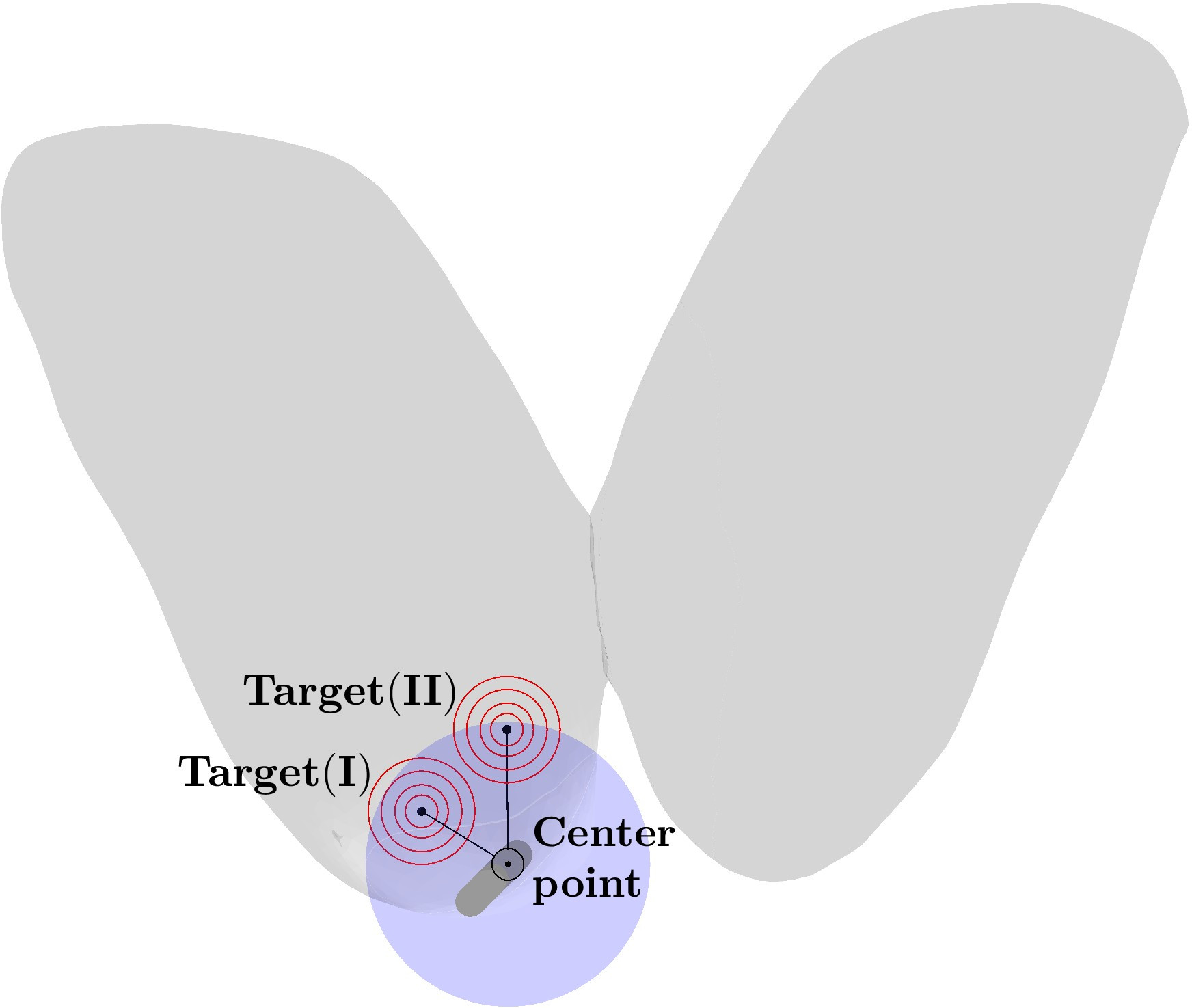}
    \caption{
    \textcolor{black}{Schematic illustration of the two synthetic stimulation targets, Target (I) and Target (II), used for noise sensitivity analysis under different PSNR levels (red circles).} Both targets are positioned within the right thalamus relative to a directional DBS lead. \textcolor{black}{The lead is shown only as a geometric reference; individual electrode contacts and stimulation patterns are not depicted in this figure.} The lead’s center point is indicated, together with a concentric spherical region of interest (ROI) of 6.0 mm radius (blue circle) representing an approximate target volume. \textcolor{black}{These targets define representative degrees of freedom (DOFs) of the lead field used to evaluate robustness of the optimization methods under noisy forward mappings.} Bilateral thalamic structures are shown in gray for anatomical reference, and distances are defined with respect to the lead center.
    }
    \label{fig:target_positions_I_and_II}
\end{figure}

\textcolor{black}{Noise sensitivity is analyzed for two target positions (I) and (II) (Fig.~\ref{fig:target_positions_I_and_II}) along the HR-LF using PSNR levels of 40.0 and 50.0~dB. The analysis considers the focused current $\Gamma$, the field ratio $\Theta$, and the targeting error, defined as the Euclidean distance between the prescribed target location and the location of the maximum current density within a region of interest (ROI) surrounding the target. The ROI is used to exclude the volume current flowing between the stimulation electrodes, thereby avoiding bias caused by the distributed current injection that is intrinsic to tES. Interpreted in terms of $\mathcal{V}\delta$, PSNR 40.0~dB (50.0~dB) corresponds to the assumption that each target in $\mathcal{V}{-30\ \mathrm{dB}}$ ($\mathcal{V}{-40\ \mathrm{dB}}$) can be stimulated with at least a 10~dB dynamic range $\varepsilon^{-1}$ for $\Theta$, which aligns with the $\varepsilon$ range of $[-10,0]$ used for L1L1(B). In this interpretation, the PSNR level determines the effective detectability threshold of the reduced lead-field rows and therefore controls the size of the feasible DOF set $\mathcal{V}\delta$ under uncertainty.}

% \begin{figure}
%     \centering
%     \begin{tikzpicture}[scale=0.5]
    
%     % Draw the thalamus (ellipsoid shape)
%     \shade[ball color=gray!40] (0,0) ellipse (2.5 and 3);
    
%     % Labels for thalamus
%     \node at (0,3.3) {\textbf{Thalamus}};
    
%     % Draw DBS probe
%     \draw[line width=1mm, black] (0,4) -- (0,-4);
    
%     % Electrode contacts on the probe
%     \fill[black] (0,2) circle (0.15);
%     \fill[black] (0,-1) circle (0.15);
    
%     % Target locations
%     \fill[red] (1,1) circle (0.2);
%     \fill[red] (-1,-2) circle (0.2);
    
%     % Label probe
%     \node[anchor=east] at (0, 4.2) {\textbf{DBS Probe}};
    
%     % Label stimulation targets
%     \node[anchor=west] at (1.2,1) {\textbf{Target (I)}};
%     \node[anchor=east] at (-1.2,-2) {\textbf{Target (II)}};
    
%     \end{tikzpicture}
    
%     \caption{\textcolor{red}{Replace with the actual illustration. This is just a placeholder.} Target positions (I) and (II) applied in analyzing the noise sensitivity.}
%     \label{fig:target_positions_I_and_II}
% \end{figure}

\begin{figure*}[ht!]
    \centering
    \begin{small}
        
    \setlength{\tabcolsep}{0.2cm} % Horizontal spacing between images
    \renewcommand{\arraystretch}{1.0} % Vertical spacing between rows
    \setlength{\arrayrulewidth}{0.1mm}

    \begin{tabular}{|>{\centering}m{0.9cm}|ccc|ccc|}
        % Column headers with small figures, all framed
                 \multicolumn{1}{c}{} & 
        \multicolumn{3}{c}{ \textbf{LR-LF}} & 
        \multicolumn{3}{c}{ \textbf{HR-LF}} \\
        \multicolumn{1}{c}{} & 
        \multicolumn{1}{c}{\textbf{Focused $\Gamma$ (A/m\textsuperscript{2})}} & 
        \multicolumn{1}{c}{\textbf{Nuisance $\Xi$ (A/m\textsuperscript{2})}} & 
        \multicolumn{1}{c}{\textbf{Ratio $\Theta$ (rel.)}} & 
        \multicolumn{1}{c}{\textbf{Focused $\Gamma$ (A/m\textsuperscript{2})}} & 
        \multicolumn{1}{c}{\textbf{Nuisance $\Xi$ (A/m\textsuperscript{2})}} & 
        \multicolumn{1}{c}{\textbf{Ratio $\Theta$ (rel.)}} \\ 
        
        \multicolumn{1}{c}{} & 
        \multicolumn{1}{c}{\includegraphics[width=2.25cm, height=0.45cm]{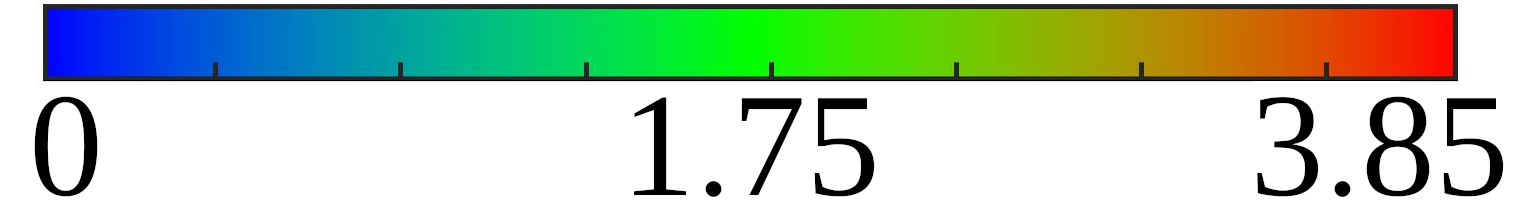}} & 
        \multicolumn{1}{c}{\includegraphics[width=2.25cm, height=0.45cm]{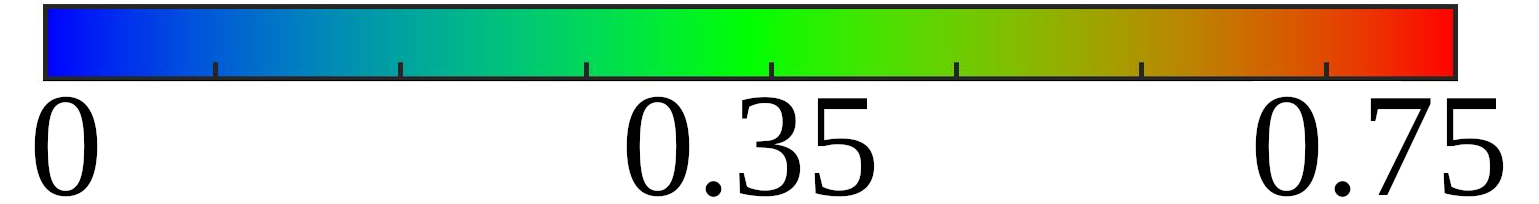}} & 
        \multicolumn{1}{c}{\includegraphics[width=2.25cm, height=0.45cm]{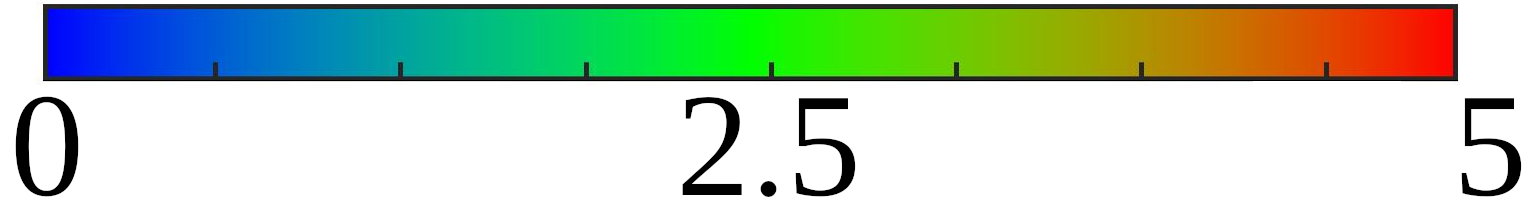}} & 
        \multicolumn{1}{c}{\includegraphics[width=2.25cm, height=0.45cm]{Figure_AB_col/colorbar_focused.png}} & 
        \multicolumn{1}{c}{\includegraphics[width=2.25cm, height=0.45cm]{Figure_AB_col/colorbar_nuisance.png}} & 
        \multicolumn{1}{c}{\includegraphics[width=2.25cm, height=0.45cm]{Figure_AB_col/colorbar_ratio.png}} \\
        \hline
        
        \cellcolor[HTML]{fff0ff} & 
        \rule{0pt}{1.80cm}
        \cellcolor[HTML]{fff0ff} \includegraphics[width=2.25cm, height=1.5cm]{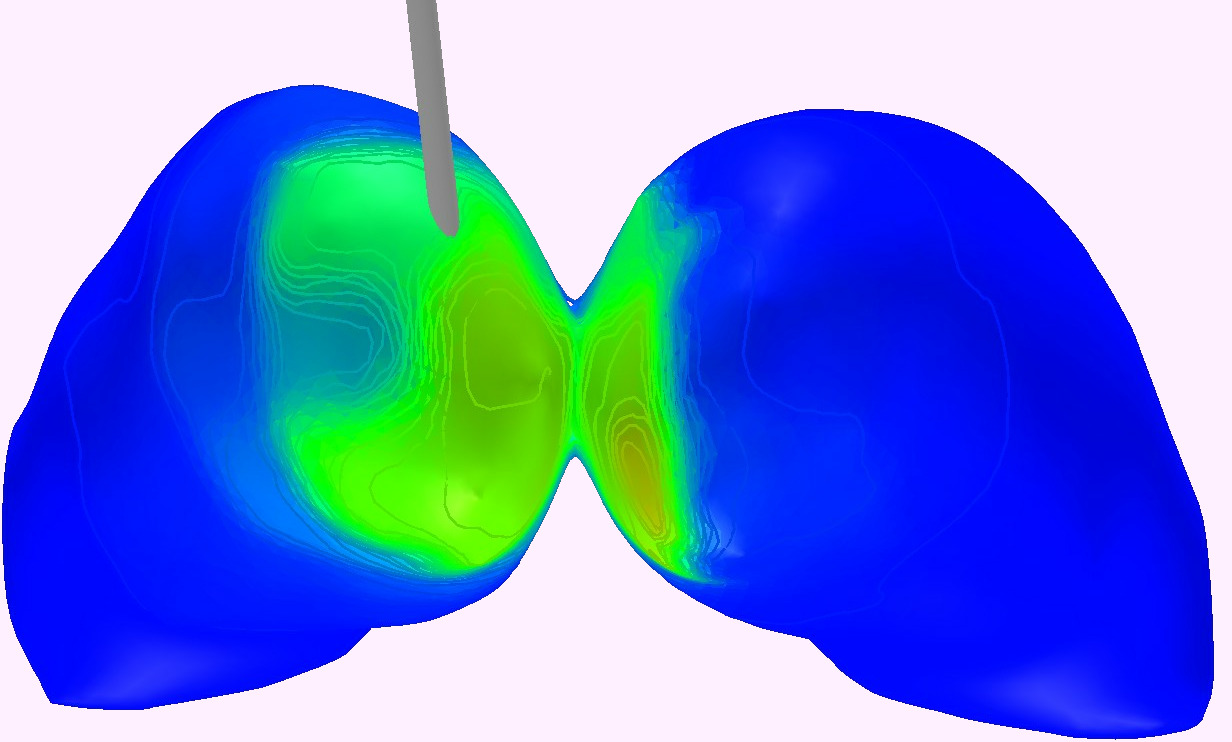} &  
        \cellcolor[HTML]{fff0ff} \includegraphics[width=2.25cm, height=1.5cm]{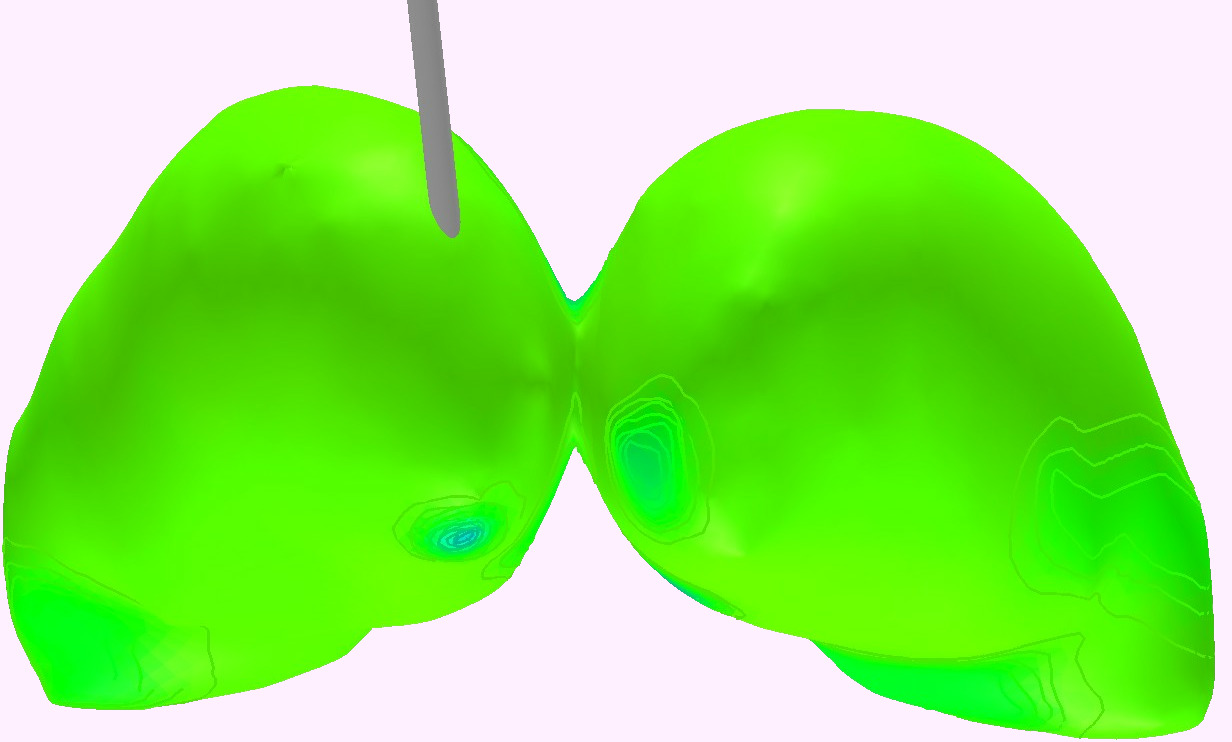} &
        \cellcolor[HTML]{fff0ff} \includegraphics[width=2.25cm, height=1.5cm]{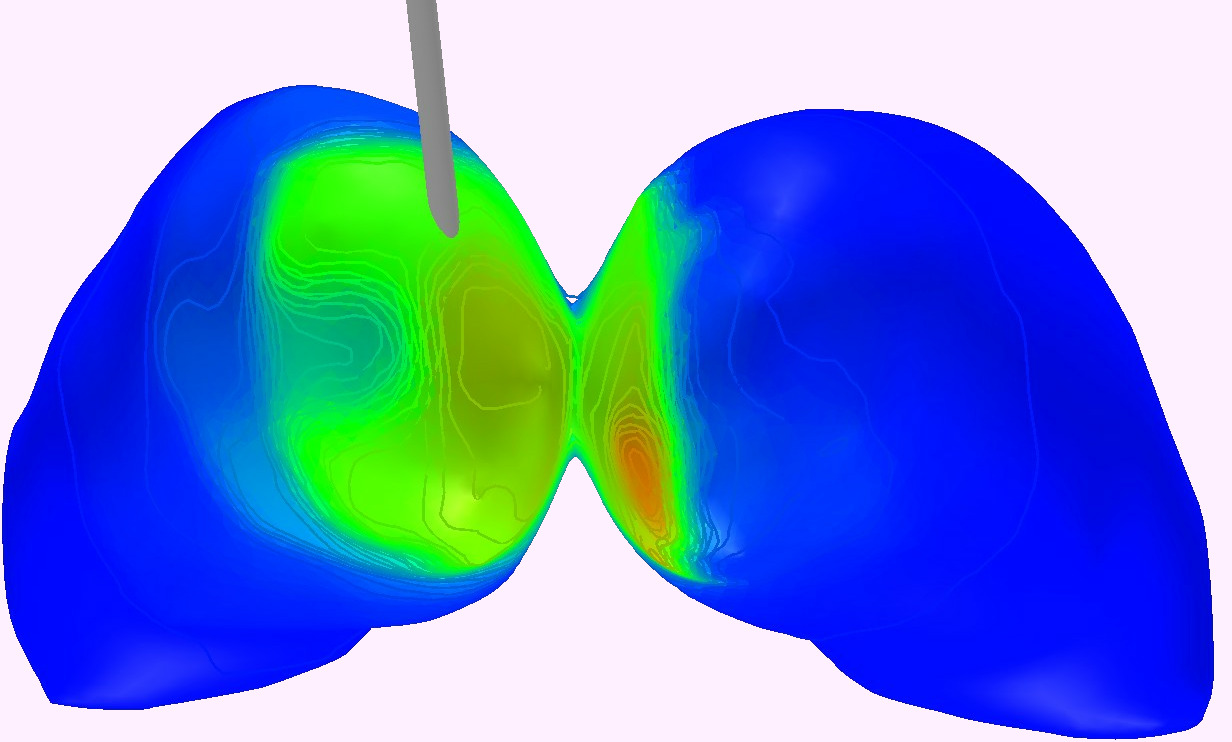} & 
        \cellcolor[HTML]{fff0ff} \includegraphics[width=2.25cm, height=1.5cm]{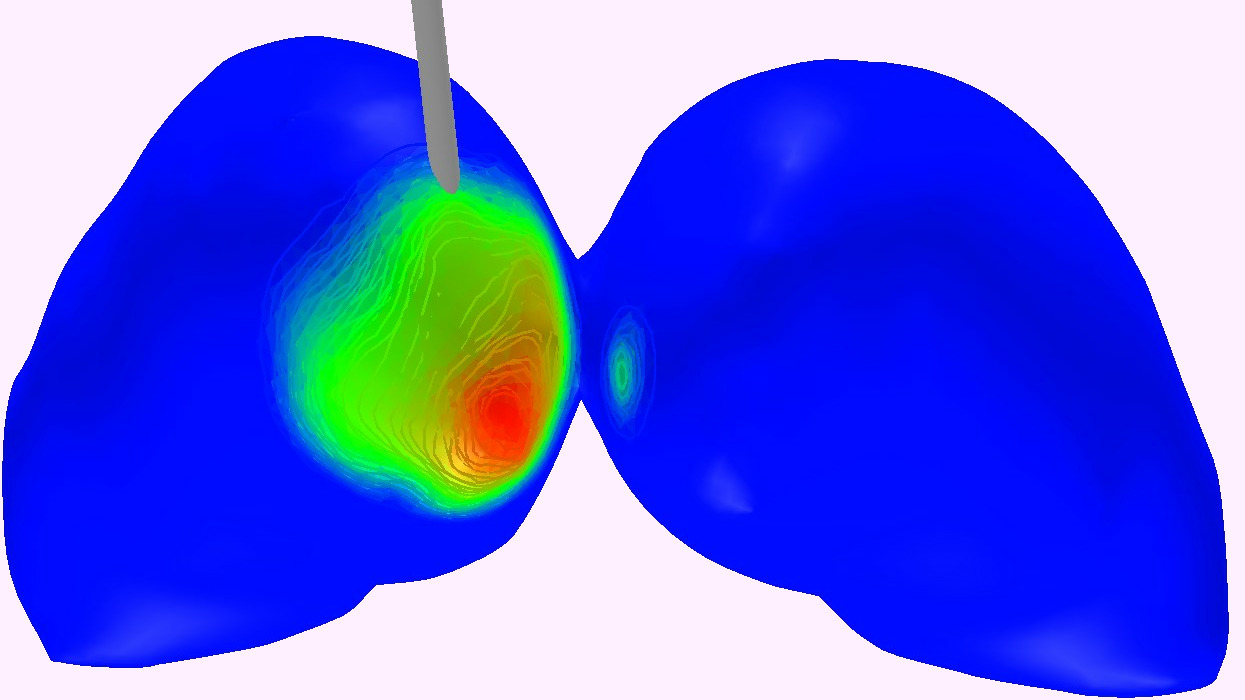} & 
        \cellcolor[HTML]{fff0ff} \includegraphics[width=2.25cm, height=1.5cm]{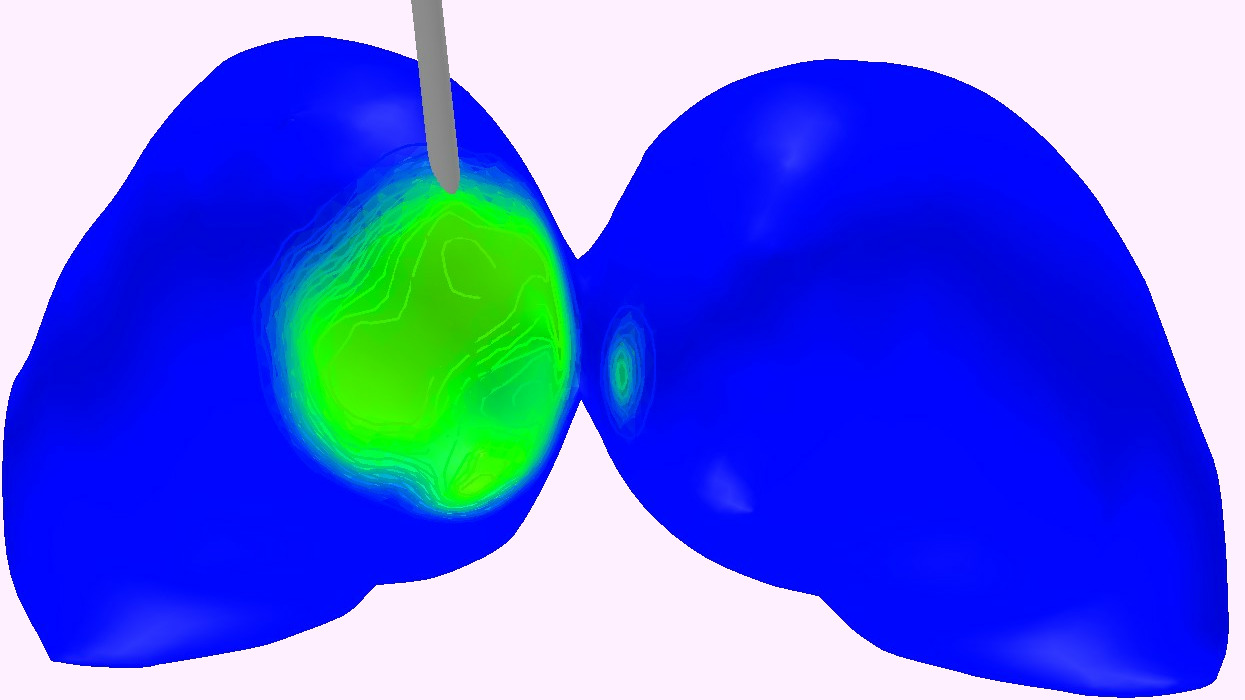} &
        \cellcolor[HTML]{fff0ff} \includegraphics[width=2.25cm, height=1.5cm]{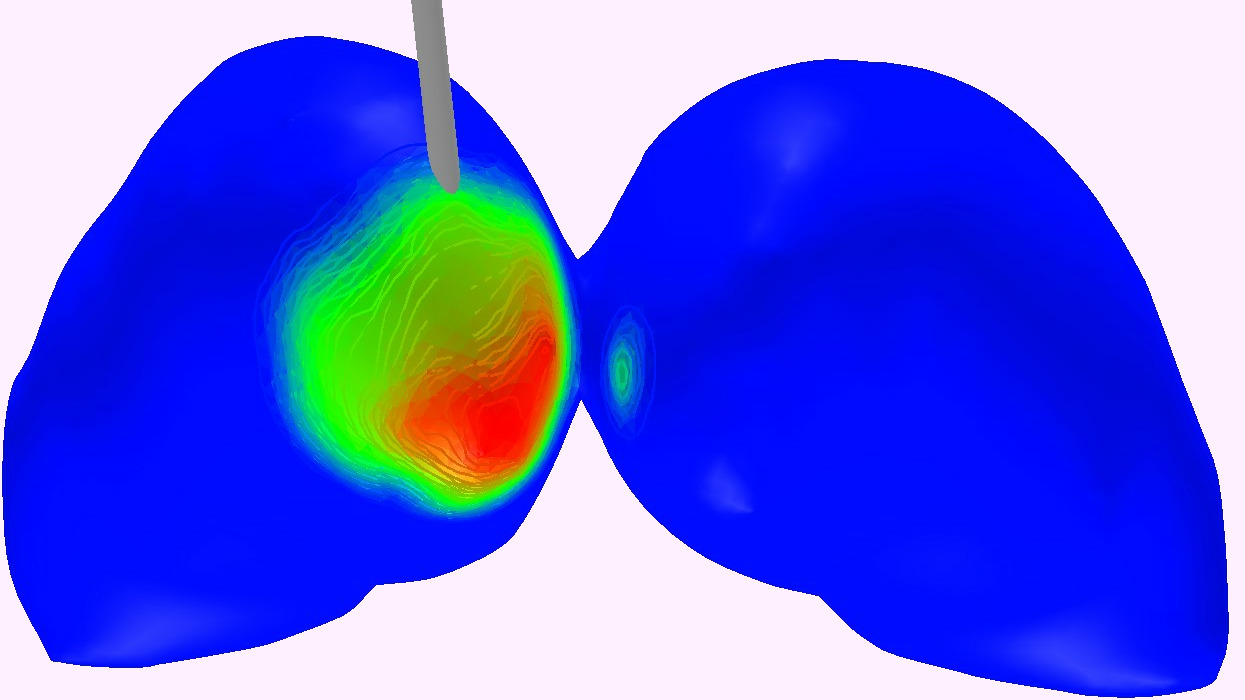} \\
        \cellcolor[HTML]{fff0ff} \raisebox{1.7cm}[0pt][0pt]{\multirow{2}{*}{\rotatebox{90}{ \textbf{RP}}}} \hskip0.1cm \raisebox{2.9cm}[0pt][0pt] {\multirow{2}{*}{\rotatebox{90}{Perpendicular \hskip0.25cm Parallel}}}&
        \rule{0pt}{1.80cm}
        \cellcolor[HTML]{fff0ff} \includegraphics[width=2.25cm, height=1.5cm]{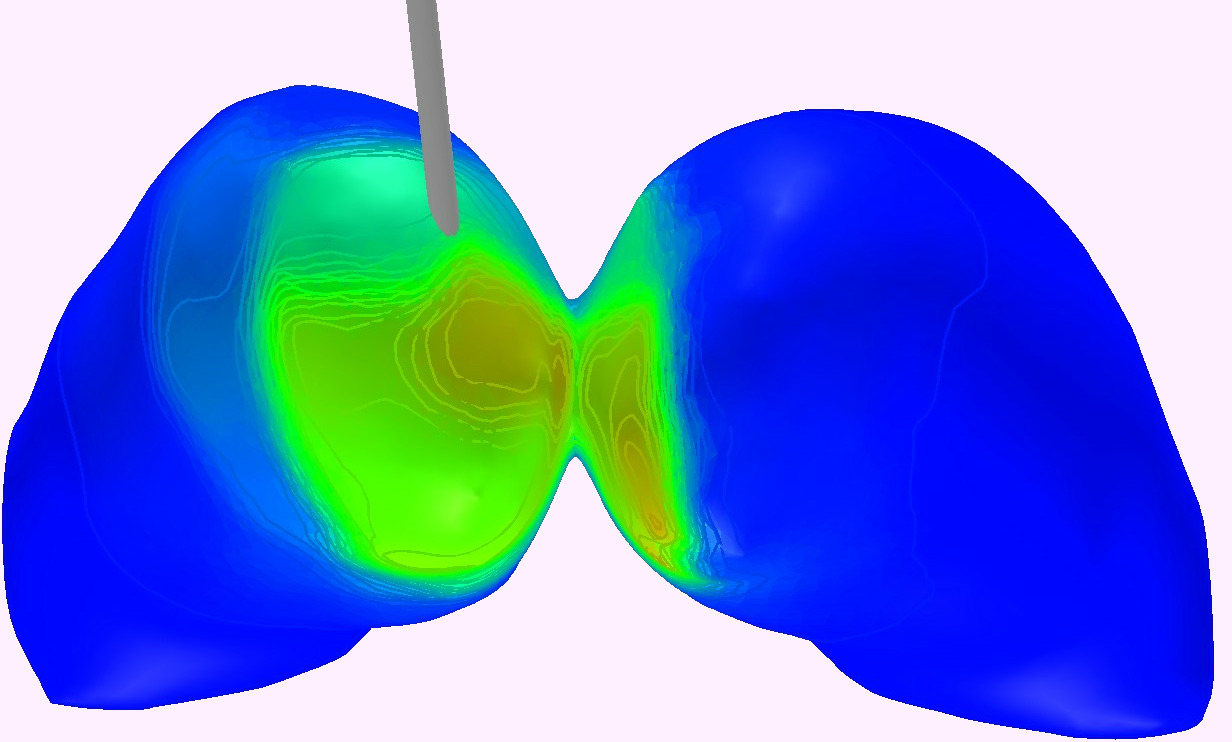} &  
        \cellcolor[HTML]{fff0ff} \includegraphics[width=2.25cm, height=1.5cm]{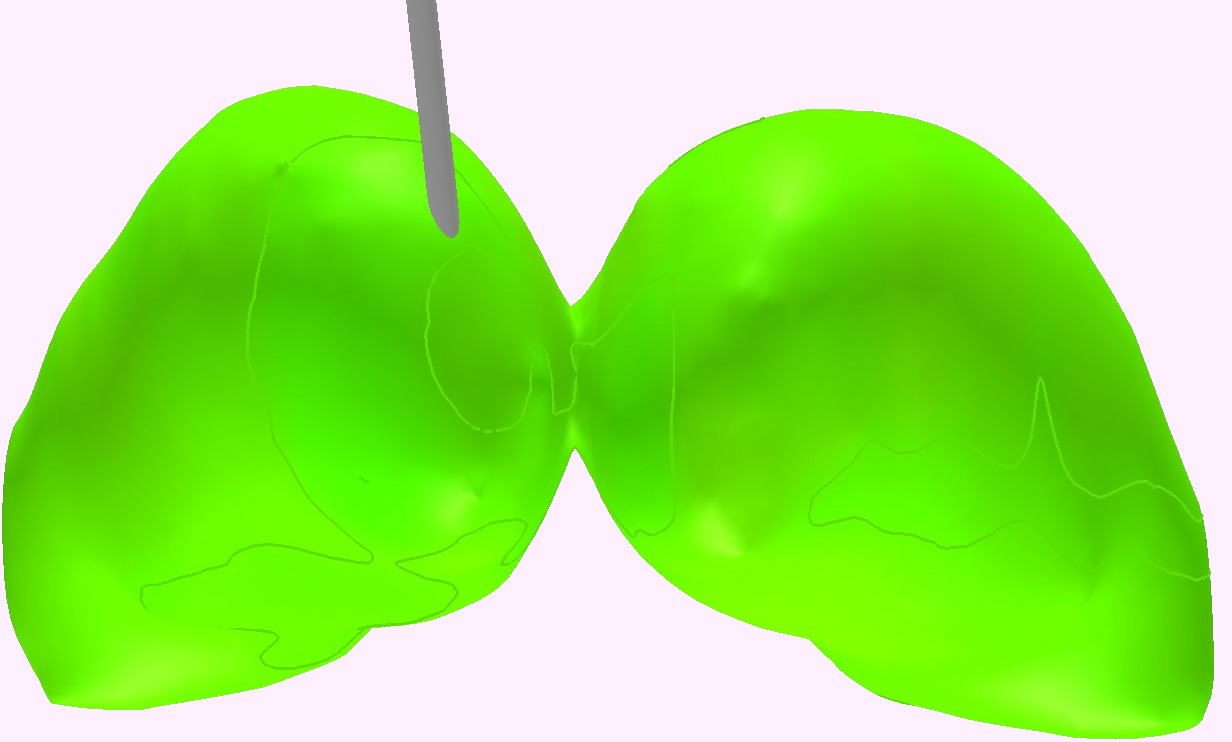} &
        \cellcolor[HTML]{fff0ff} \includegraphics[width=2.25cm, height=1.5cm]{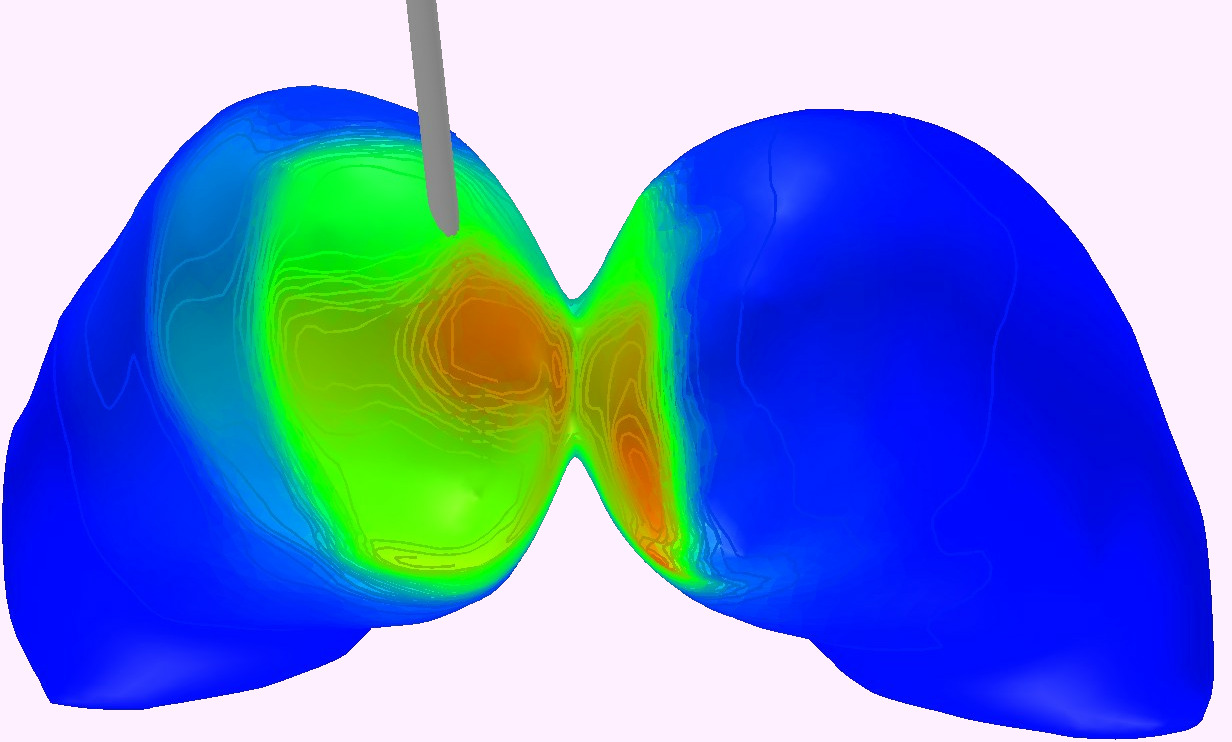} & 
        \cellcolor[HTML]{fff0ff} \includegraphics[width=2.25cm, height=1.5cm]{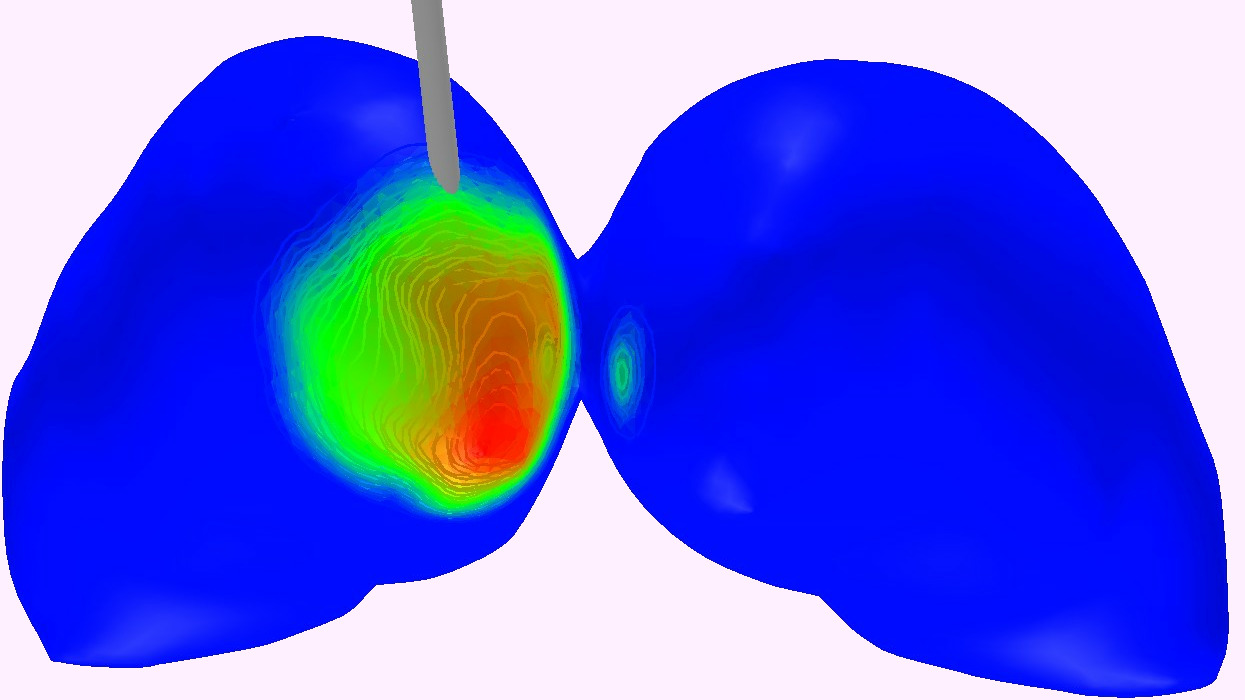} & 
        \cellcolor[HTML]{fff0ff} \includegraphics[width=2.25cm, height=1.5cm]{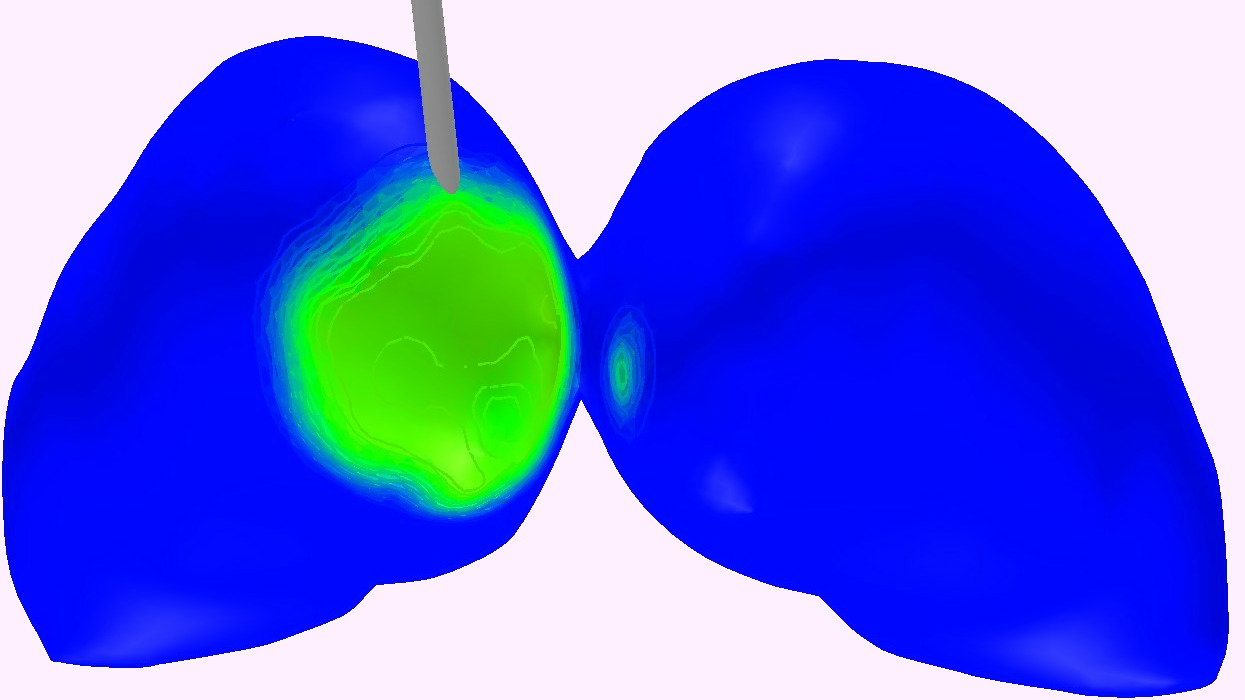} &
        \cellcolor[HTML]{fff0ff} \includegraphics[width=2.25cm, height=1.5cm]{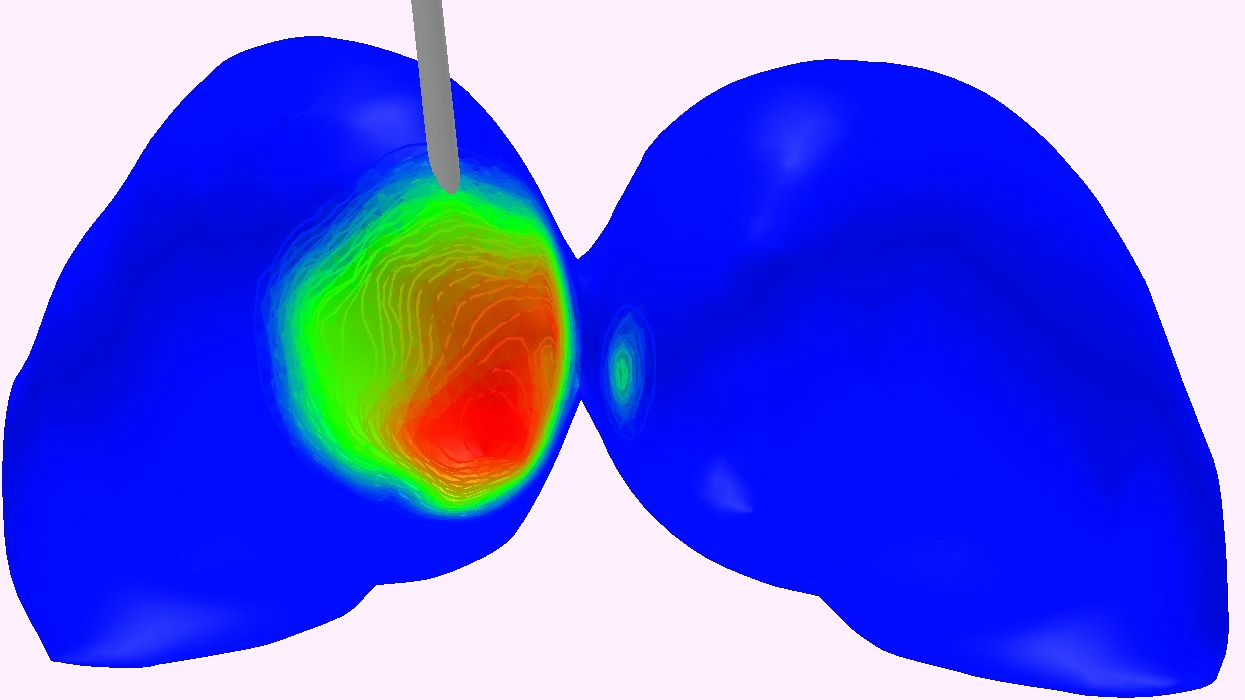} \\
        \hline
        
        \cellcolor[HTML]{ccffff} & 
        \rule{0pt}{1.80cm}
        \cellcolor[HTML]{ccffff} \includegraphics[width=2.25cm, height=1.5cm]{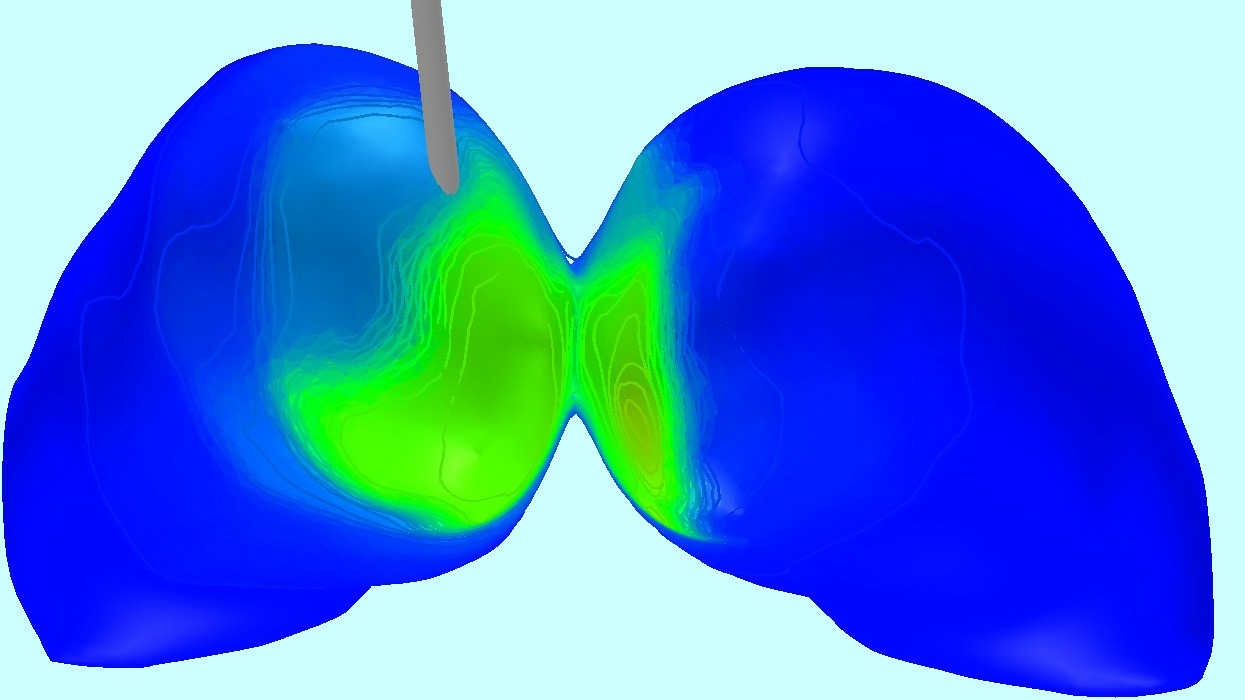} &  
        \cellcolor[HTML]{ccffff} \includegraphics[width=2.25cm, height=1.5cm]{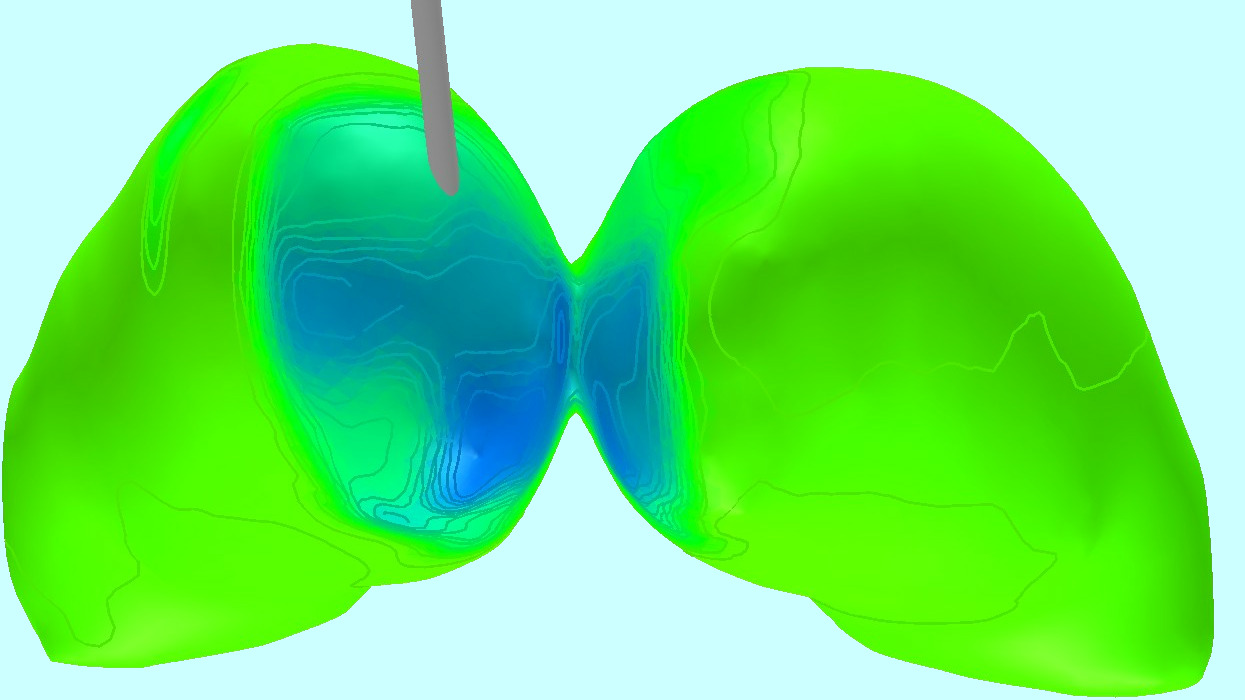} &
        \cellcolor[HTML]{ccffff} \includegraphics[width=2.25cm, height=1.5cm]{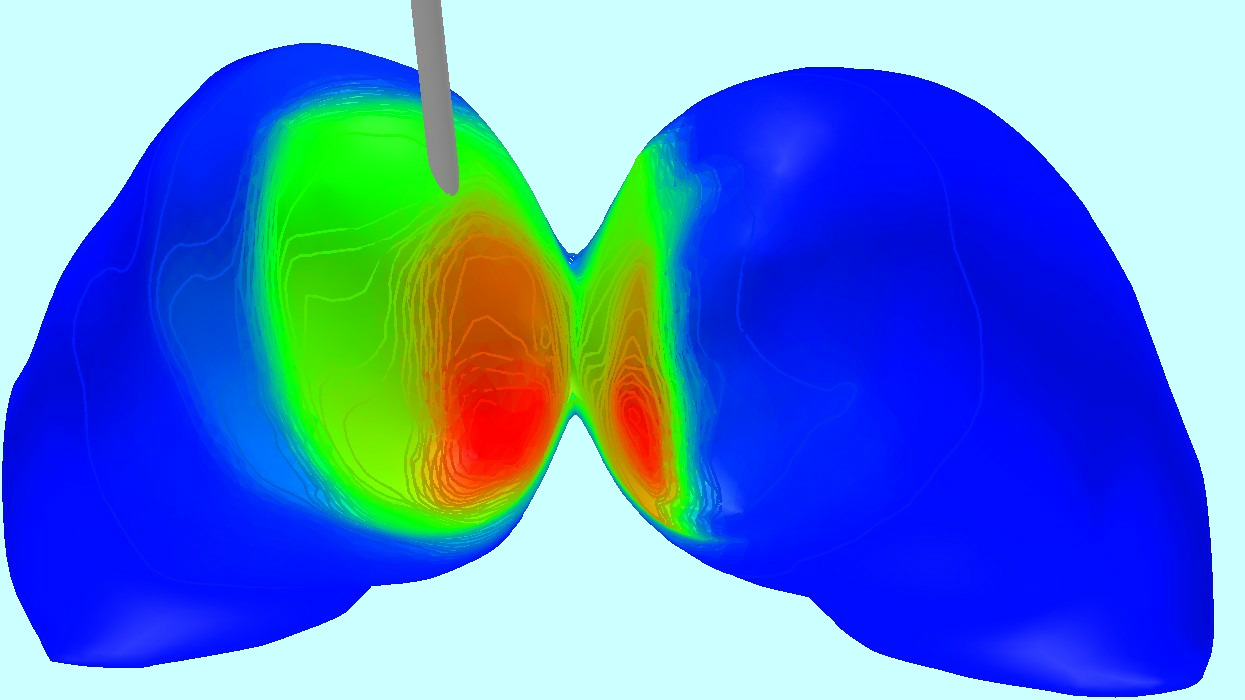} & 
        \cellcolor[HTML]{ccffff} \includegraphics[width=2.25cm, height=1.5cm]{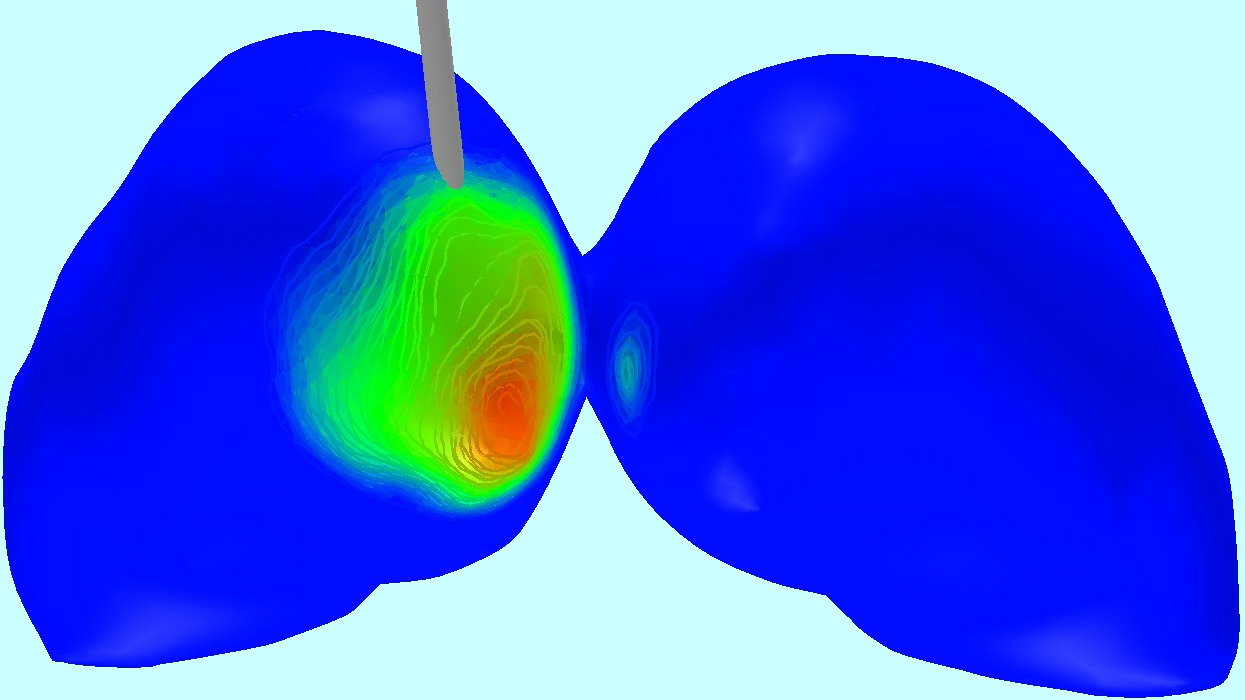} & 
        \cellcolor[HTML]{ccffff} \includegraphics[width=2.25cm, height=1.5cm]{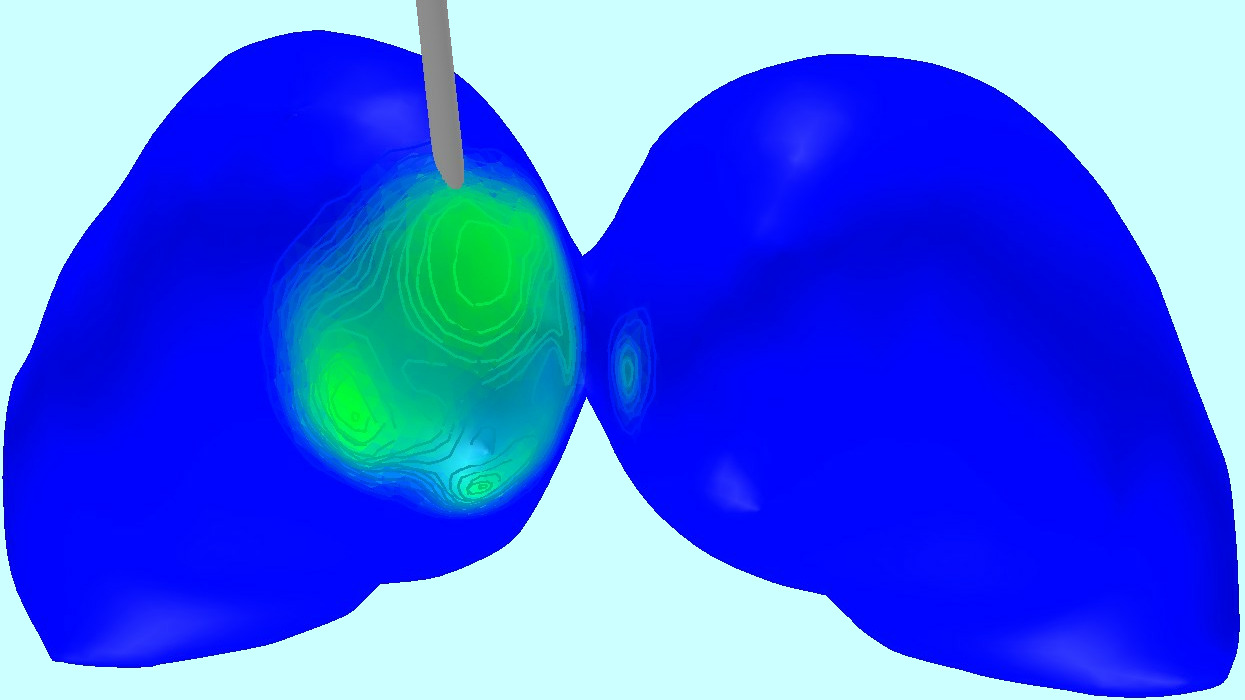} &
        \cellcolor[HTML]{ccffff} \includegraphics[width=2.25cm, height=1.5cm]{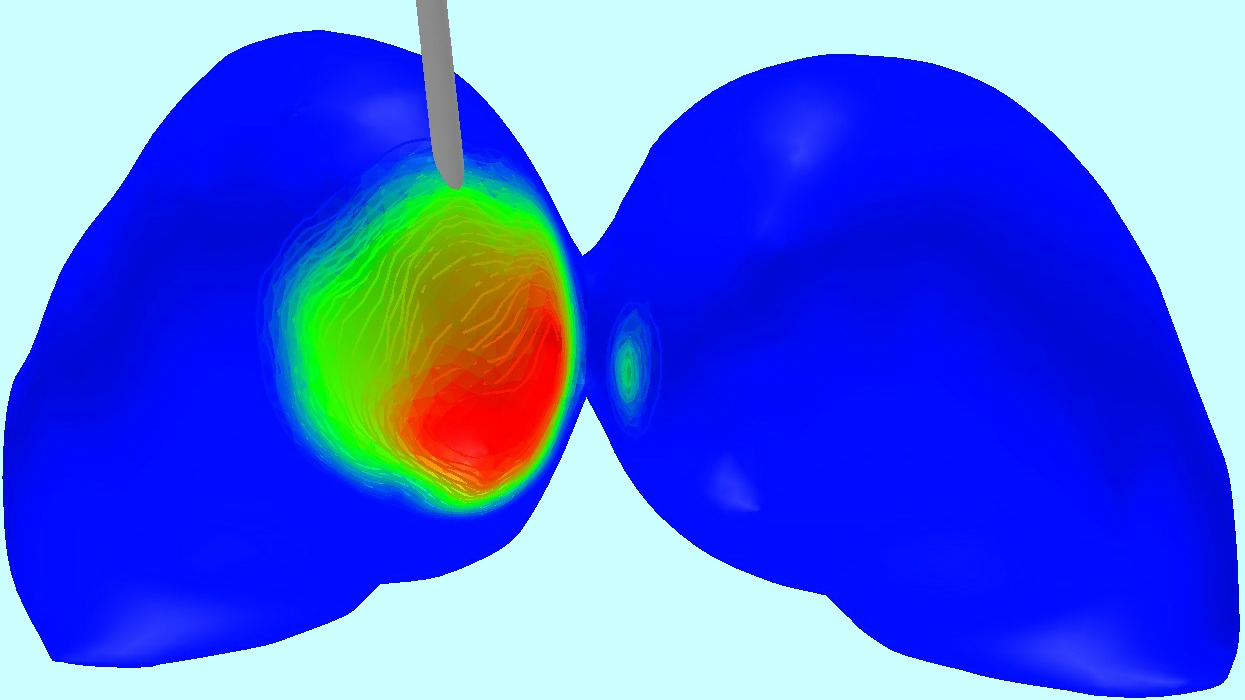} \\
        \cellcolor[HTML]{ccffff} \raisebox{3.15cm}[0pt][0pt]{\multirow{2}{*}{\rotatebox{90}{\textbf{L1L1(B)}, $\varepsilon \in  [ -10, 0] $ dB}}} \hskip0.1cm \raisebox{2.9cm}[0pt][0pt] {\multirow{2}{*}{\rotatebox{90}{Perpendicular \hskip0.25cm Parallel}}}& 
        \rule{0pt}{1.80cm}
        \cellcolor[HTML]{ccffff} \includegraphics[width=2.25cm, height=1.5cm]{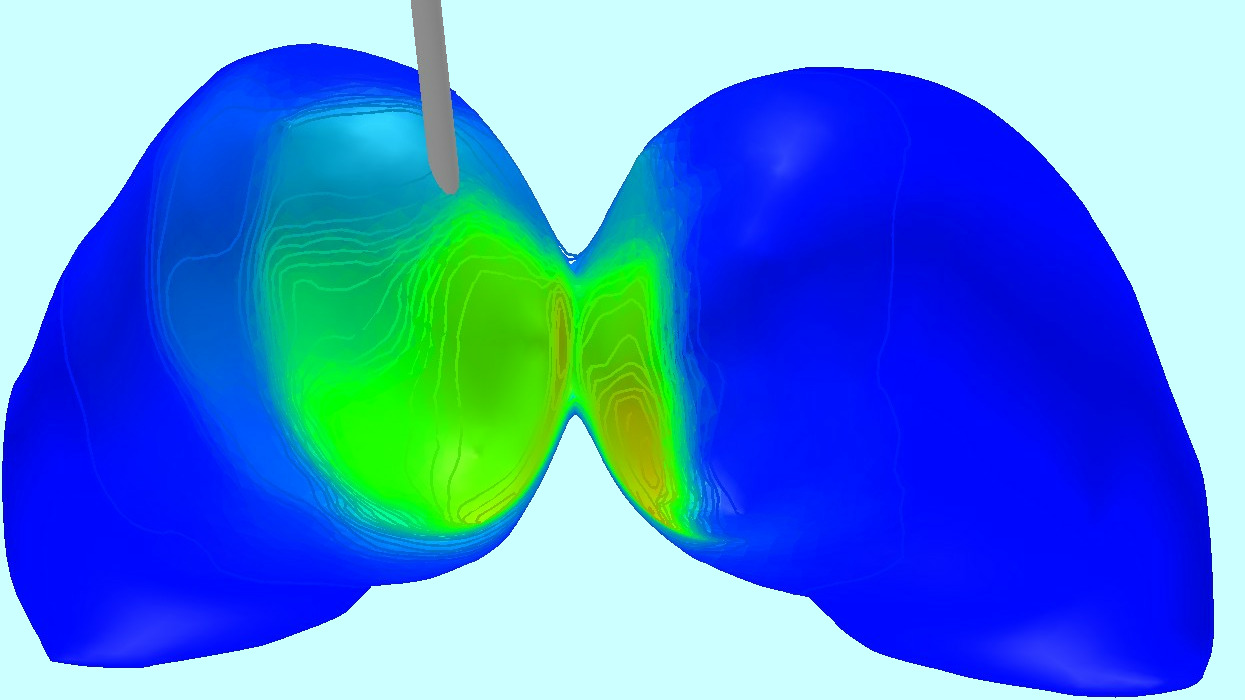} &  
        \cellcolor[HTML]{ccffff} \includegraphics[width=2.25cm, height=1.5cm]{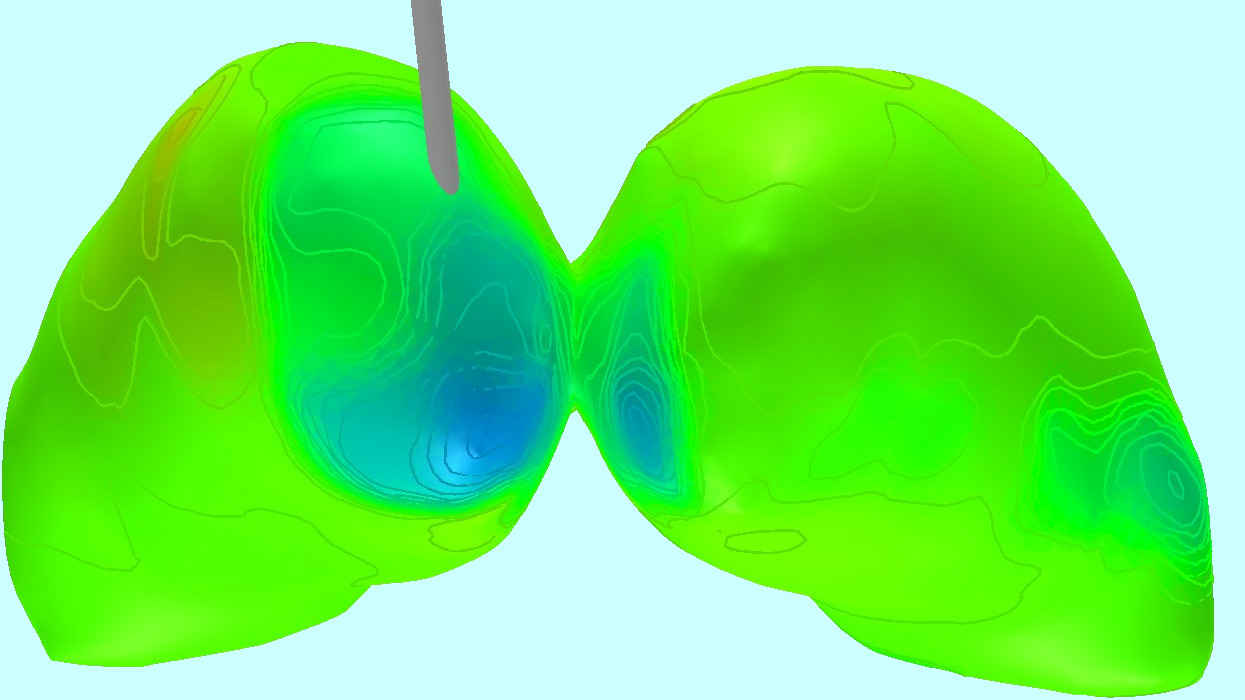} &
        \cellcolor[HTML]{ccffff} \includegraphics[width=2.25cm, height=1.5cm]{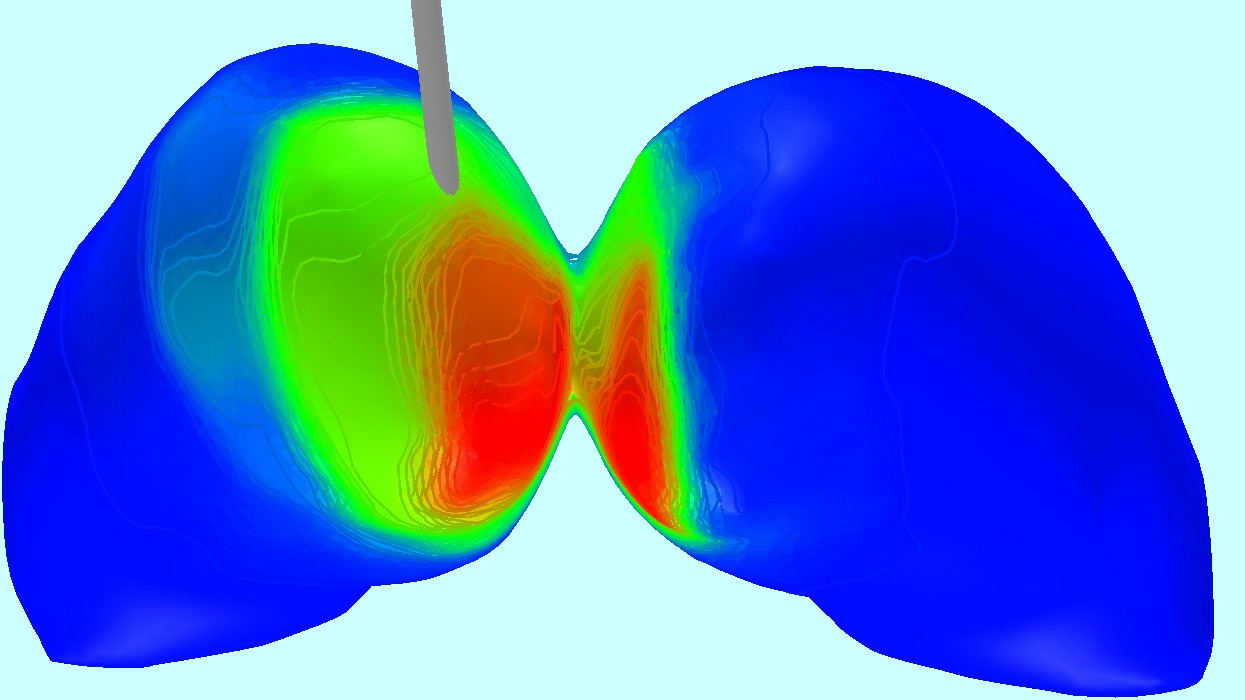} & 
        \cellcolor[HTML]{ccffff} \includegraphics[width=2.25cm, height=1.5cm]{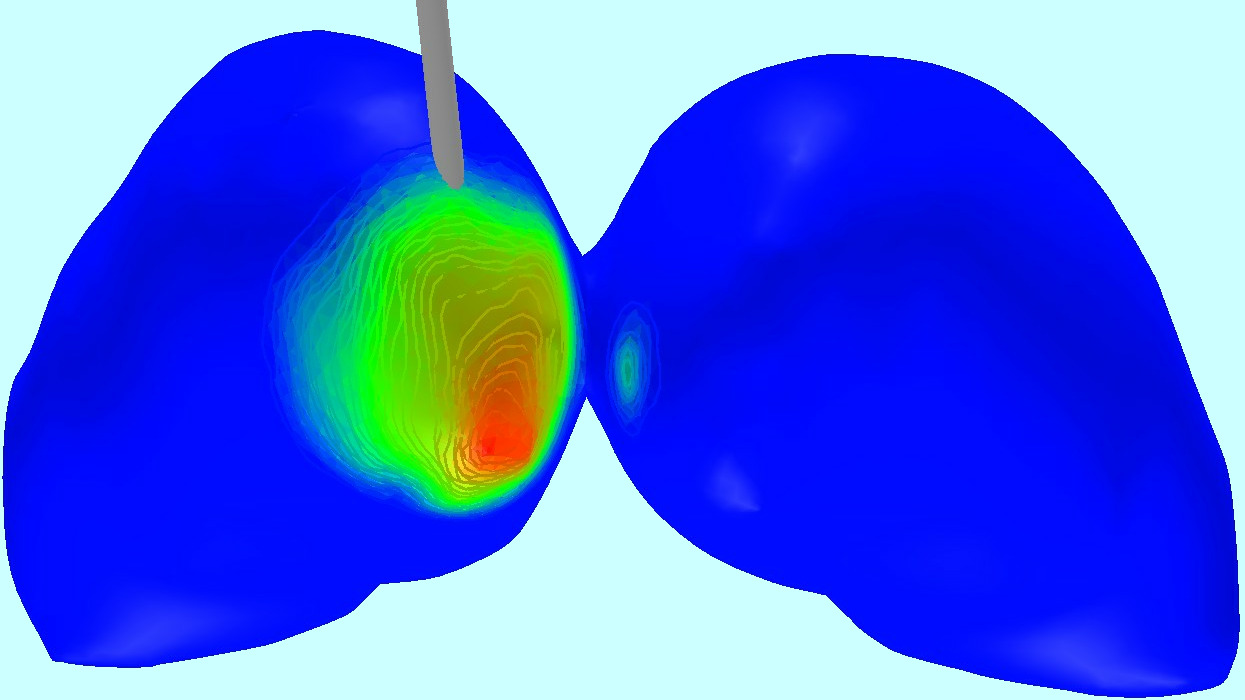} & 
        \cellcolor[HTML]{ccffff} \includegraphics[width=2.25cm, height=1.5cm]{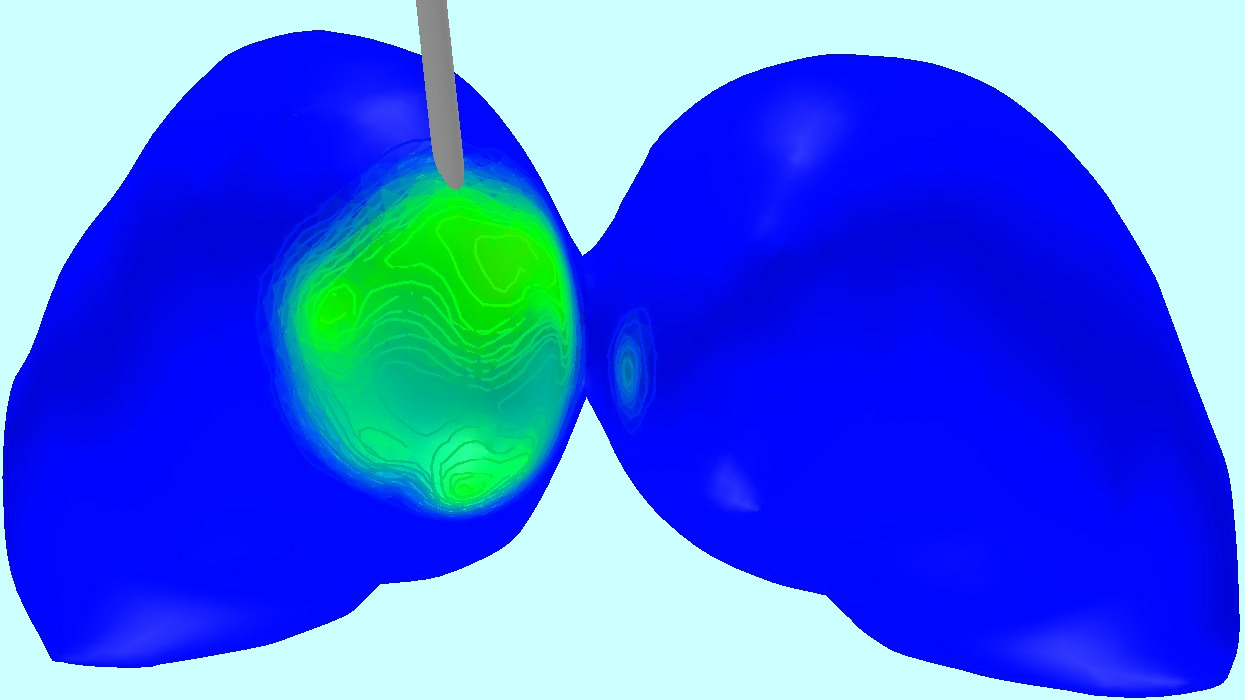} &
        \cellcolor[HTML]{ccffff} \includegraphics[width=2.25cm, height=1.5cm]{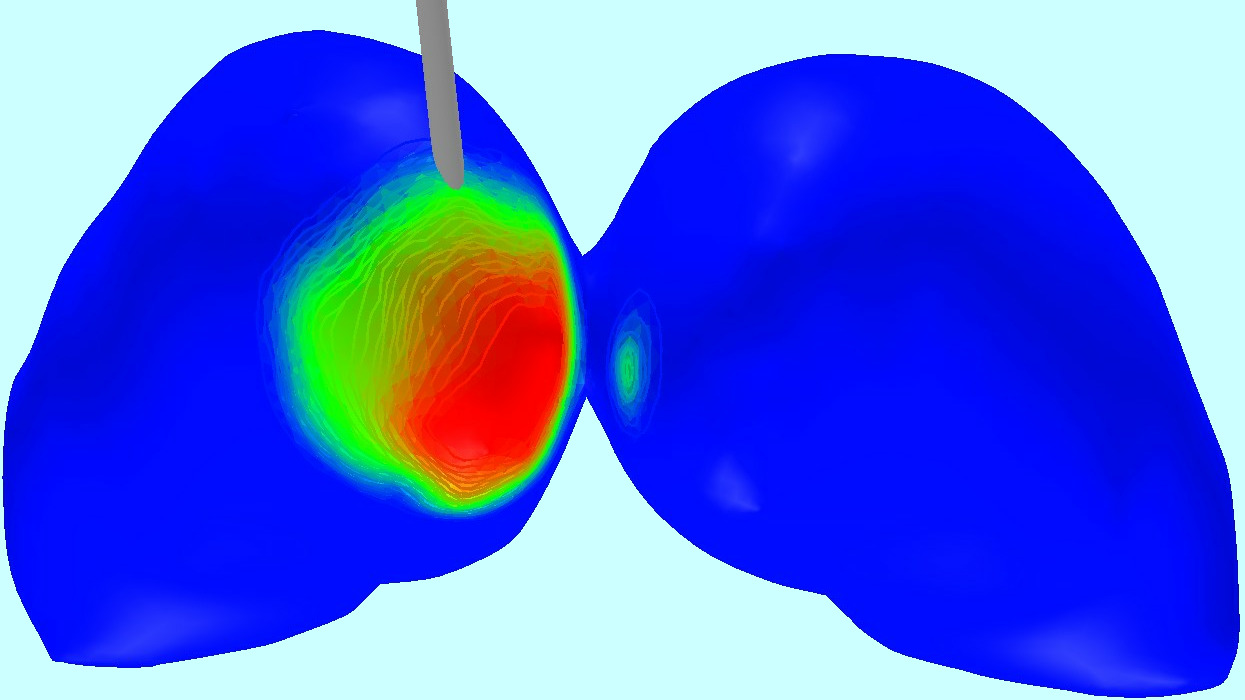} \\
        \hline
        
        \cellcolor[HTML]{ccffff} & 
        \rule{0pt}{1.80cm}
        \cellcolor[HTML]{ccffff} \includegraphics[width=2.25cm, height=1.5cm]{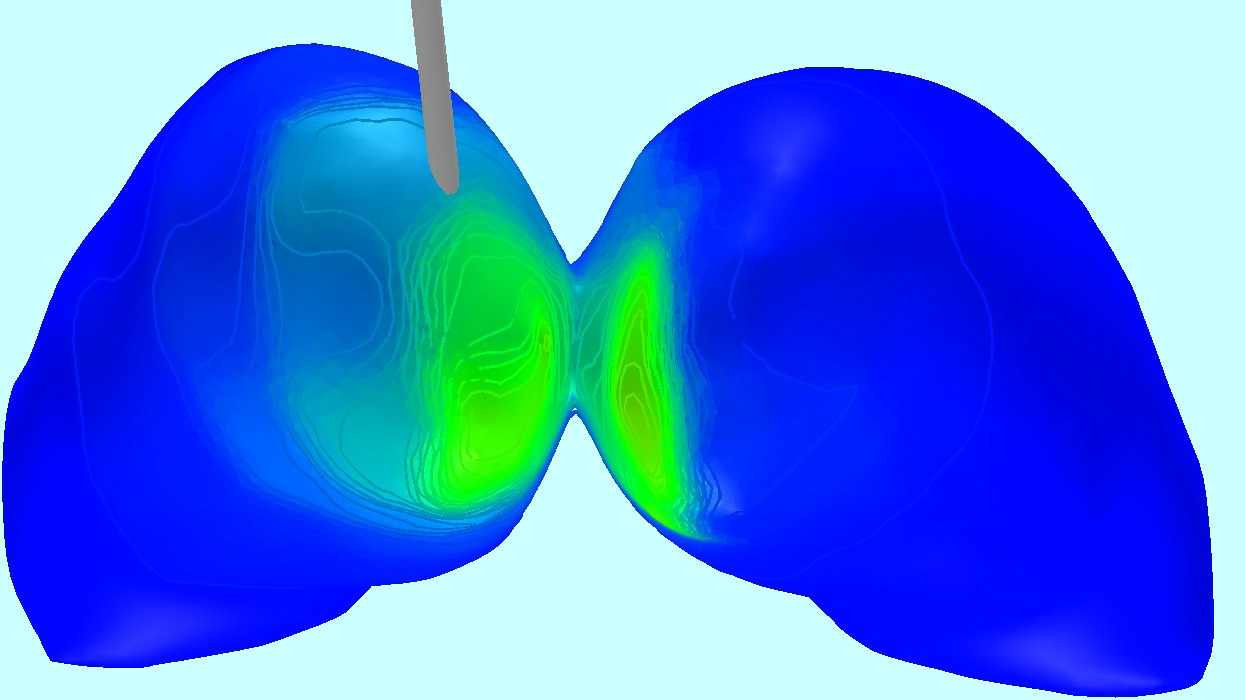} &  
        \cellcolor[HTML]{ccffff} \includegraphics[width=2.25cm, height=1.5cm]{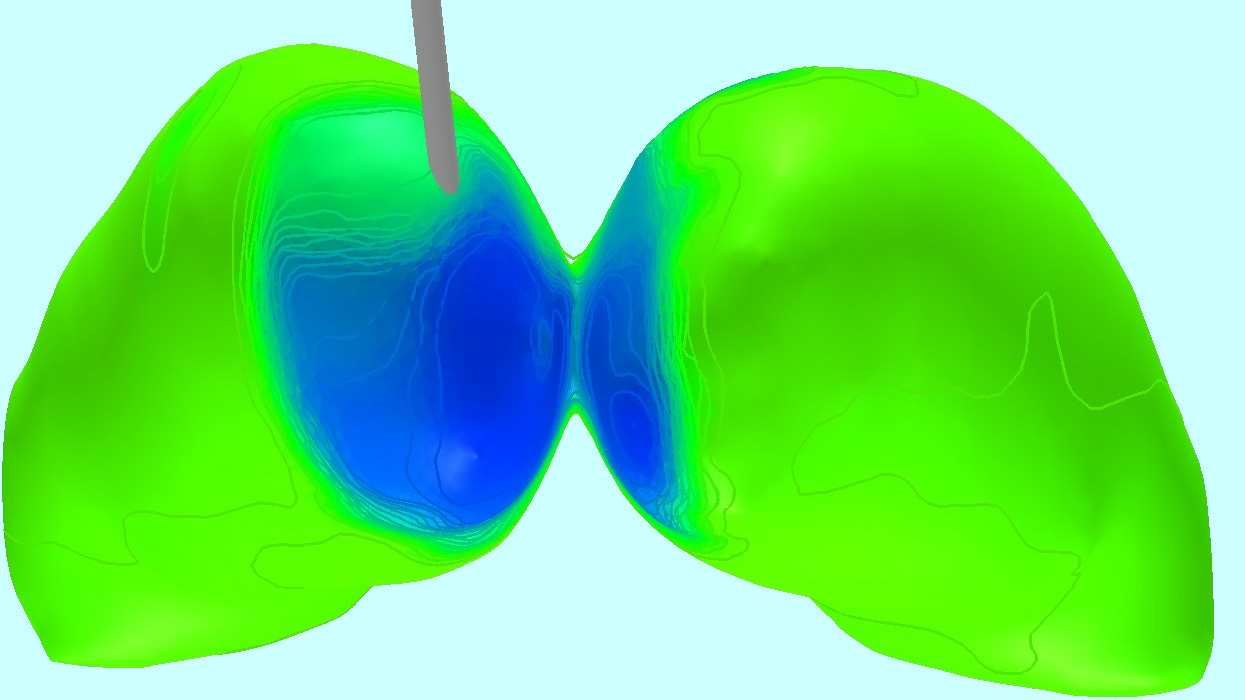} &
        \cellcolor[HTML]{ccffff} \includegraphics[width=2.25cm, height=1.5cm]{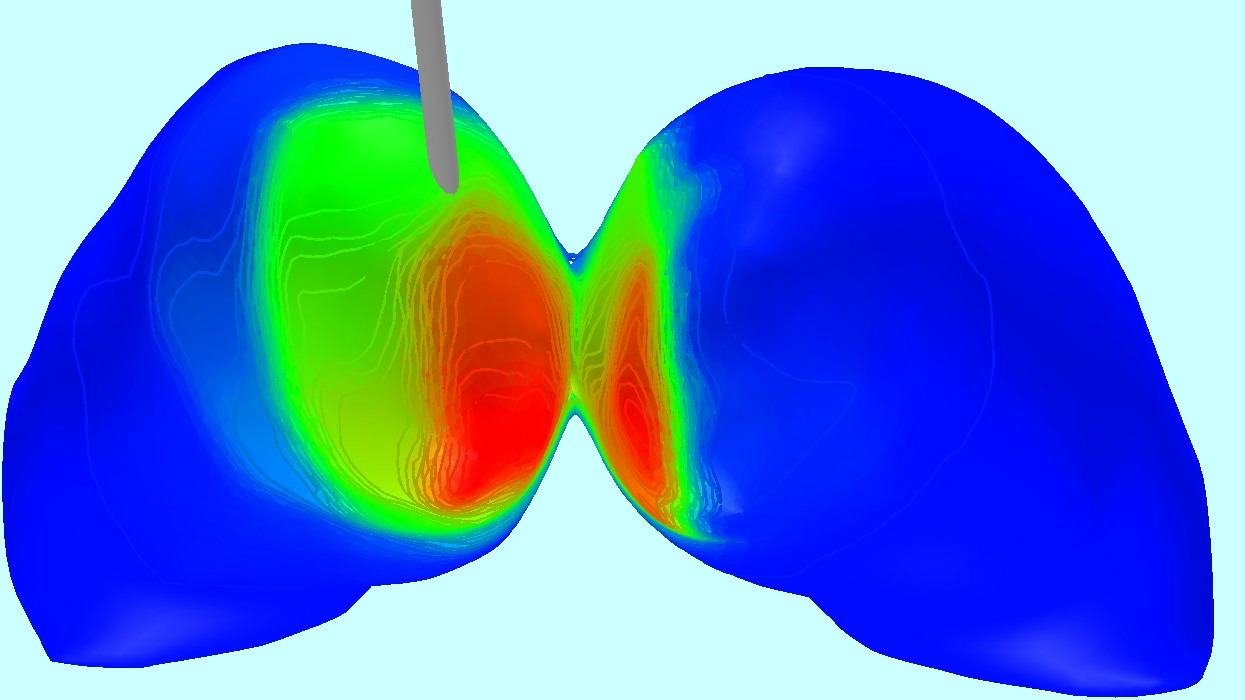} & 
        \cellcolor[HTML]{ccffff} \includegraphics[width=2.25cm, height=1.5cm]{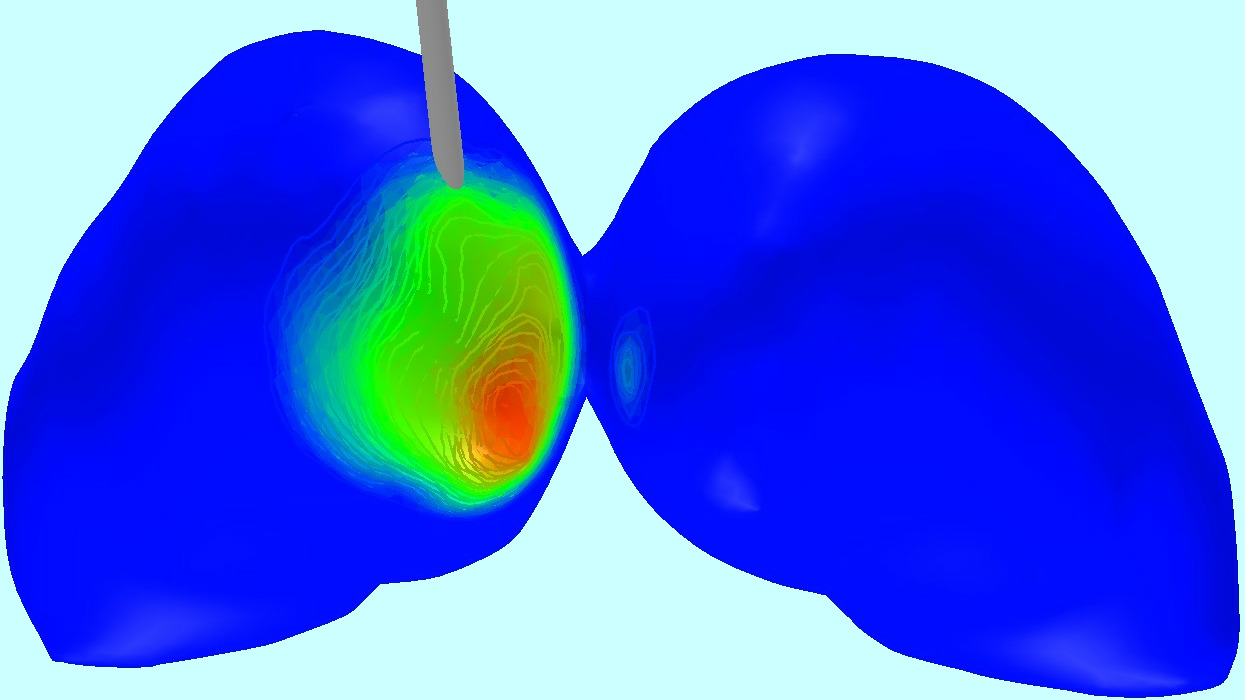} & 
        \cellcolor[HTML]{ccffff} \includegraphics[width=2.25cm, height=1.5cm]{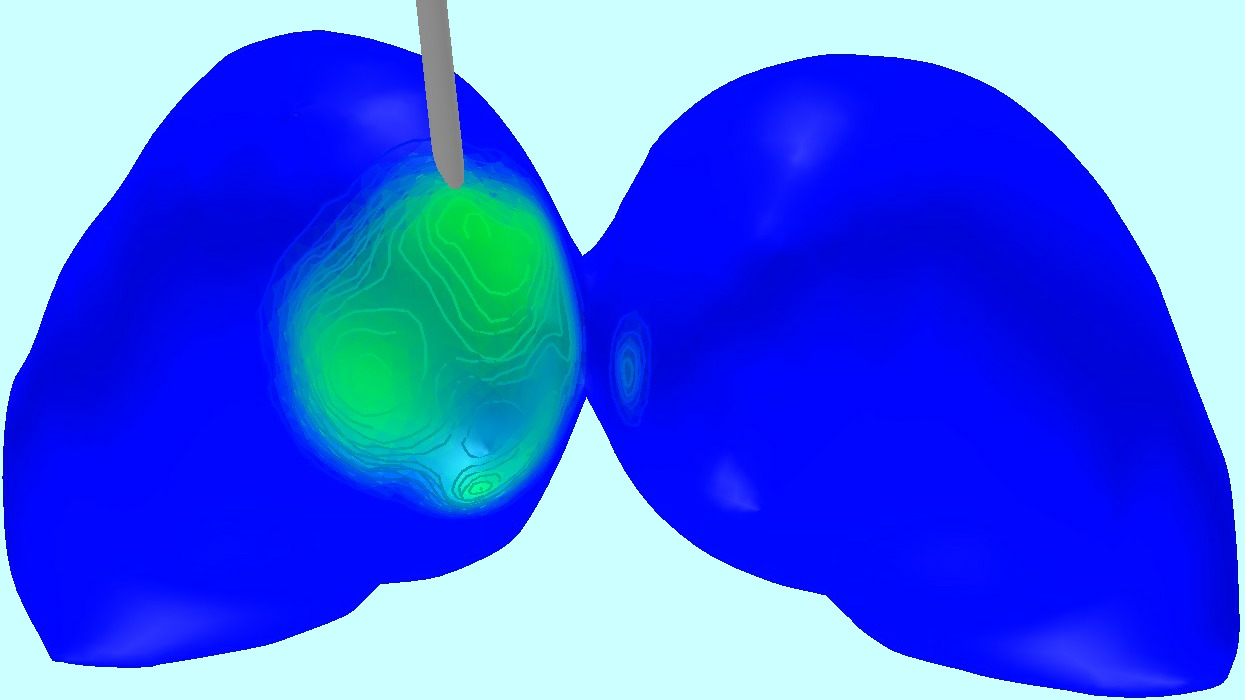} &
        \cellcolor[HTML]{ccffff} \includegraphics[width=2.25cm, height=1.5cm]{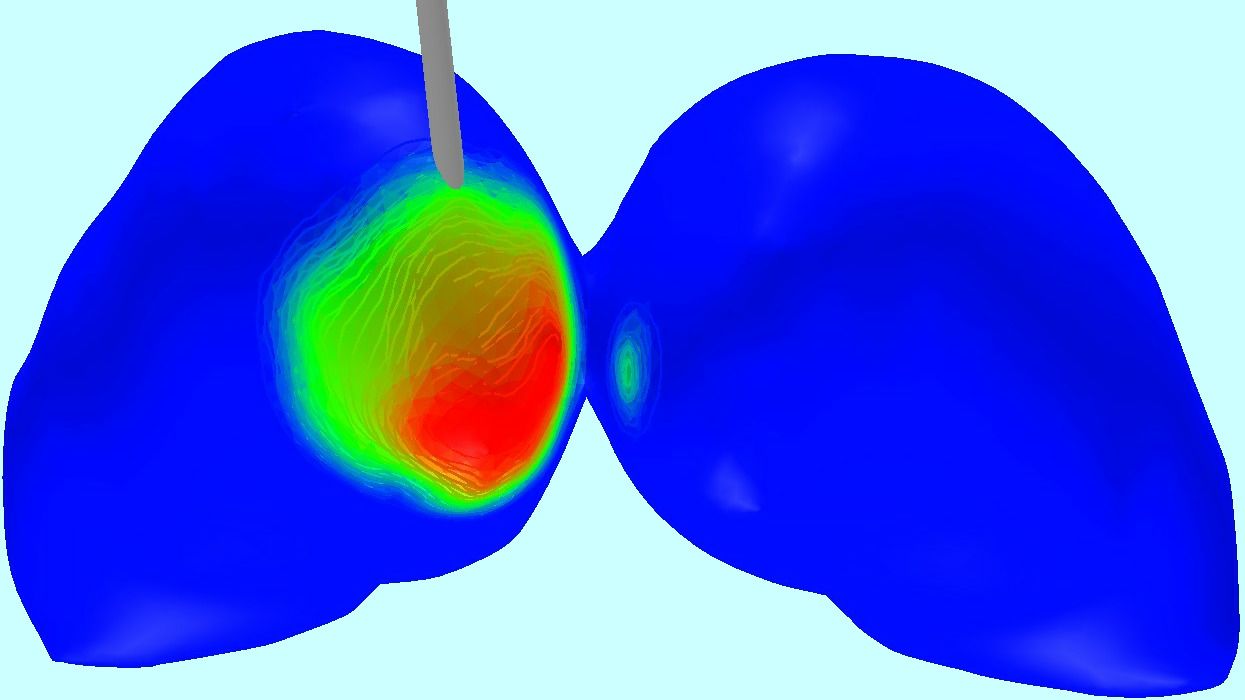} \\
        \cellcolor[HTML]{ccffff} \raisebox{3.25cm}[0pt][0pt]{\multirow{2}{*}{\rotatebox{90}{\textbf{L1L1(A)}, $\varepsilon \in  [ -160, 0] $ dB}}} \hskip0.1cm \raisebox{2.9cm}[0pt][0pt] {\multirow{2}{*}{\rotatebox{90}{Perpendicular \hskip0.25cm Parallel}}} &
        \rule{0pt}{1.80cm}
        \cellcolor[HTML]{ccffff} \includegraphics[width=2.25cm, height=1.5cm]{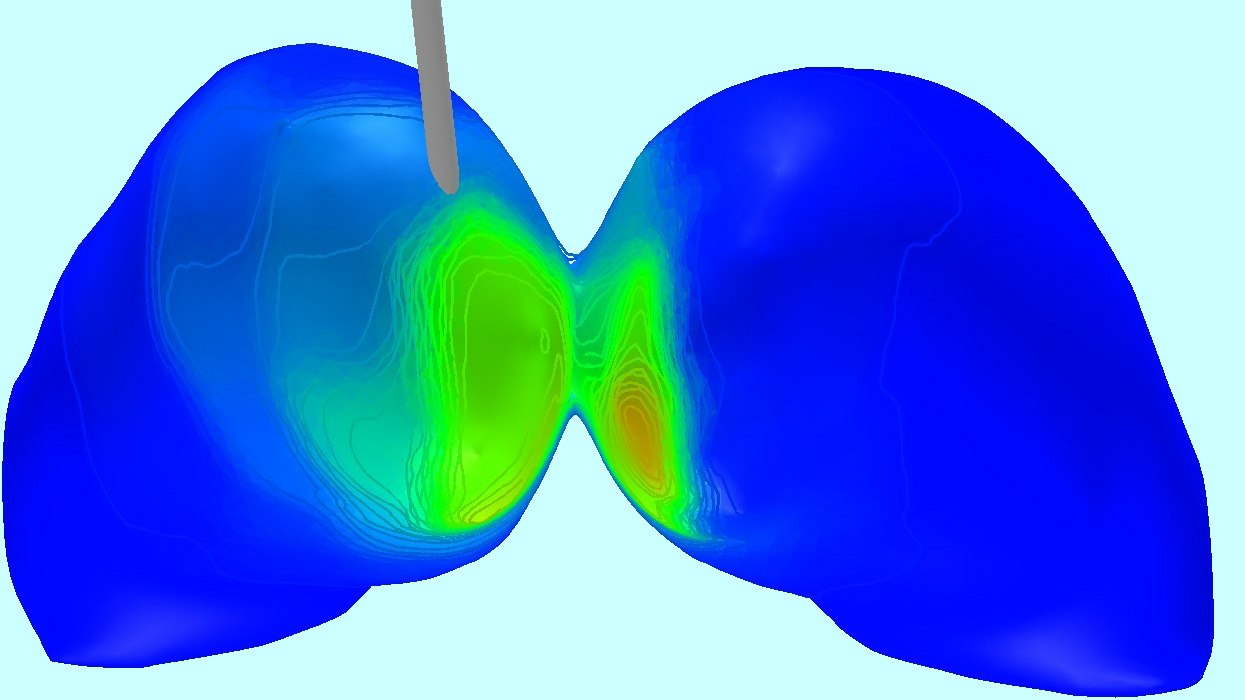} &  
        \cellcolor[HTML]{ccffff} \includegraphics[width=2.25cm, height=1.5cm]{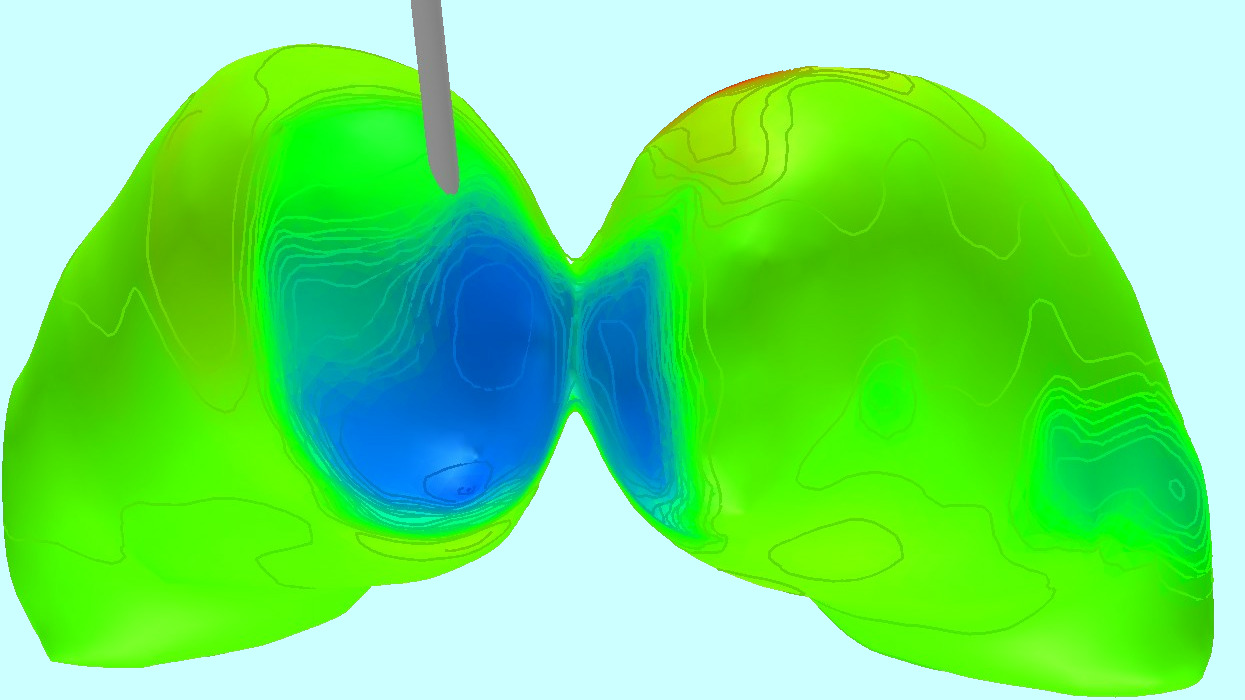} &
        \cellcolor[HTML]{ccffff} \includegraphics[width=2.25cm, height=1.5cm]{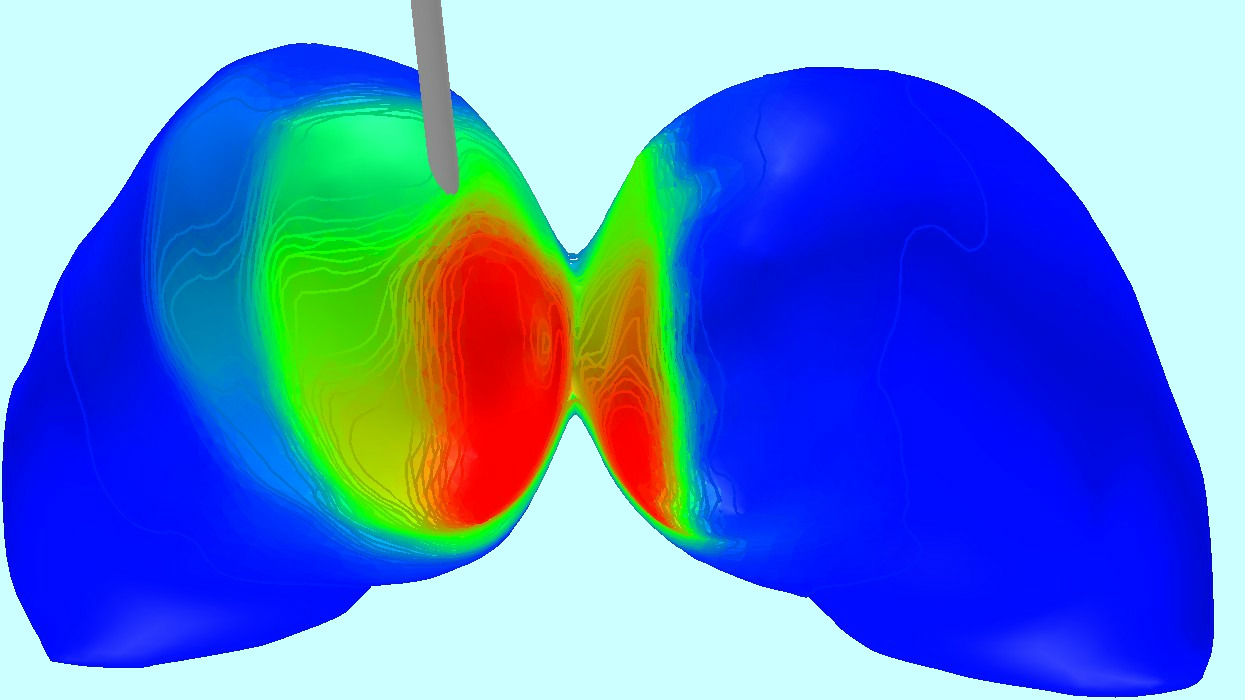} & 
        \cellcolor[HTML]{ccffff} \includegraphics[width=2.25cm, height=1.5cm]{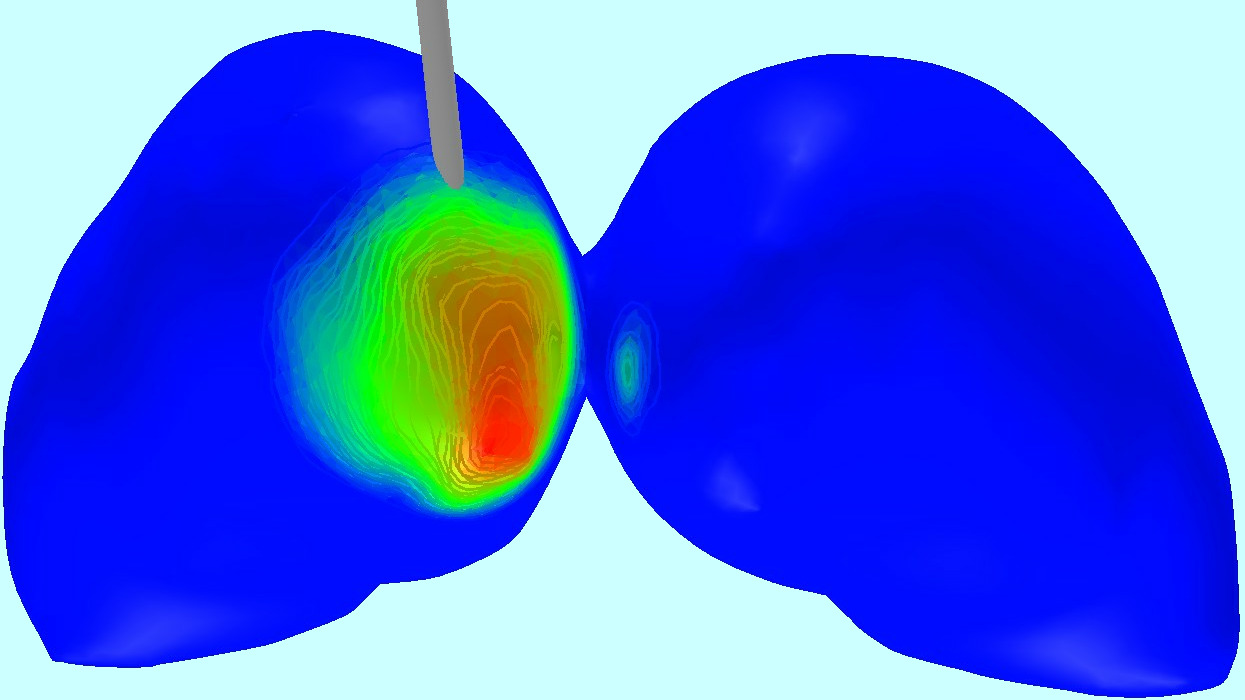} & 
        \cellcolor[HTML]{ccffff} \includegraphics[width=2.25cm, height=1.5cm]{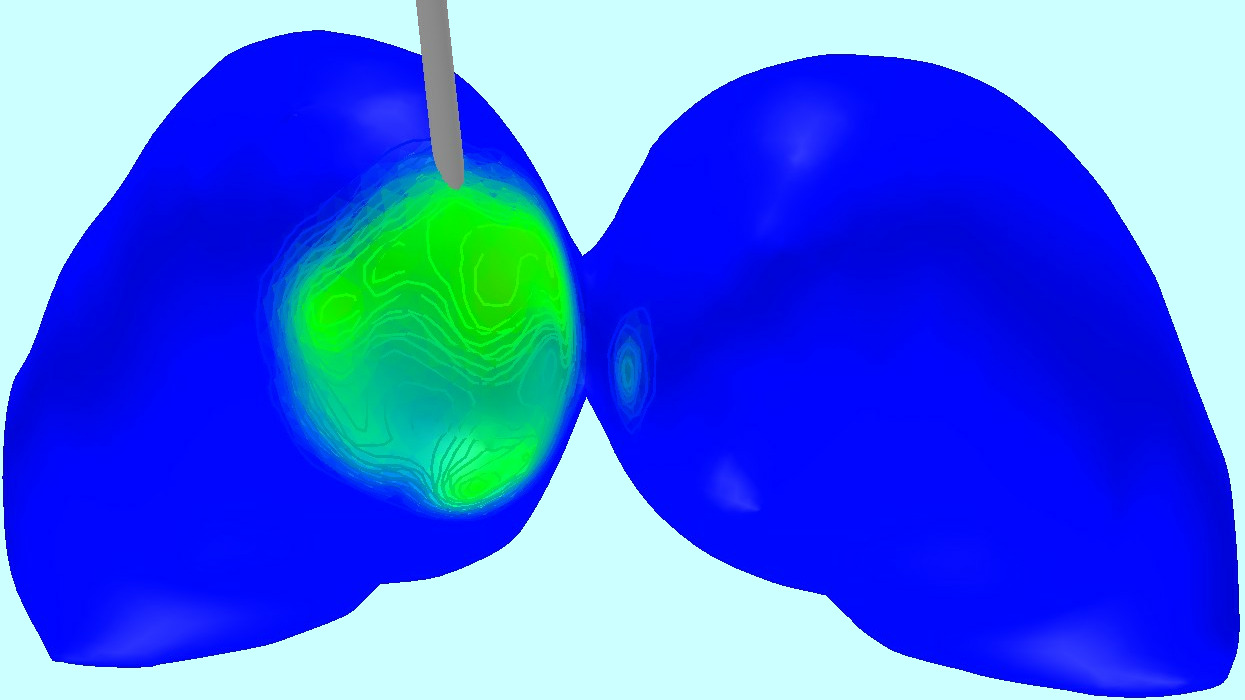} &
        \cellcolor[HTML]{ccffff} \includegraphics[width=2.25cm, height=1.5cm]{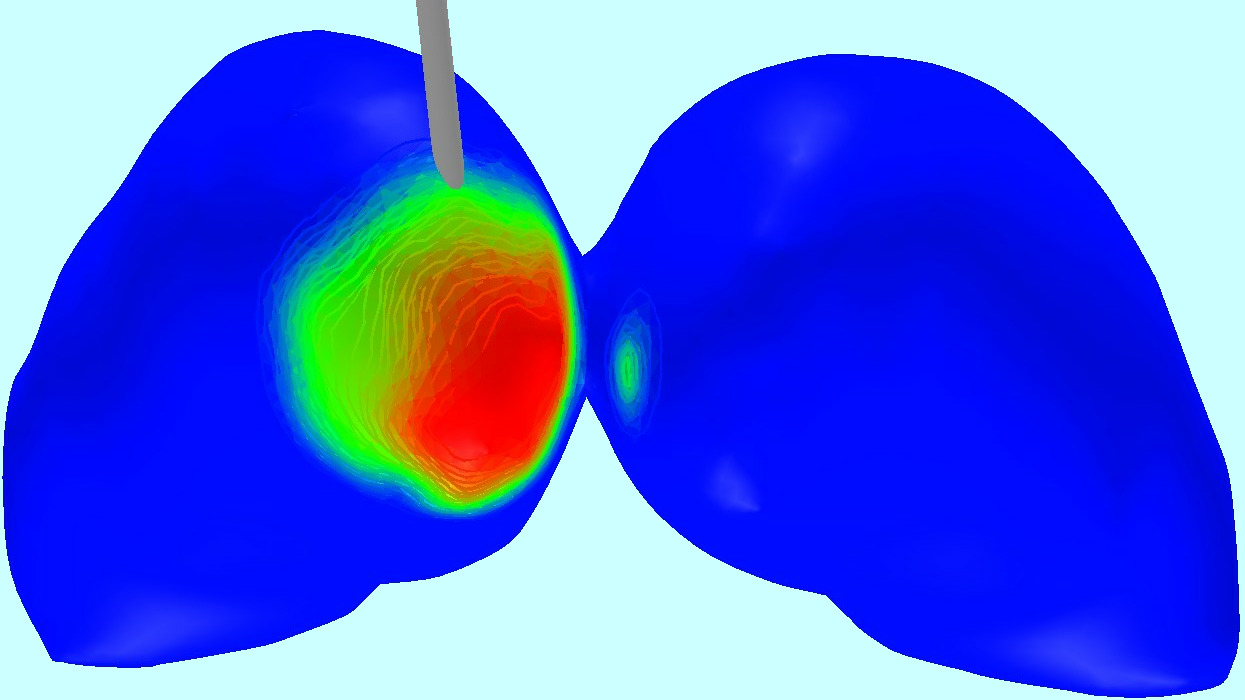} \\ 
        \hline
        
        \cellcolor[HTML]{e5ffd4}  &
        \rule{0pt}{1.80cm}
        \cellcolor[HTML]{e5ffd4} \includegraphics[width=2.25cm, height=1.5cm]{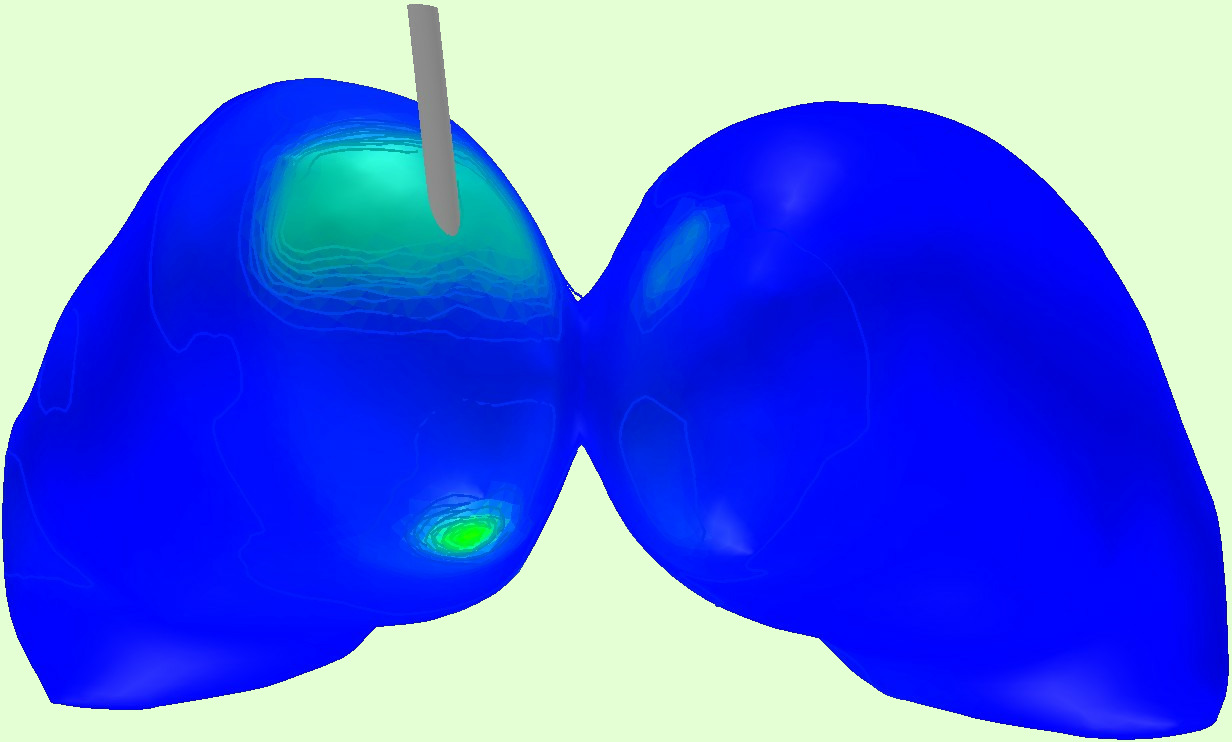} & 
        \cellcolor[HTML]{e5ffd4} \includegraphics[width=2.25cm, height=1.5cm]{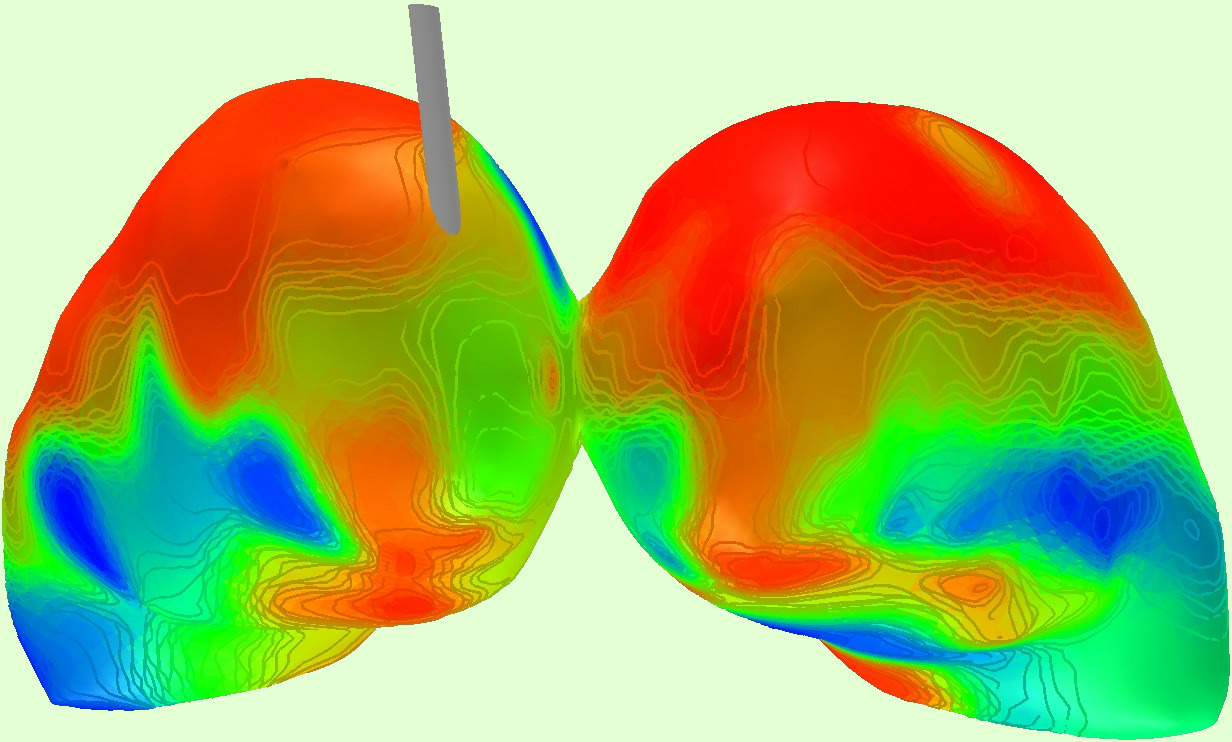} & 
        \cellcolor[HTML]{e5ffd4} \includegraphics[width=2.25cm, height=1.5cm]{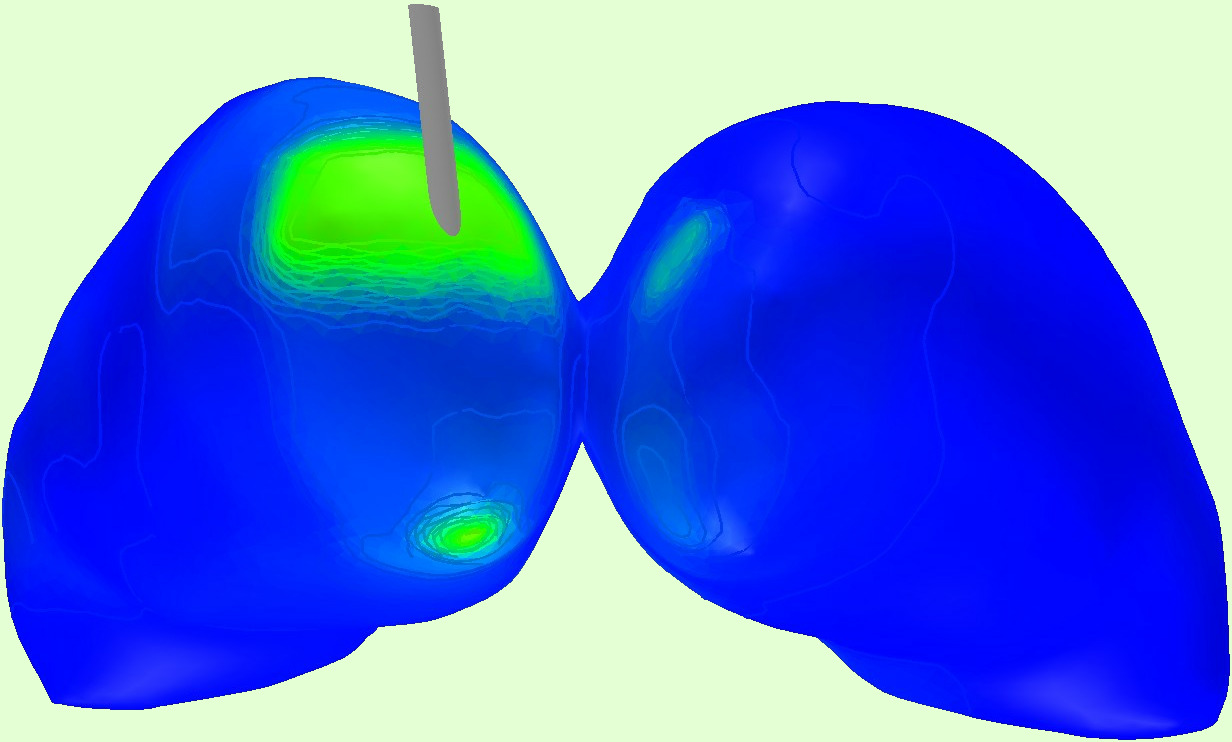} & 
        \cellcolor[HTML]{e5ffd4} \includegraphics[width=2.25cm, height=1.5cm]{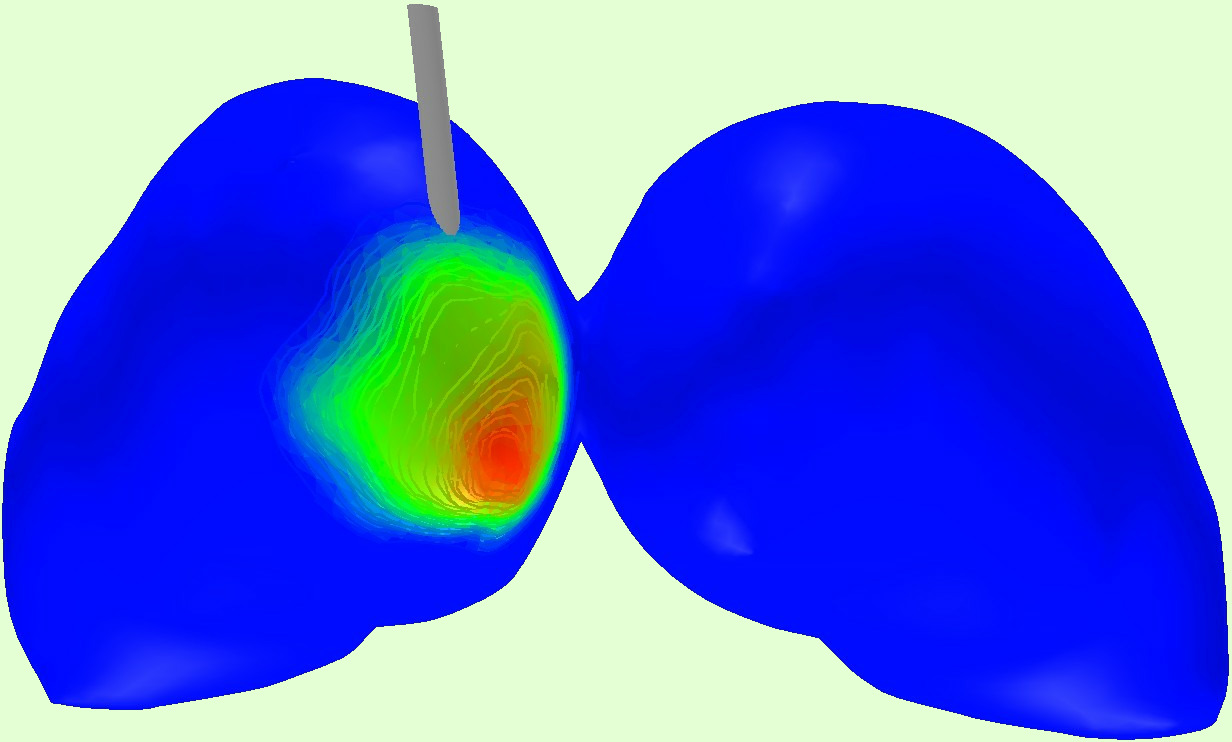} & 
        \cellcolor[HTML]{e5ffd4} \includegraphics[width=2.25cm, height=1.5cm]{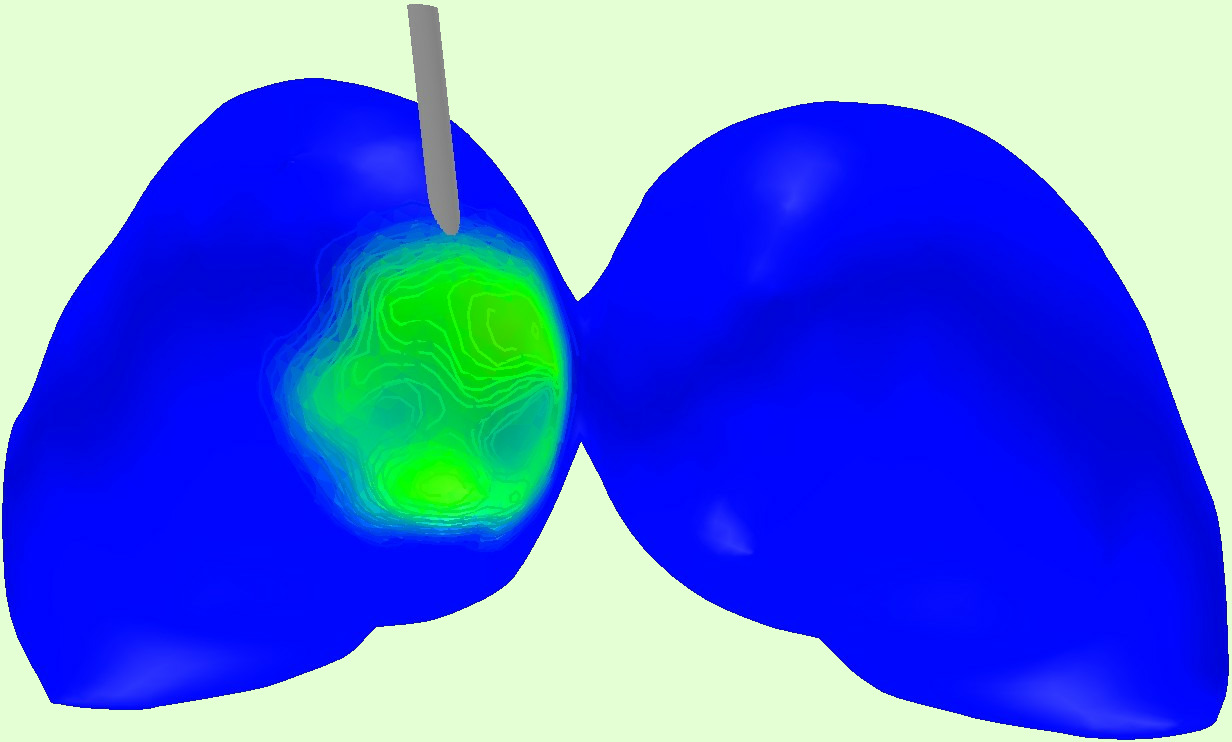} & 
        \cellcolor[HTML]{e5ffd4} \includegraphics[width=2.25cm, height=1.5cm]{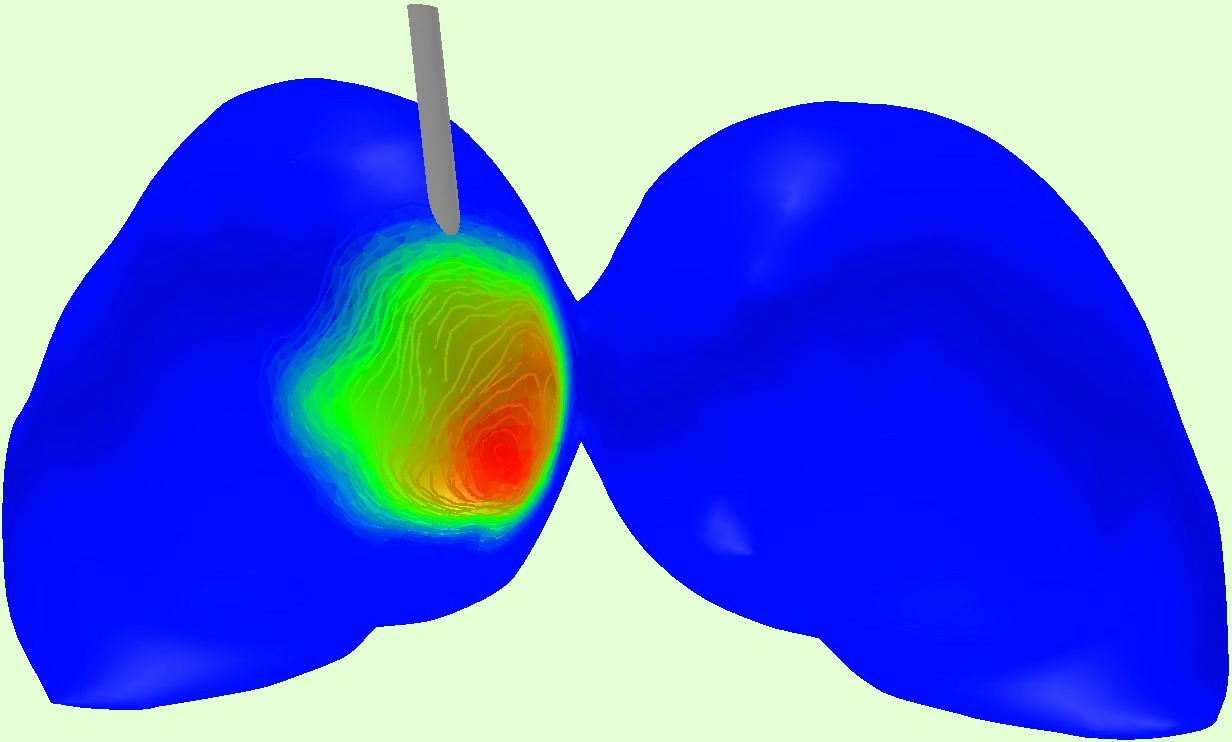} \\ 
        \cellcolor[HTML]{e5ffd4} \raisebox{2.90cm}[0pt][0pt]{\multirow{2}{*}{\rotatebox{90}{\textbf{TLS},  $\beta \in [-50, 40] $ dB}}} \hskip0.1cm \raisebox{2.9cm}[0pt][0pt] {\multirow{2}{*}{\rotatebox{90}{Perpendicular \hskip0.25cm Parallel}}}&
        \rule{0pt}{1.80cm}
        \cellcolor[HTML]{e5ffd4} \includegraphics[width=2.25cm, height=1.5cm]{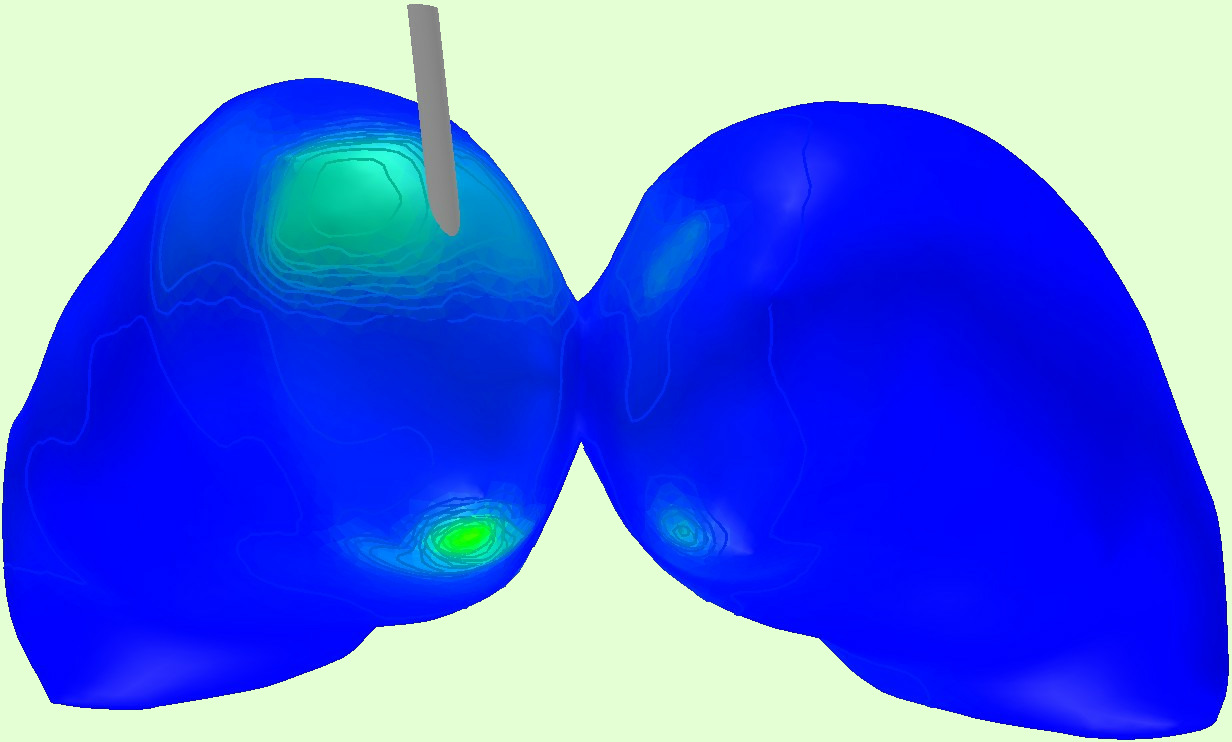} & 
        \cellcolor[HTML]{e5ffd4} \includegraphics[width=2.25cm, height=1.5cm]{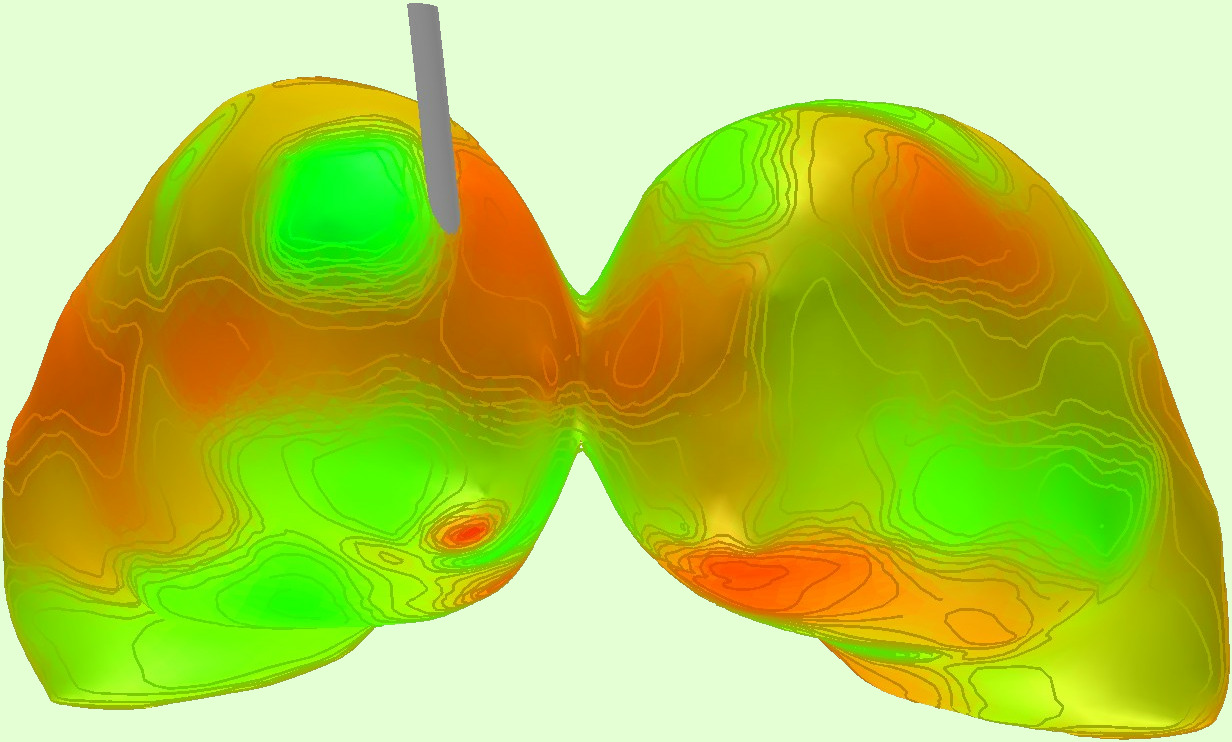} & 
        \cellcolor[HTML]{e5ffd4} \includegraphics[width=2.25cm, height=1.5cm]{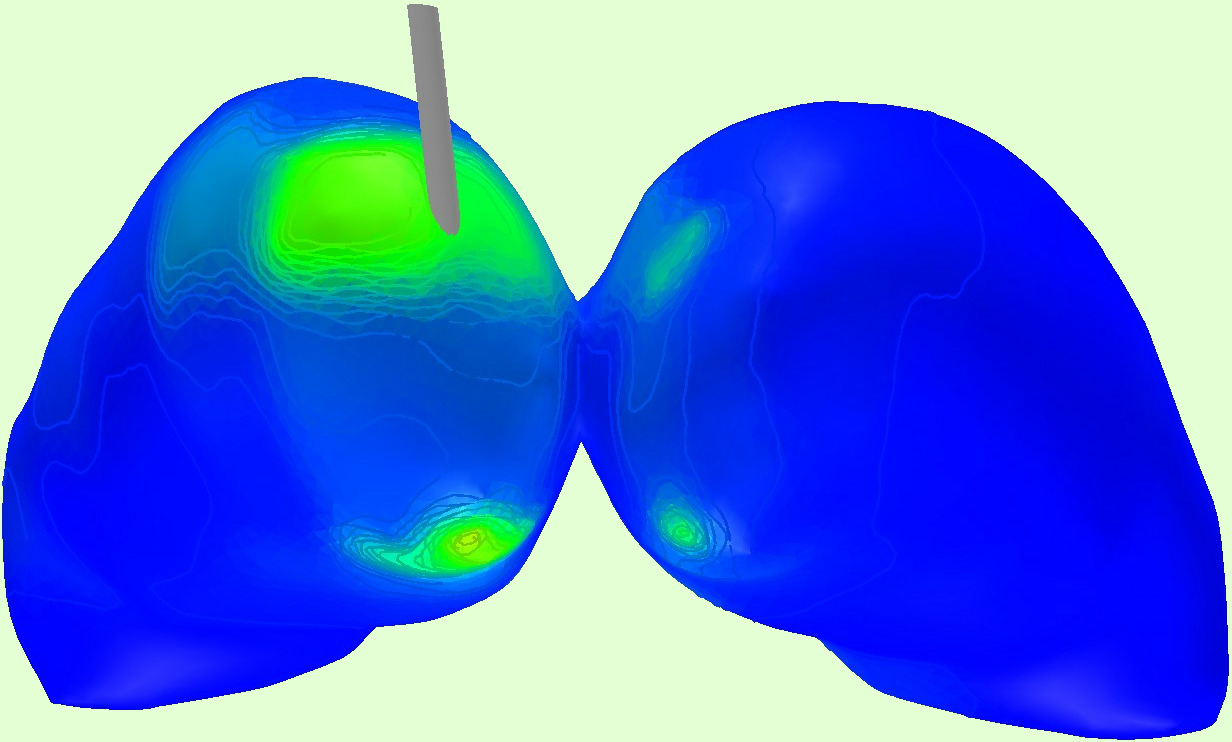} & 
        \cellcolor[HTML]{e5ffd4} \includegraphics[width=2.25cm, height=1.5cm]{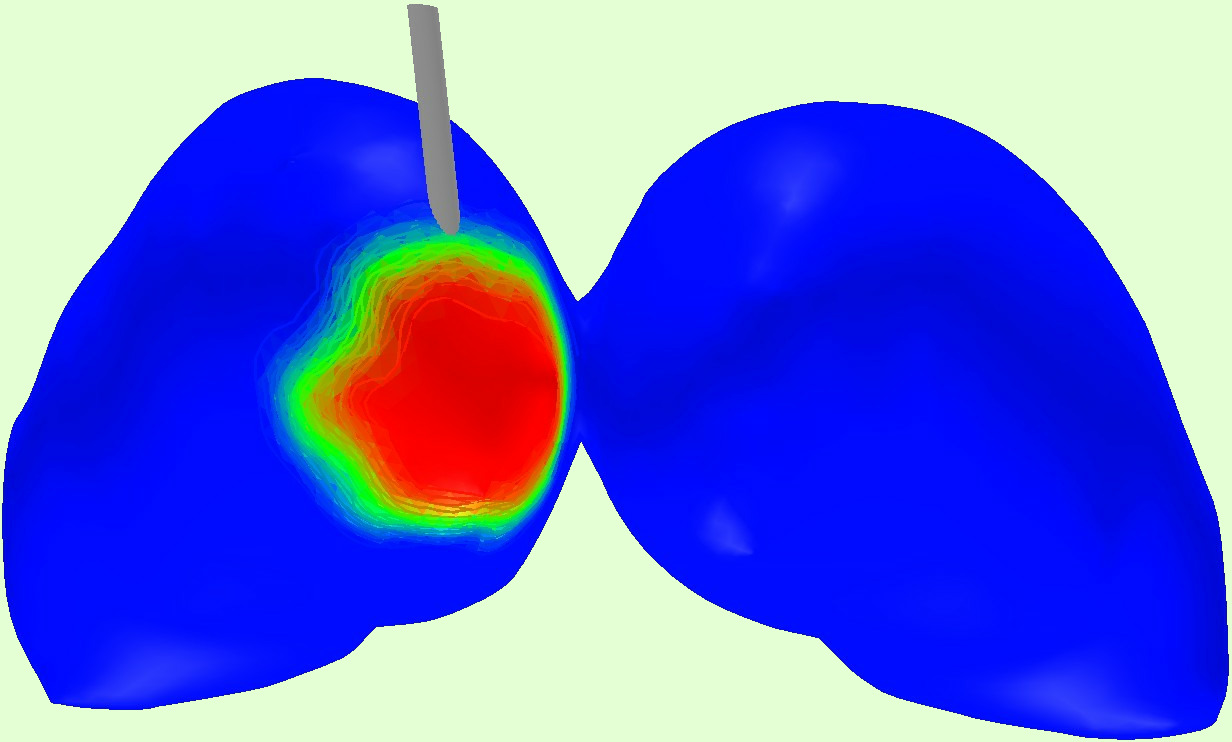} & 
        \cellcolor[HTML]{e5ffd4} \includegraphics[width=2.25cm, height=1.5cm]{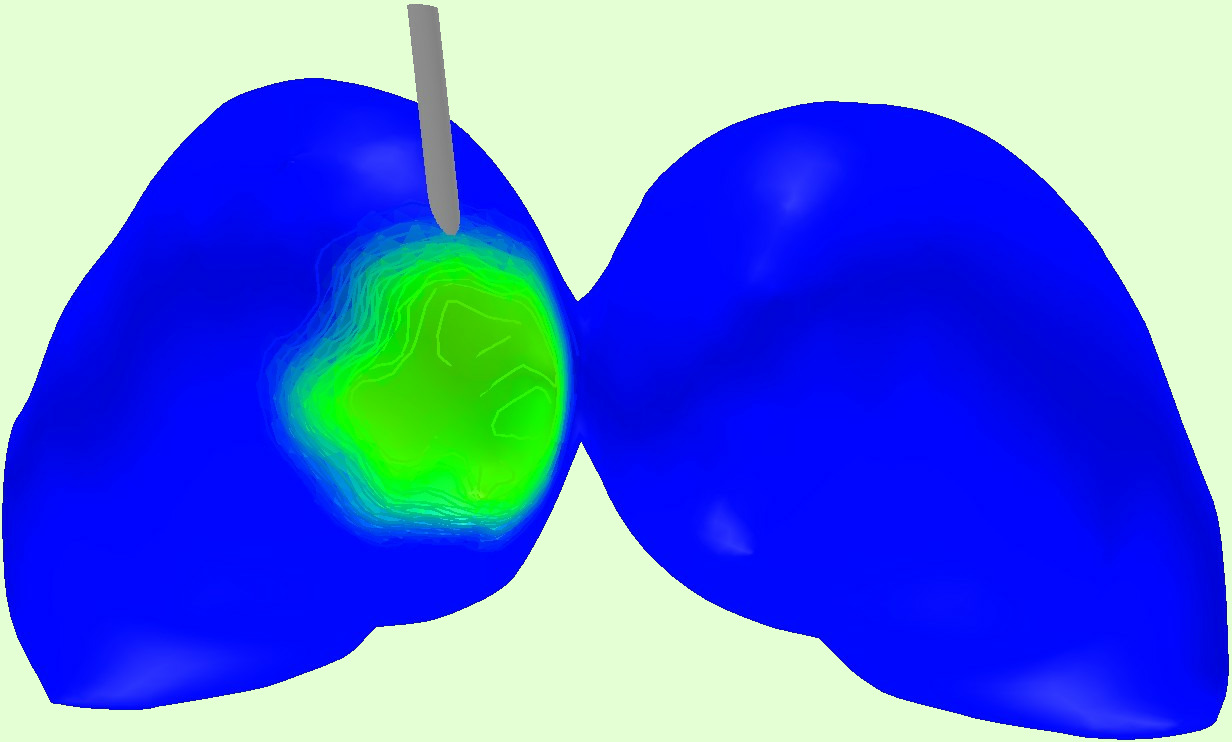} & 
        \cellcolor[HTML]{e5ffd4} \includegraphics[width=2.25cm, height=1.5cm]{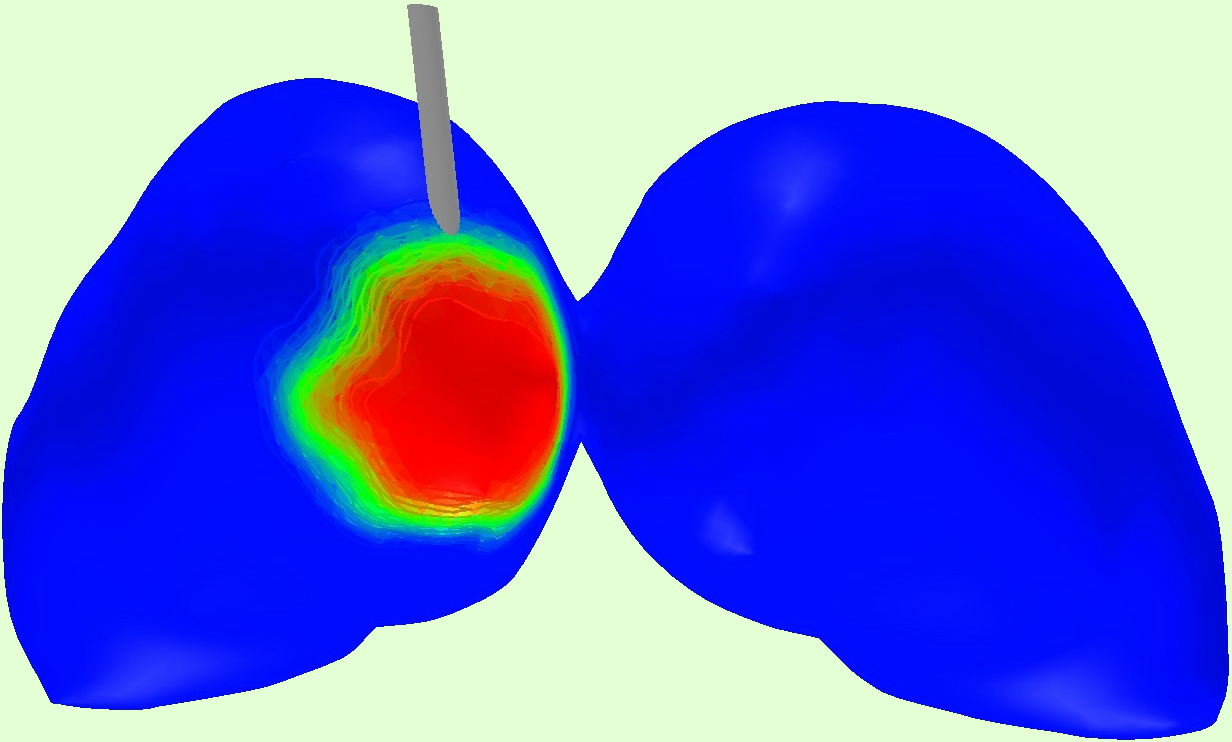} \\ 
        \hline
        
    \end{tabular}
    \end{small}
    \caption{
Volumetric optimization results for the 8-contact probe and noiseless lead fields. Decision variable $\Gamma$, $\Xi$ and $\Theta$ for optimizers found via Reciprocity Principle (RP), L1-norm regularized L1-norm fitting (L1L1) evaluated as L1L1(A) and L1L1(B) with feasible $\varepsilon$ threshold range: [-160, 0]~dB and [-10, 0]~dB, respectively, and Ti\-kho\-nov Regularized Least Squares (TLS) methods illustrated as a function of the target position and orientation (parallel or perpendicular w.r.t.\ the probe).  \textbf{Left:} distributions for a Low-Resolution Lead Field (\textbf{LR-LF}) where the spatial degrees of freedom (DOFs) cover both thalamus lobes. \textbf{Right:} distribution for a High-Resolution Lead Field (\textbf{HR-LF}) where the DOFs is confined to a region of interest within a radius of 6.0 mm (millimeters) from the center of the lead. The optimal electrode configuration found by the L1L1 method maximizes the field ratio $\Theta$ subject to $\Gamma >= \Gamma_0$. The $\Theta$ distribution found using L1L1(A) is elevated as compared to the case of RP and TLS. The limited dynamic range ($\varepsilon$ range) of L1L1(B) is reflected in an overall slightly lower $\Theta$.
\label{Fig:Abbott_results}}
\end{figure*}
\begin{figure*}[ht!]
    \centering
    \begin{small}
        
    \setlength{\tabcolsep}{0.2cm} % Horizontal spacing between images
    \renewcommand{\arraystretch}{1.0} % Vertical spacing between rows
    \setlength{\arrayrulewidth}{0.1mm}

    \begin{tabular}{|>{\centering}m{0.9cm}|ccc|ccc|}

        \multicolumn{1}{c}{} & 
        \multicolumn{3}{c}{\textbf{LR-LF}} & 
        \multicolumn{3}{c}{\textbf{HR-LF}} \\
        %\hline
        % Column headers with small figures, all framed
        \multicolumn{1}{c}{} & 
        \multicolumn{1}{c}{\textbf{Focused $\Gamma$ (A/m\textsuperscript{2})}} & 
        \multicolumn{1}{c}{\textbf{Nuisance $\Xi$ (A/m\textsuperscript{2})}} & 
        \multicolumn{1}{c}{\textbf{Ratio $\Theta$ (rel.)}} & 
        \multicolumn{1}{c}{\textbf{Focused $\Gamma$ (A/m\textsuperscript{2})}} & 
        \multicolumn{1}{c}{\textbf{Nuisance $\Xi$ (A/m\textsuperscript{2})}} & 
        \multicolumn{1}{c}{\textbf{Ratio $\Theta$ (rel.)}} \\ 
        
        \multicolumn{1}{c}{} & 
        \multicolumn{1}{c}{\includegraphics[width=2.25cm, height=0.45cm]{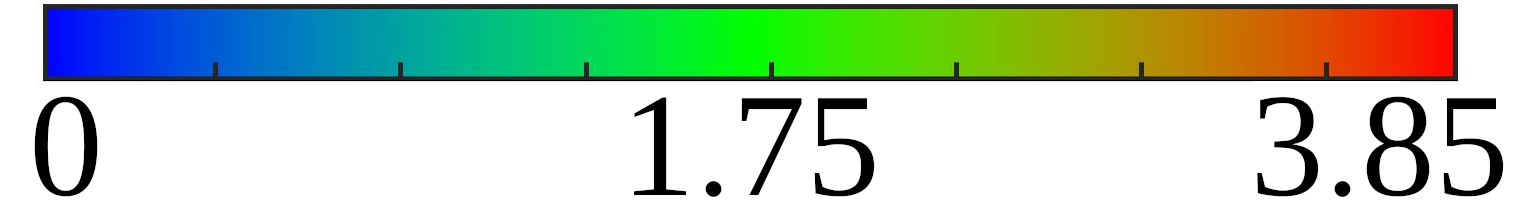}} & 
        \multicolumn{1}{c}{\includegraphics[width=2.25cm, height=0.45cm]{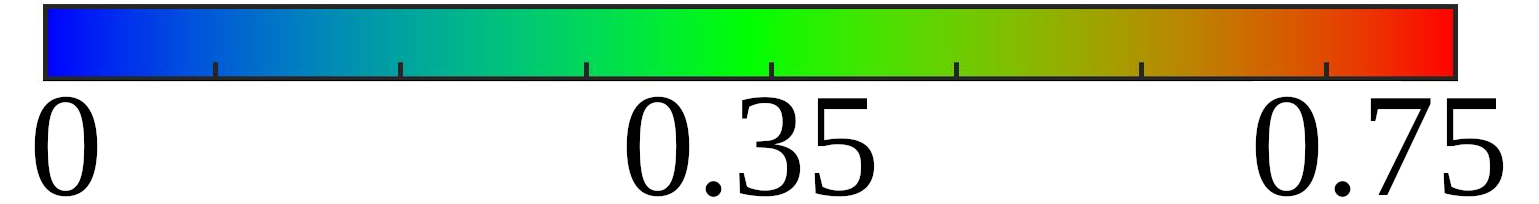}} & 
        \multicolumn{1}{c}{\includegraphics[width=2.25cm, height=0.45cm]{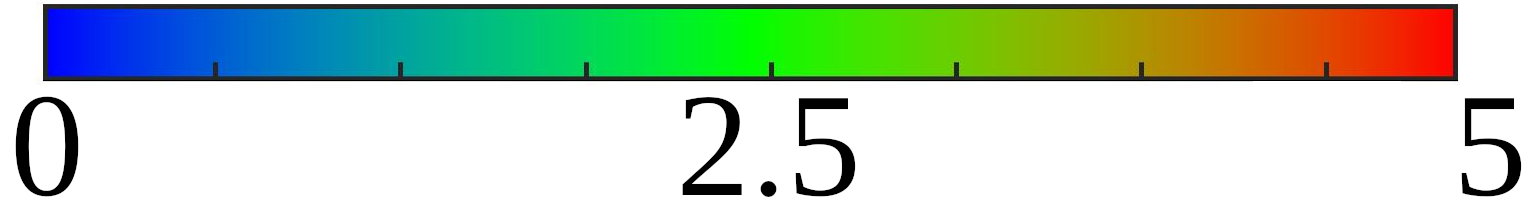}} & 
        \multicolumn{1}{c}{\includegraphics[width=2.25cm, height=0.45cm]{Figure_MED_col/colorbar_focused.png}} & 
        \multicolumn{1}{c}{\includegraphics[width=2.25cm, height=0.45cm]{Figure_MED_col/colorbar_nuisance.png}} & 
        \multicolumn{1}{c}{\includegraphics[width=2.25cm, height=0.45cm]{Figure_MED_col/colorbar_ratio.png}} \\
        \hline
        
        \cellcolor[HTML]{fff0ff} & 
        \rule{0pt}{1.80cm}
        \cellcolor[HTML]{fff0ff} \includegraphics[width=2.25cm, height=1.5cm]{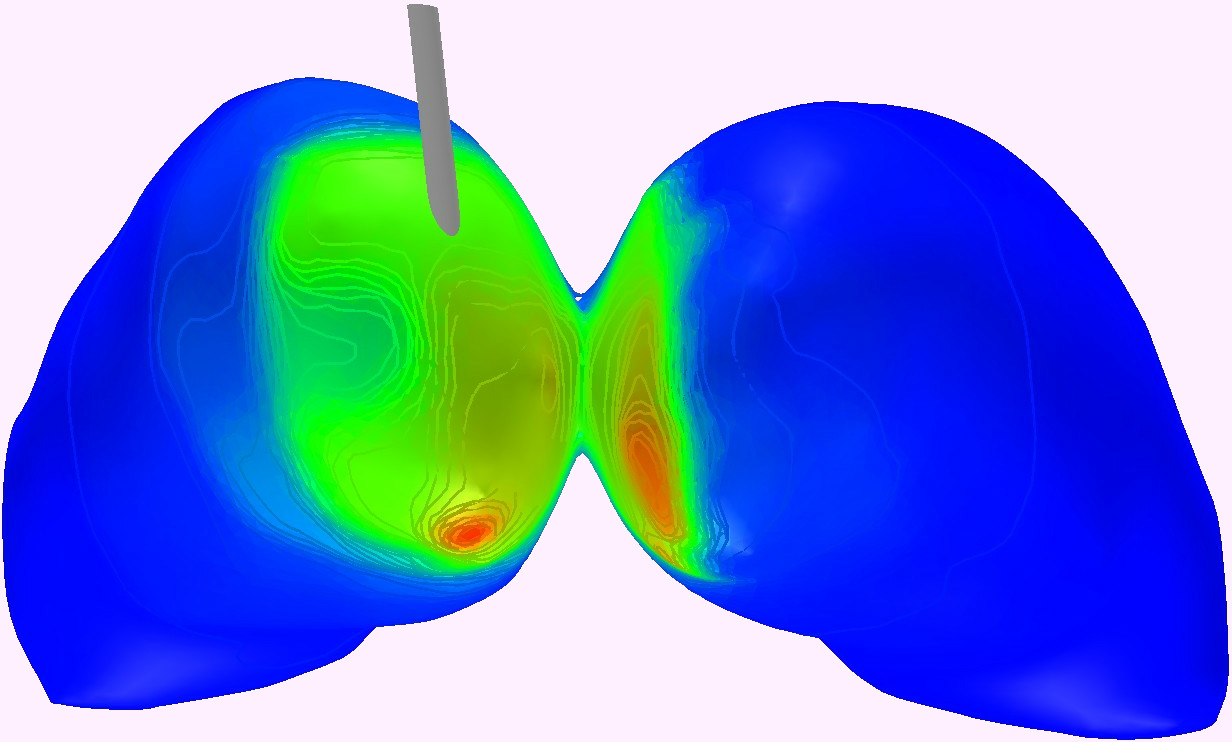} &  
        \cellcolor[HTML]{fff0ff} \includegraphics[width=2.25cm, height=1.5cm]{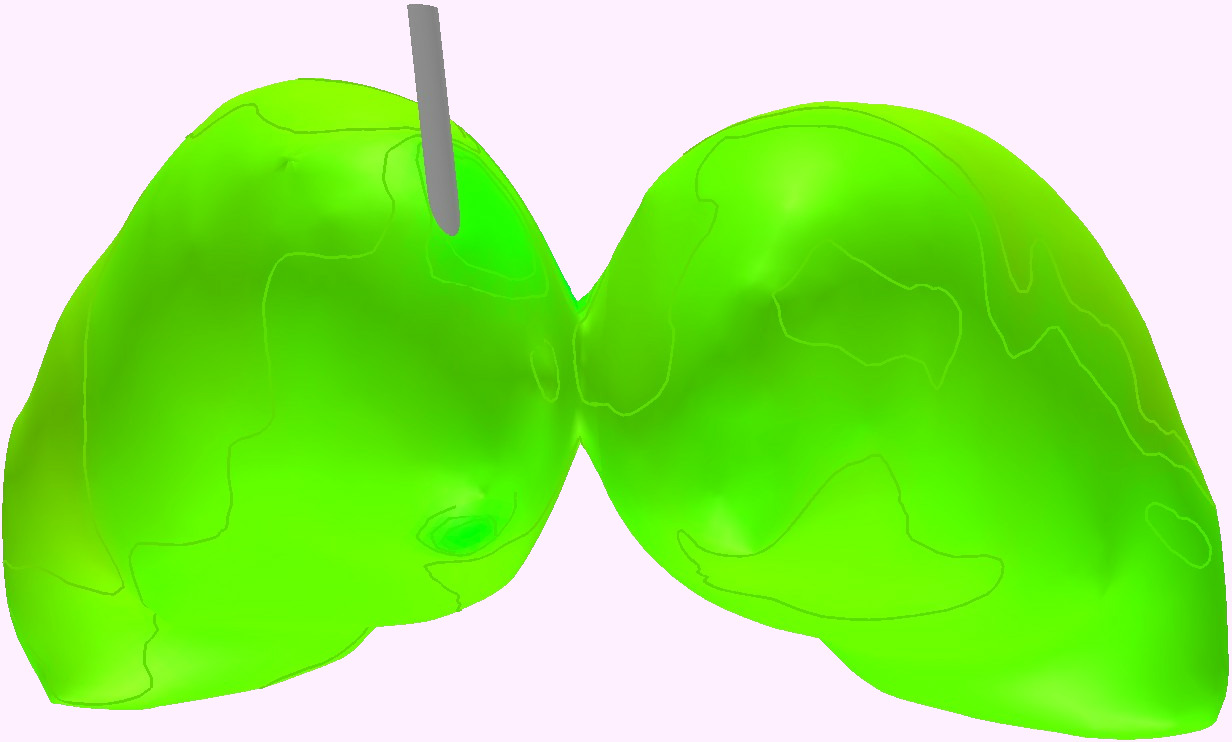} &
        \cellcolor[HTML]{fff0ff} \includegraphics[width=2.25cm, height=1.5cm]{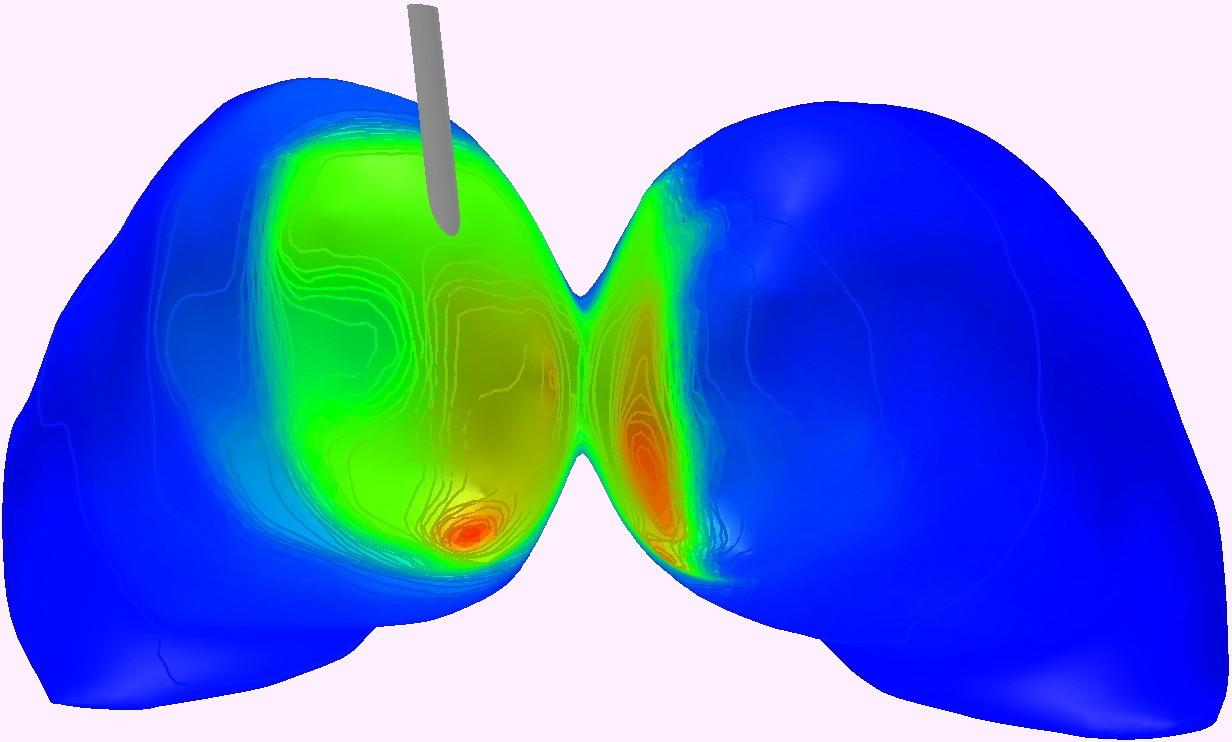} & 
        \cellcolor[HTML]{fff0ff} \includegraphics[width=2.25cm, height=1.5cm]{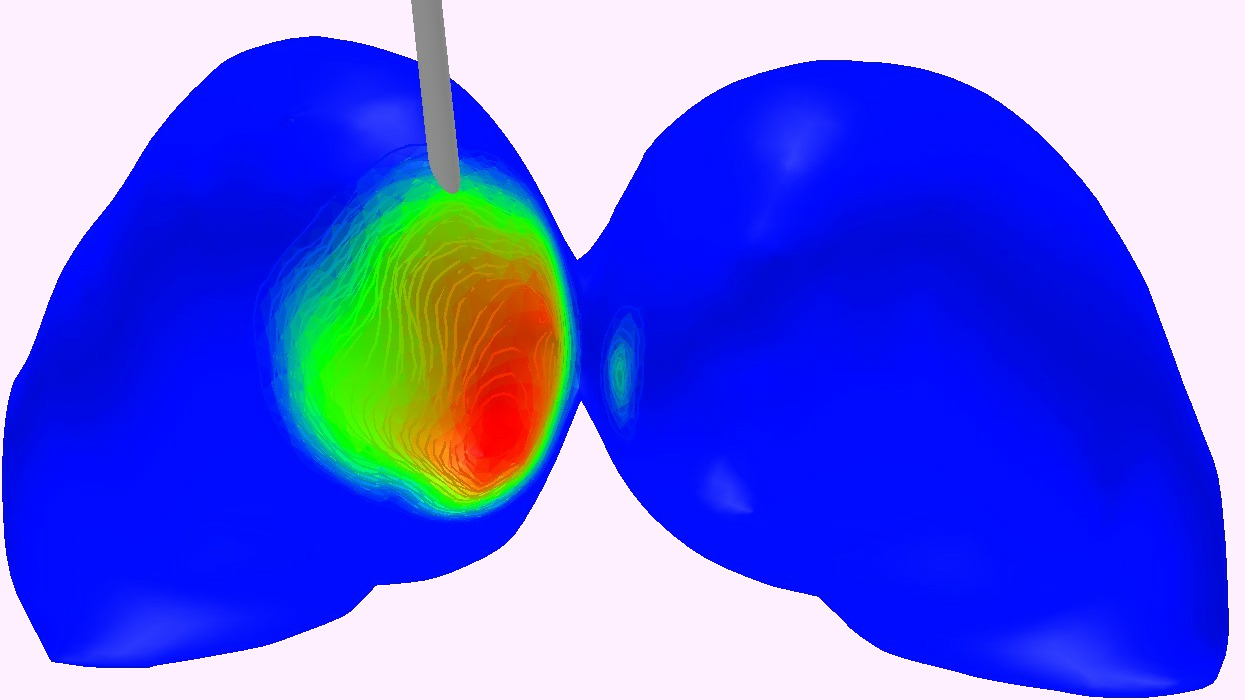} & 
        \cellcolor[HTML]{fff0ff} \includegraphics[width=2.25cm, height=1.5cm]{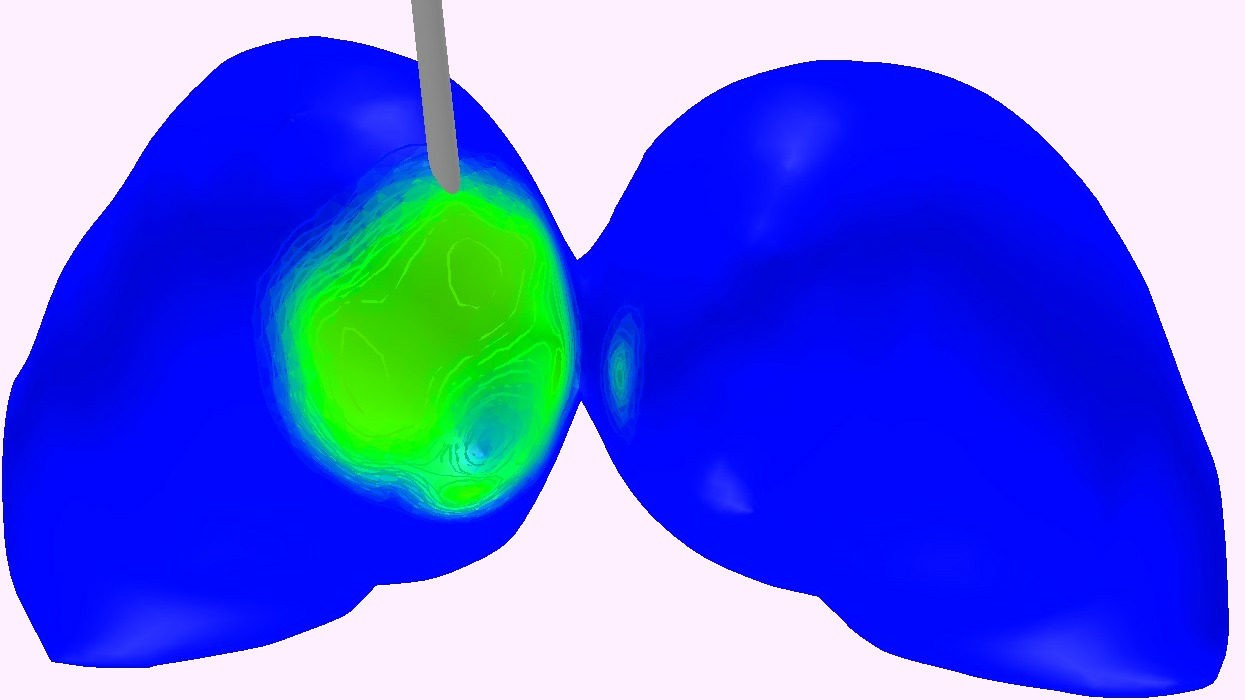} &
        \cellcolor[HTML]{fff0ff} \includegraphics[width=2.25cm, height=1.5cm]{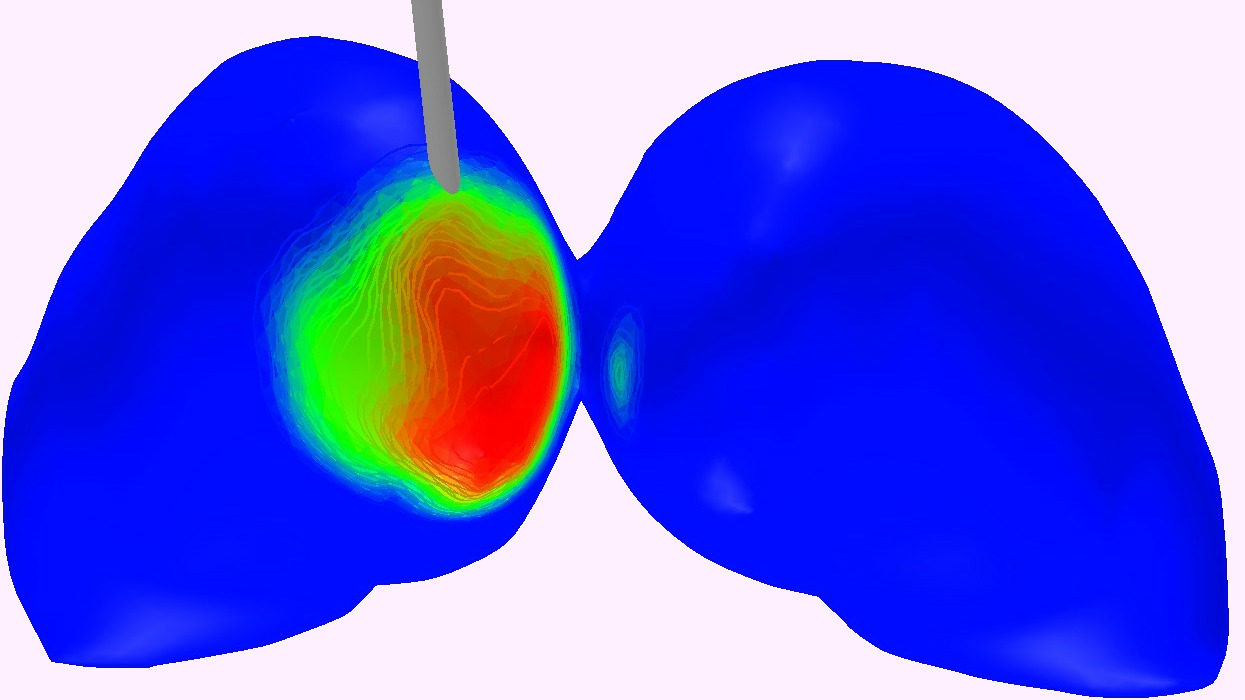} \\
        \cellcolor[HTML]{fff0ff} \raisebox{1.7cm}[0pt][0pt]
        {\multirow{2}{*}{\rotatebox{90}{\textbf{RP}} }} \hskip0.1cm
        \raisebox{2.9cm}[0pt][0pt] {\multirow{2}{*}{\rotatebox{90}{ Perpendicular \hskip0.25cm Parallel}}}
        &    
        \rule{0pt}{1.80cm}
        \cellcolor[HTML]{fff0ff} \includegraphics[width=2.25cm, height=1.5cm]{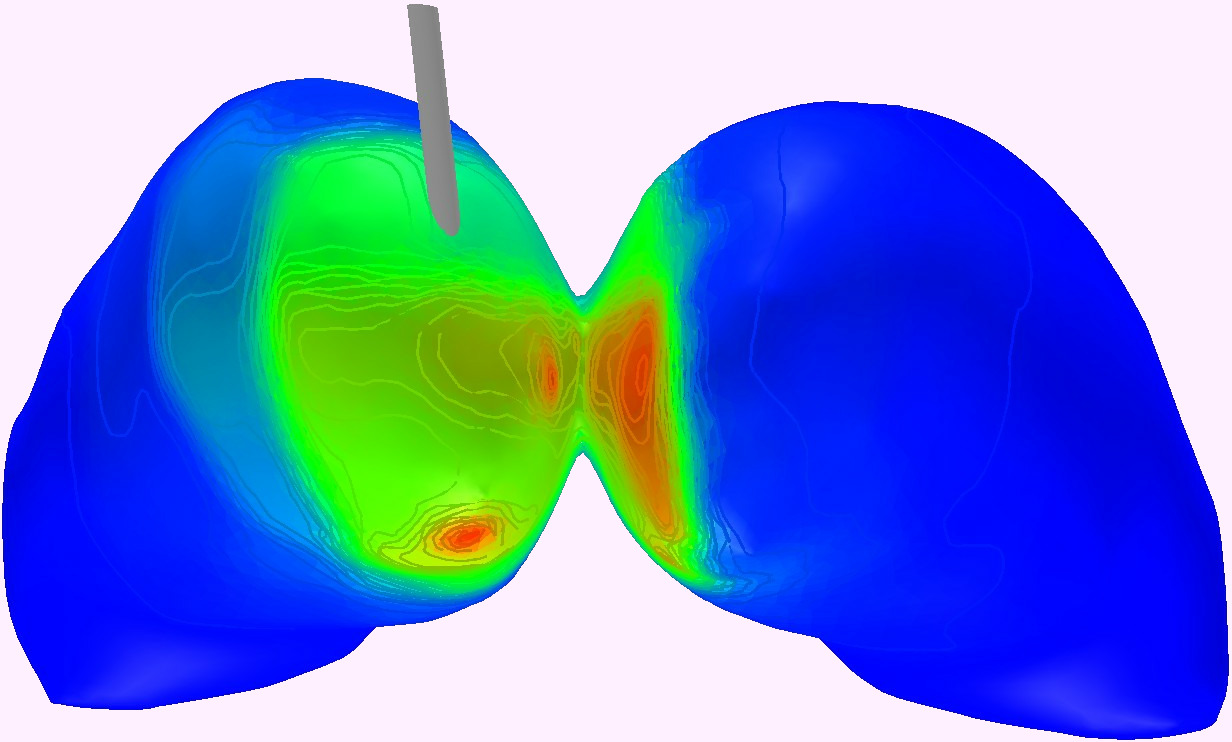} &  
        \cellcolor[HTML]{fff0ff} \includegraphics[width=2.25cm, height=1.5cm]{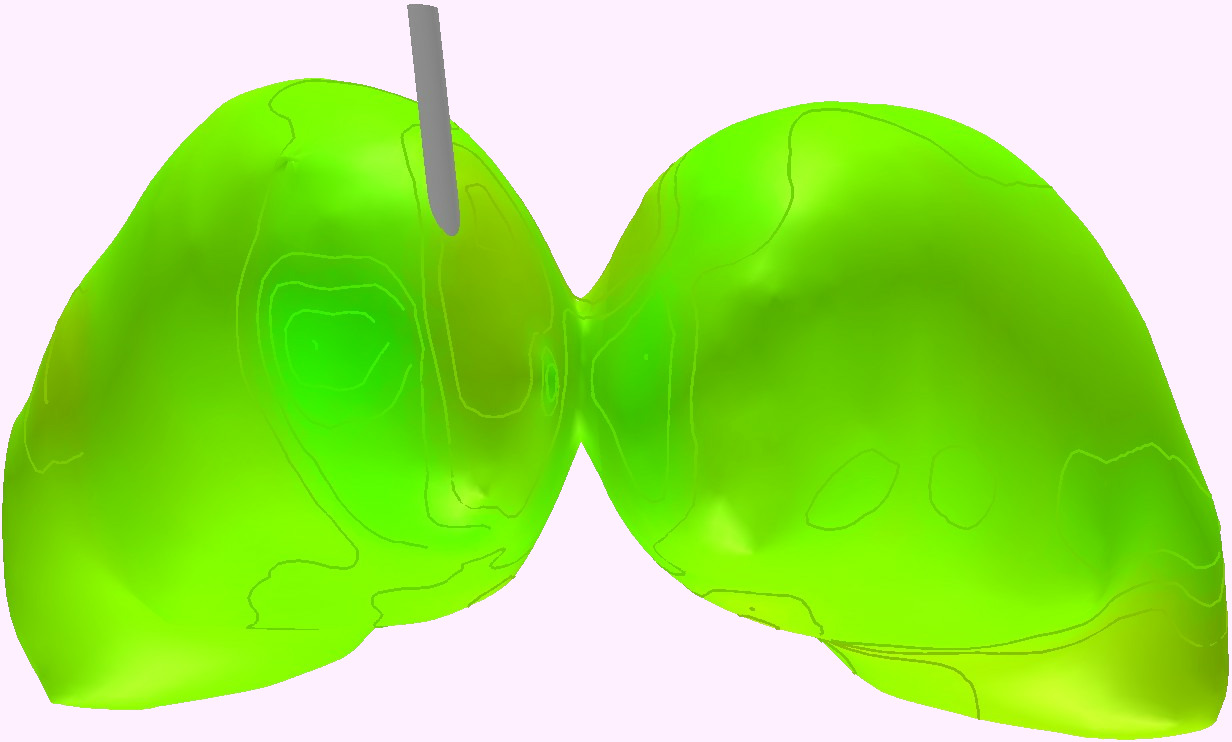} &
        \cellcolor[HTML]{fff0ff} \includegraphics[width=2.25cm, height=1.5cm]{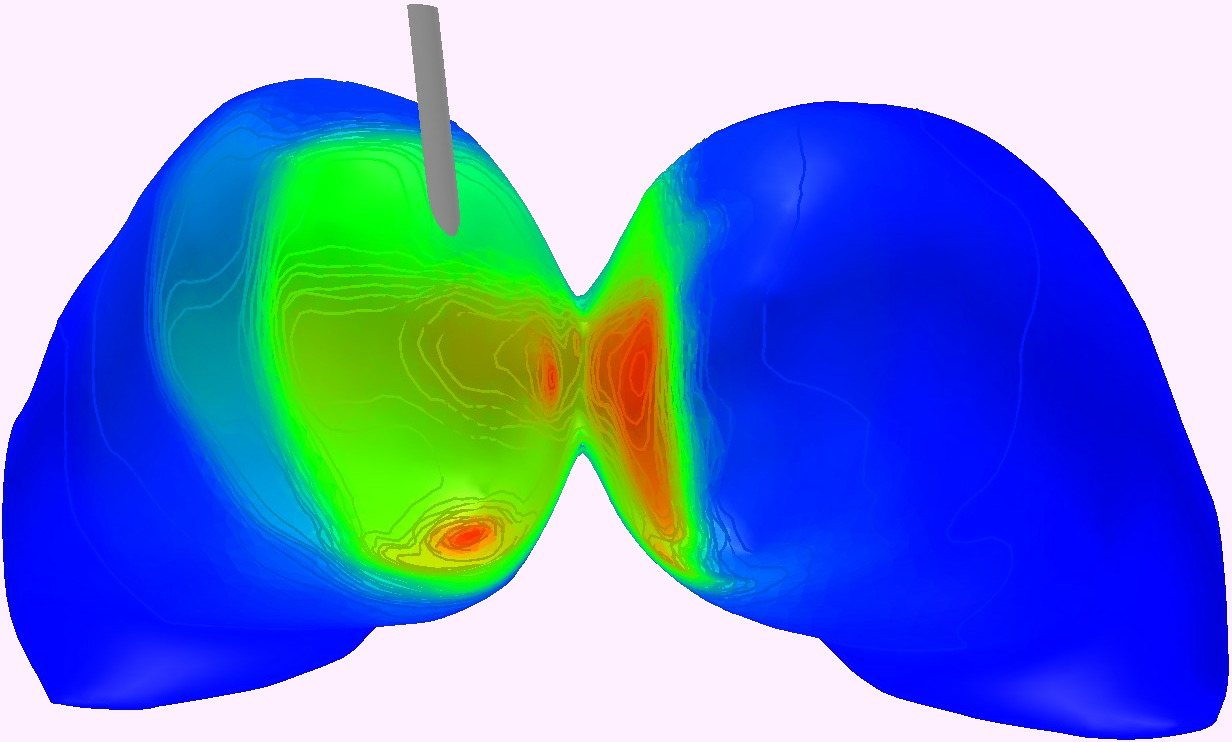} & 
        \cellcolor[HTML]{fff0ff} \includegraphics[width=2.25cm, height=1.5cm]{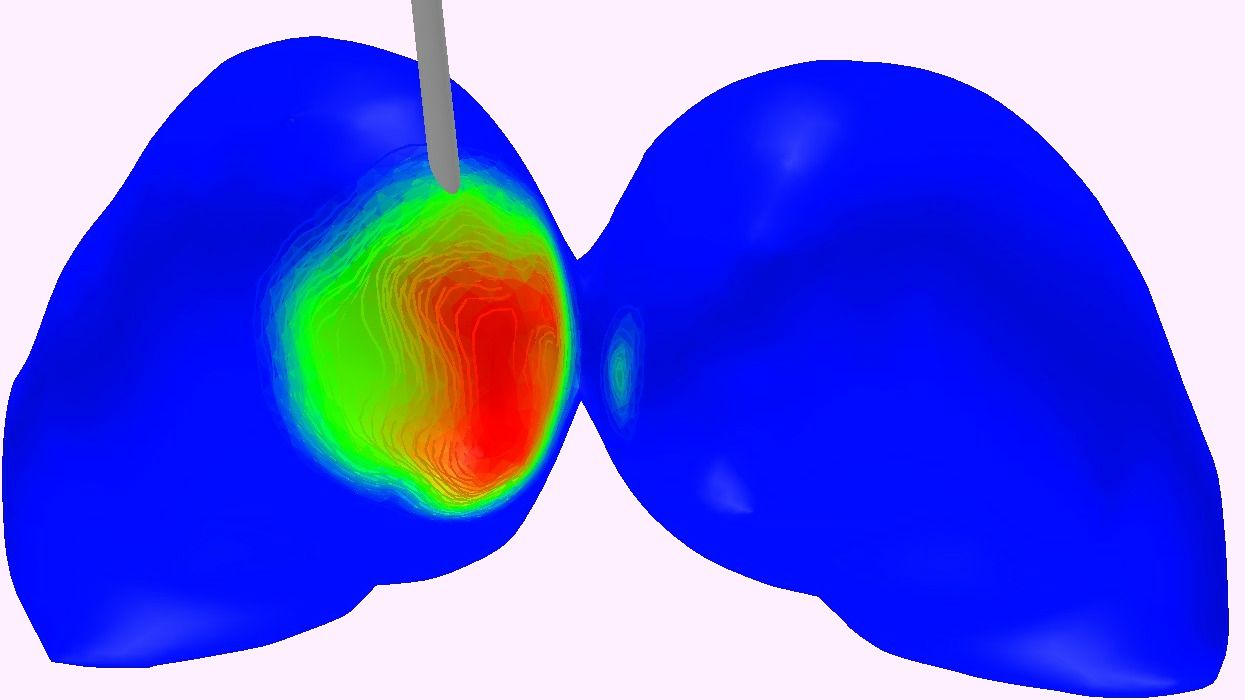} & 
        \cellcolor[HTML]{fff0ff} \includegraphics[width=2.25cm, height=1.5cm]{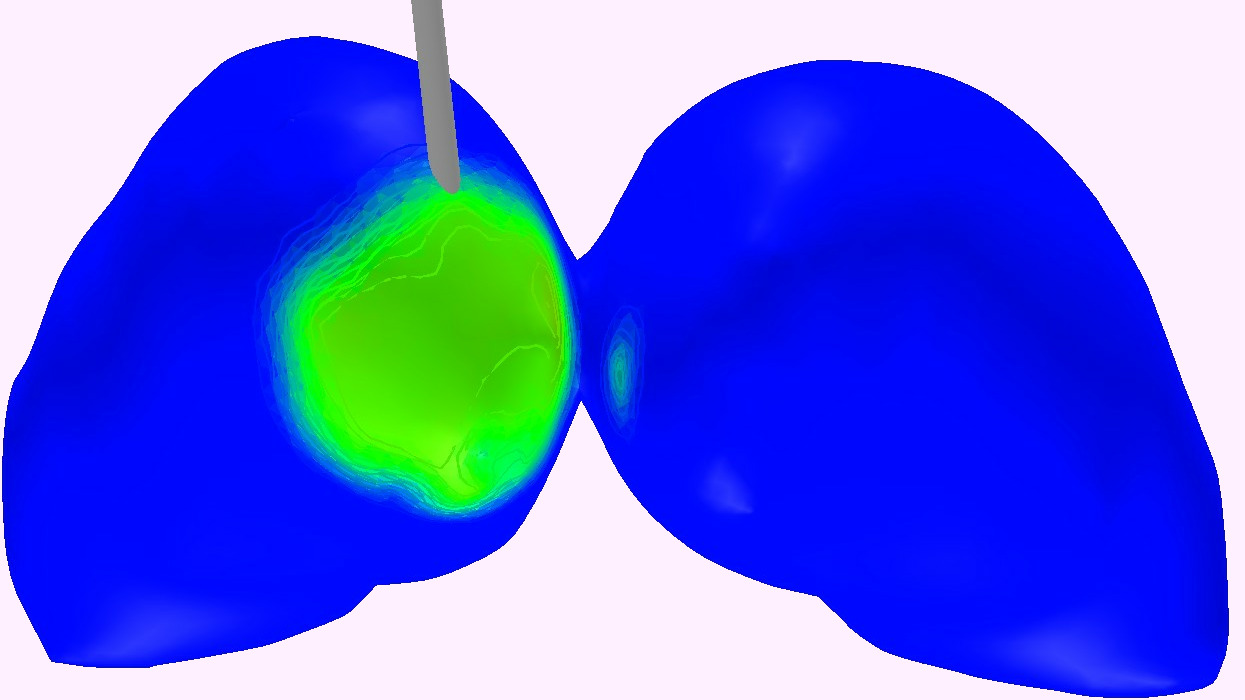} &
        \cellcolor[HTML]{fff0ff} \includegraphics[width=2.25cm, height=1.5cm]{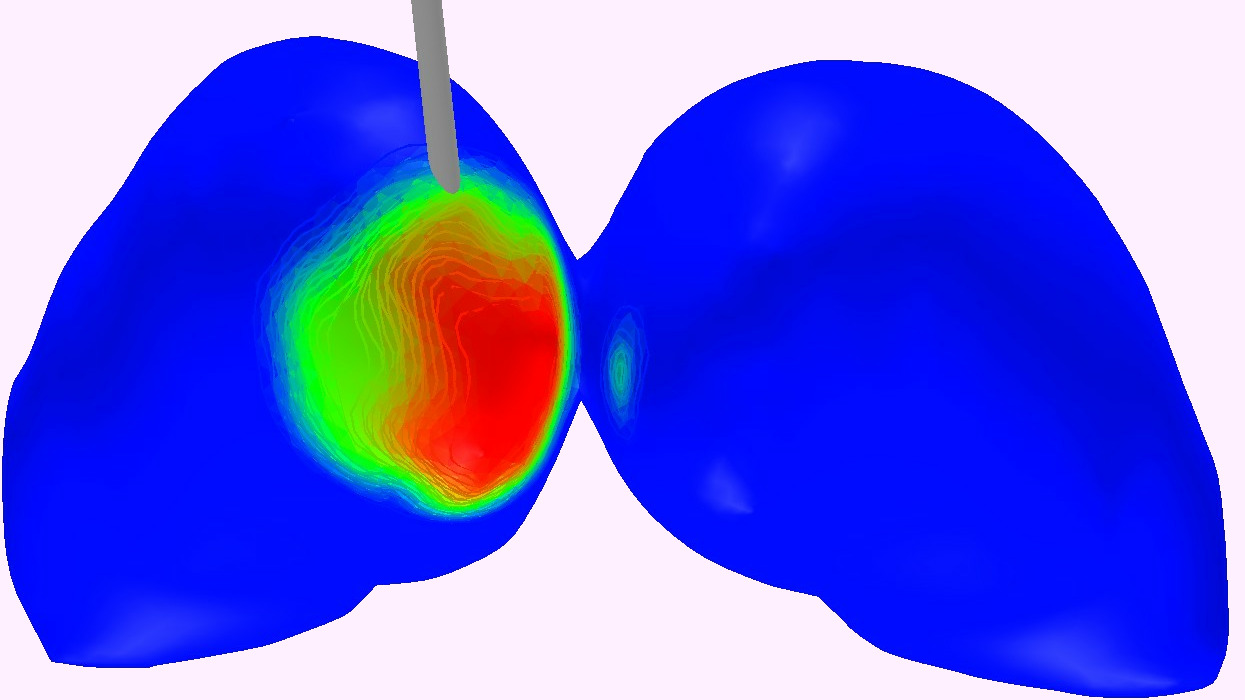} \\
        \hline
        
        \cellcolor[HTML]{ccffff} & 
        \rule{0pt}{1.80cm}
        \cellcolor[HTML]{ccffff} \includegraphics[width=2.25cm, height=1.5cm]{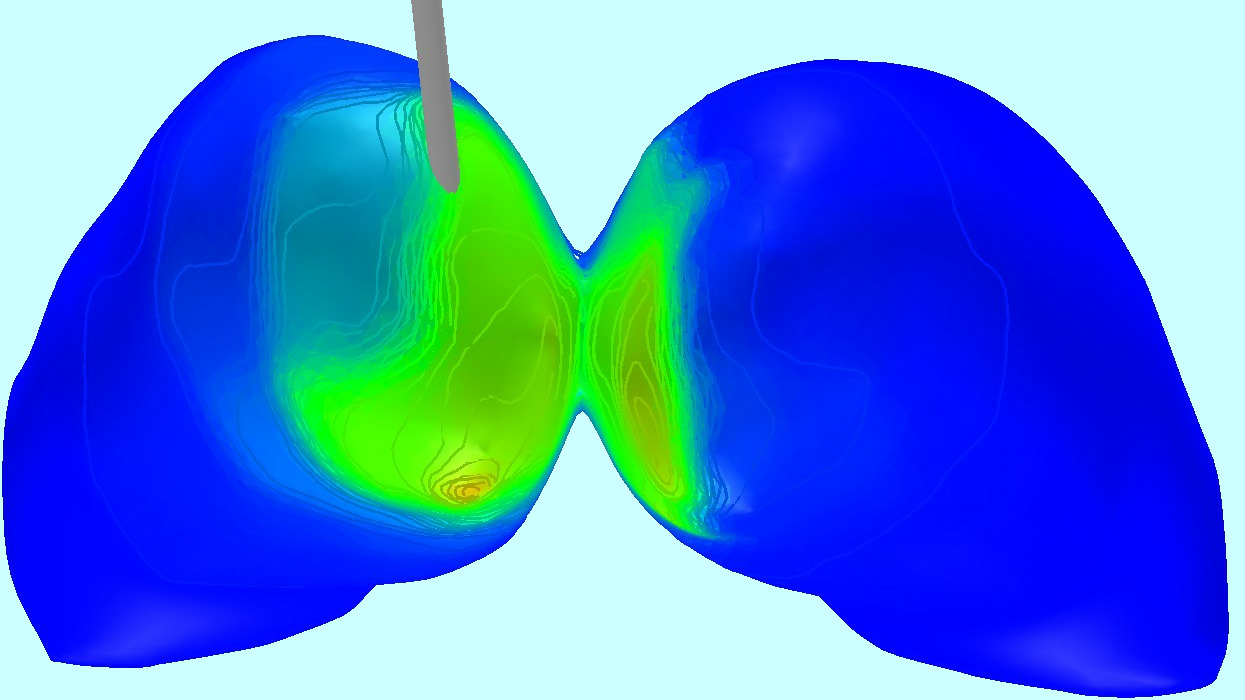} &  
        \cellcolor[HTML]{ccffff} \includegraphics[width=2.25cm, height=1.5cm]{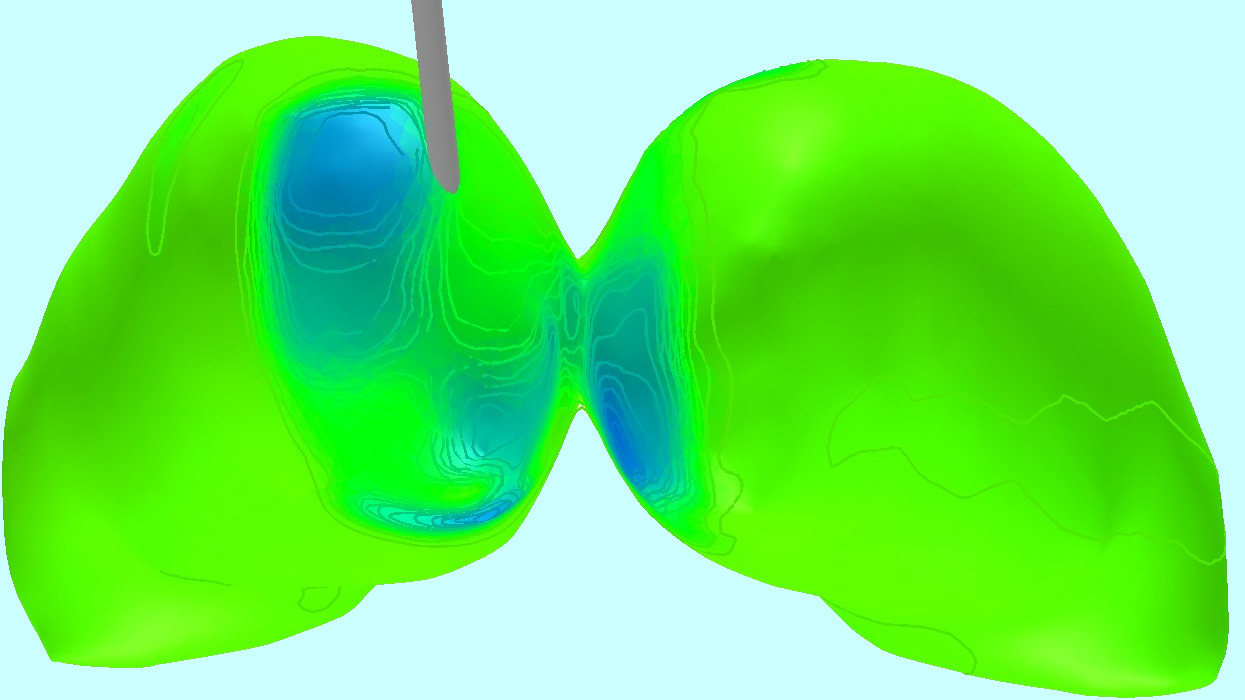} &
        \cellcolor[HTML]{ccffff} \includegraphics[width=2.25cm, height=1.5cm]{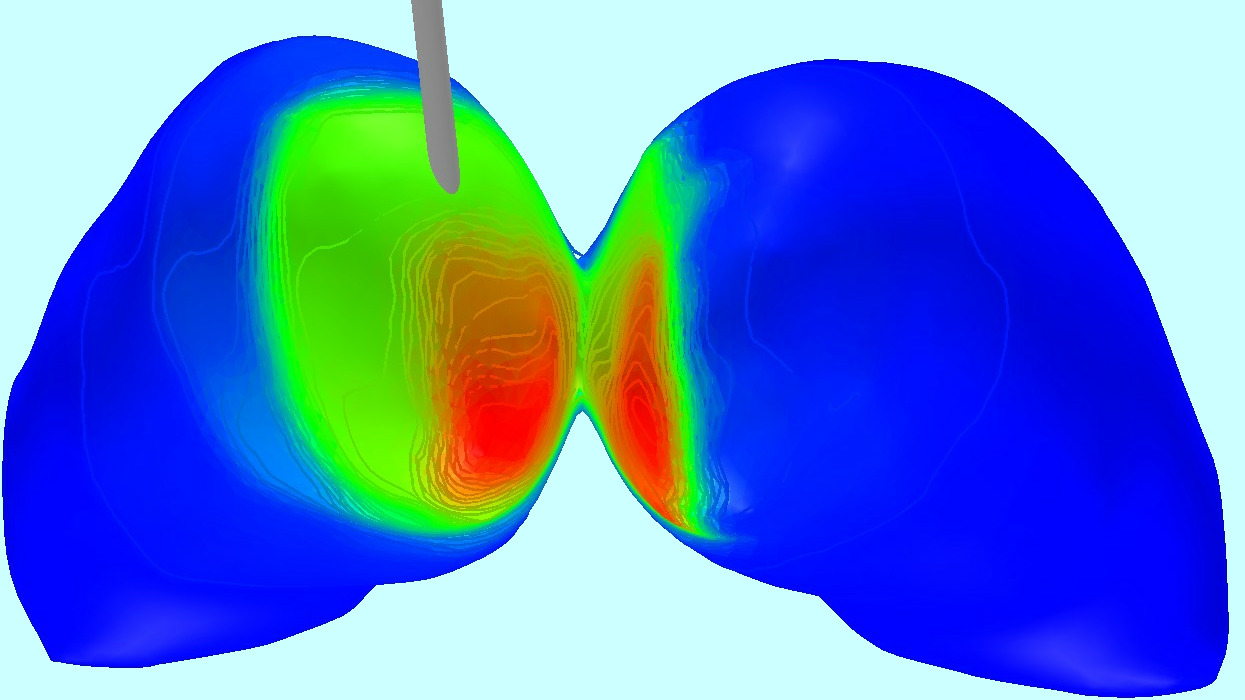} & 
        \cellcolor[HTML]{ccffff} \includegraphics[width=2.25cm, height=1.5cm]{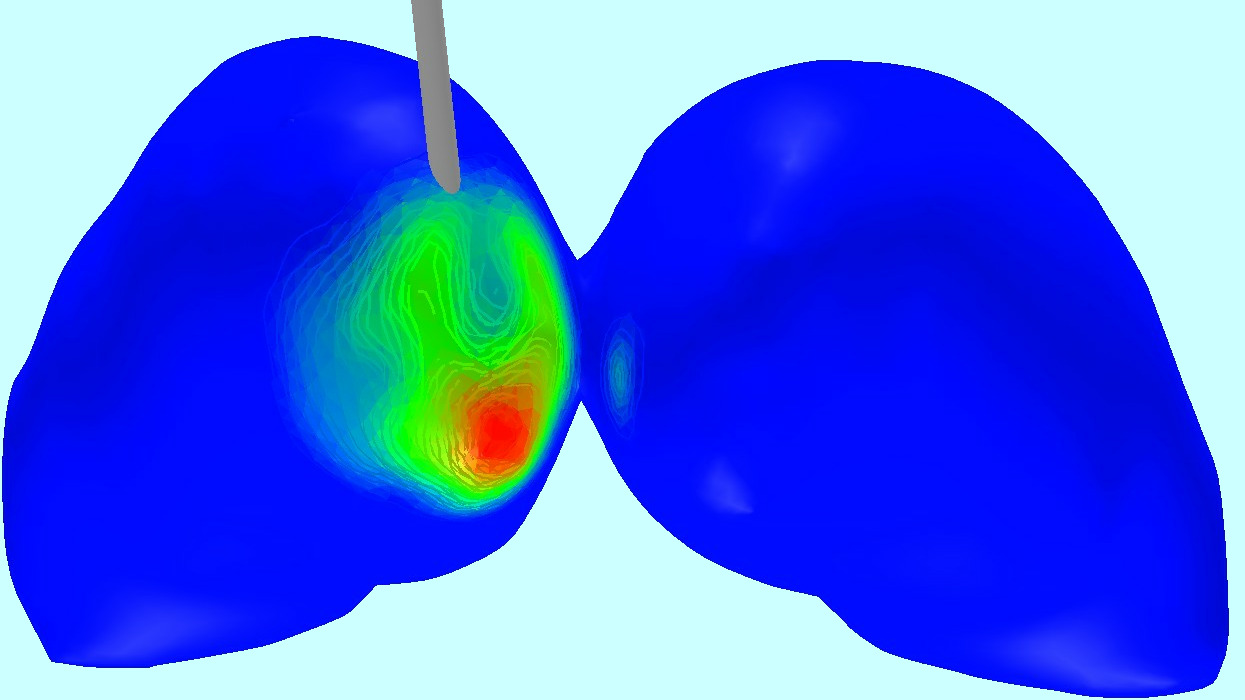} & 
        \cellcolor[HTML]{ccffff} \includegraphics[width=2.25cm, height=1.5cm]{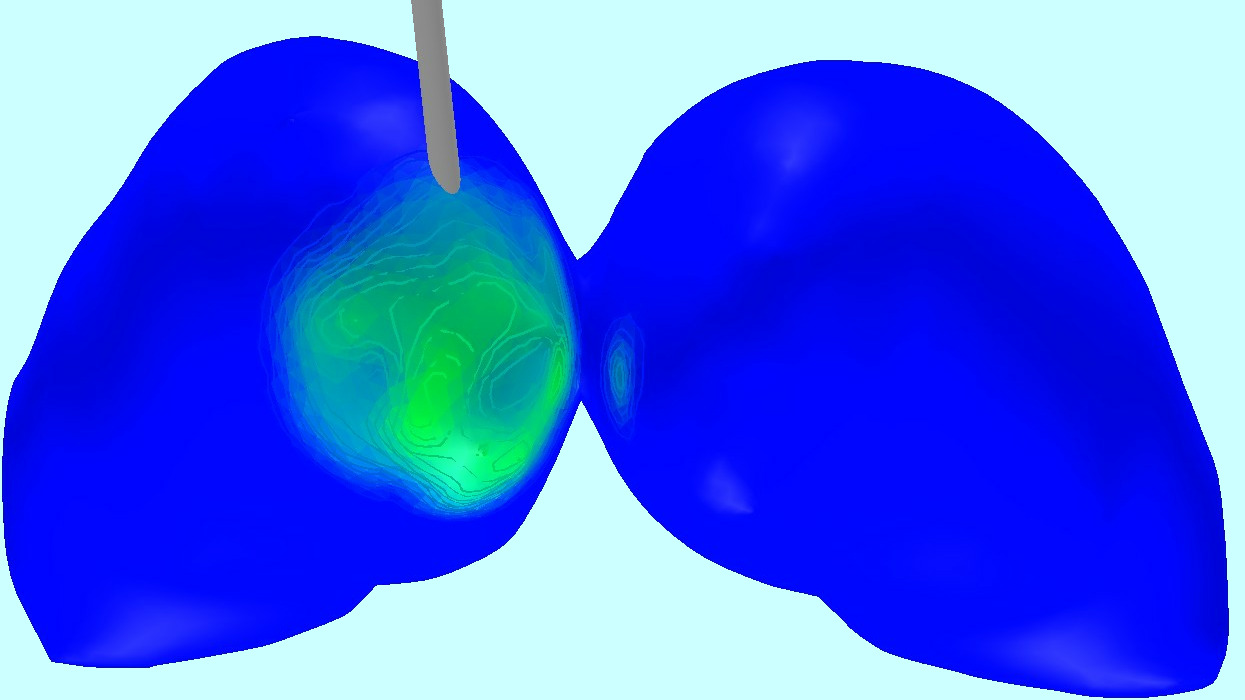} &
        \cellcolor[HTML]{ccffff} \includegraphics[width=2.25cm, height=1.5cm]{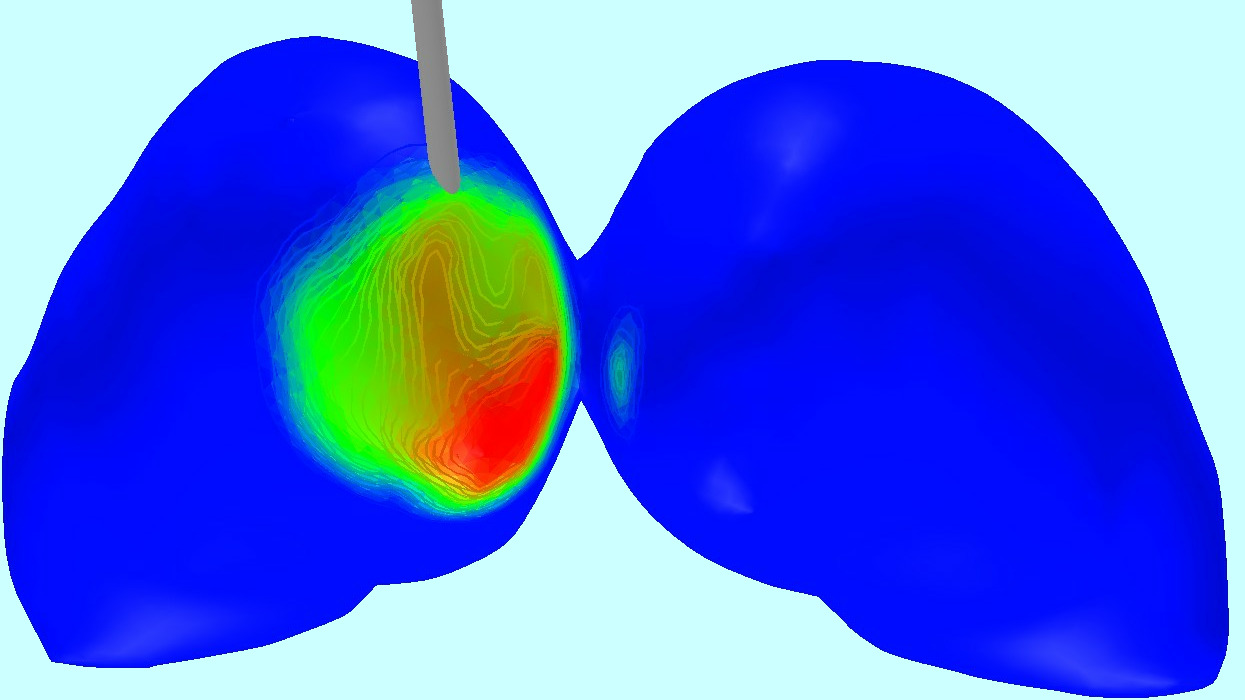} \\
        
        \cellcolor[HTML]{ccffff} \raisebox{3.15cm}[0pt][0pt]{\multirow{2}{*}{\rotatebox{90}{\textbf{L1L1(B)},  $\varepsilon \in  [-10, 0] $ dB}}} \hskip0.1cm        \raisebox{2.9cm}[0pt][0pt] {\multirow{2}{*}{\rotatebox{90}{Perpendicular \hskip0.25cm Parallel}}}&
        \rule{0pt}{1.80cm}
        \cellcolor[HTML]{ccffff} \includegraphics[width=2.25cm, height=1.5cm]{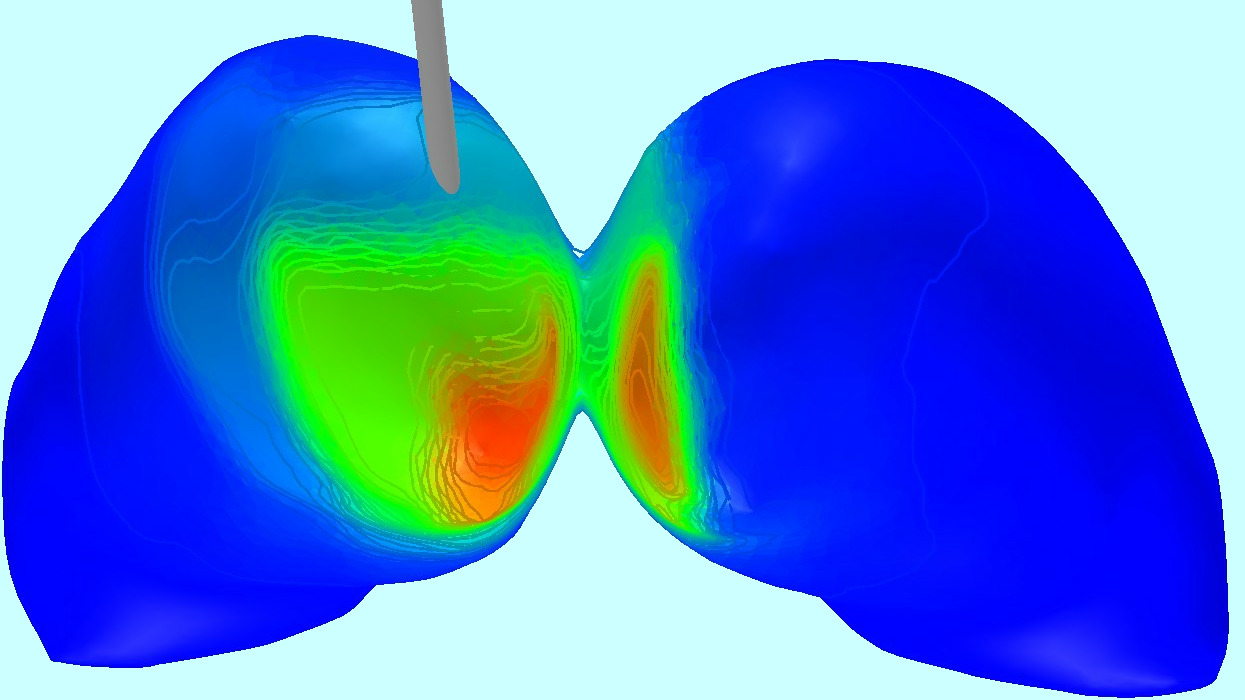} &  
        \cellcolor[HTML]{ccffff} \includegraphics[width=2.25cm, height=1.5cm]{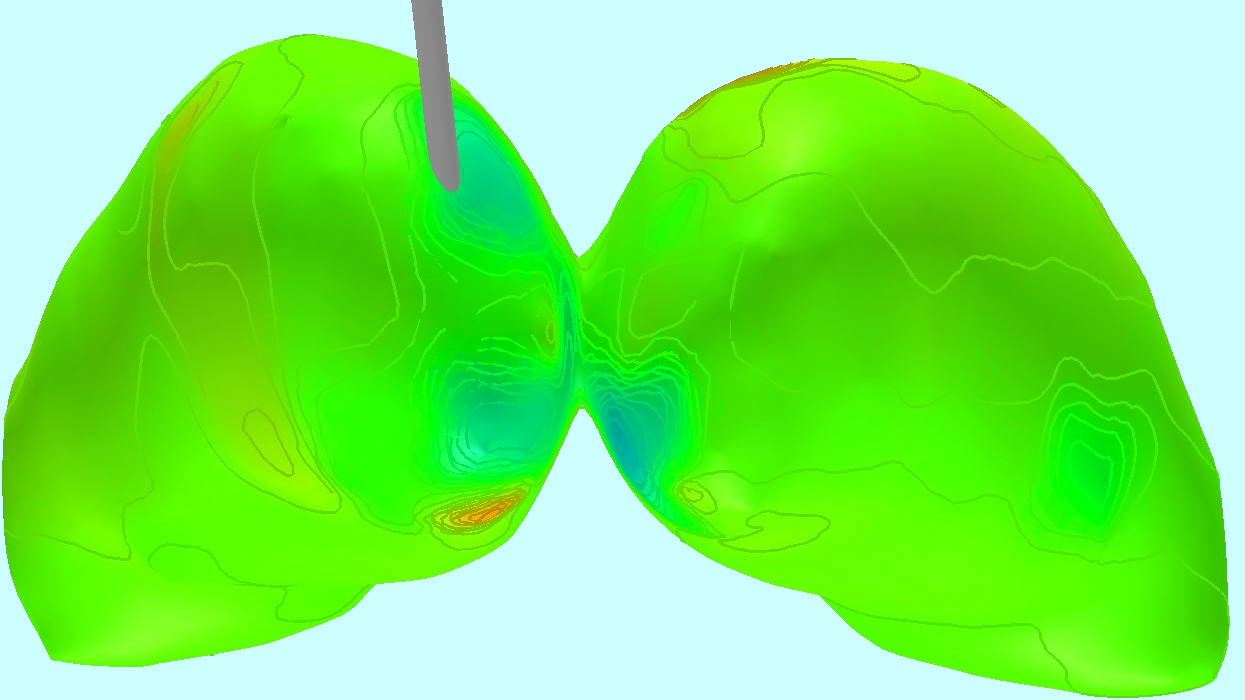} &
        \cellcolor[HTML]{ccffff} \includegraphics[width=2.25cm, height=1.5cm]{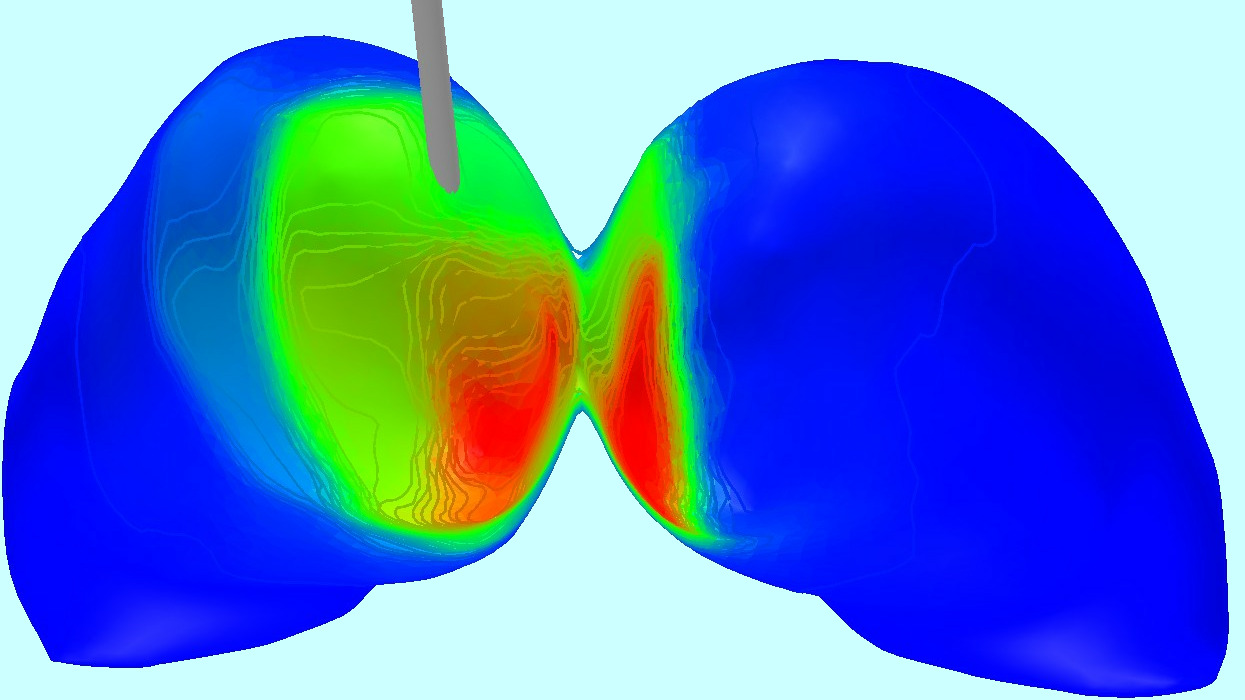} & 
        \cellcolor[HTML]{ccffff} \includegraphics[width=2.25cm, height=1.5cm]{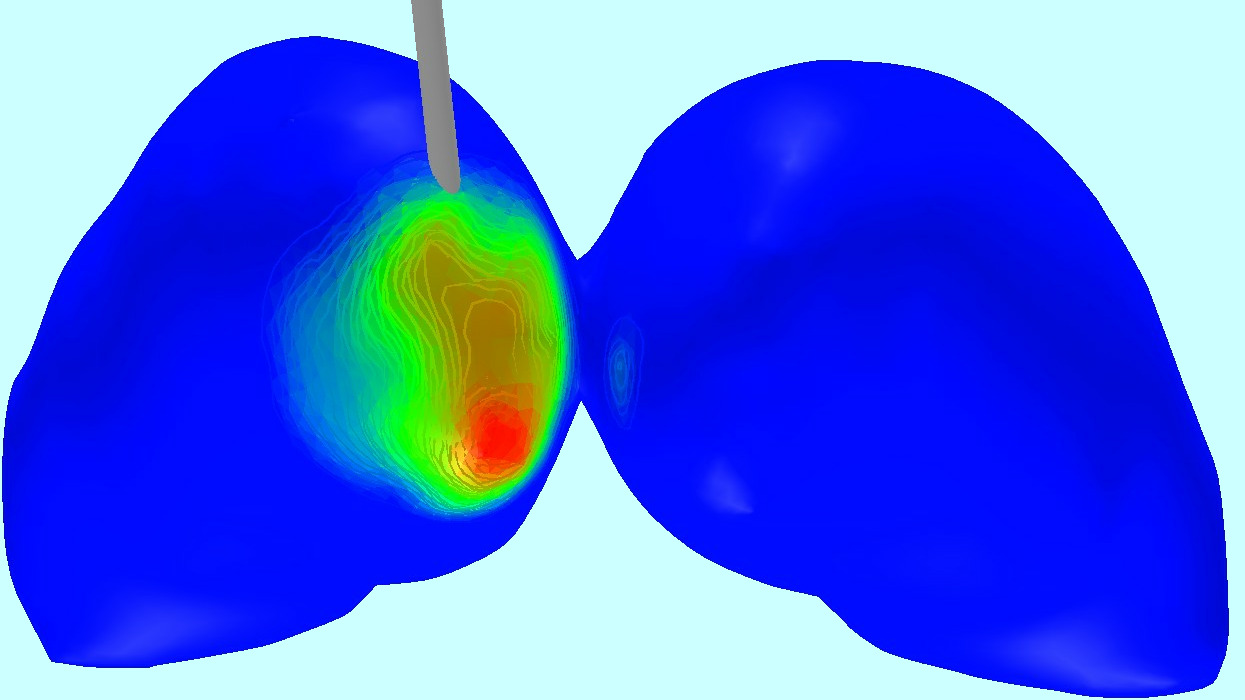} & 
        \cellcolor[HTML]{ccffff} \includegraphics[width=2.25cm, height=1.5cm]{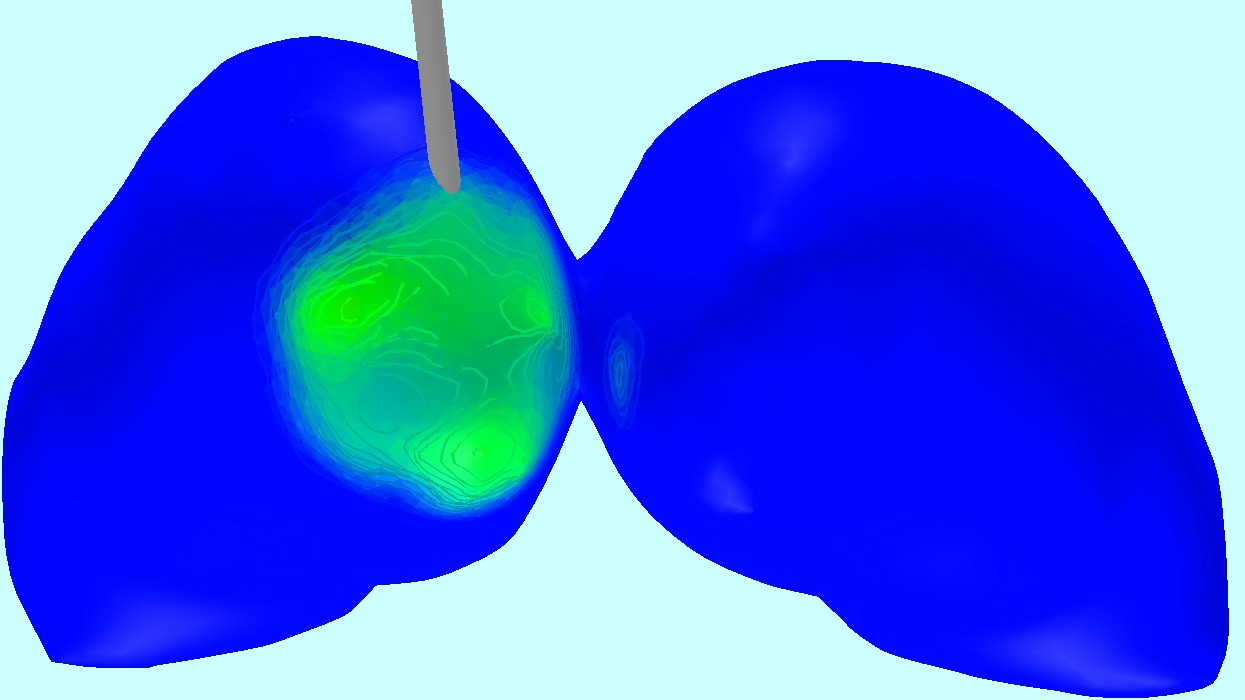} &
        \cellcolor[HTML]{ccffff} \includegraphics[width=2.25cm, height=1.5cm]{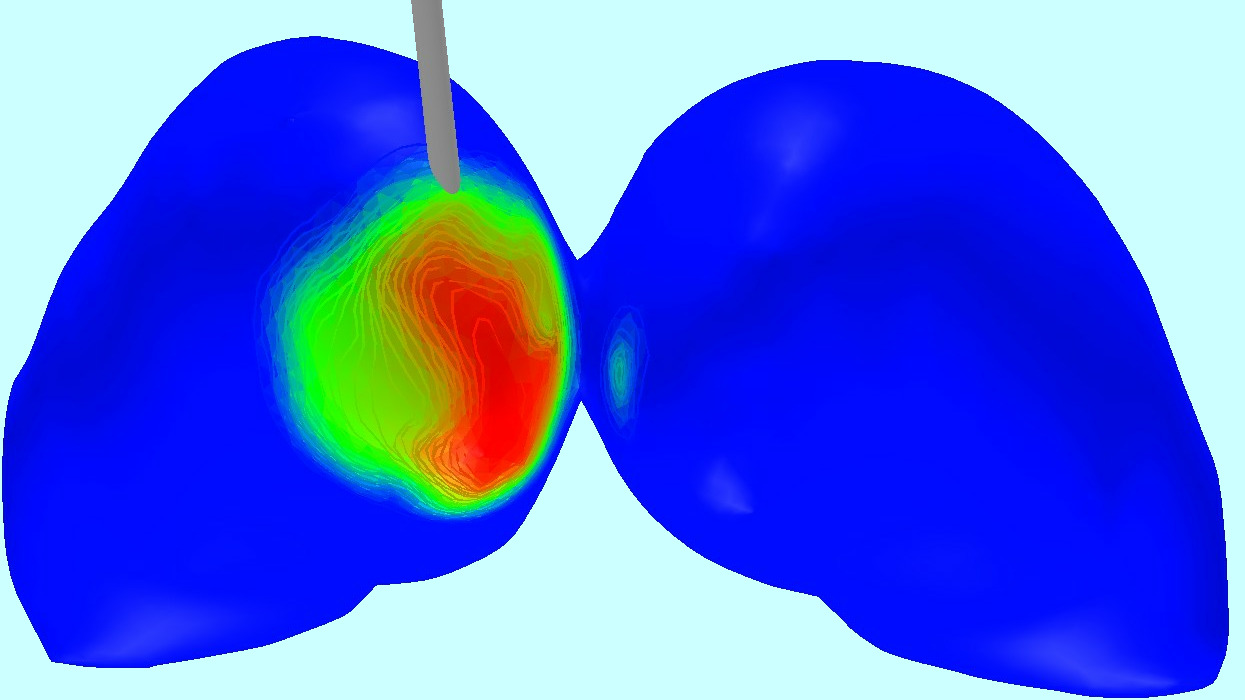} \\
        \hline
        
        \cellcolor[HTML]{ccffff} & 
        \rule{0pt}{1.80cm}
        \cellcolor[HTML]{ccffff} \includegraphics[width=2.25cm, height=1.5cm]{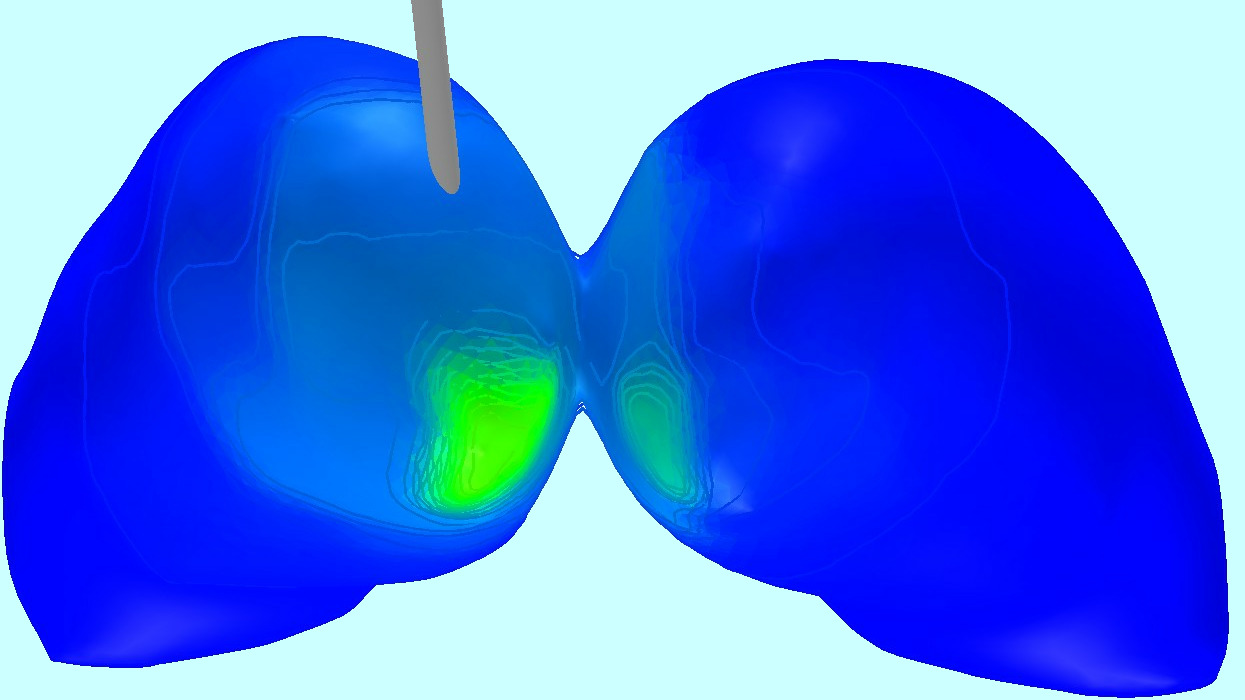} &  
        \cellcolor[HTML]{ccffff} \includegraphics[width=2.25cm, height=1.5cm]{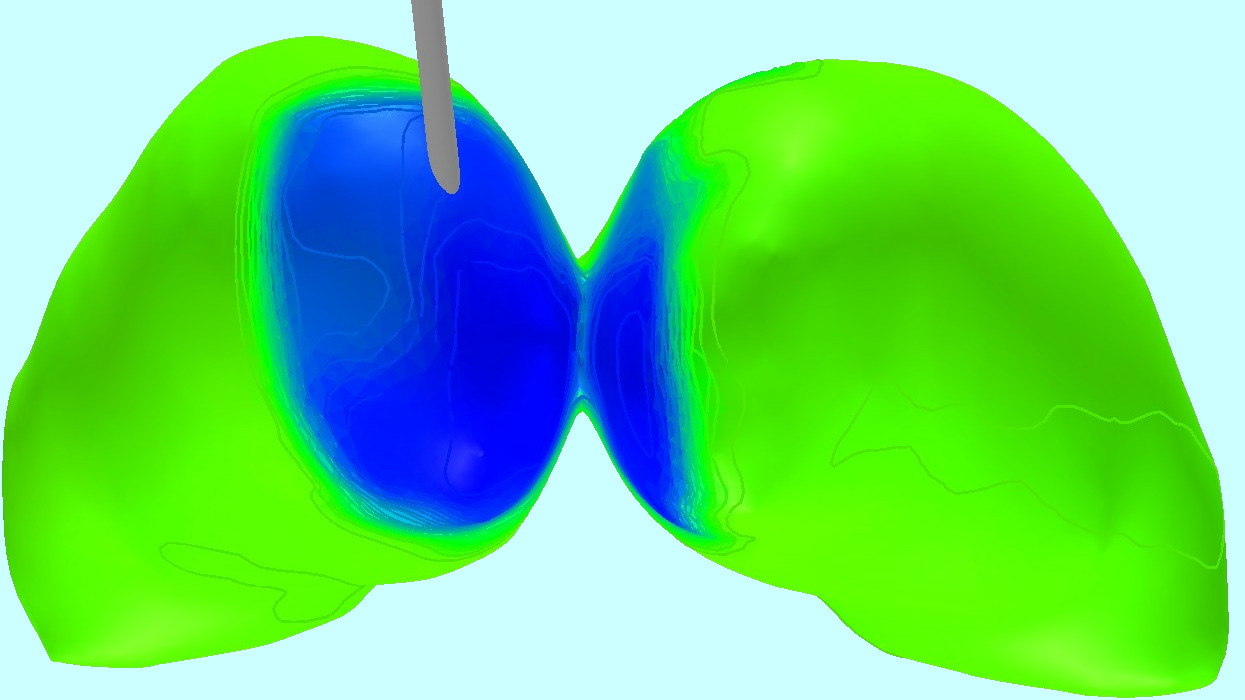} &
        \cellcolor[HTML]{ccffff} \includegraphics[width=2.25cm, height=1.5cm]{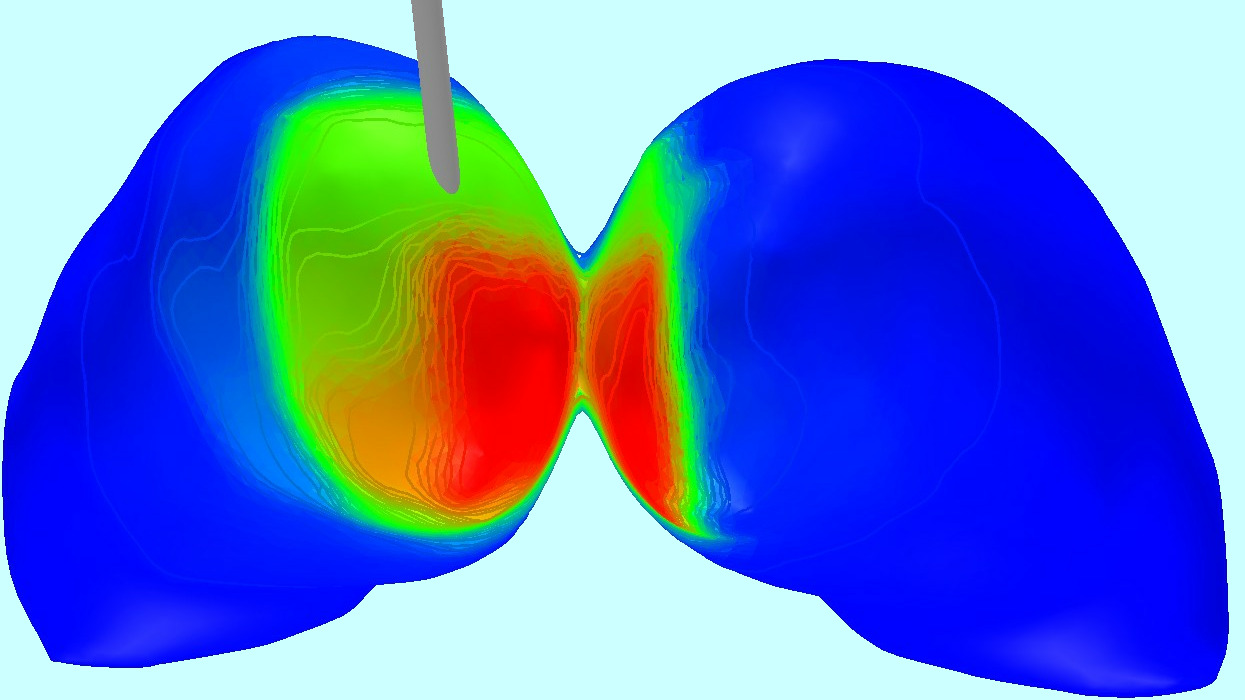} & 
        \cellcolor[HTML]{ccffff} \includegraphics[width=2.25cm, height=1.5cm]{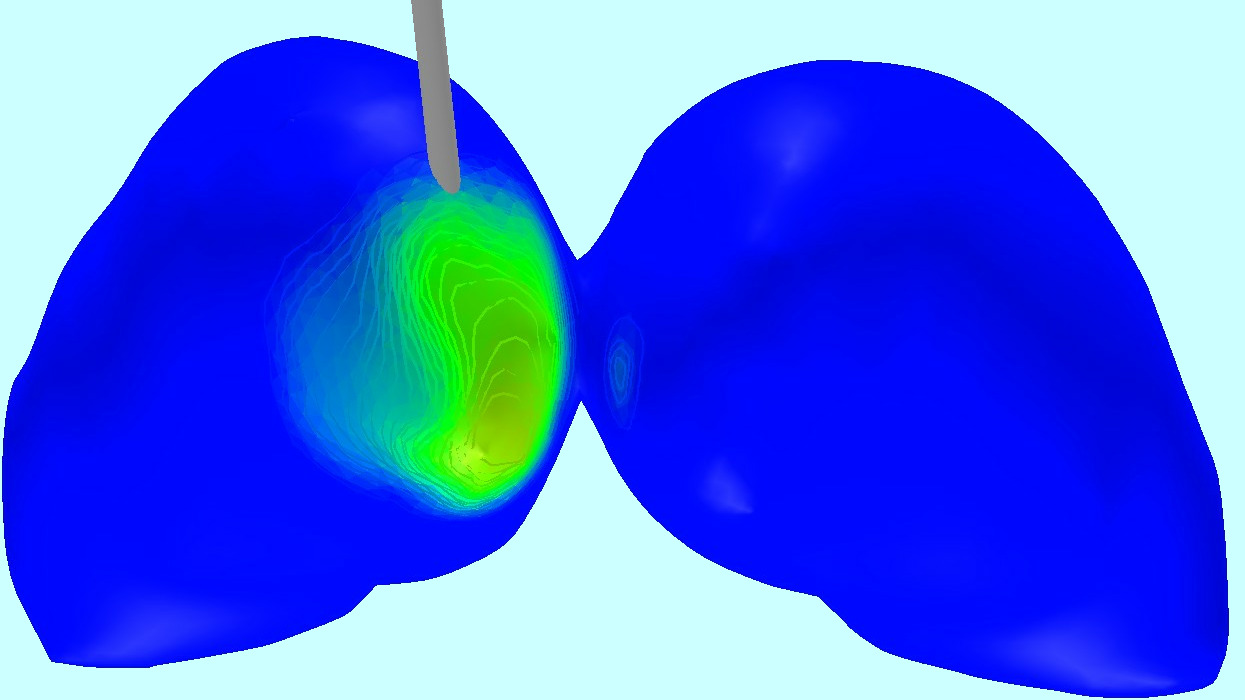} & 
        \cellcolor[HTML]{ccffff} \includegraphics[width=2.25cm, height=1.5cm]{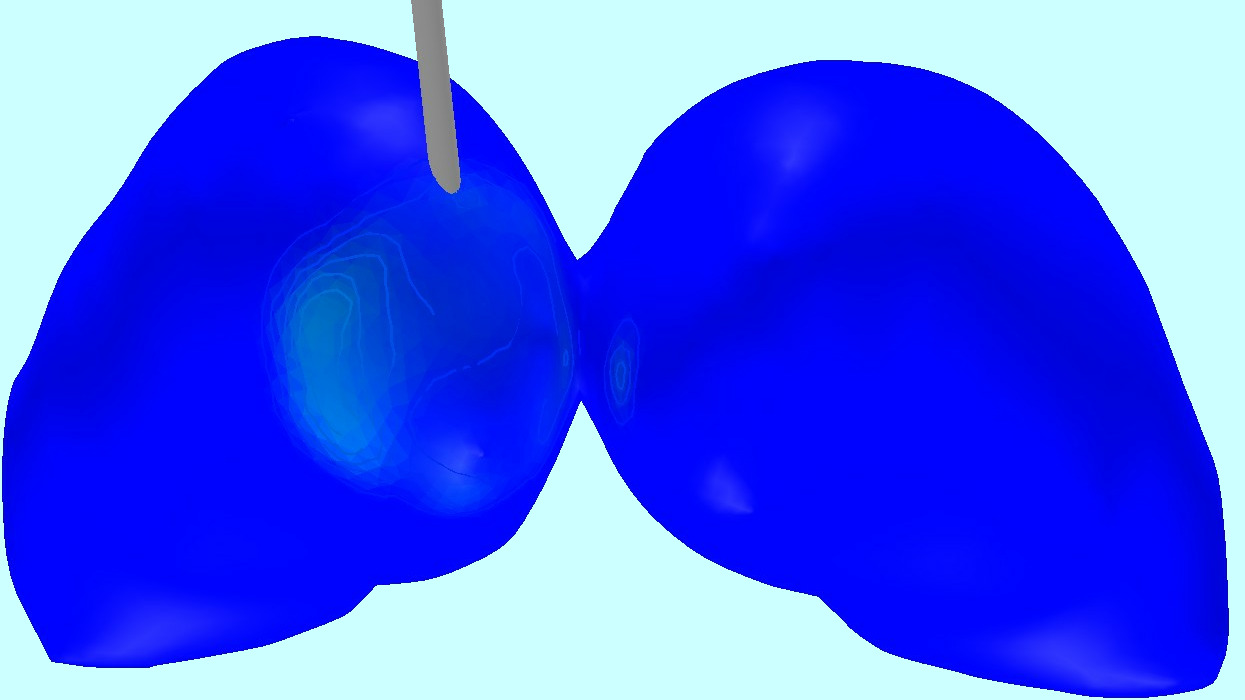} &
        \cellcolor[HTML]{ccffff} \includegraphics[width=2.25cm, height=1.5cm]{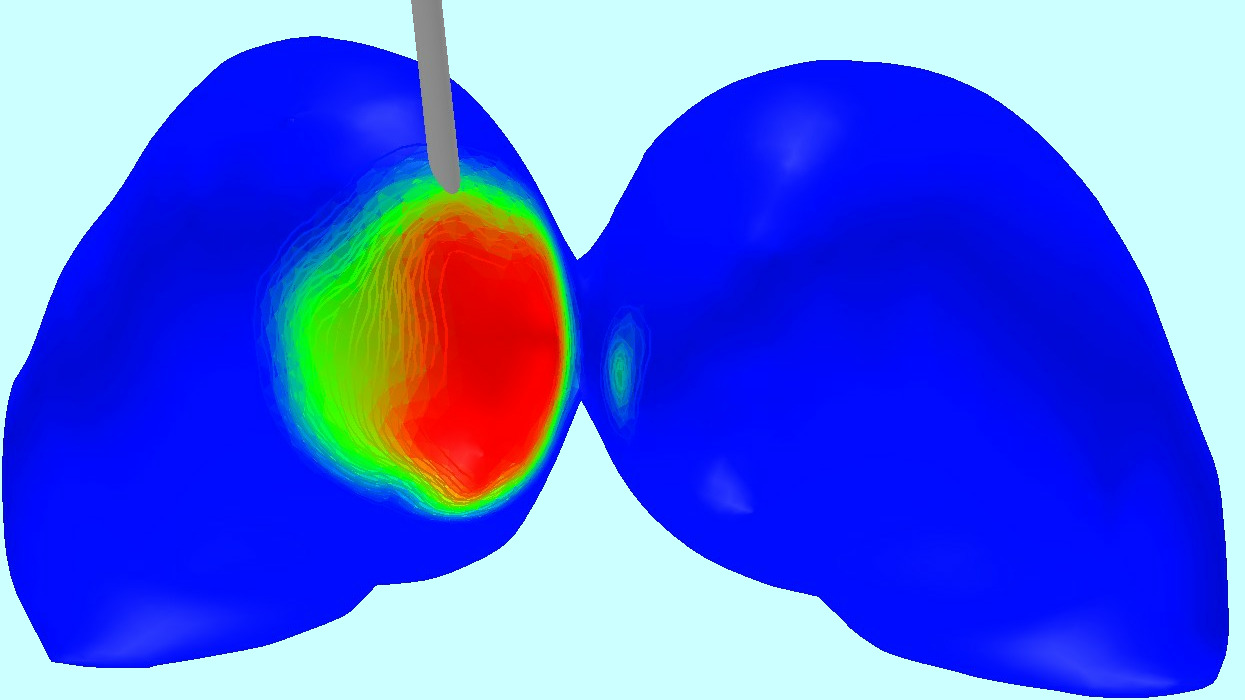} \\
       
        \cellcolor[HTML]{ccffff} \raisebox{3.25cm}[0pt][0pt]{\multirow{2}{*}{\rotatebox{90}{\textbf{L1L1(A)},  $\varepsilon \in  [-160, 0] $ dB}}}  \hskip0.1cm \raisebox{2.9cm}[0pt][0pt] {\multirow{2}{*}{\rotatebox{90}{Perpendicular \hskip0.25cm Parallel}}} &
        \rule{0pt}{1.80cm}
        \cellcolor[HTML]{ccffff} \includegraphics[width=2.25cm, height=1.5cm]{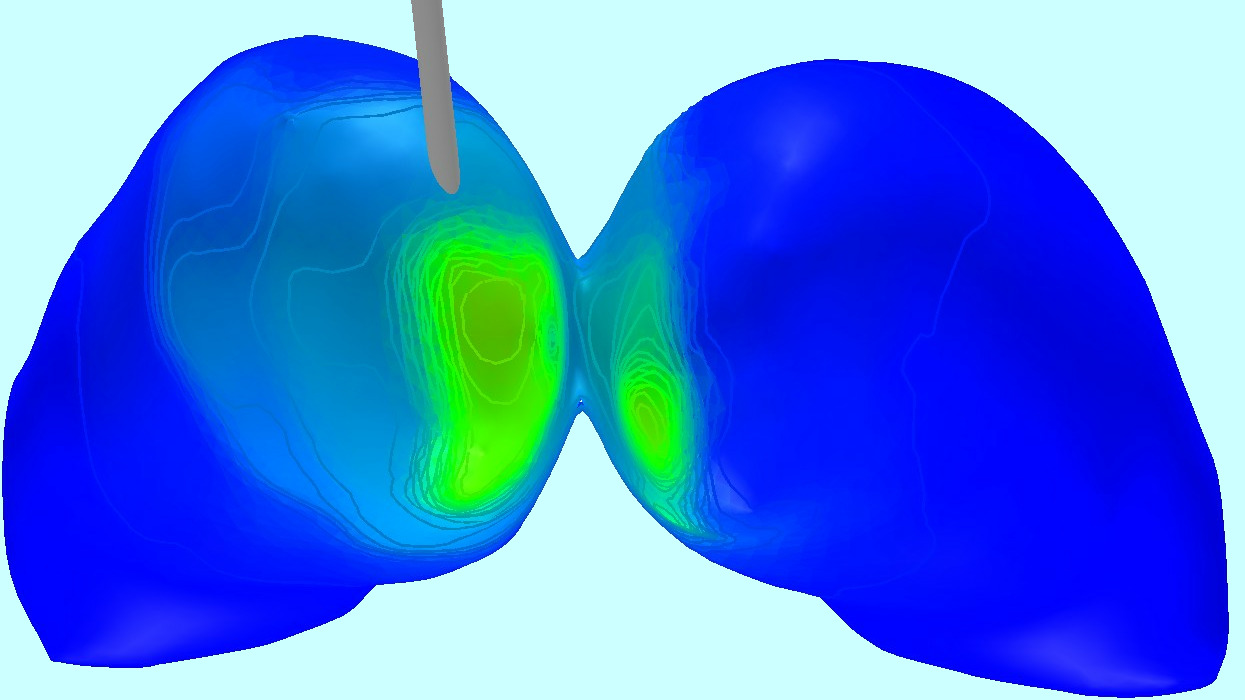} &  
        \cellcolor[HTML]{ccffff} \includegraphics[width=2.25cm, height=1.5cm]{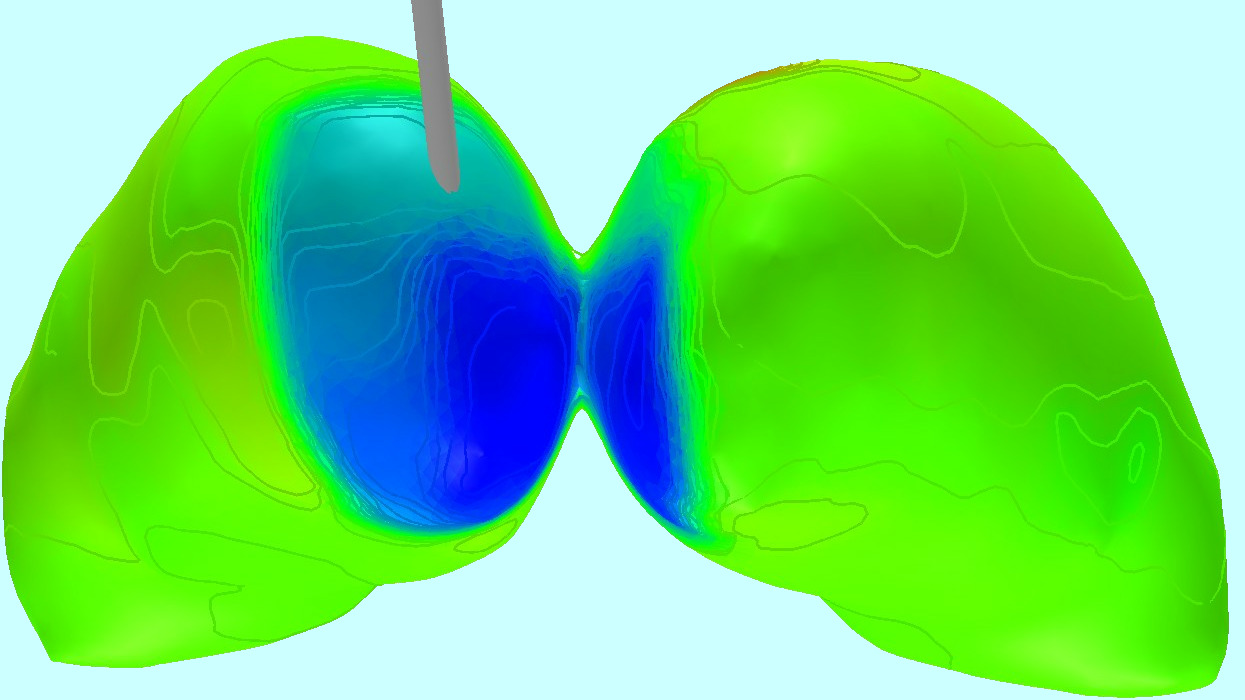} &
        \cellcolor[HTML]{ccffff} \includegraphics[width=2.25cm, height=1.5cm]{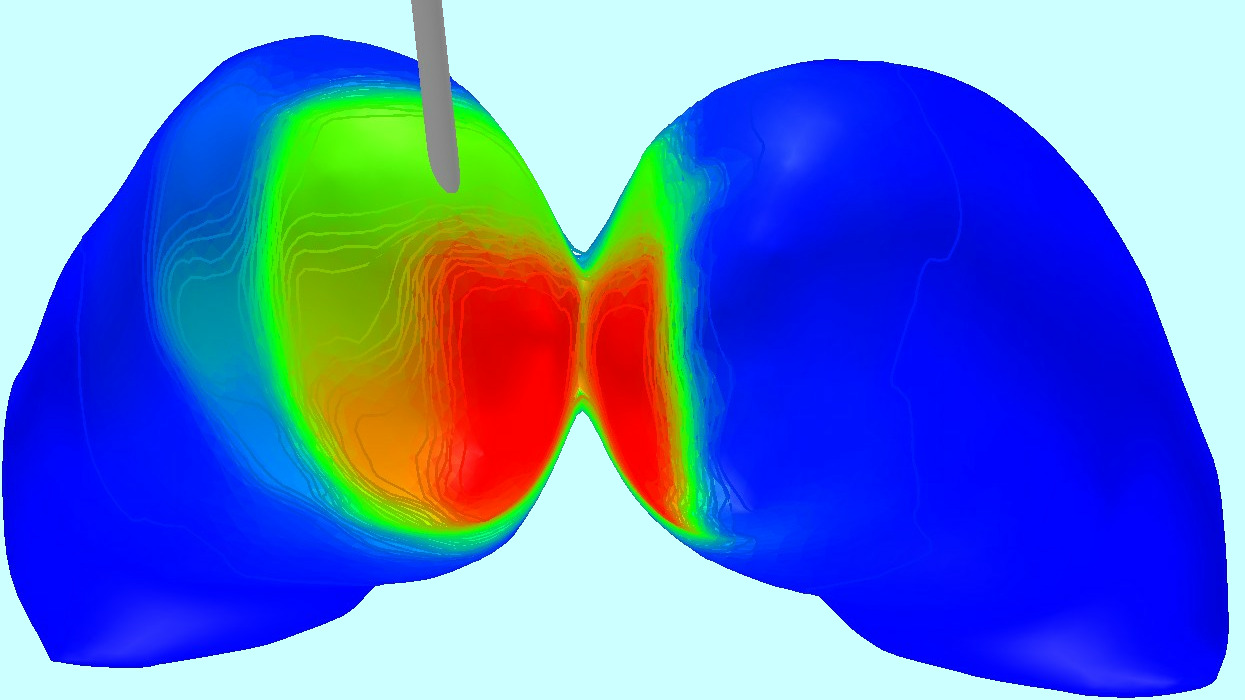} & 
        \cellcolor[HTML]{ccffff} \includegraphics[width=2.25cm, height=1.5cm]{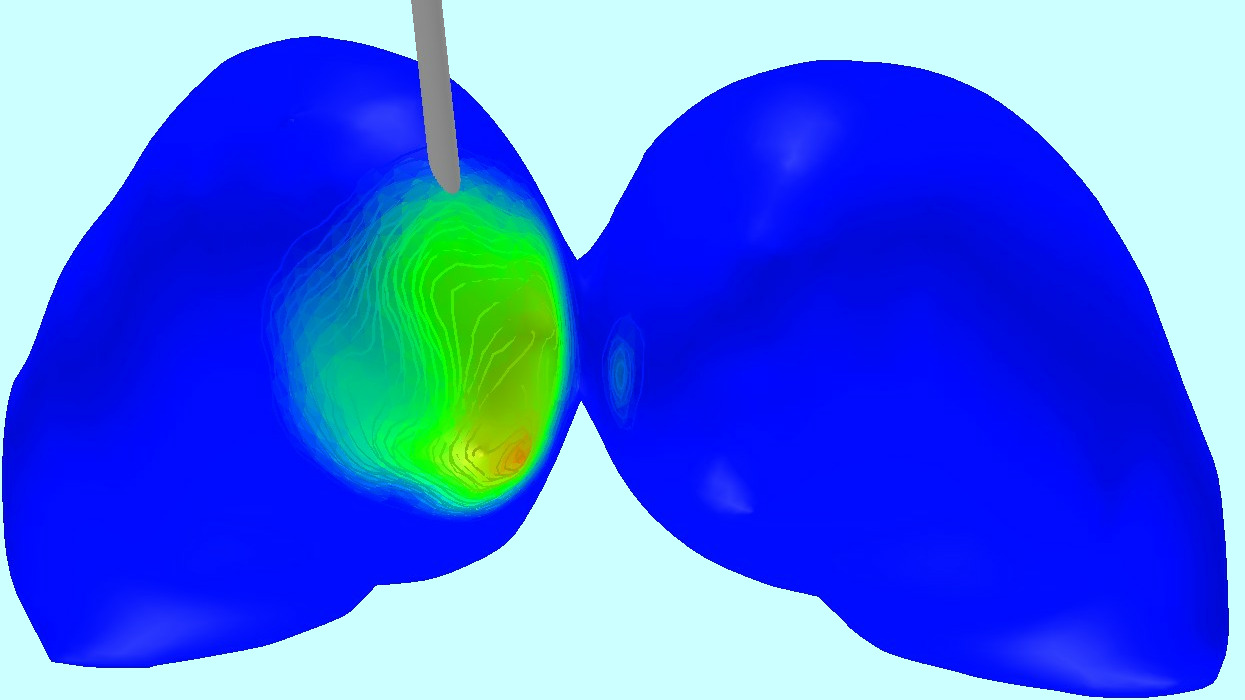} & 
        \cellcolor[HTML]{ccffff} \includegraphics[width=2.25cm, height=1.5cm]{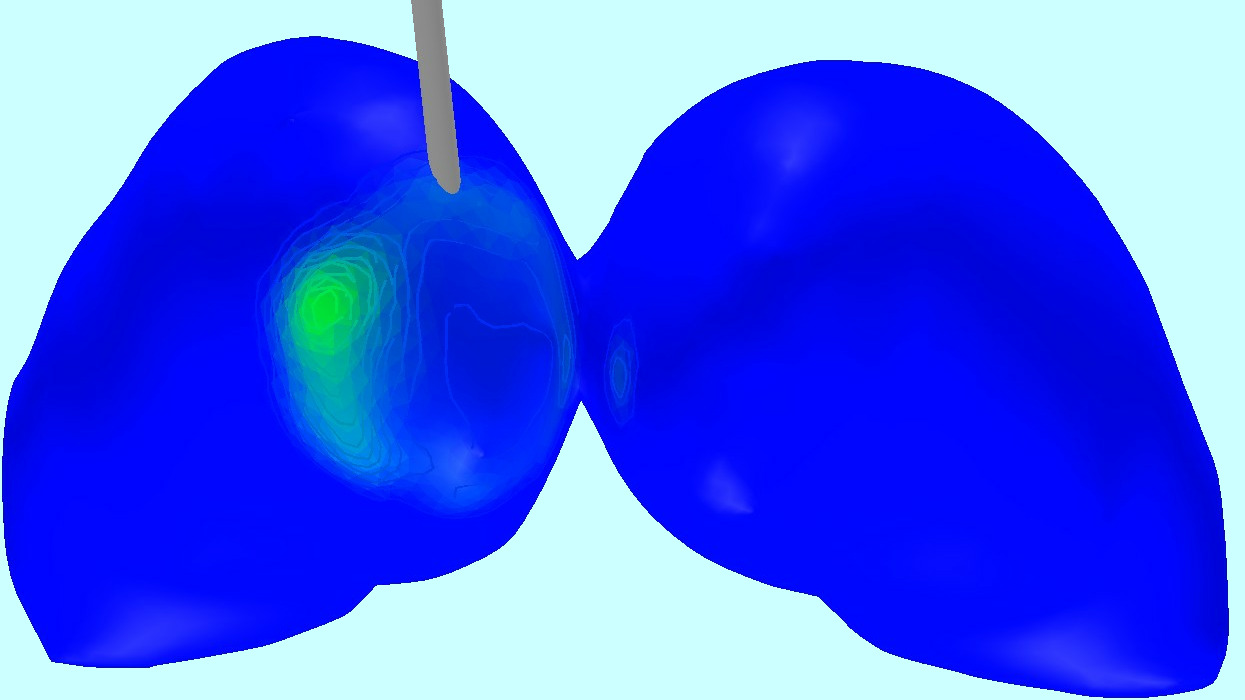} &
        \cellcolor[HTML]{ccffff} \includegraphics[width=2.25cm, height=1.5cm]{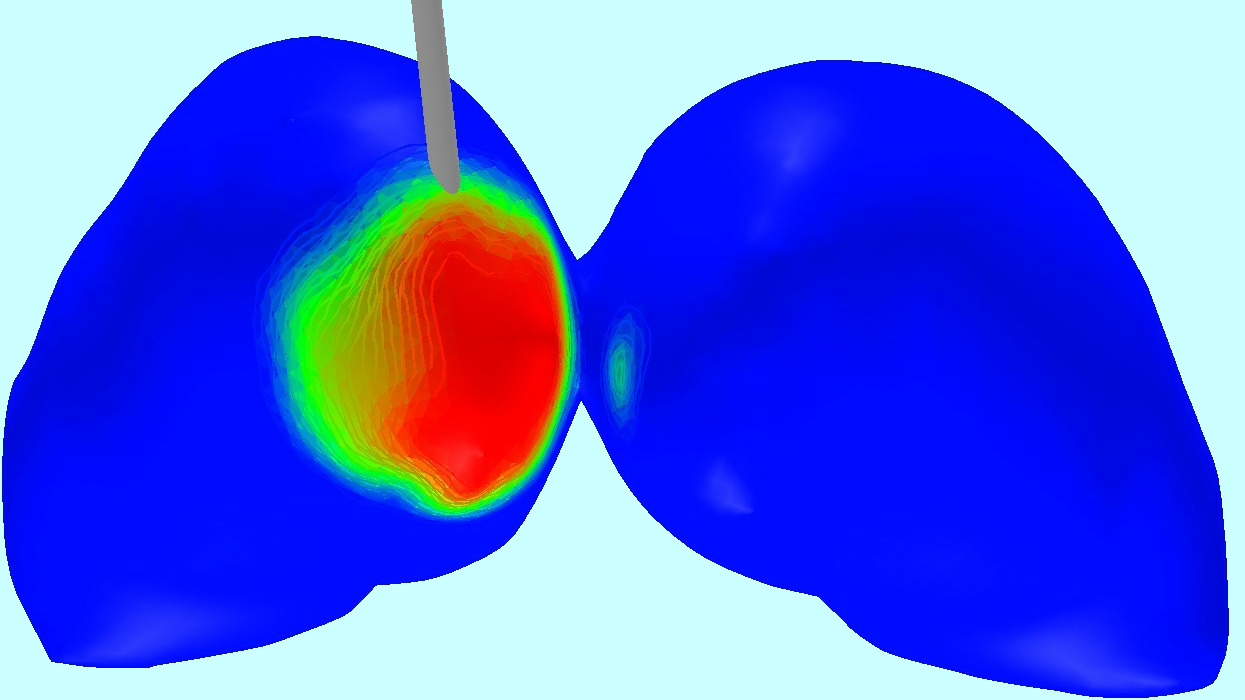} \\ 
        \hline
        
        \cellcolor[HTML]{e5ffd4}  &
        \rule{0pt}{1.80cm}
        \cellcolor[HTML]{e5ffd4} \includegraphics[width=2.25cm, height=1.5cm]{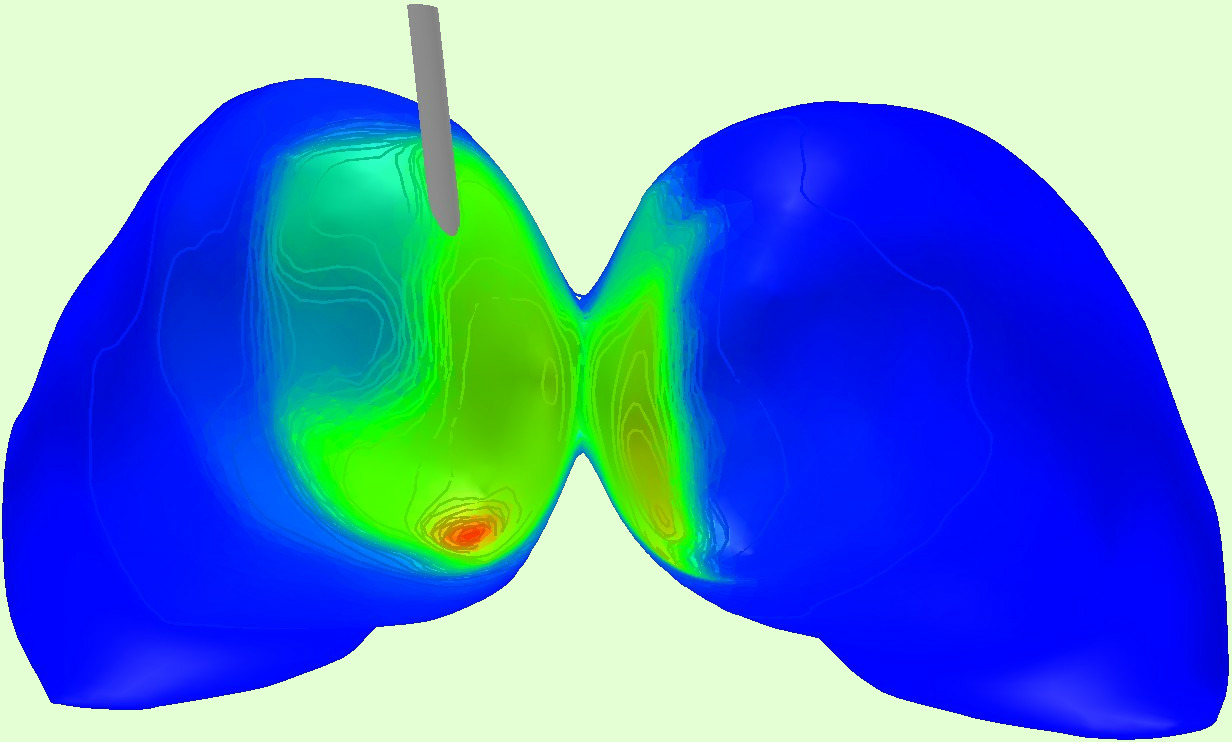} & 
        \cellcolor[HTML]{e5ffd4} \includegraphics[width=2.25cm, height=1.5cm]{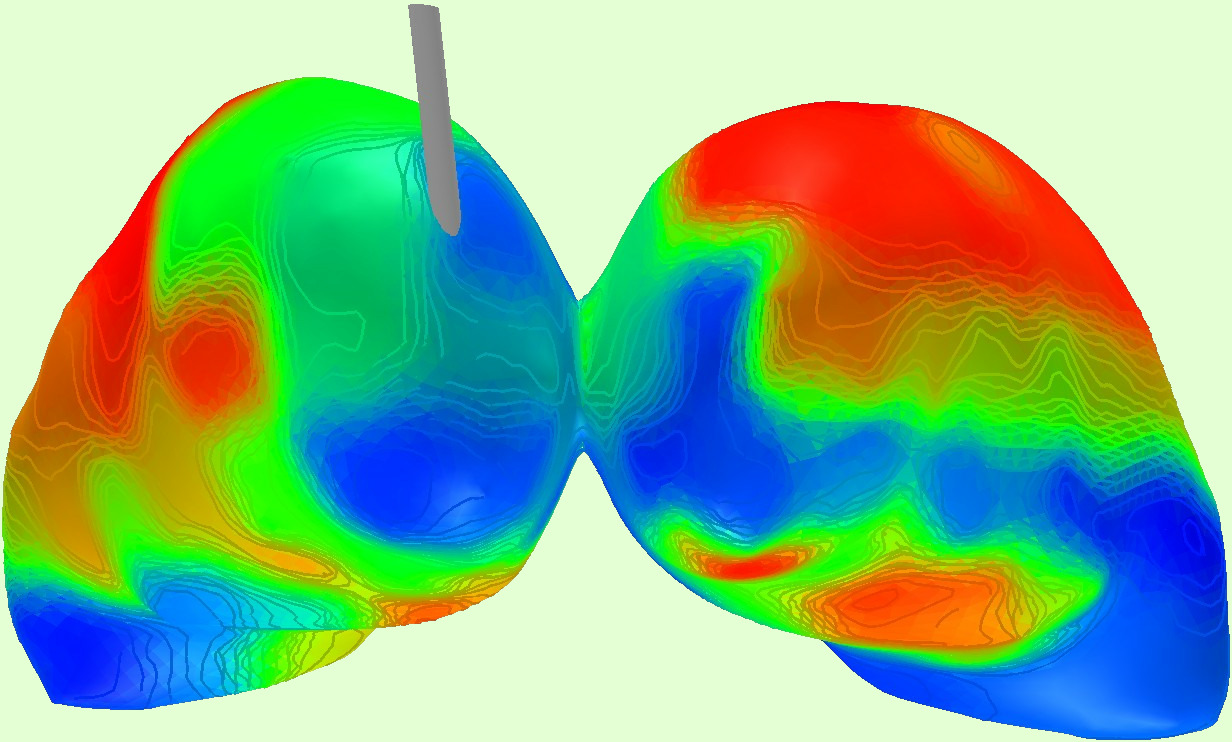} & 
        \cellcolor[HTML]{e5ffd4} \includegraphics[width=2.25cm, height=1.5cm]{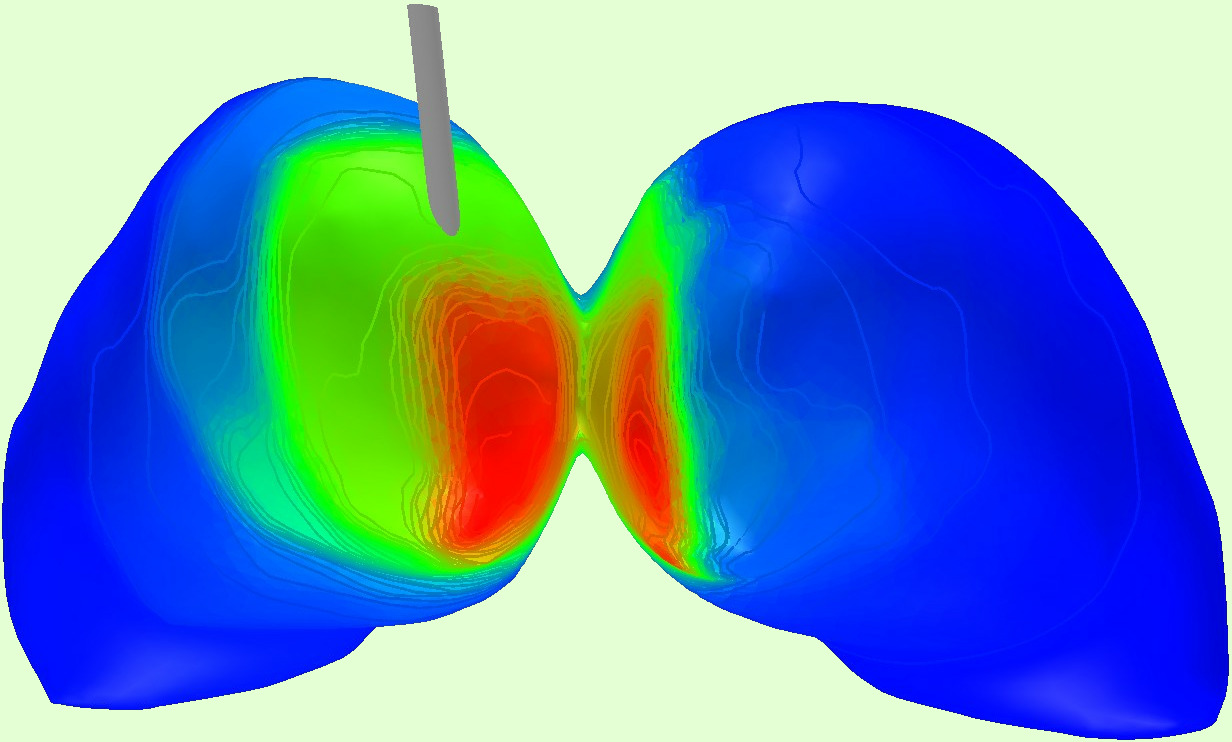} & 
        \cellcolor[HTML]{e5ffd4} \includegraphics[width=2.25cm, height=1.5cm]{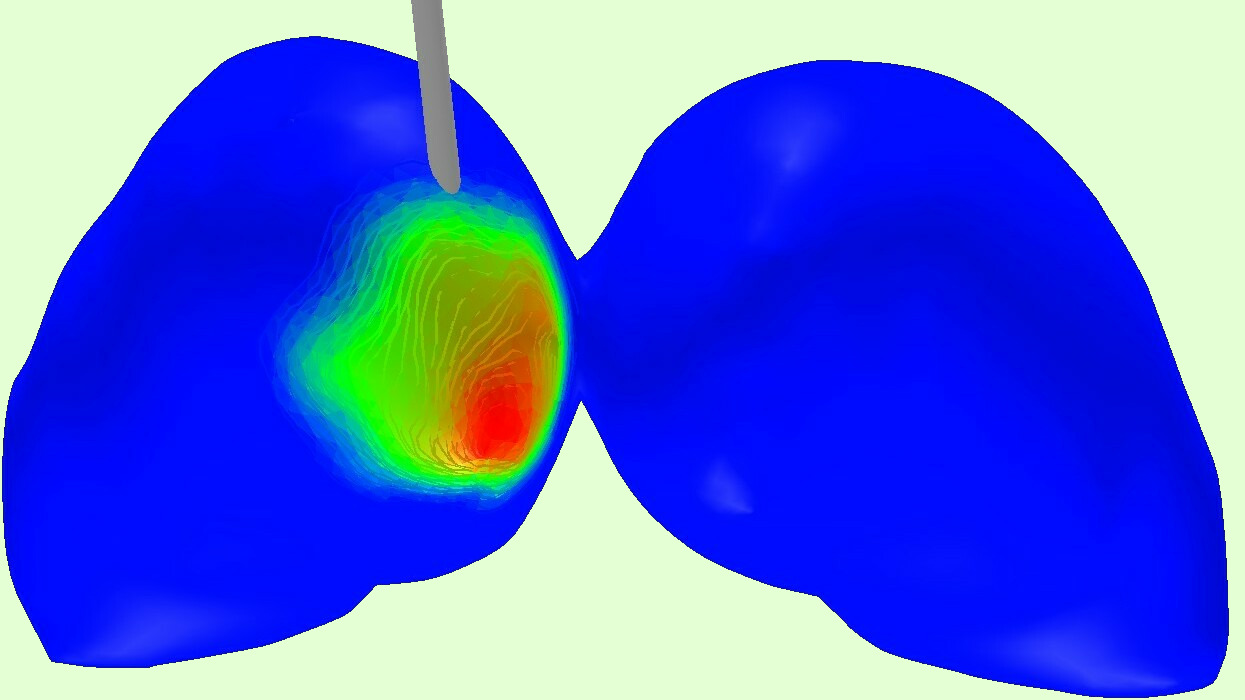} & 
        \cellcolor[HTML]{e5ffd4} \includegraphics[width=2.25cm, height=1.5cm]{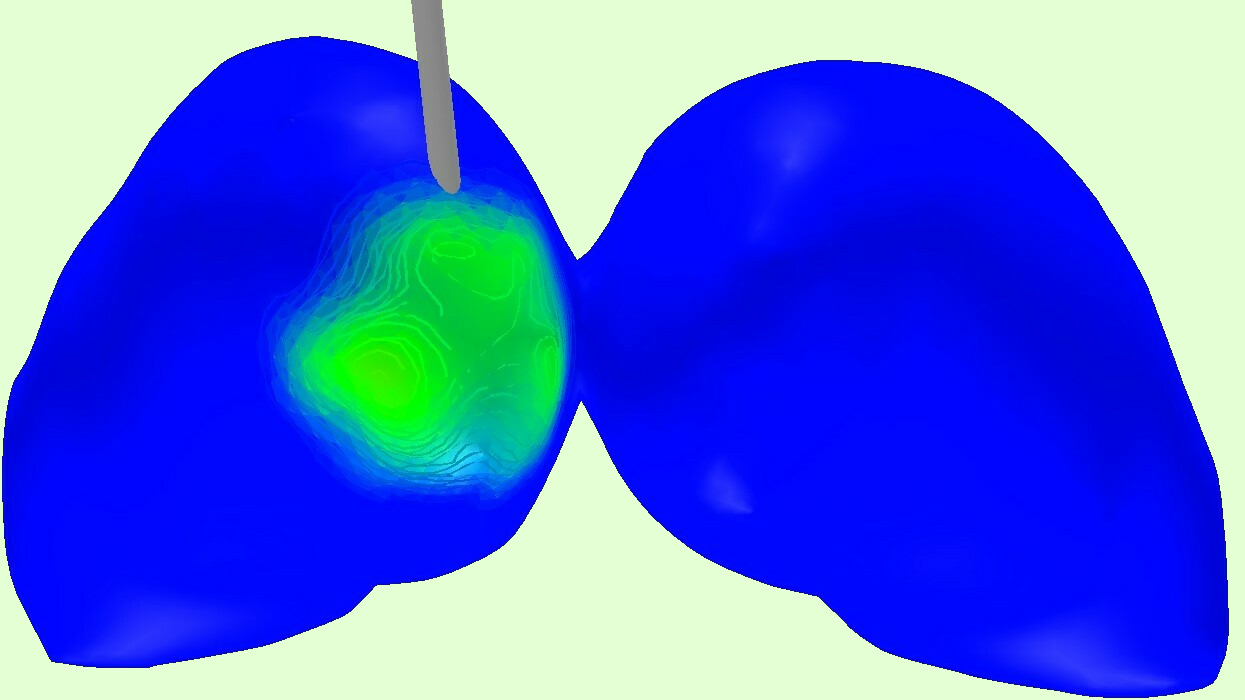} & 
        \cellcolor[HTML]{e5ffd4} \includegraphics[width=2.25cm, height=1.5cm]{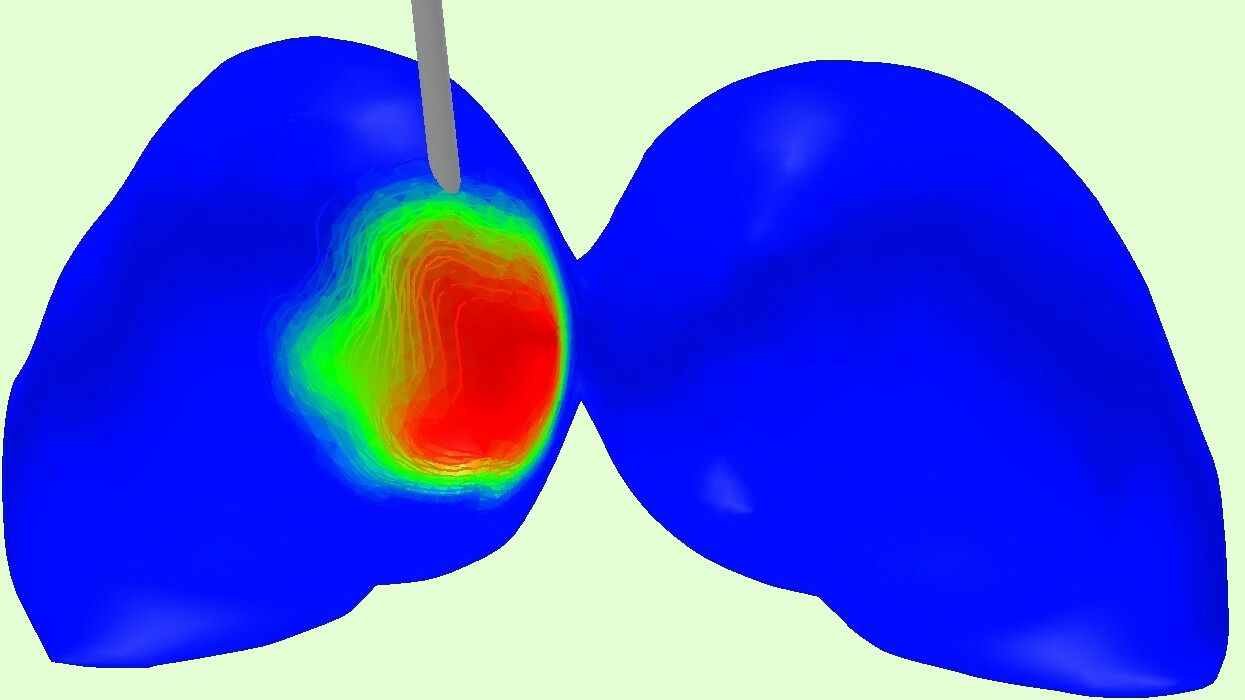} \\ 
        \cellcolor[HTML]{e5ffd4} \raisebox{2.9cm}[0pt][0pt]{\multirow{2}{*}{\rotatebox{90}{\textbf{TLS}, $\beta = [-50, 40 ] $ dB}}}  \hskip0.1cm \raisebox{2.9cm}[0pt][0pt] {\multirow{2}{*}{\rotatebox{90}{Perpendicular \hskip0.25cm Parallel}}}&
        \rule{0pt}{1.80cm}
        \cellcolor[HTML]{e5ffd4} \includegraphics[width=2.25cm, height=1.5cm]{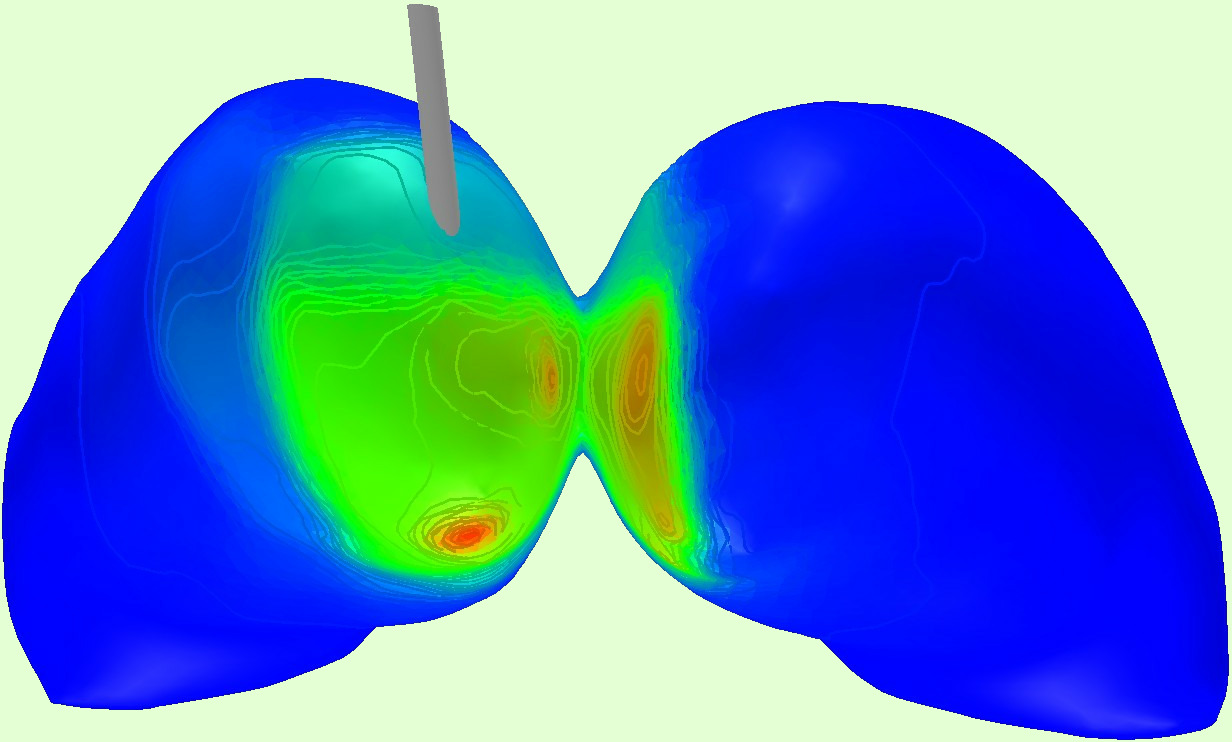} & 
        \cellcolor[HTML]{e5ffd4} \includegraphics[width=2.25cm, height=1.5cm]{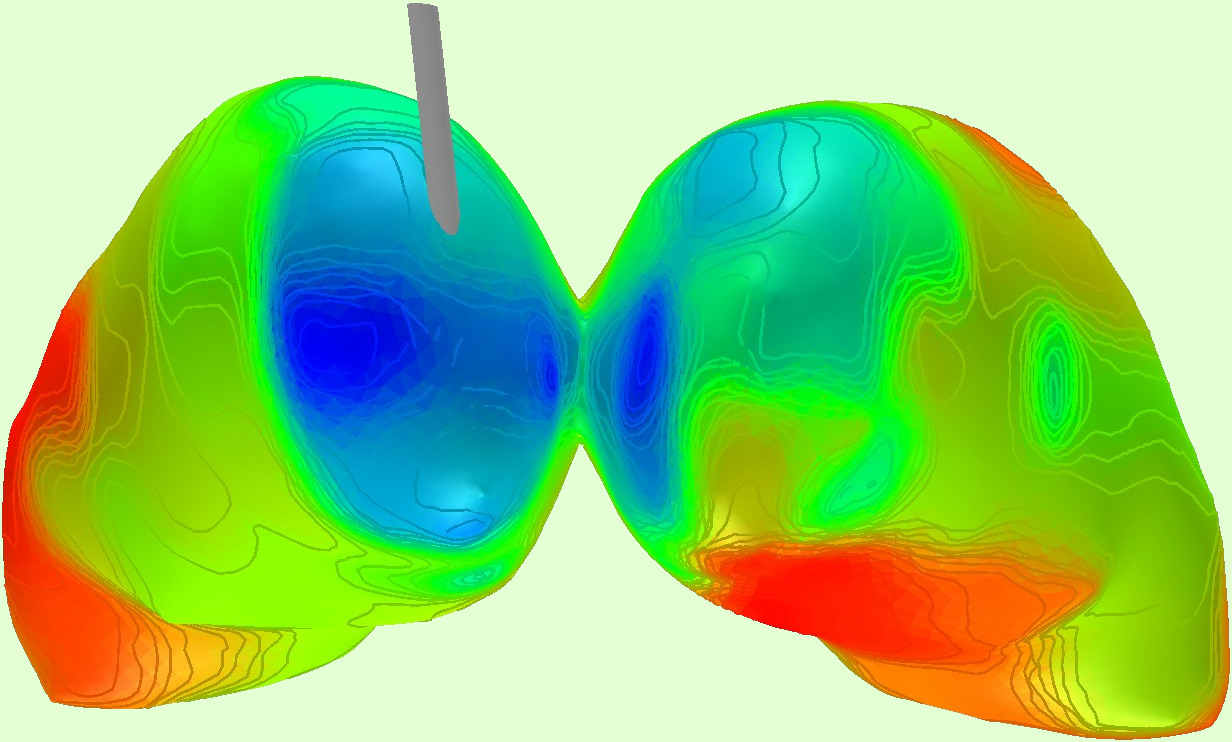} & 
        \cellcolor[HTML]{e5ffd4} \includegraphics[width=2.25cm, height=1.5cm]{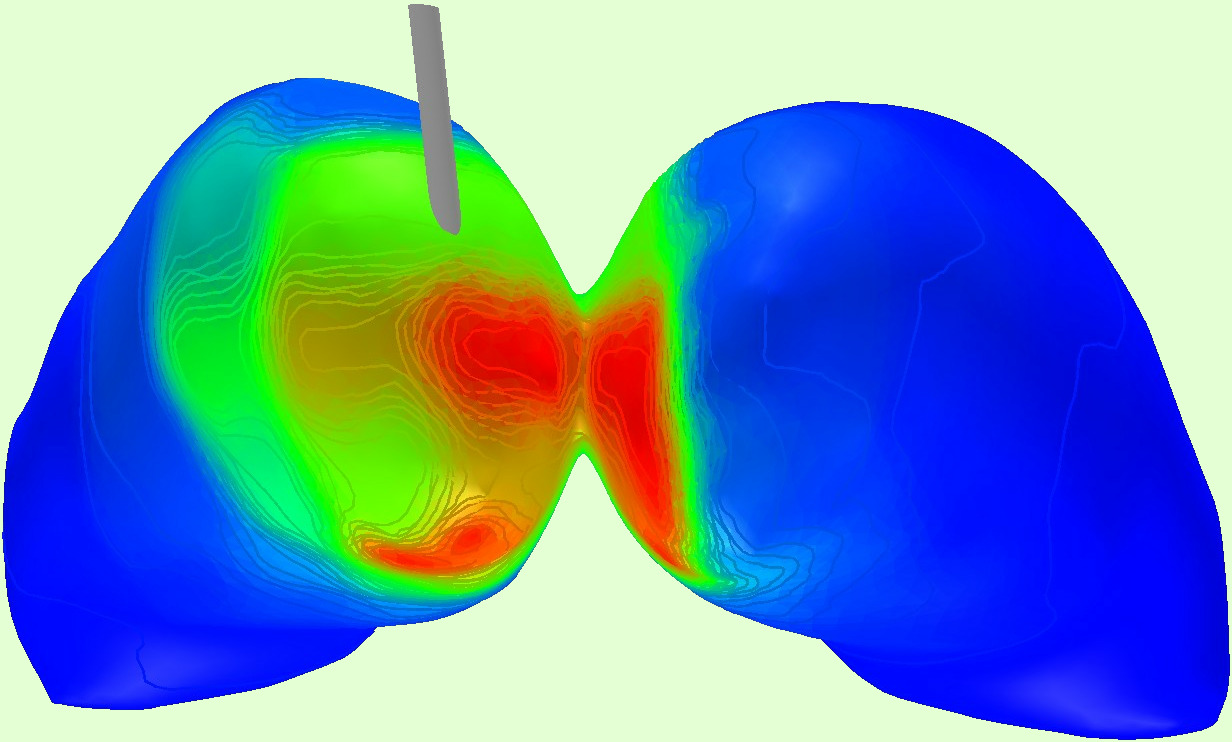} & 
        \cellcolor[HTML]{e5ffd4} \includegraphics[width=2.25cm, height=1.5cm]{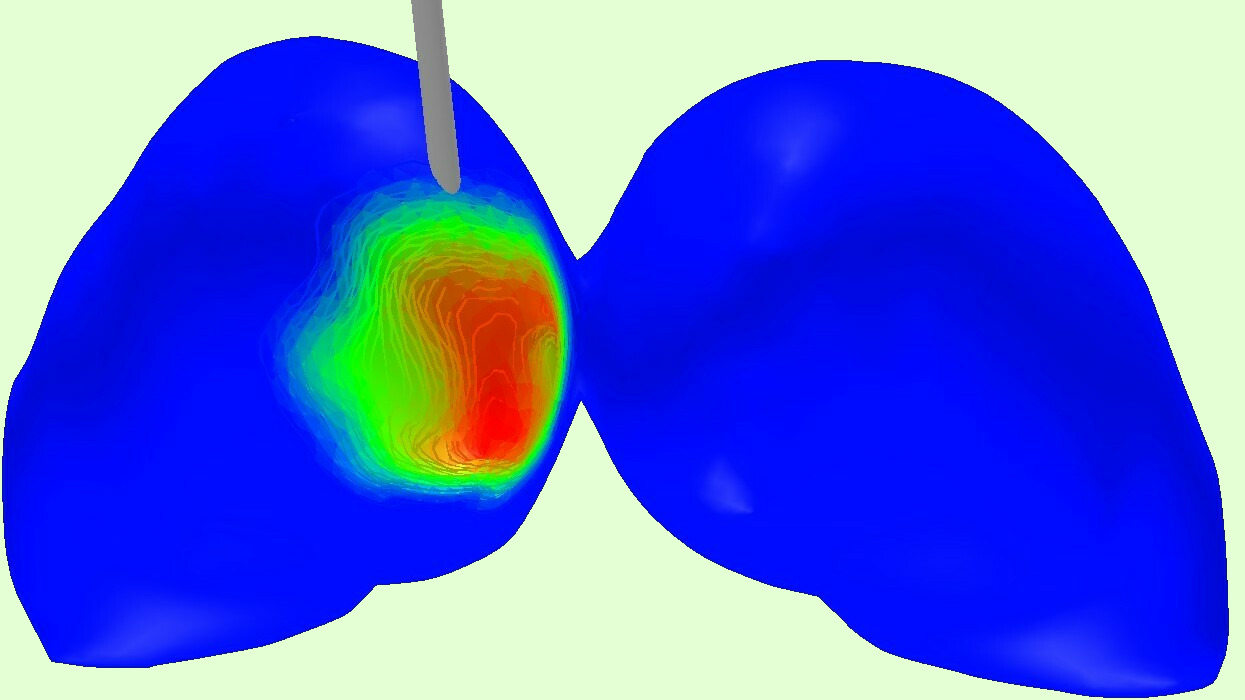} & 
        \cellcolor[HTML]{e5ffd4} \includegraphics[width=2.25cm, height=1.5cm]{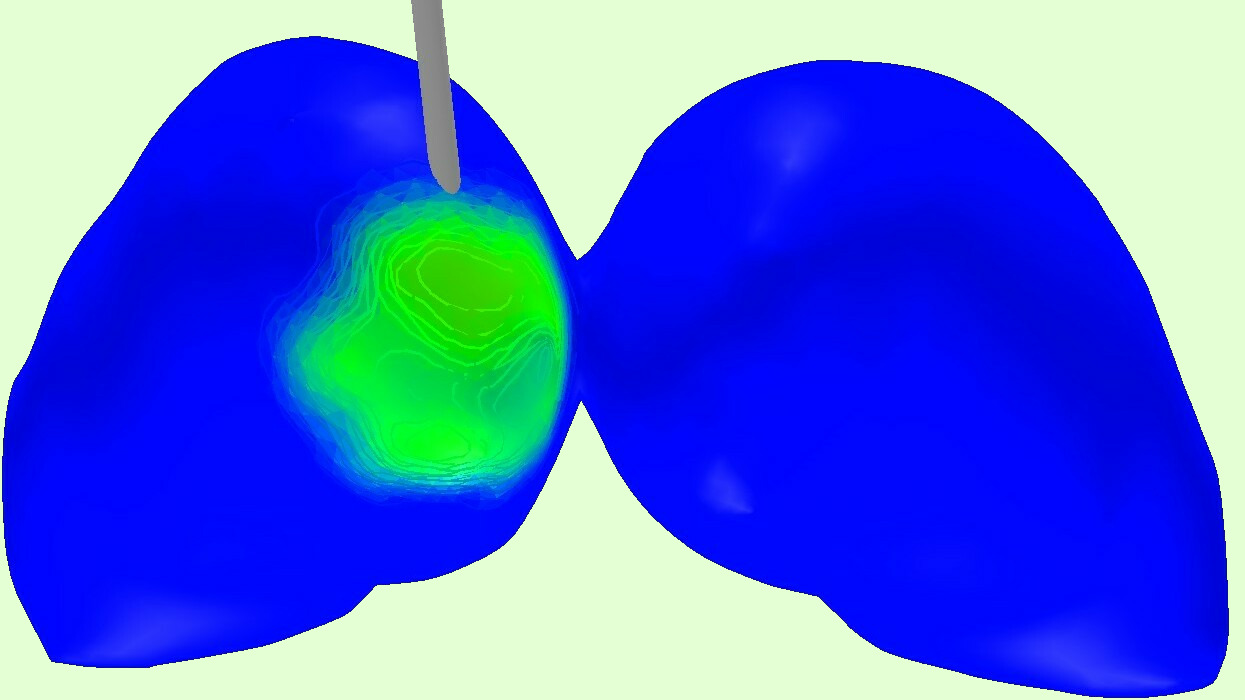} & 
        \cellcolor[HTML]{e5ffd4} \includegraphics[width=2.25cm, height=1.5cm]{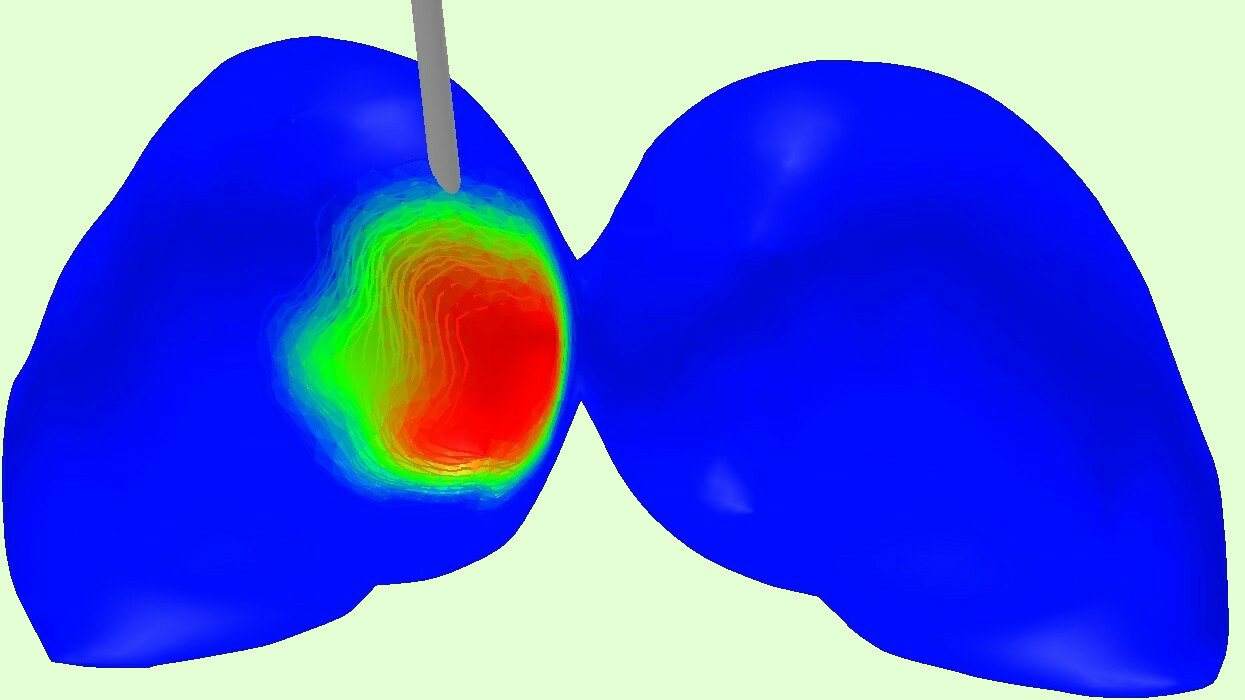} \\ 
        \hline
        %\hline
        
    \end{tabular}
    \end{small}
    \caption{Volumetric optimization results for the 40-contact probe and noiseless lead fields. As compared to the results obtained with 8-contact probe the overall level of $\Theta$ is greater for L1L1(A), TLS and RP, while for L1L1(B) with limited dynamic range ($\varepsilon$ range) there is no clear increase.}
    \label{Fig:Medtronic_results}
\end{figure*}

%% The comparison of works 
%\input{tables/Table_Comparison_Methods_L1L1}

\section{Results}
\label{sec:results}

The applied \textcolor{black}{optimization} algorithms generate distinct distributions of the variables $\Gamma$, $\Xi$, and $\Theta$, which can be differentiated based on how they are spatially organized across the target domain under the optimal electrode configurations for the 8- (Fig.~\ref{Fig:Abbott_results}) and 40-contact (Fig.~\ref{Fig:Medtronic_results}) leads. The direction of the current dipole plays the role in shaping electric field targeting. The anode–cathode pair exhibiting the strongest positive and negative values may activate electrodes located as close as \( 0.189 \pm 0.009\,\text{mm} \) or as far as \( 6.375 \pm 0.375\,\text{mm} \) apart across the lead. With the parallel orientation, the targeting is highly concentrated near the tip of the probe, whereas the perpendicular orientation results in a broader spatial spread, with activity shifted medially relative to the lead.

Bipolar configurations obtained via RP yield the highest focused current density \( \Gamma \); however, this comes at the expense of focality, as indicated by elevated nuisance current \( \Xi \) and a reduced field ratio \( \Theta \). This trade-off is understandable, as the multipolar configurations are designed to maximize the field ratio \( \Theta \) subject to the constraint \( \Gamma \geq \Gamma_0 > 0 \), with \( \Gamma_0 = 0.8\,\text{mA} \), thereby allowing activation of more than two contacts.

\begin{figure}[ht!]
\centering
    \begin{footnotesize}
    \centering
    \begin{minipage}{8.5cm}
        \centering
        {\bf 8-contact Electrode Configuration} \\
        \hrule
        \vskip0.2cm
        %%% A, B, C (Text)
        \begin{minipage}{0.35cm}
        \rotatebox{90}{{\bf Noiseless} lead field}
        \end{minipage}
        \begin{minipage}{3.3cm}
            \centering
            \includegraphics[width = 3.3cm]{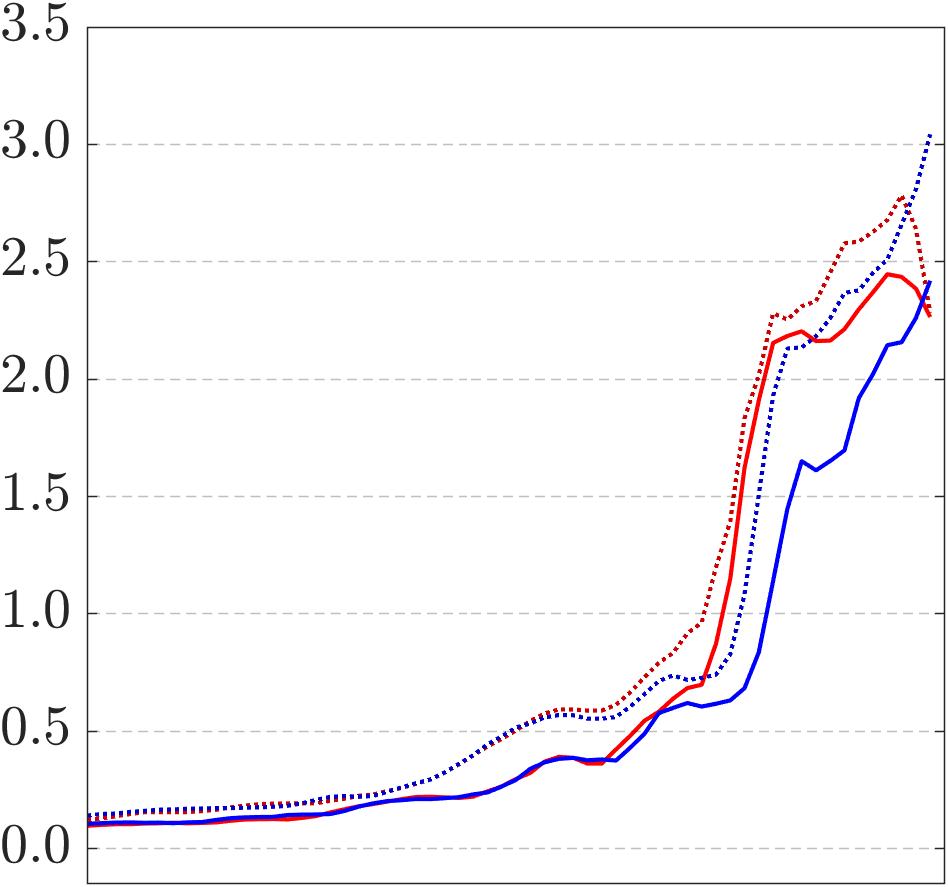} \\
            \textbf{A1.} {Focus} $\Gamma$ (A/m\textsuperscript{2})
        \end{minipage}
        \hskip0.25cm
        \begin{minipage}{3.3cm}
            \centering
            \includegraphics[width = 3.3cm]{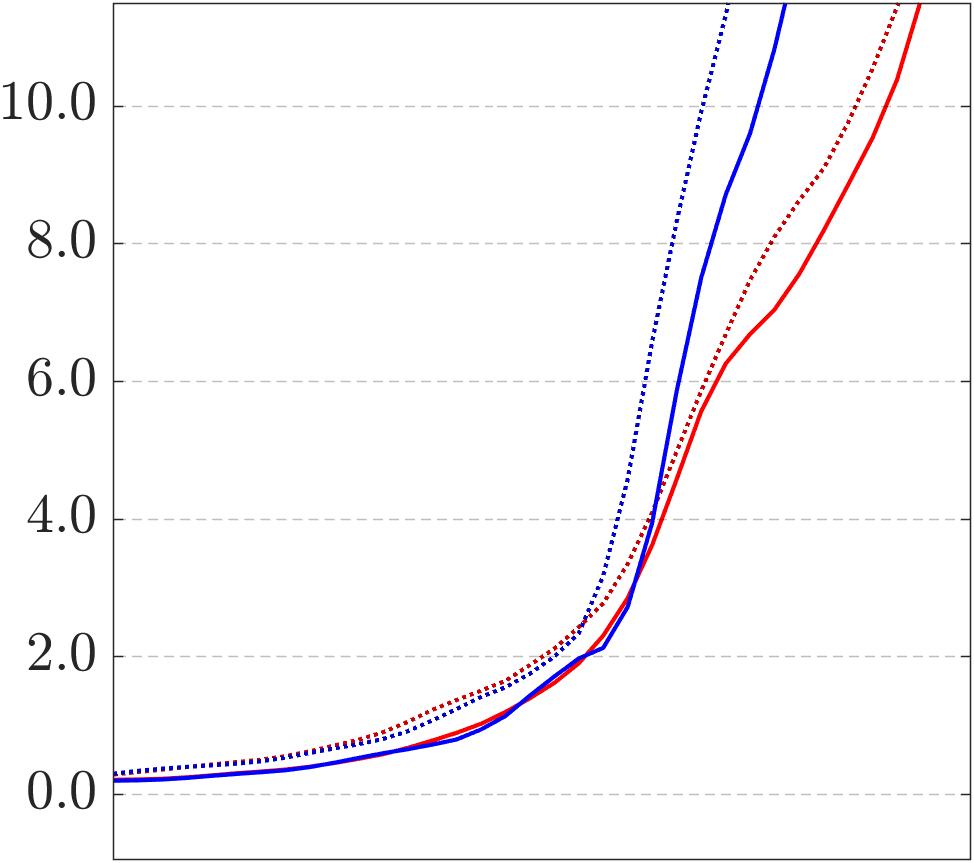} \\
            \textbf{B1.} {Ratio} $\Theta$ (rel.)
        \end{minipage}
        %%%%
        \vskip0.2cm
        %%%% A, B, C
        \begin{minipage}{0.35cm}
        \rotatebox{90}{{\bf Noisy} lead field}
        \end{minipage}
        \begin{minipage}{3.3cm}
            \centering
            \includegraphics[width = 3.3cm]{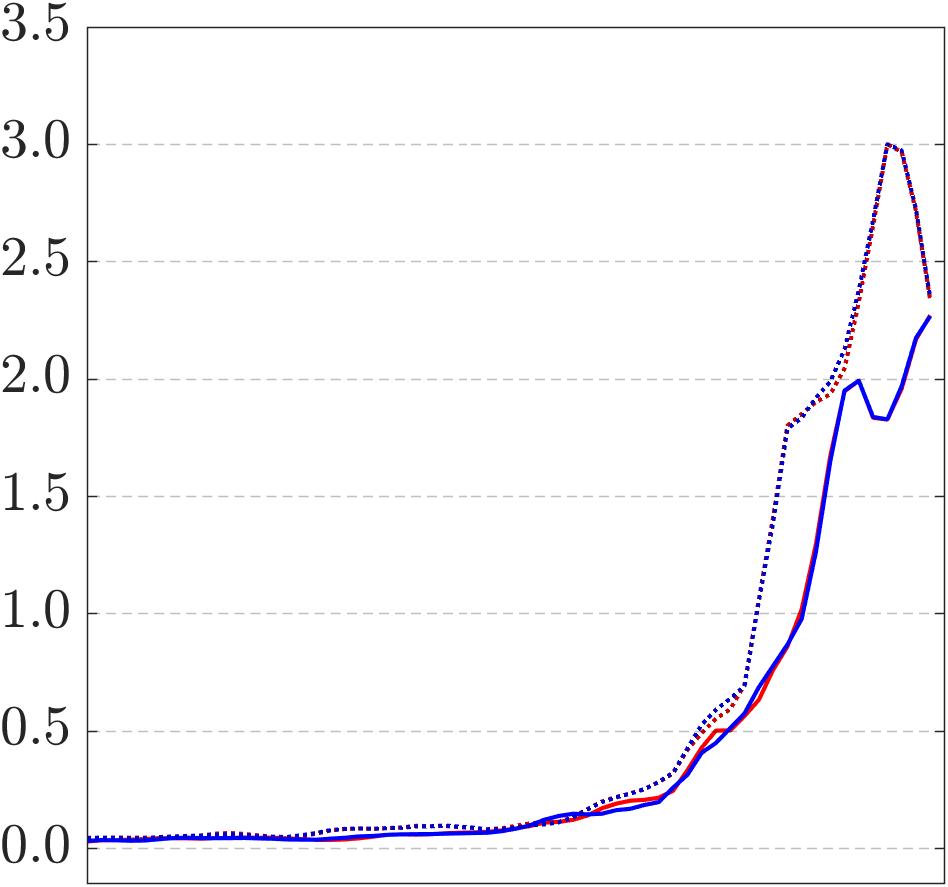} \\
            \textbf{C1.} {Focus} $\Gamma$ (A/m\textsuperscript{2})
        \end{minipage}
        \hskip0.25cm
        \begin{minipage}{3.3cm}
            \centering
            \includegraphics[width = 3.3cm]{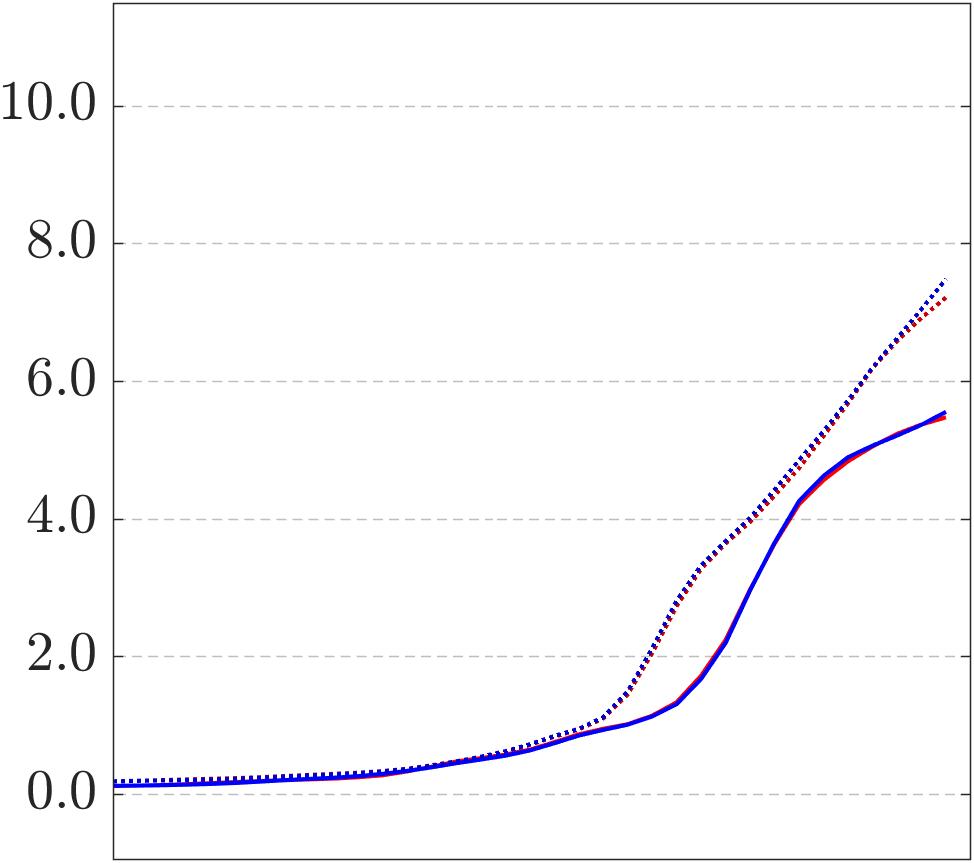} \\
            \textbf{D1.} {Ratio} $\Theta$ (rel.)
        \end{minipage}
        %%%% D, F, G
        %%%%
    \end{minipage} \\
    \vskip0.5cm
    %%%%%%%%%%%%%%%%%%%%%%%%%%%%%%%%%%%%%%%%%%%%%%%%%%%%
    \begin{minipage}{8.5cm}
        \centering
        {\bf 40-contact Electrode Configuration} \\
        \hrule
        \vskip0.2cm
        %%% A, B, C (Text)
        \begin{minipage}{0.35cm}
        \rotatebox{90}{{\bf Noiseless} lead field}
        \end{minipage}
        \begin{minipage}{3.3cm}
            \centering
            \includegraphics[width = 3.3cm]{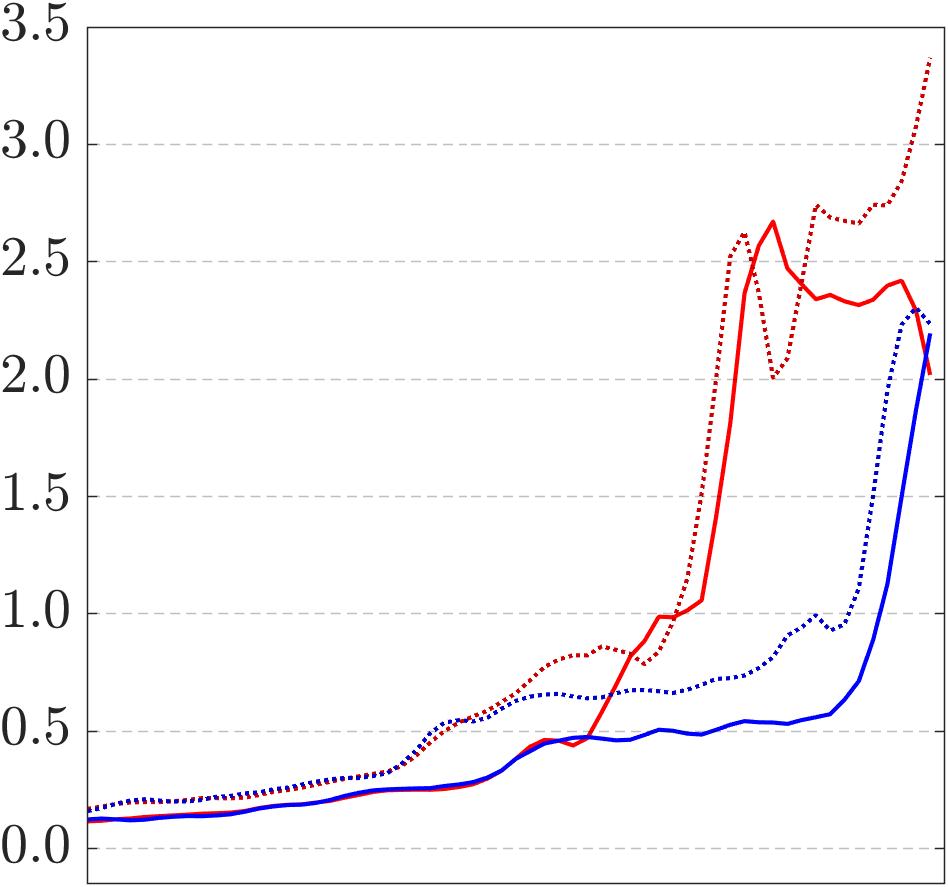} \\
            \textbf{A2.} {Focus} $\Gamma$ (A/m\textsuperscript{2})
        \end{minipage}
        \hskip0.25cm
        \begin{minipage}{3.3cm}
            \centering
            \includegraphics[width = 3.3cm]{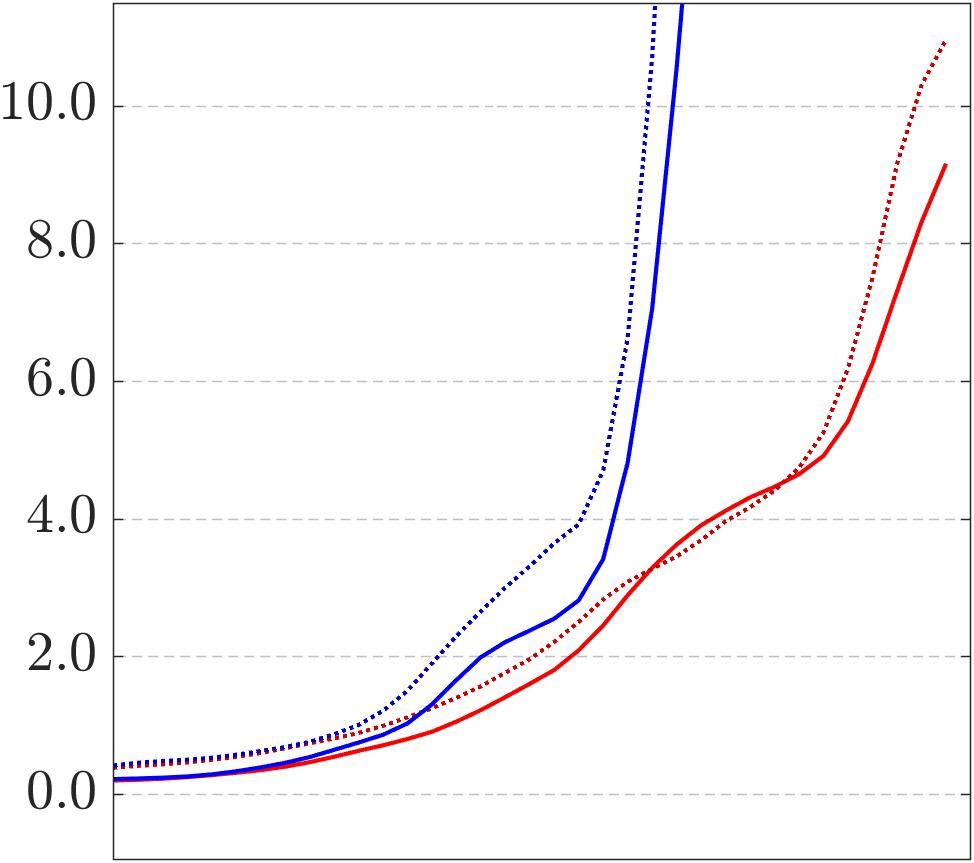} \\
            \textbf{B2.} {Ratio} $\Theta$ (rel.)
        \end{minipage}
        %%%%
        \vskip0.2cm
        %%%% D, E, F
        %%%% A, B, C
        \begin{minipage}{0.35cm}
        \rotatebox{90}{{\bf Noisy} lead field}
        \end{minipage}
        \begin{minipage}{3.3cm}
            \centering
            \includegraphics[width = 3.3cm]{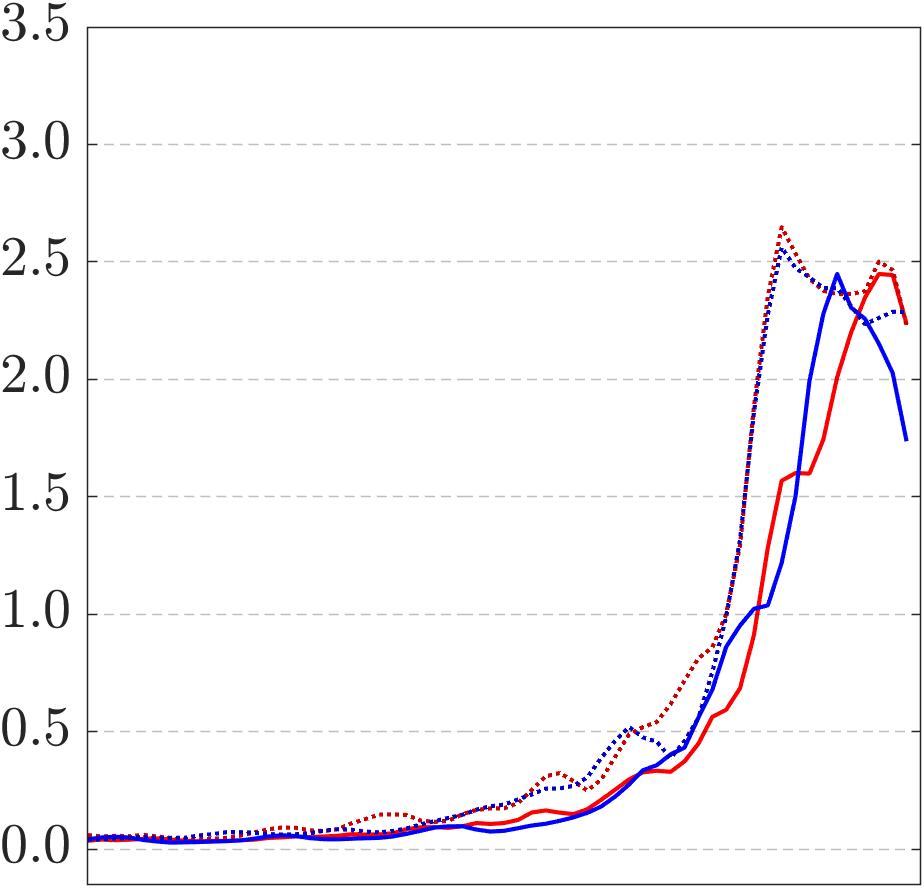} \\
            \textbf{C2.} {Focus} $\Gamma$ (A/m\textsuperscript{2})
        \end{minipage}
        \hskip0.25cm
        \begin{minipage}{3.3cm}
            \centering
            \includegraphics[width = 3.3cm]{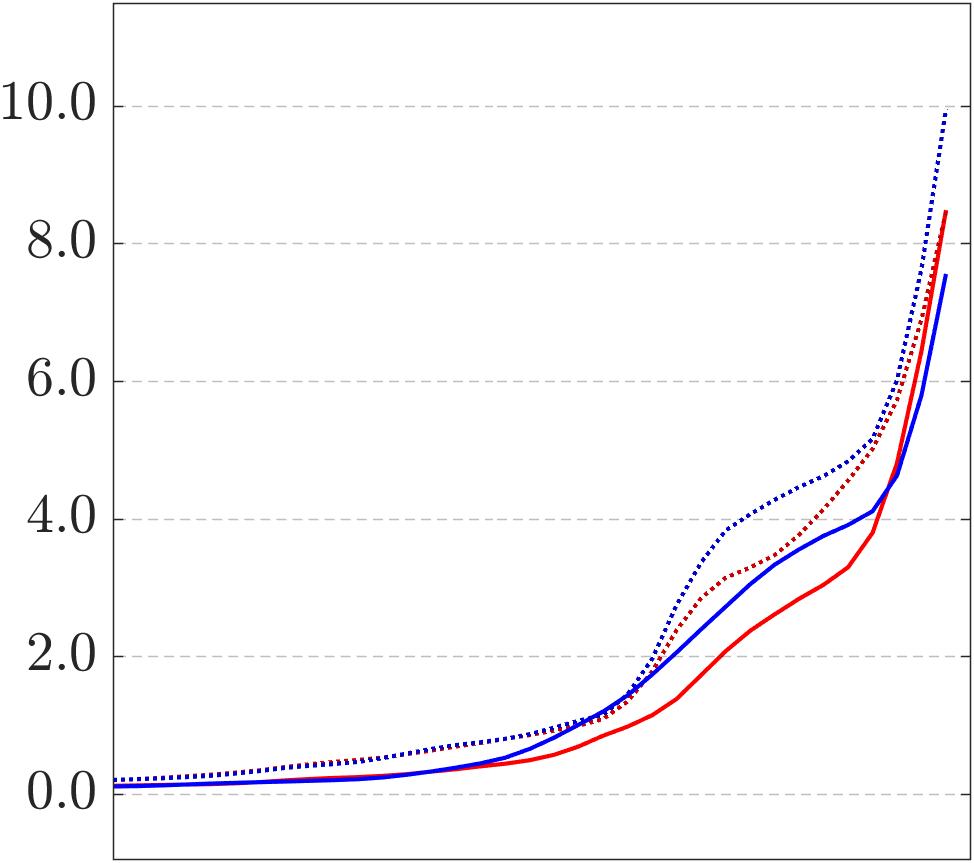} \\
            \textbf{D2.} {Ratio} $\Theta$ (rel.) \\
        \end{minipage}
        \end{minipage}     
        \vskip0.2cm
               \begin{tikzpicture}
    \node[draw=none] at (0,0) {
        \begin{tabular}{ll}
                    \tikz \draw[blue, thick] (0,0) -- (0.6,0); & L1L1(A), $\varepsilon \in [-160, 0]$, Parallel \\ 
            \tikz \draw[red, thick] (0,0) -- (0.6,0); & L1L1(B), $\varepsilon \in [-10, 0]$, Parallel \\
                        \tikz \draw[blue, dotted, thick] (0,0) -- (0.6,0); & L1L1(A), $\varepsilon \in [-160, 0]$, Perpendicular \\
            \tikz \draw[red, dotted,thick] (0,0) -- (0.6,0); & L1L1(B), $\varepsilon \in [-10, 0]$, Perpendicular 
        \end{tabular}
    };
\end{tikzpicture} 
\end{footnotesize}
\caption{Performance of L1L1(A) and L1L1(B) for optimizing the field ratio $\Theta$ subject to criterion $\Gamma \geq \Gamma_0$ with noiseless and noisy (PSNR 40 dB) lead field. L1L1(A) and L1L1(B) suppress the nuisance current density~$\Xi$ pattern using two different dynamic ranges, namely, $\varepsilon \in [-160, 0]$ (blue) and $\varepsilon \in [-10,  0]$ (red), respectively. The results are sub-grouped by the direction of the current dipole, whether its alignment is parallel (solid line) or perpendicular (dotted line) with the applied lead. The degrees of freedom have been sorted according to ascending field ratio $\Theta$. The results show that the overly optimistic dynamic range applied in L1L1(A) results in very high field values of $\Theta$ in the noiseless case, which are, however, absent in the noisy one. In contrast, L1L1(B) is only marginally affected by the noise due to more realistic assumption of the maximum obtainable dynamic range under uncertainty. }
\label{Fig:Norm_differences_plot}
\end{figure}

Figure~\ref{Fig:Norm_differences_plot} and Table~\ref{Tab:Together} compare noiseless and noisy (\( \text{PSNR} = 40\,\text{dB} \)) results for \( \Gamma \) and \( \Theta \) using the L1L1 method over a LR-LF. In the noiseless setting, the L1L1(B) reconstructions with \( \varepsilon \in [-10, 0]\,\text{dB} \) yield slightly higher \( \Gamma \) and \( \Xi \) values, and moderate \( \Theta \) compared to L1L1(A) with \( \varepsilon \in [-160, 0]\,\text{dB} \). L1L1(A) provided unrealistically high values of the field ratio \( \Theta \): up to \( 27.07 \) in the 8-contact directional DBS leads and as high as \( 536.05 \) in the 40-contact case indicating a case of overfitting.

\begin{figure*}[ht!]
    \centering
    \begin{scriptsize}

   {\bf 8-contact Electrode Configuration} \\
    \vskip0.1cm
    \hrule
    \vskip0.1cm

    \begin{minipage}{9.12cm}
        \centering
        \begin{minipage}{0.15cm}
            \centering
            \rotatebox{90}{Parallel}
        \end{minipage}
        \begin{minipage}{2.87cm}
            \centering
            \includegraphics[width=2.86cm]{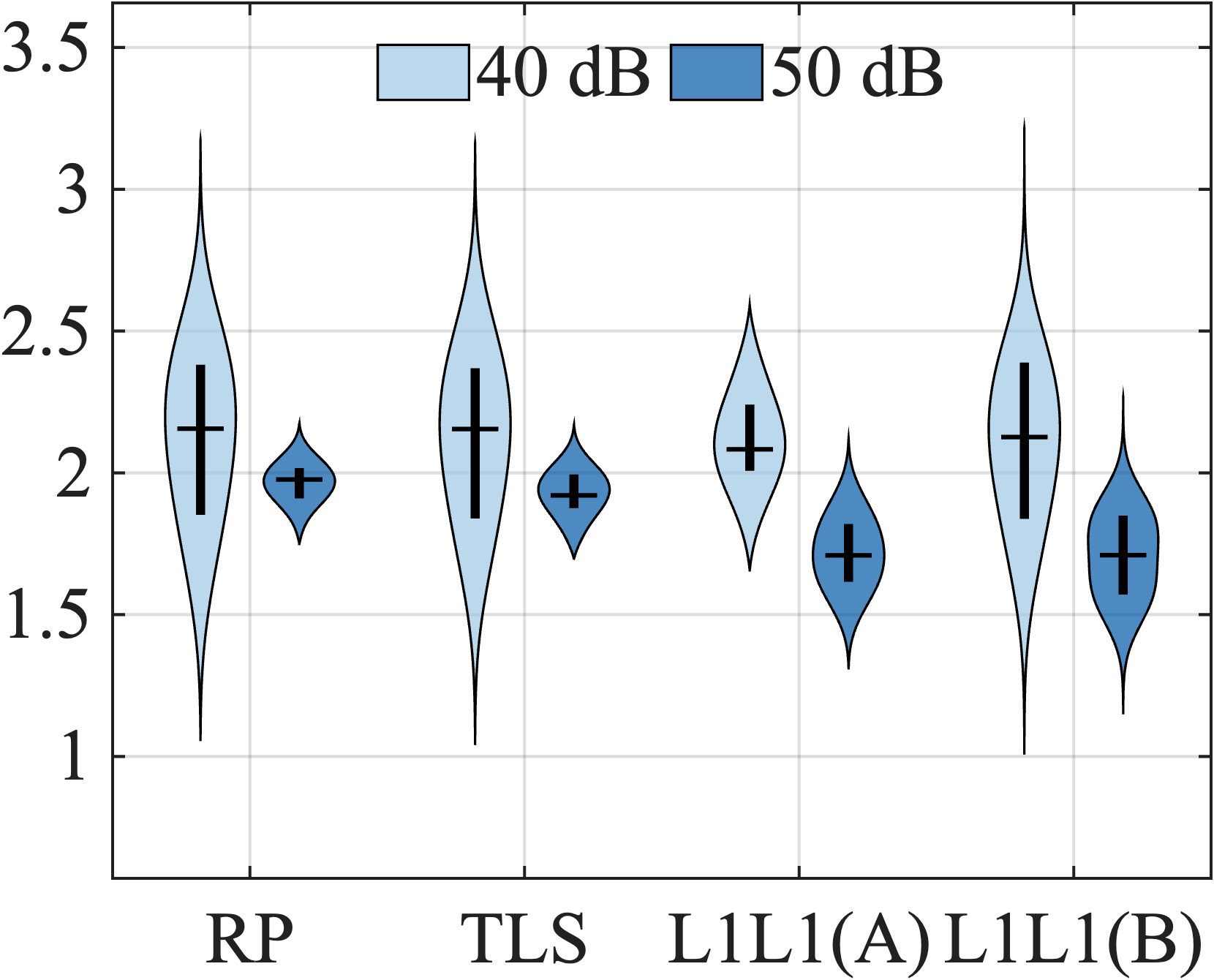} \\
            Focus $\Gamma$ (A/m\textsuperscript{2})
        \end{minipage}
        \begin{minipage}{2.87cm}
            \centering
            \includegraphics[width=2.86cm]{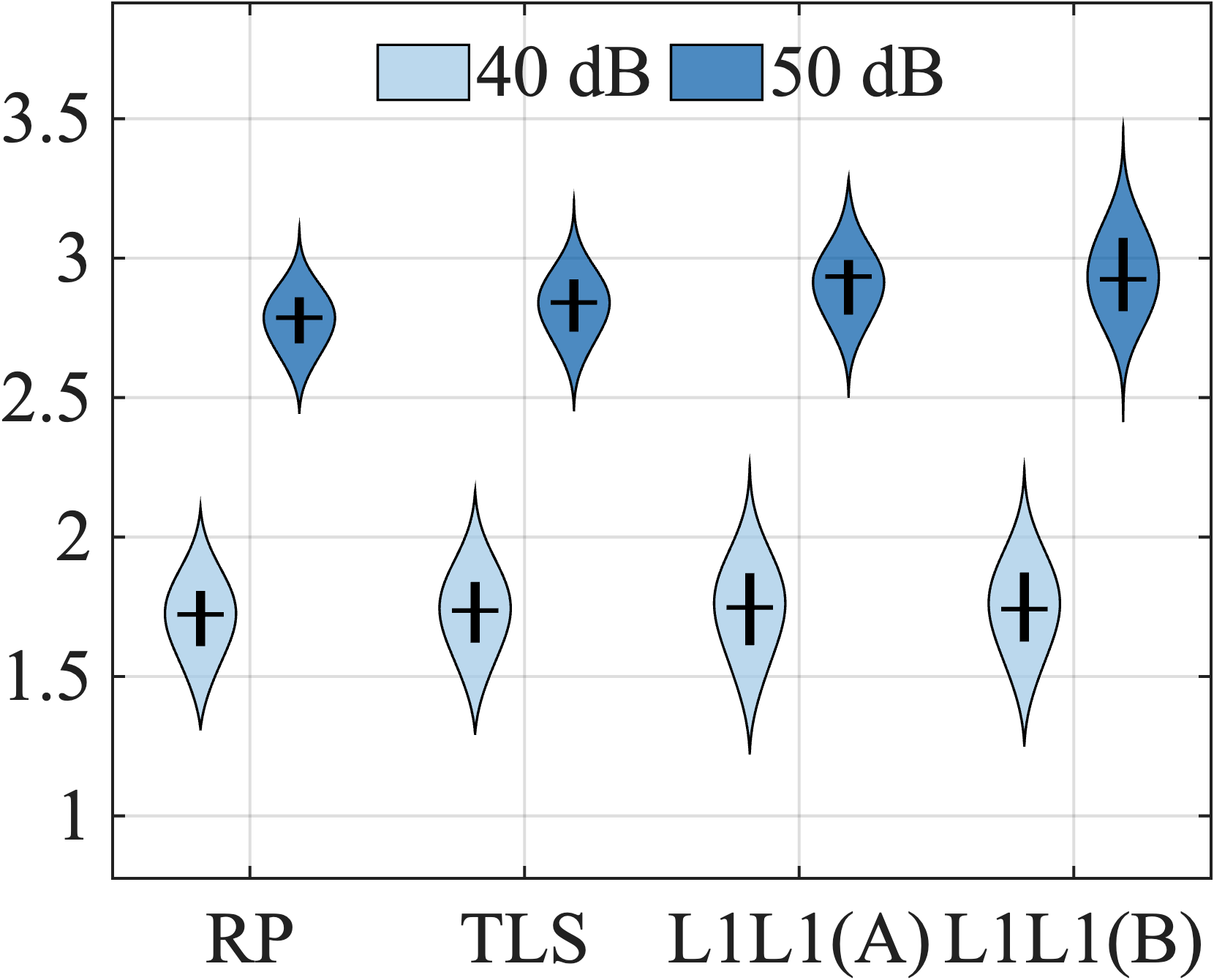} \\
            Ratio $\Theta$ (rel.)
        \end{minipage}
        \begin{minipage}{2.87cm}
            \centering
            \includegraphics[width=2.86cm]{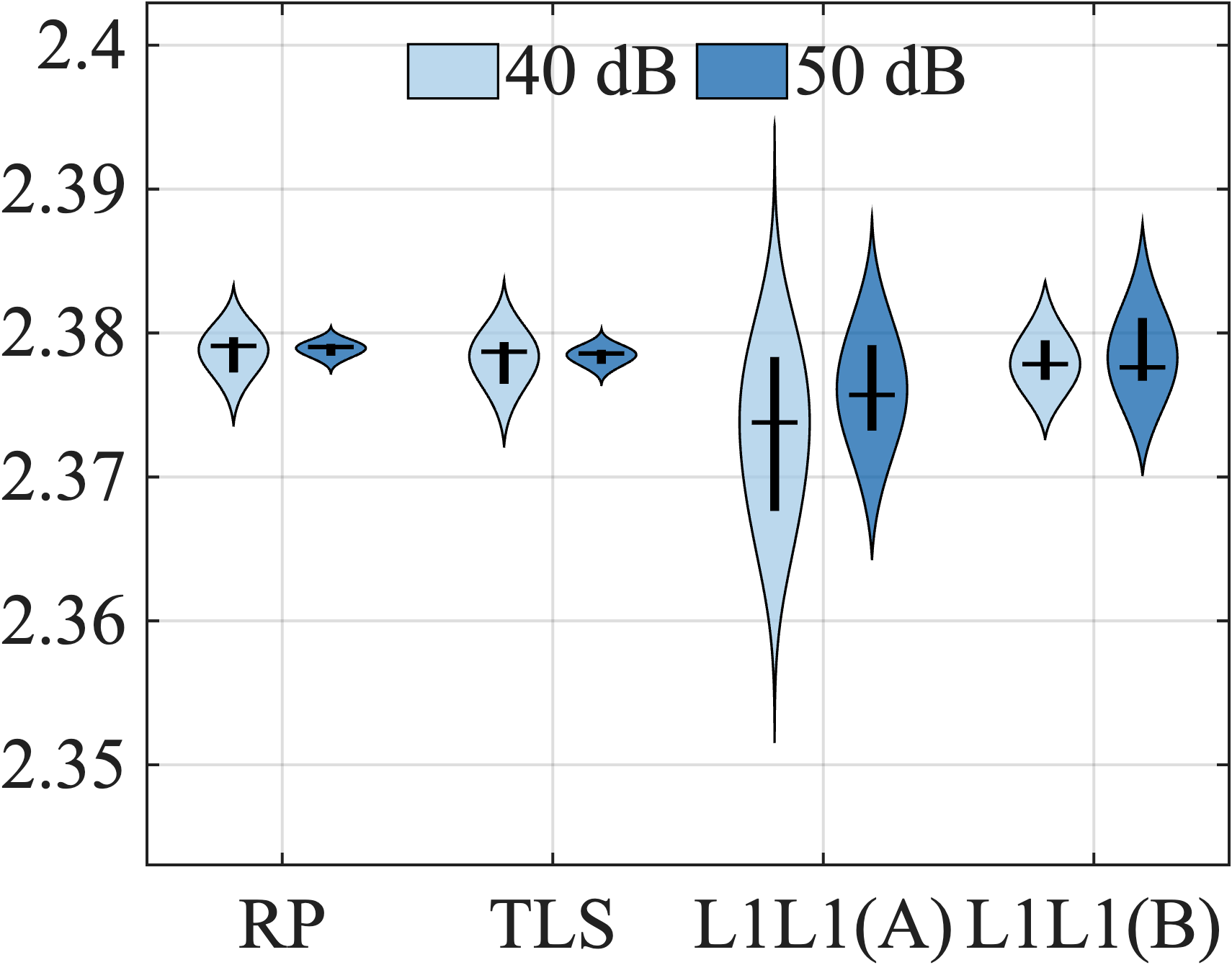} \\
            Targeting error (mm)
        \end{minipage} \\

        \vskip0.1cm

        \begin{minipage}{0.15cm}
            \centering
            \rotatebox{90}{Perpendicular}
        \end{minipage}
        \begin{minipage}{2.87cm}
            \centering
            \includegraphics[width=2.86cm]{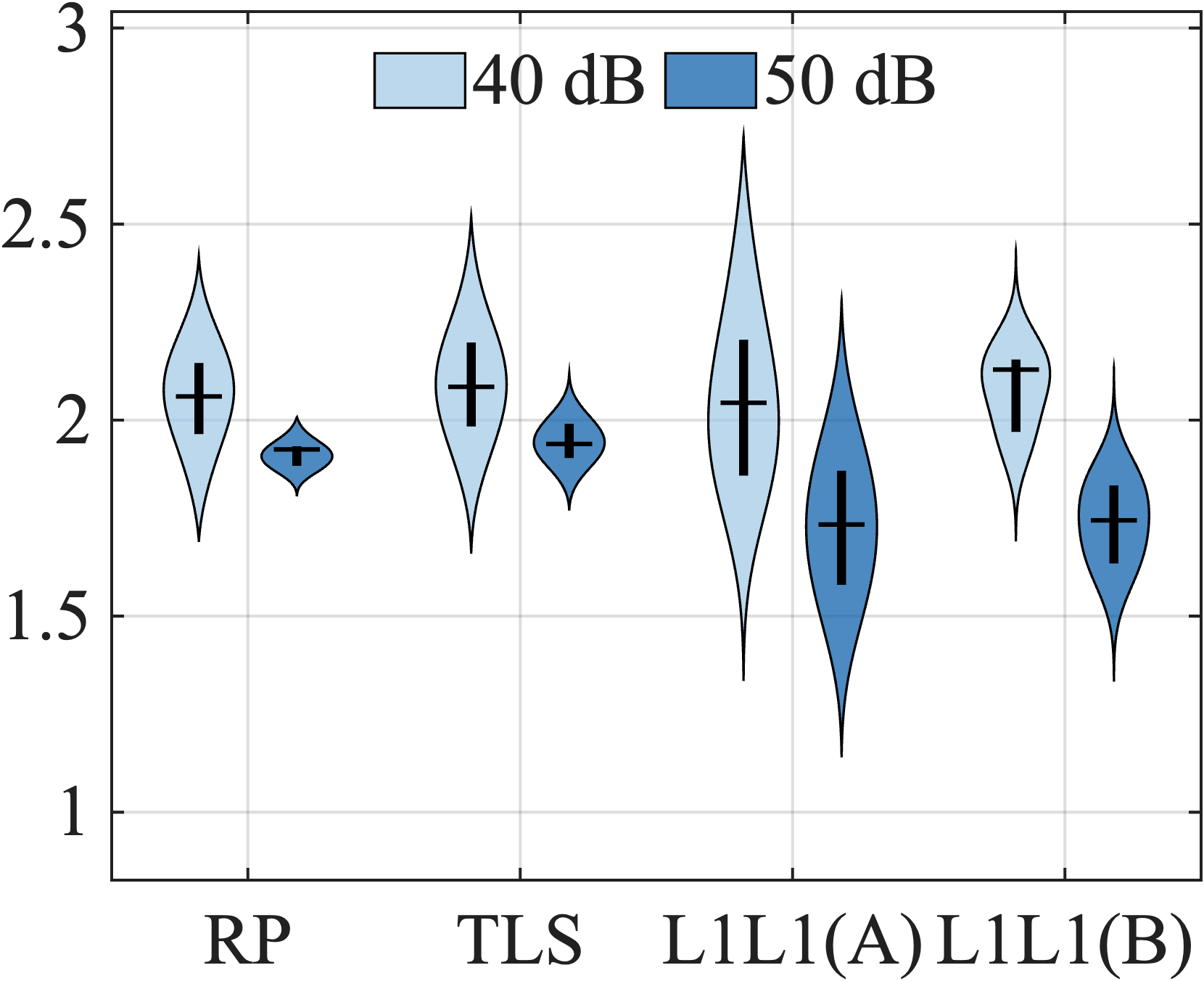} \\
            Focus $\Gamma$ (A/m\textsuperscript{2})
        \end{minipage}
        \begin{minipage}{2.87cm}
            \centering
            \includegraphics[width=2.86cm]{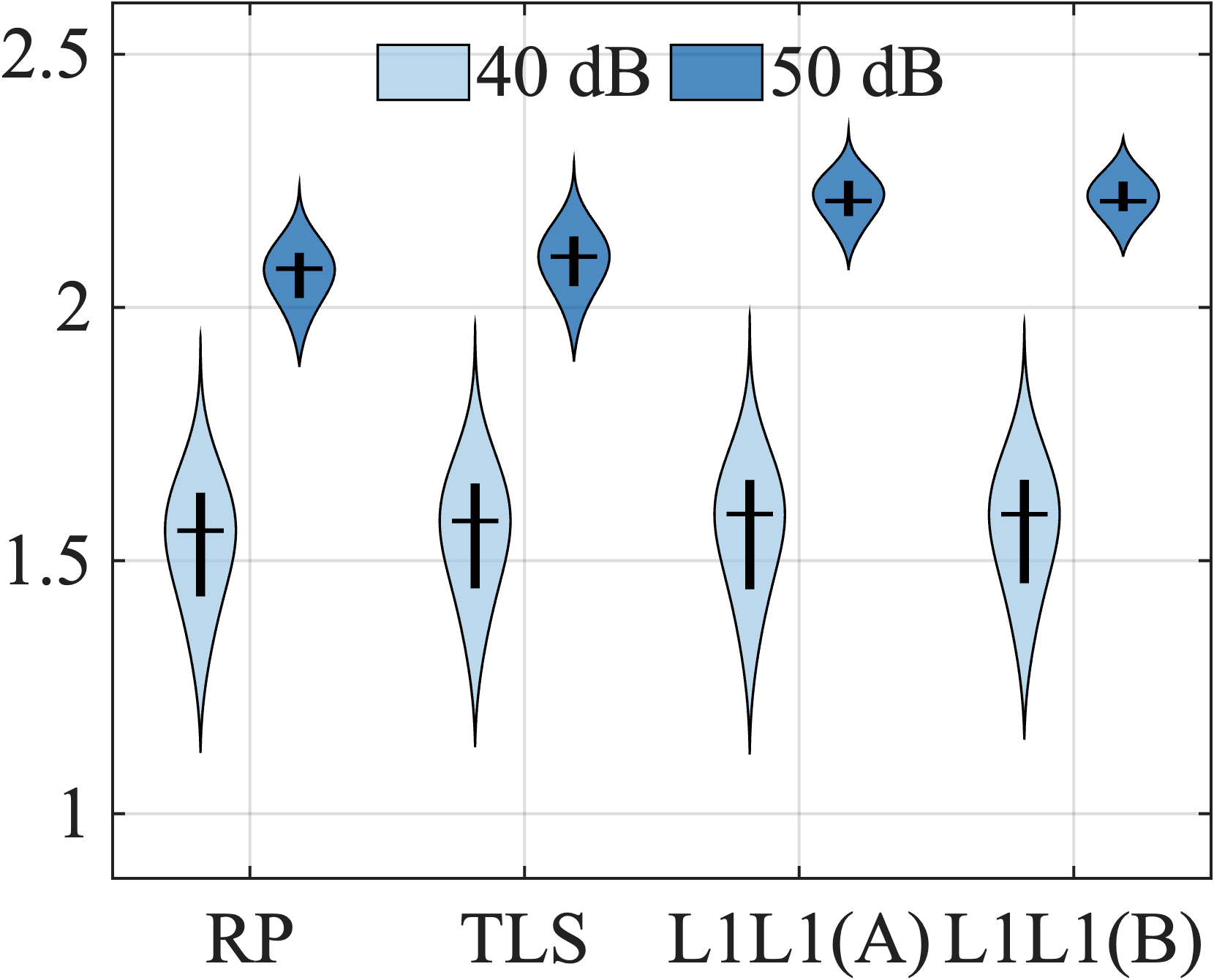} \\
            Ratio $\Theta$ (rel.)
        \end{minipage}
        \begin{minipage}{2.87cm}
            \centering
            \includegraphics[width=2.87cm]{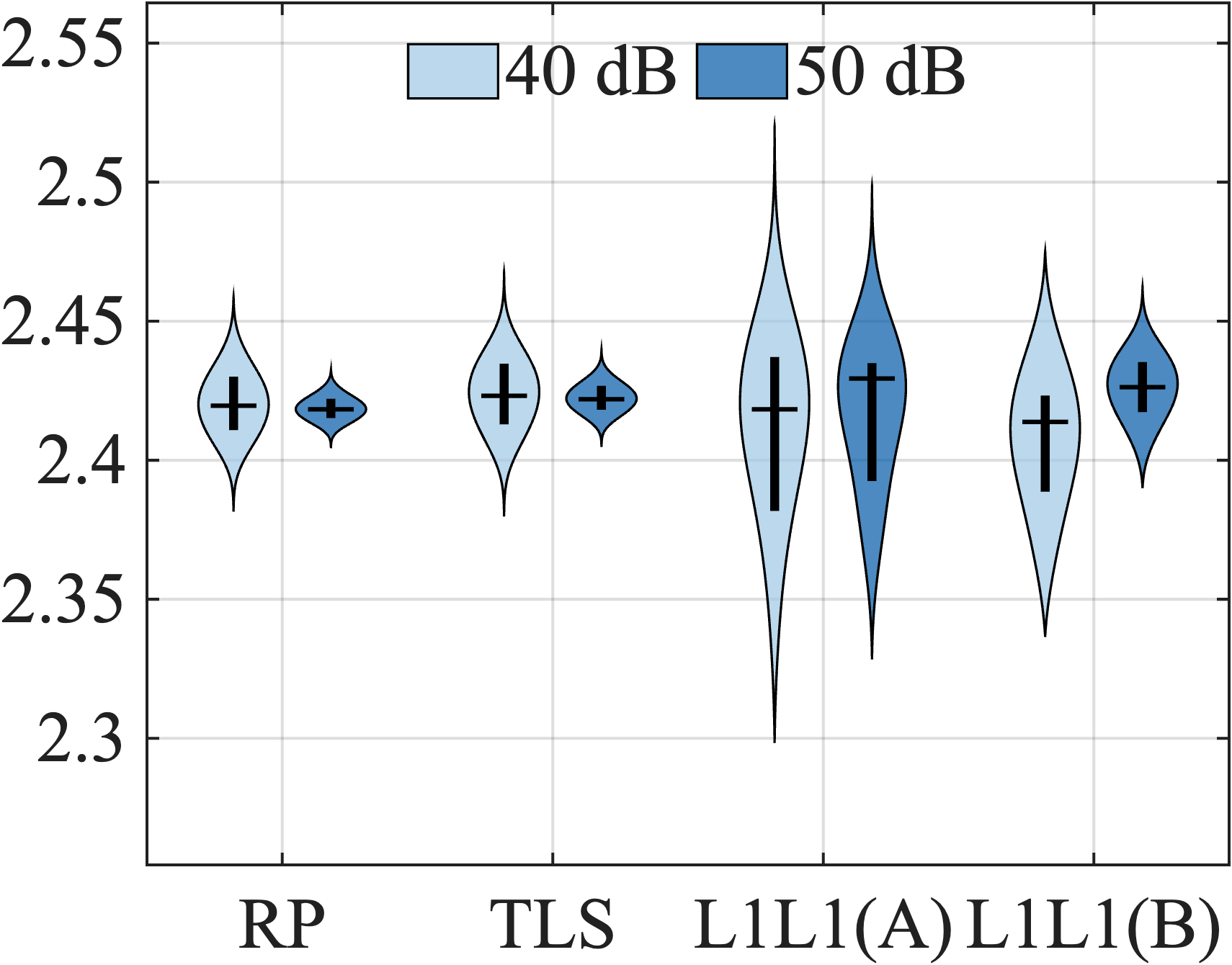} \\
            Targeting error (mm)
        \end{minipage} \\

        \vskip0.1cm
        Target (I)
    \end{minipage}
    \begin{minipage}{9.12cm}
        \centering
        \begin{minipage}{0.15cm}
            \centering
            \rotatebox{90}{Parallel}
        \end{minipage}
        \begin{minipage}{2.87cm}
            \centering
            \includegraphics[width=2.86cm]{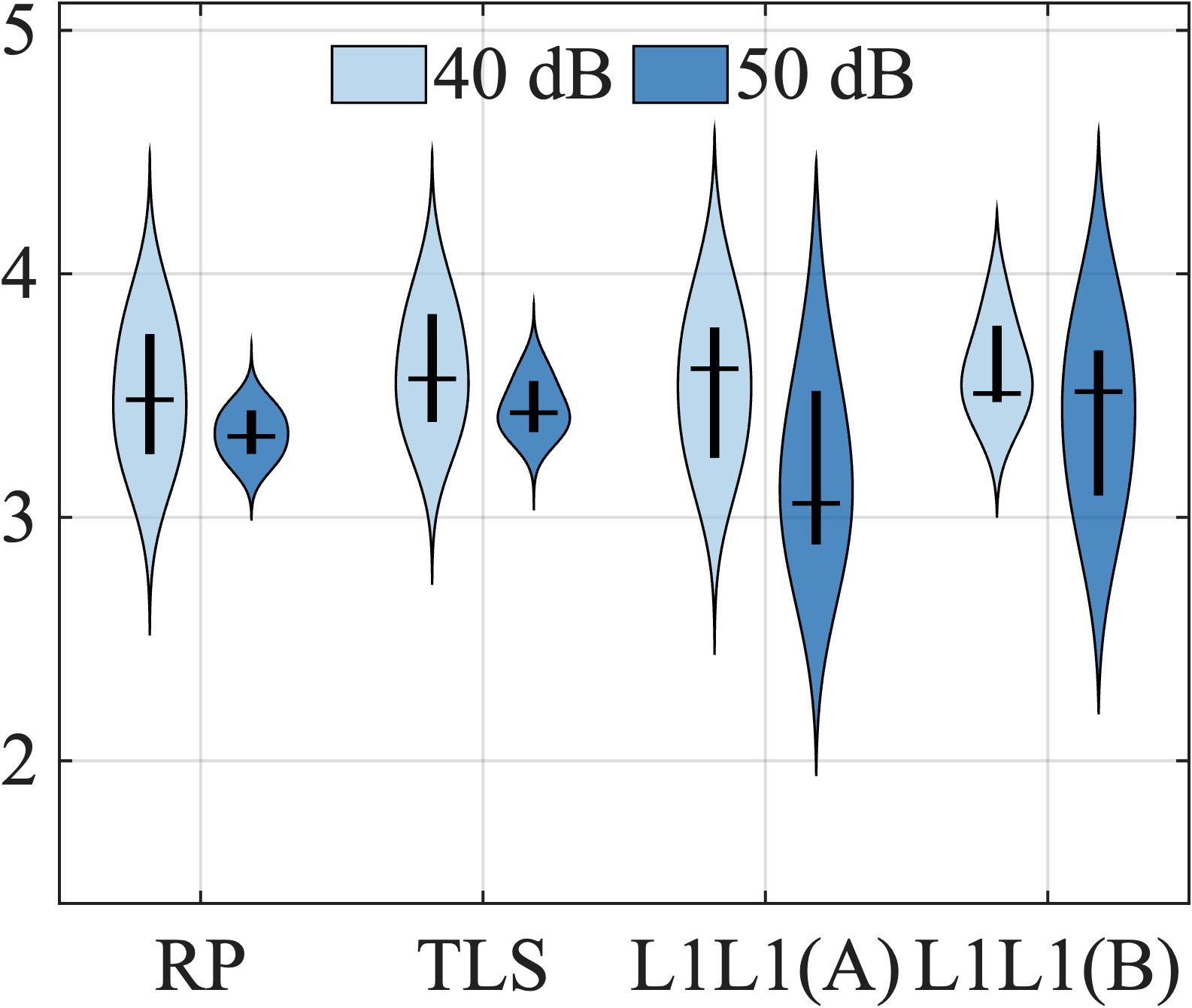} \\
            Focus $\Gamma$ (A/m\textsuperscript{2})
        \end{minipage}
        \begin{minipage}{2.87cm}
            \centering
            \includegraphics[width=2.86cm]{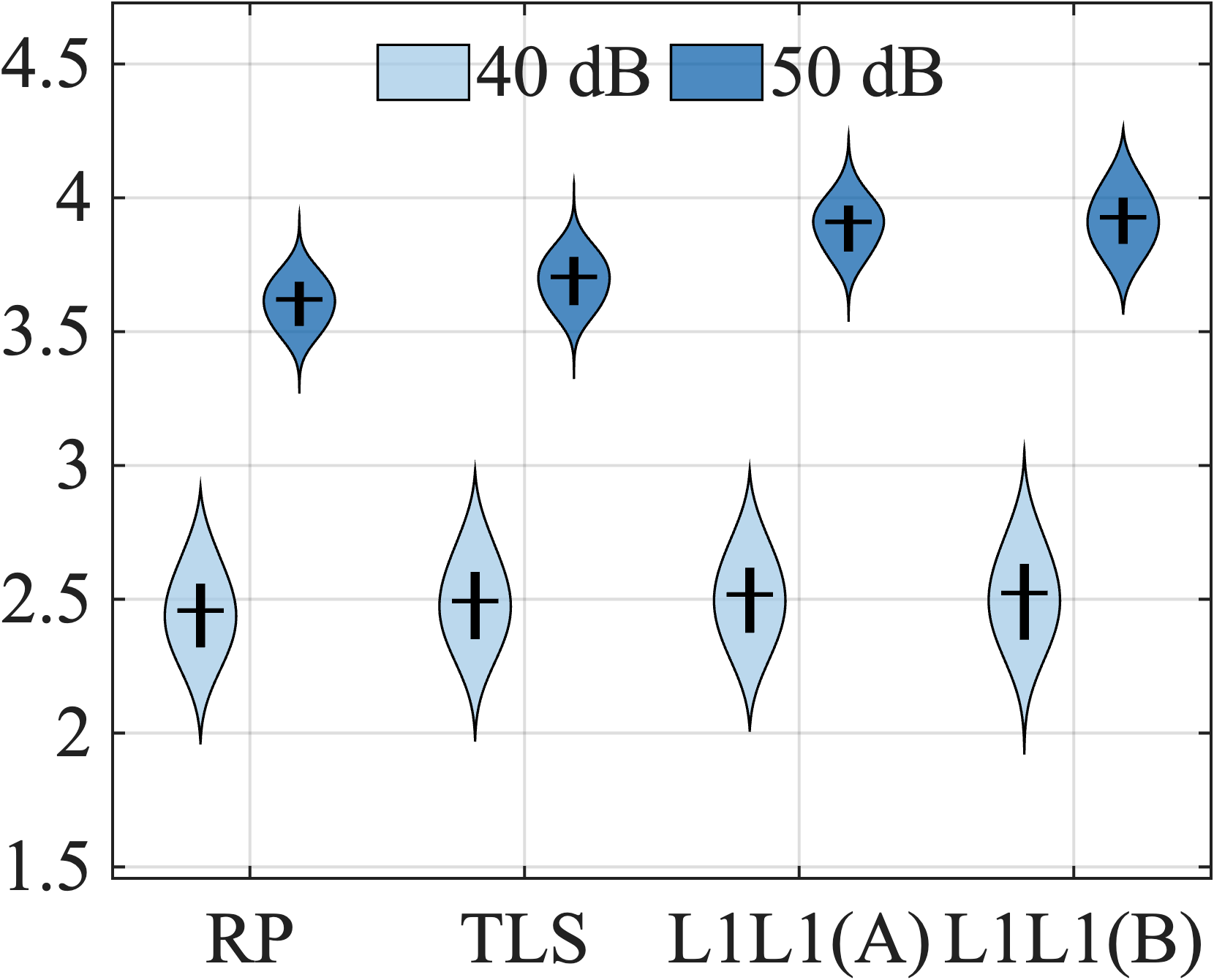} \\
            Ratio $\Theta$ (rel.)
        \end{minipage}
        \begin{minipage}{2.87cm}
            \centering
            \includegraphics[width=2.86cm]{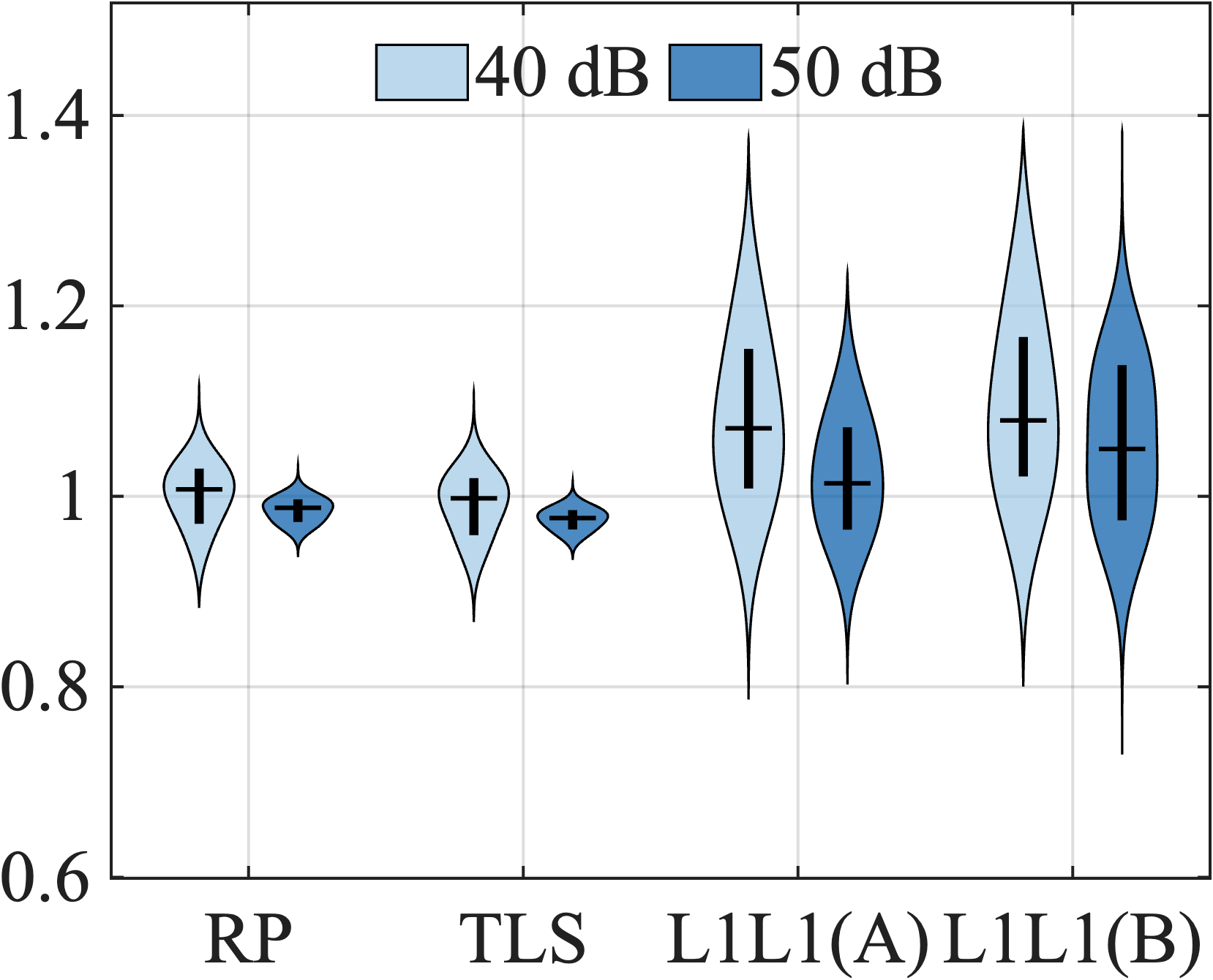} \\
            Targeting error (mm)
        \end{minipage} \\

        \vskip0.1cm

        \begin{minipage}{0.15cm}
            \centering
            \rotatebox{90}{Perpendicular}
        \end{minipage}
        \begin{minipage}{2.87cm}
            \centering
            \includegraphics[width=2.86cm]{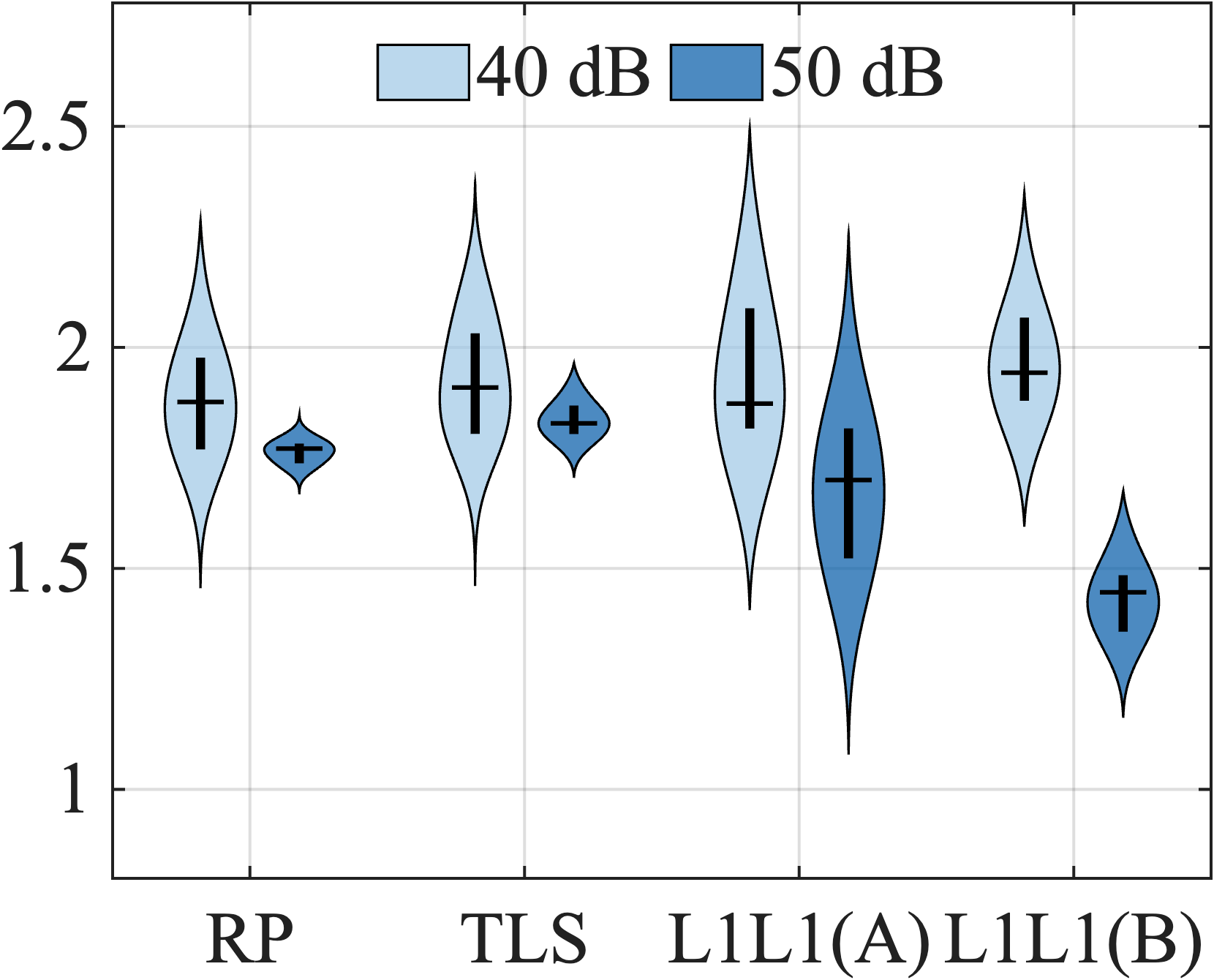} \\
            Focus $\Gamma$ (A/m\textsuperscript{2})
        \end{minipage}
        \begin{minipage}{2.87cm}
            \centering
            \includegraphics[width=2.86cm]{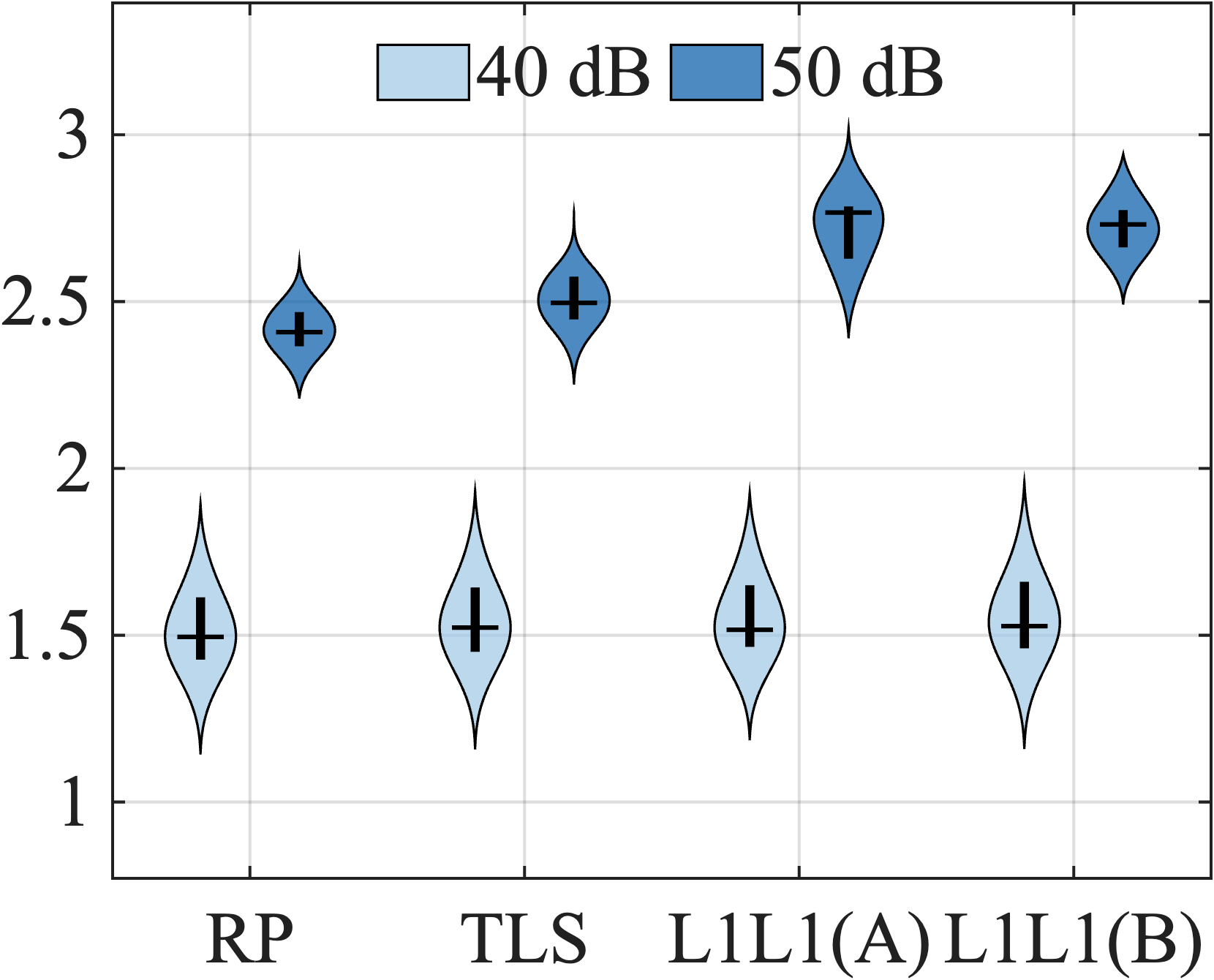} \\
            Ratio $\Theta$ (rel.)
        \end{minipage}
        \begin{minipage}{2.87cm}
            \centering
            \includegraphics[width=2.86cm]{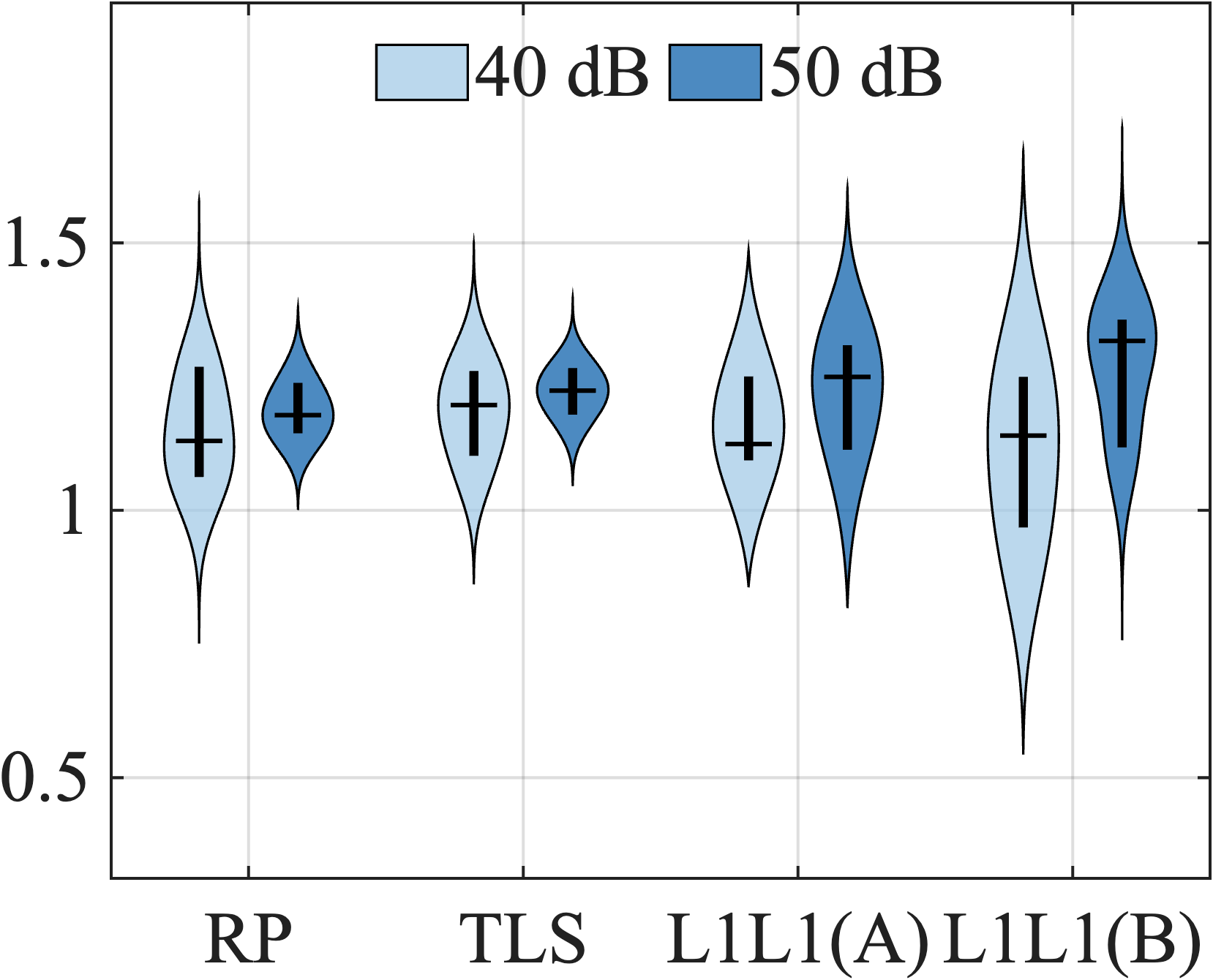} \\
            Targeting error (mm)
        \end{minipage} \\

        \vskip0.1cm
        Target (II)
    \end{minipage}

    \vskip0.25cm

    {\bf 40-contact Electrode Configuration} \\
    \vskip0.1cm
    \hrule
    \vskip0.1cm

    \begin{minipage}{9.12cm}
        \centering
        \begin{minipage}{0.15cm}
            \centering
            \rotatebox{90}{Parallel}
        \end{minipage}
        \begin{minipage}{2.87cm}
            \centering
            \includegraphics[width=2.86cm]{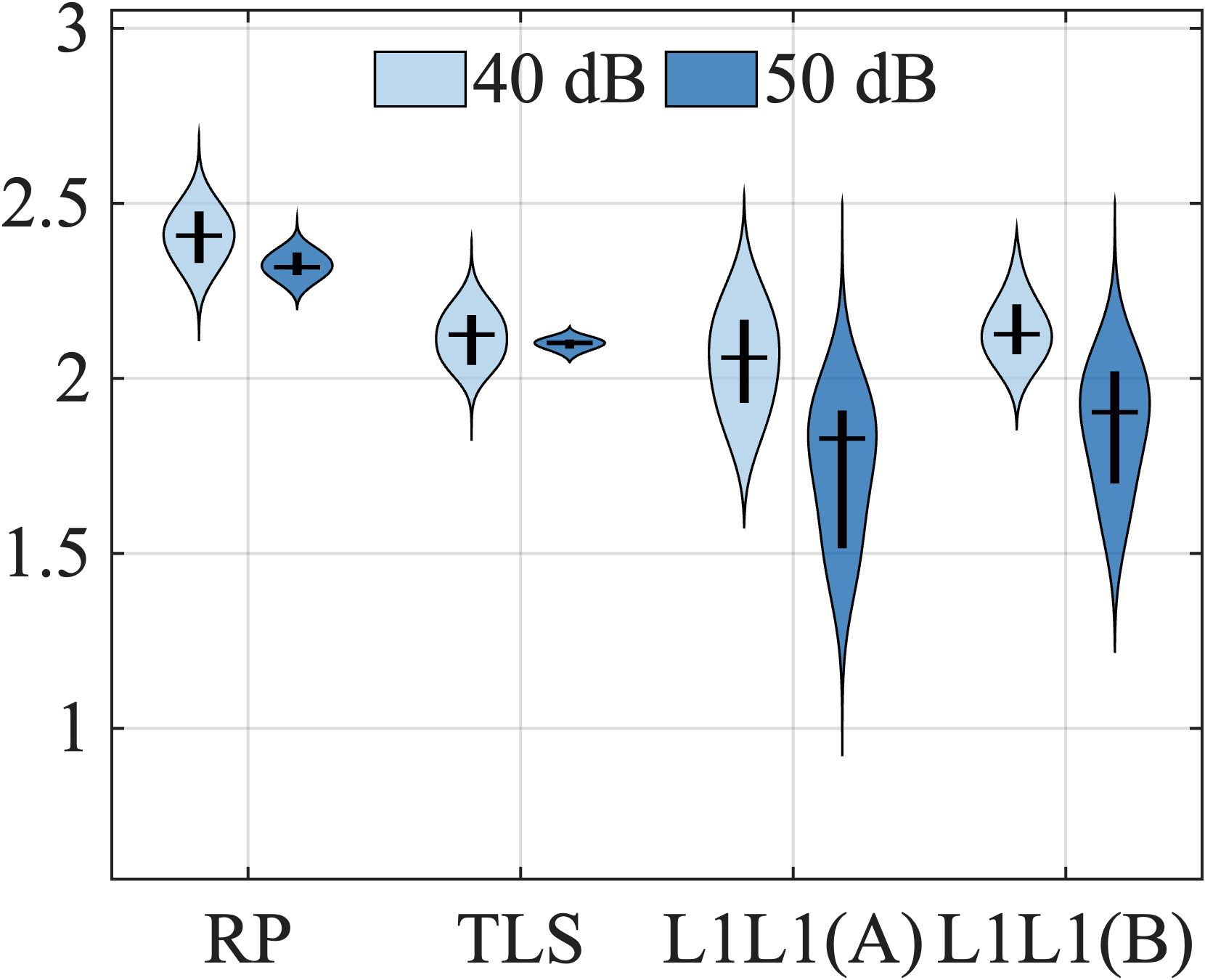} \\
            Focus $\Gamma$ (A/m\textsuperscript{2})
        \end{minipage}
        \begin{minipage}{2.87cm}
            \centering
            \includegraphics[width=2.86cm]{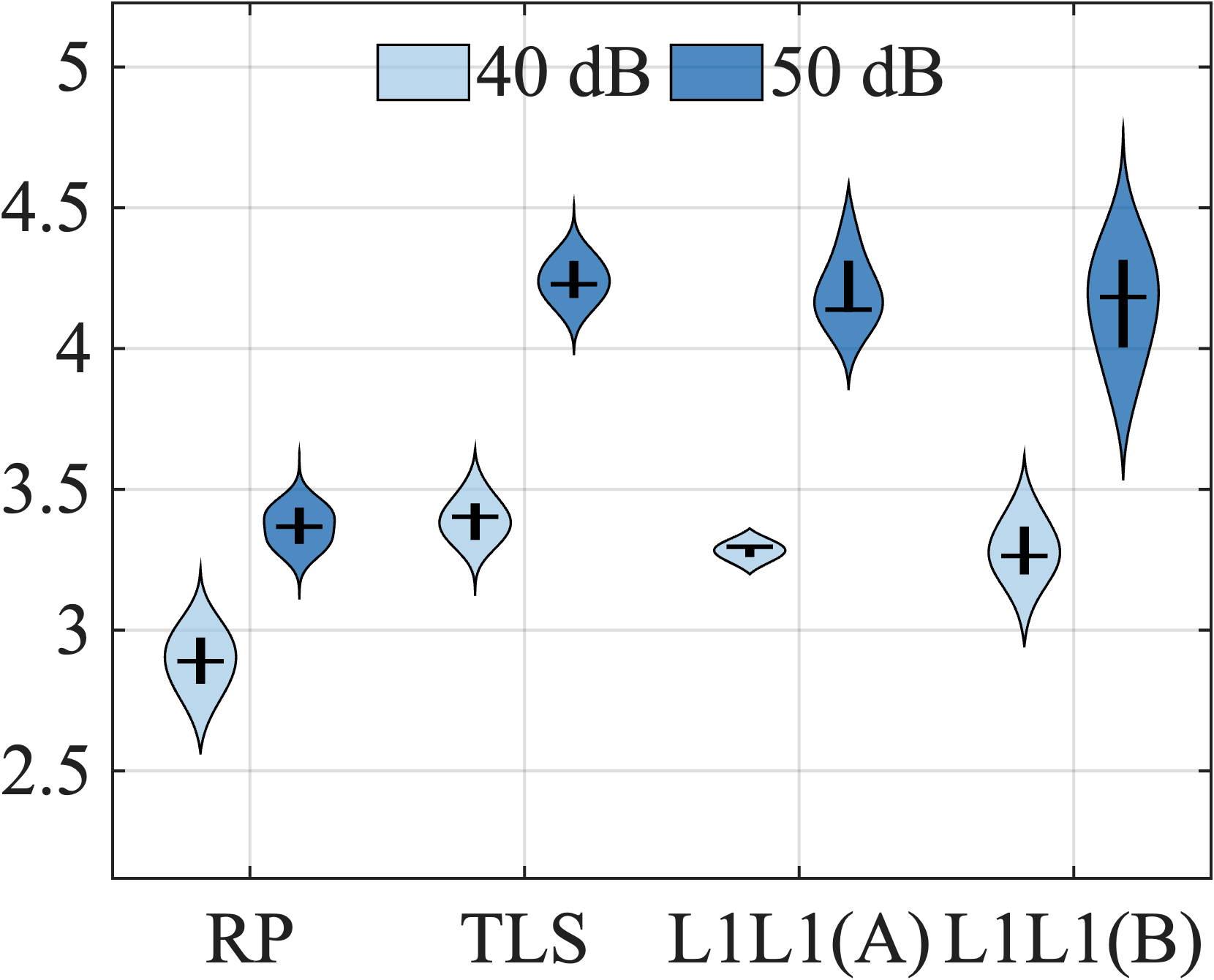} \\
            Ratio $\Theta$ (rel.)
        \end{minipage}
        \begin{minipage}{2.87cm}
            \centering
            \includegraphics[width=2.86cm]{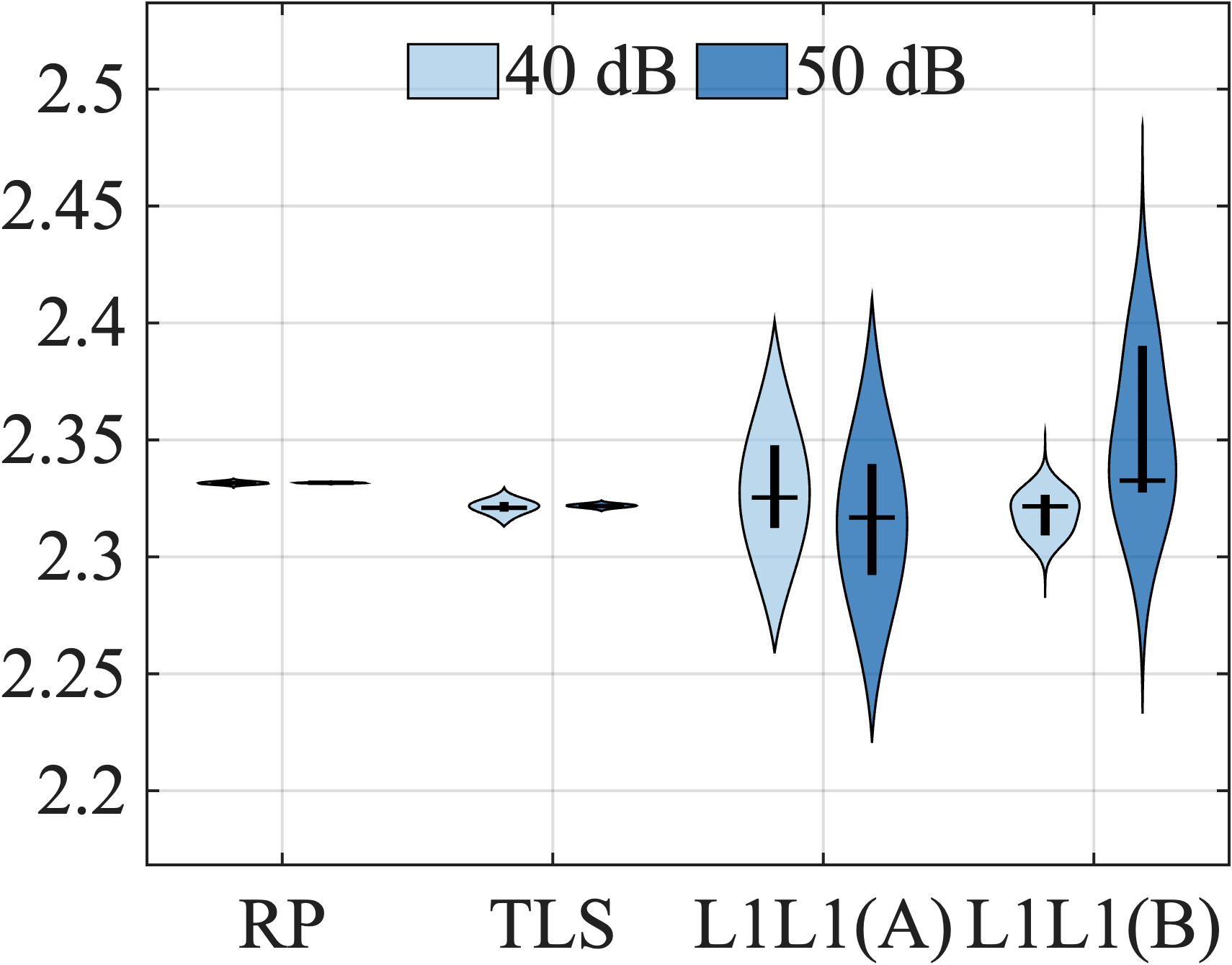} \\
            Targeting error (mm)
        \end{minipage} \\

        \vskip0.1cm

        \begin{minipage}{0.15cm}
            \centering
            \rotatebox{90}{Perpendicular}
        \end{minipage}
        \begin{minipage}{2.87cm}
            \centering
            \includegraphics[width=2.86cm]{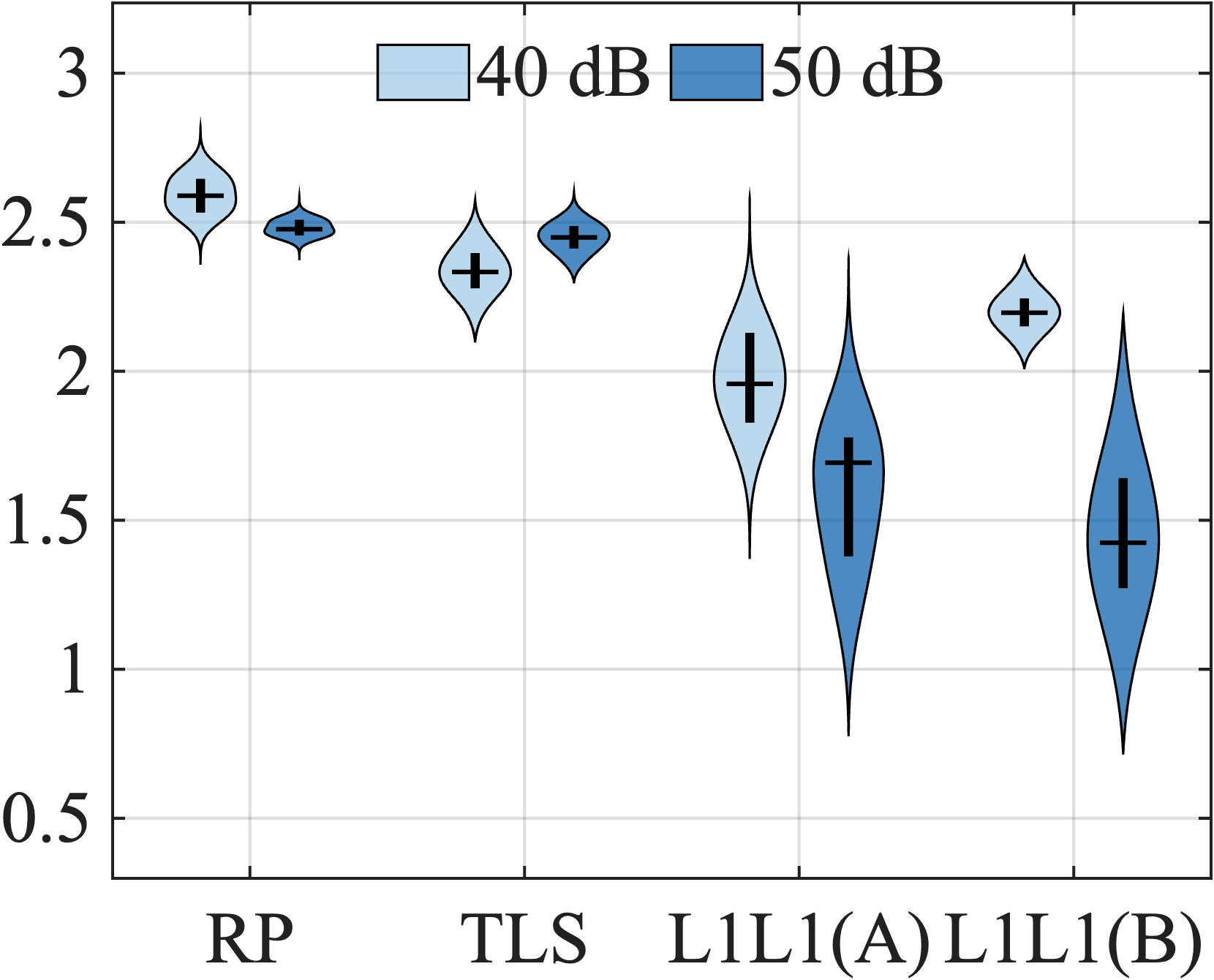} \\
            Focus $\Gamma$ (A/m\textsuperscript{2})
        \end{minipage}
        \begin{minipage}{2.87cm}
            \centering
            \includegraphics[width=2.86cm]{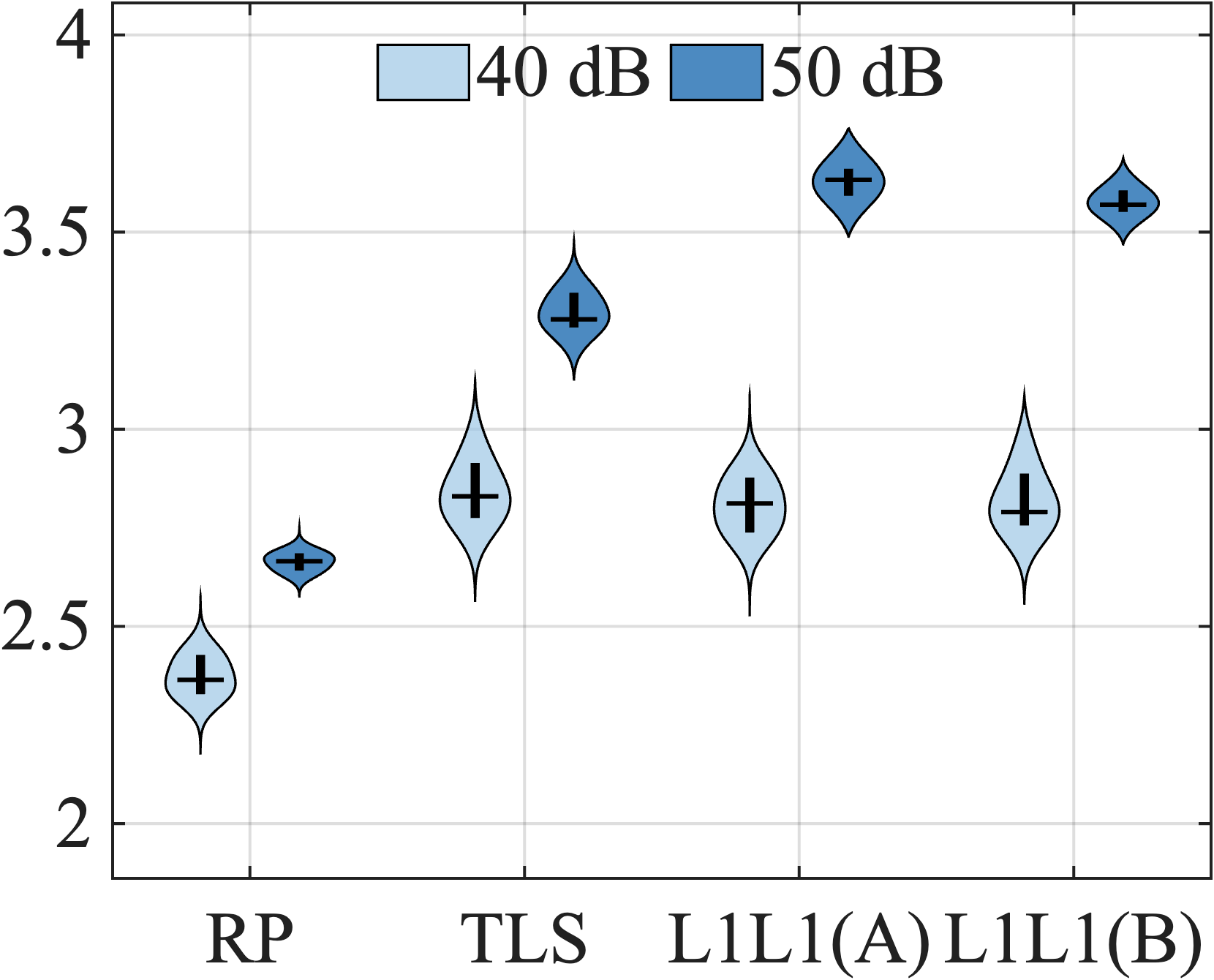} \\
            Ratio $\Theta$ (rel.)
        \end{minipage}
        \begin{minipage}{2.87cm}
            \centering
            \includegraphics[width=2.86cm]{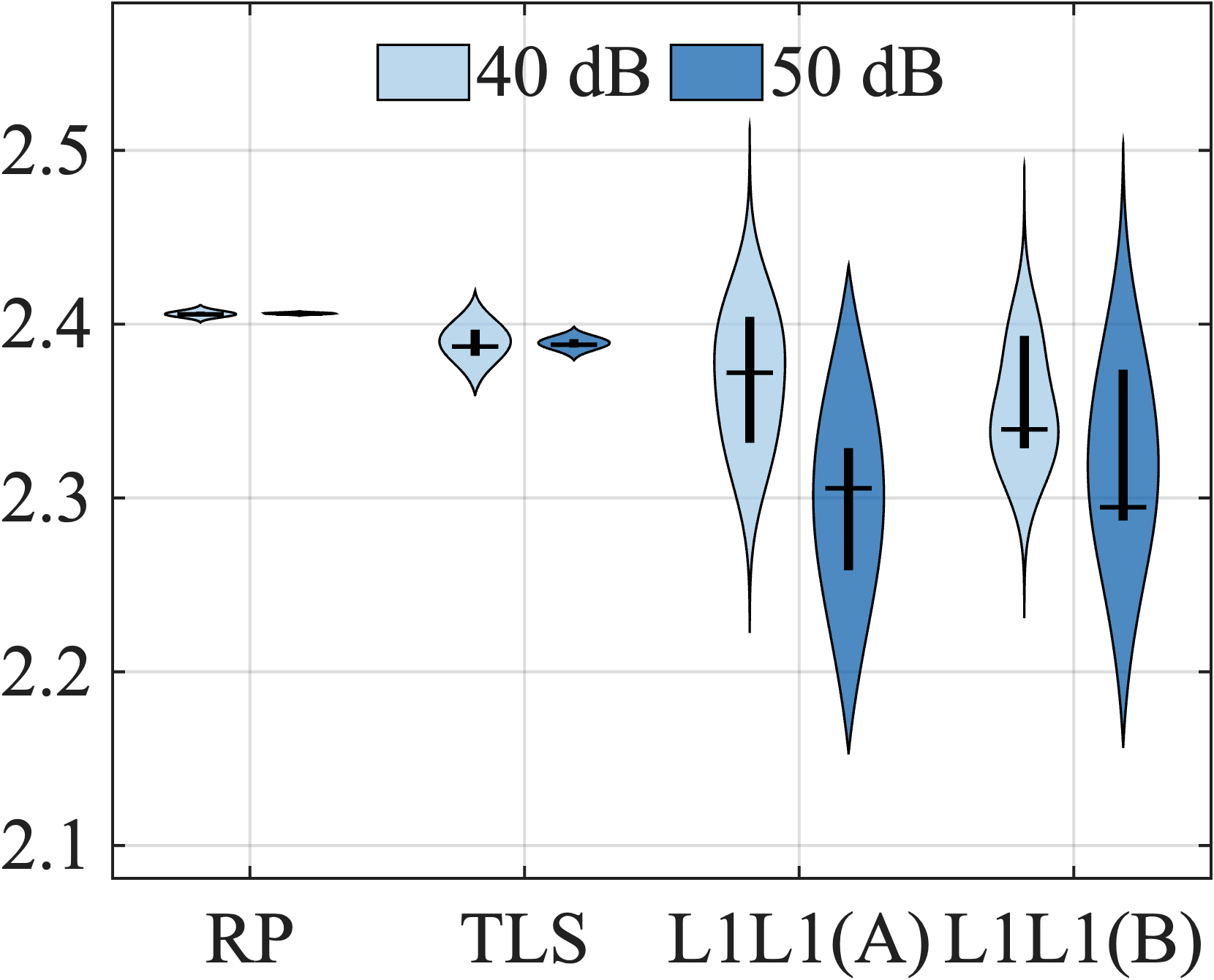} \\
            Targeting error (mm)
        \end{minipage} \\

        \vskip0.1cm
        Target (I)
    \end{minipage}
    \begin{minipage}{9.12cm}
        \centering
        \begin{minipage}{0.15cm}
            \centering
            \rotatebox{90}{Parallel}
        \end{minipage}
        \begin{minipage}{2.87cm}
            \centering
            \includegraphics[width=2.86cm]{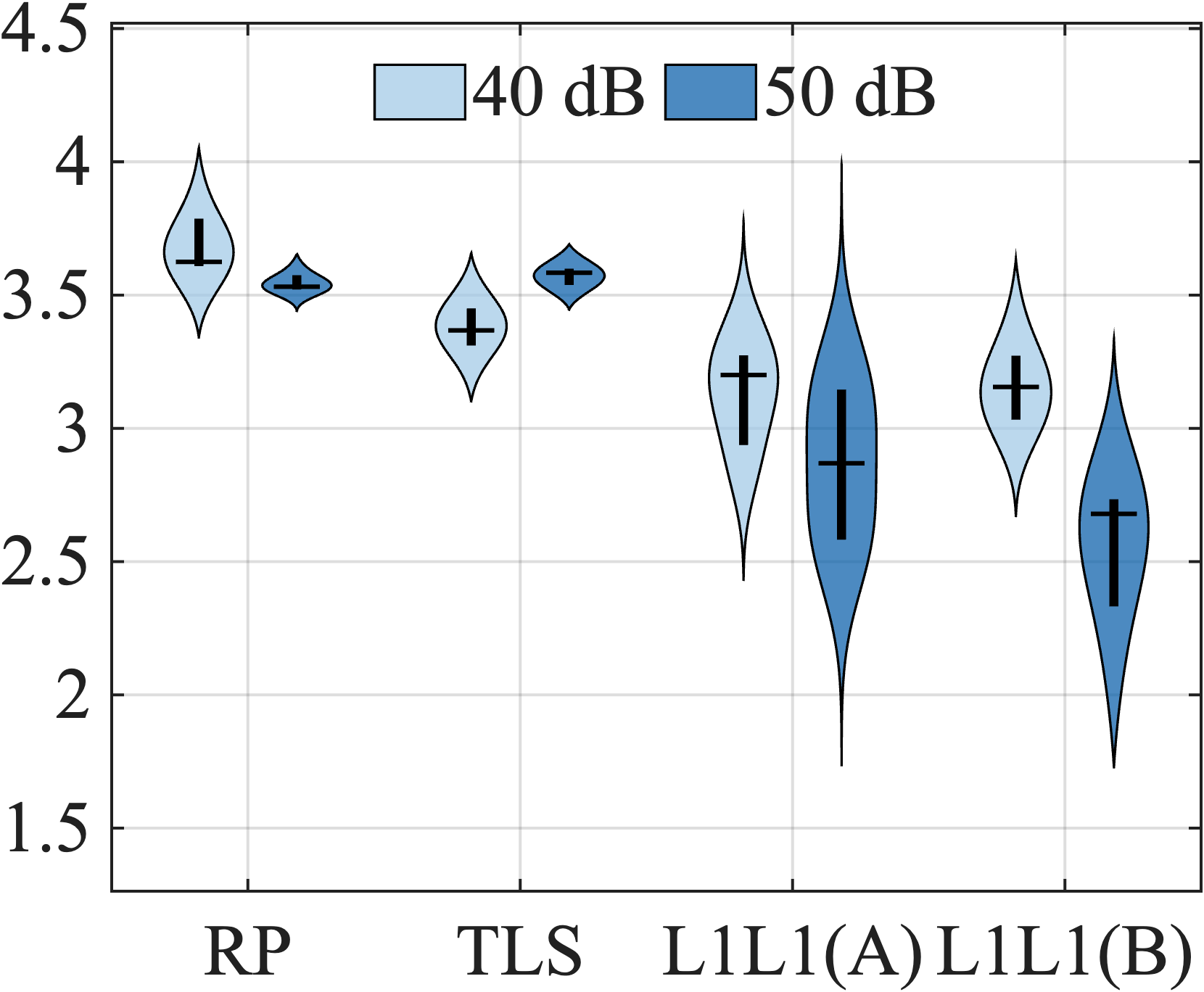} \\
            Focus $\Gamma$ (A/m\textsuperscript{2})
        \end{minipage}
        \begin{minipage}{2.87cm}
            \centering
            \includegraphics[width=2.86cm]{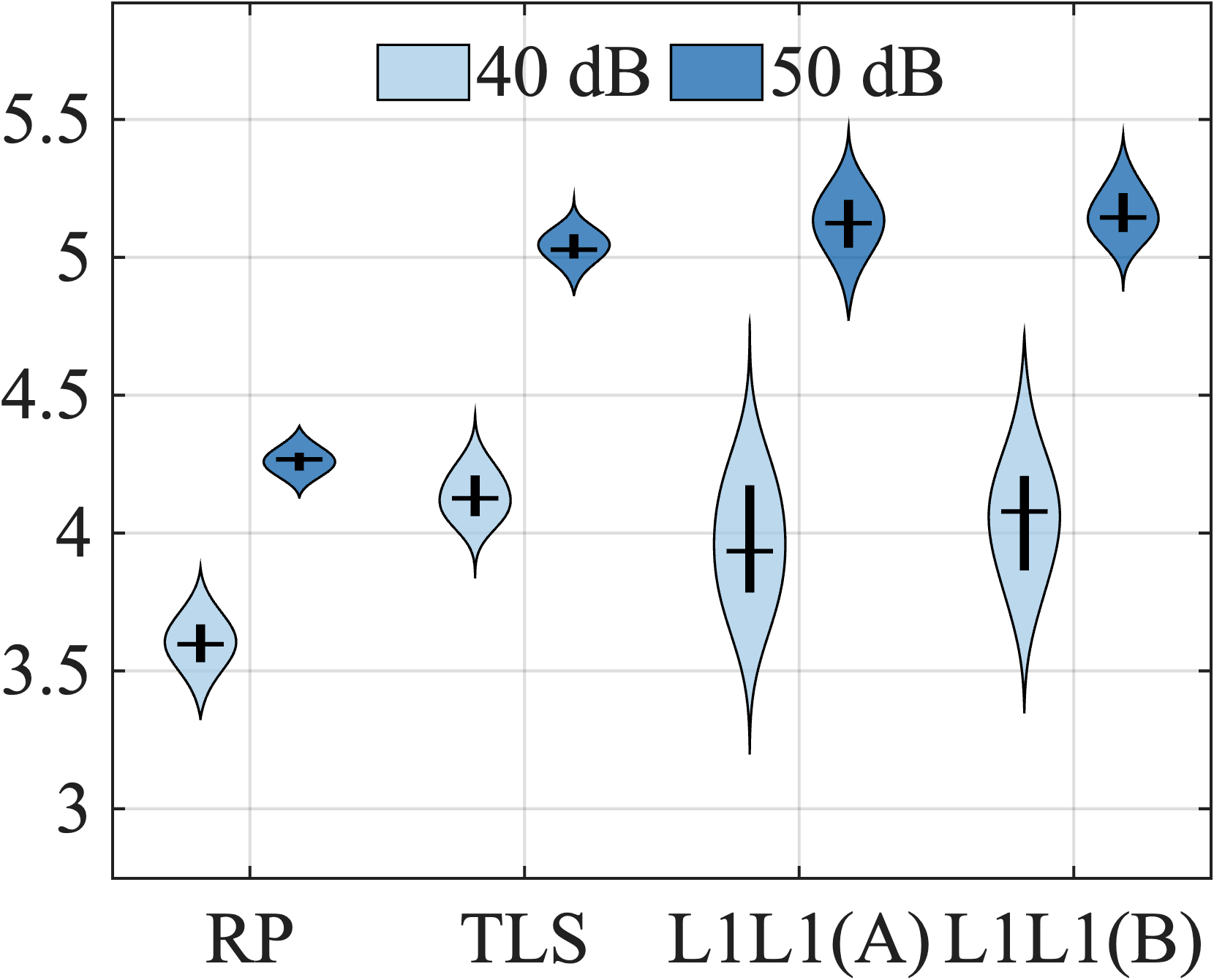} \\
            Ratio $\Theta$ (rel.)
        \end{minipage}
        \begin{minipage}{2.87cm}
            \centering
            \includegraphics[width=2.86cm]{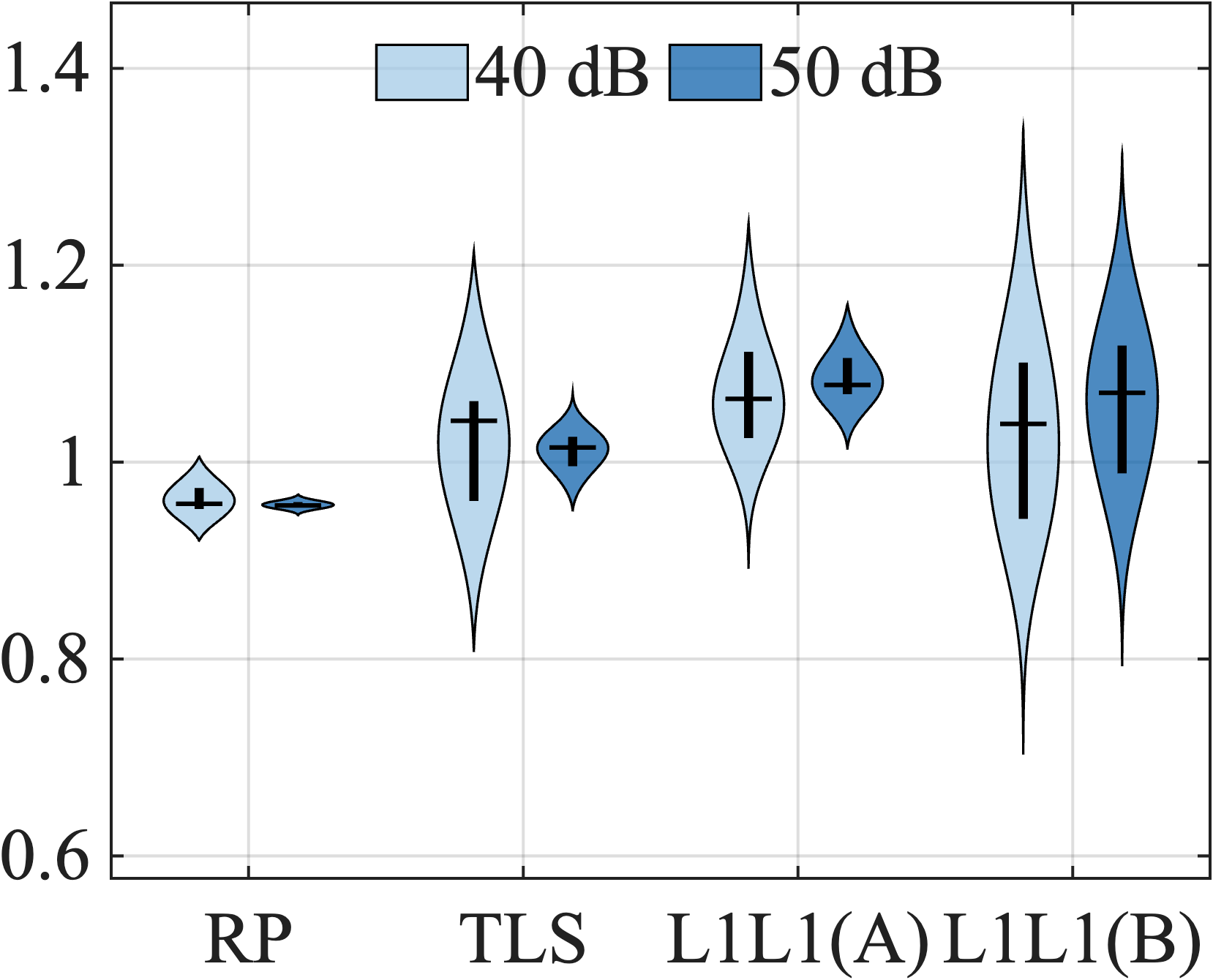} \\
            Targeting error (mm)
        \end{minipage} \\

        \vskip0.1cm

        \begin{minipage}{0.15cm}
            \centering
            \rotatebox{90}{Perpendicular}
        \end{minipage}
        \begin{minipage}{2.87cm}
            \centering
            \includegraphics[width=2.86cm]{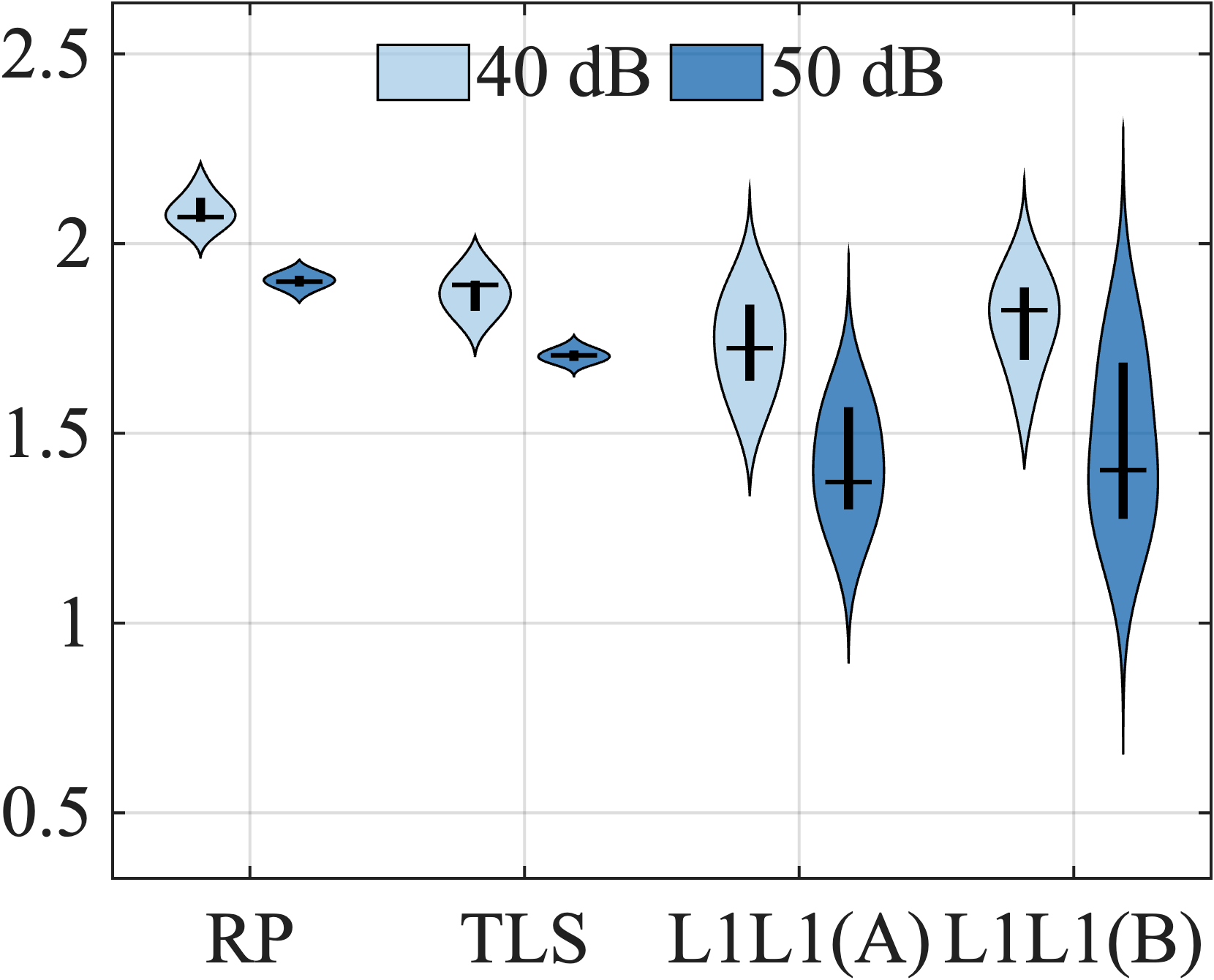} \\
            Focus $\Gamma$ (A/m\textsuperscript{2})
        \end{minipage}
        \begin{minipage}{2.87cm}
            \centering
            \includegraphics[width=2.86cm]{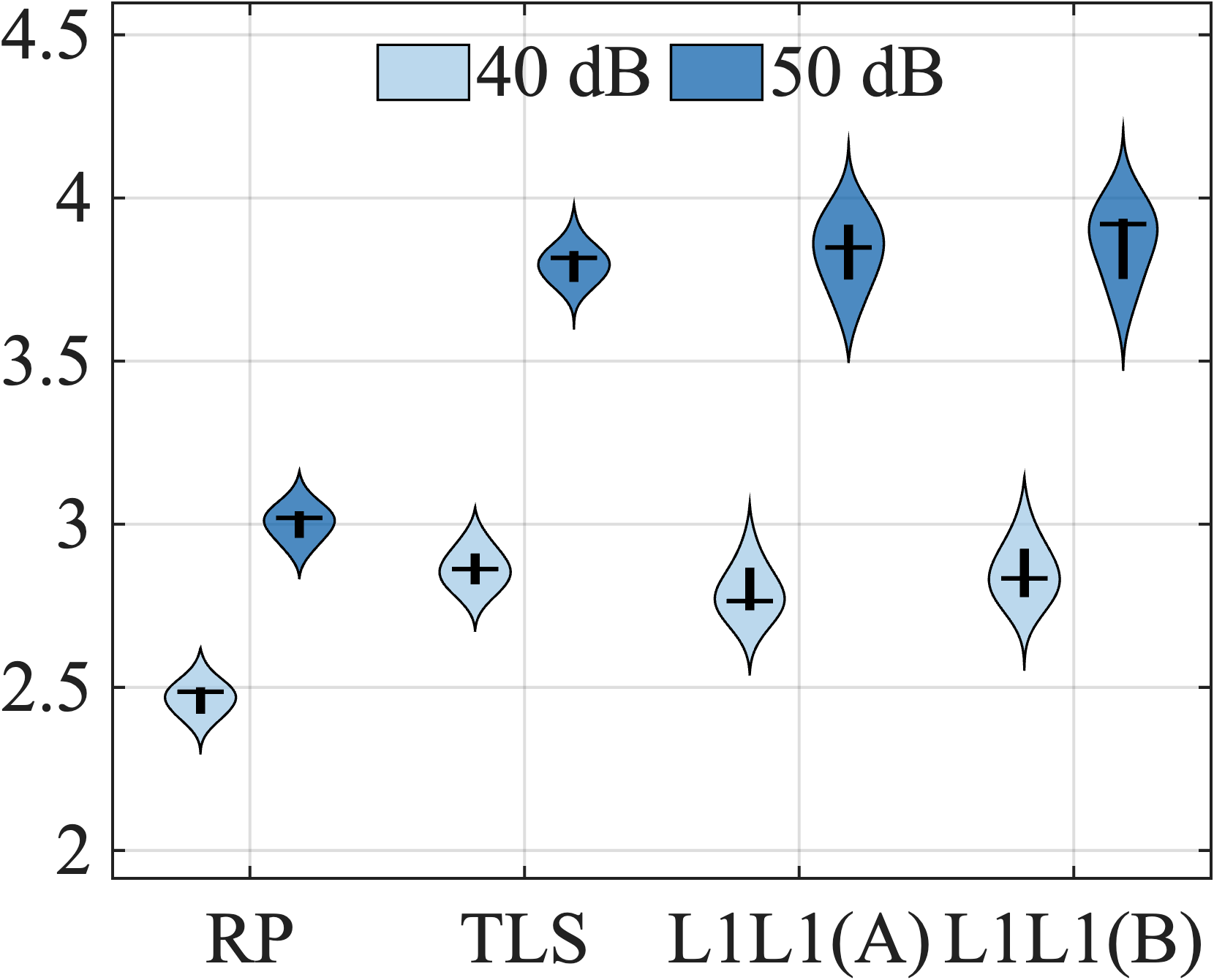} \\
            Ratio $\Theta$ (rel.)
        \end{minipage}
        \begin{minipage}{2.87cm}
            \centering
            \includegraphics[width=2.86cm]{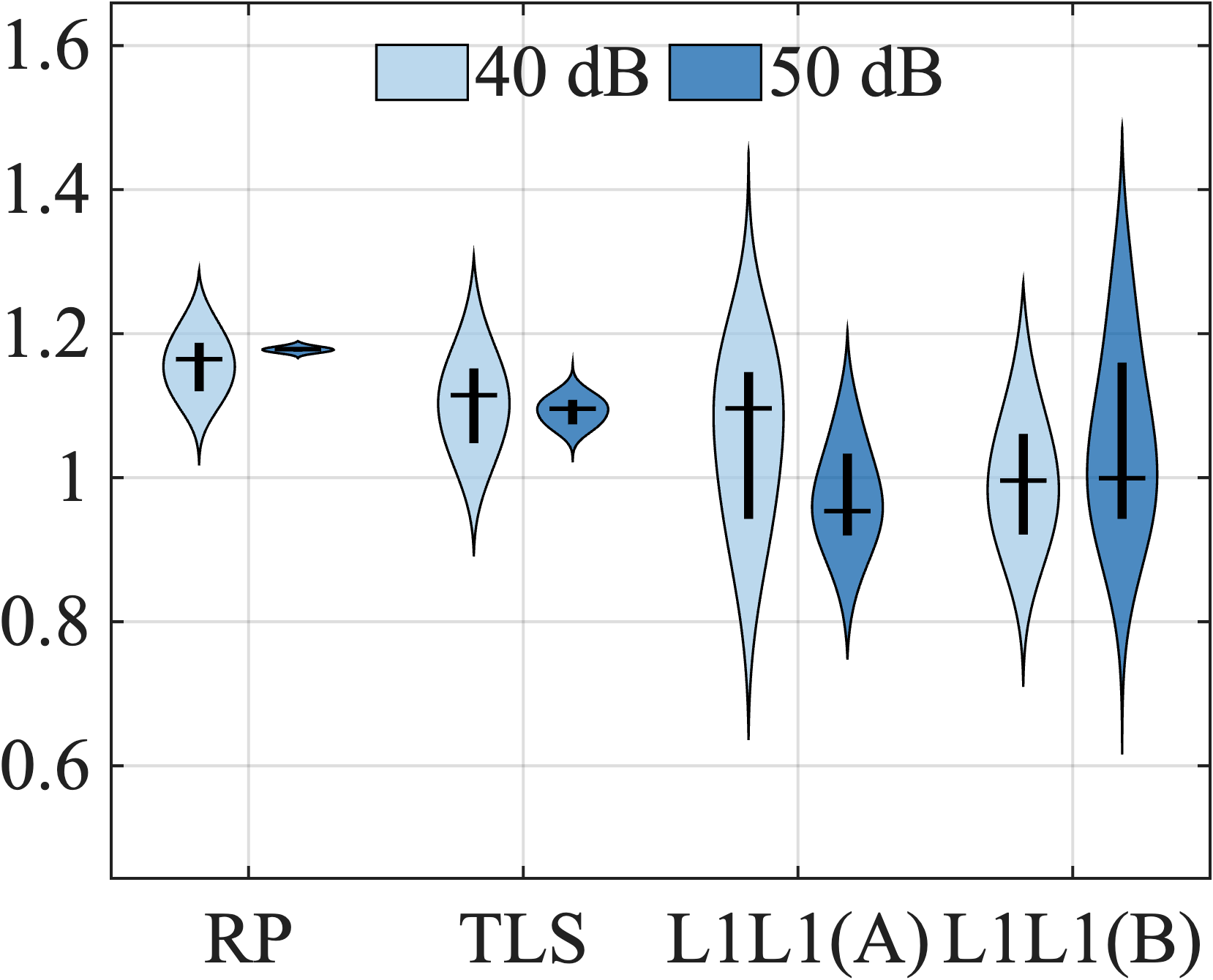} \\
            Targeting error (mm)
        \end{minipage} \\

        \vskip0.1cm
        Target (II)
    \end{minipage}

    \end{scriptsize}

    \caption{\textcolor{black}{Violin plots showing the focused current density $\Gamma$, field ratio $\Theta$, and targeting error (the distance between the current density maximizer and the actual target position) obtained with a noise-corrupted HR-LF and the 8- and 40-channel electrode configurations for test targets (I) and (II). The optimizers have been found through RP, L1L1(A), L1L1(B), and TLS methods. Each bar collates a sample of twenty decision variable values corresponding to optimizers obtained with as many independent lead field noise realizations (PSNR 40 or 50 dB). The violin width represents a kernel density estimate of the data distribution based on observations within the interval $[Q_1 - 1.5 \times \mathrm{IQR},, Q_3 + 1.5 \times \mathrm{IQR}]$, where $Q_1$ and $Q_3$ denote the first (25th percentile) and third (75th percentile) quartiles, respectively, and $\mathrm{IQR}=Q_3-Q_1$. Observations outside this interval are treated as outliers and excluded from the density estimate. The vertical black bar spans the interquartile range from $Q_1$ to $Q_3$, and the horizontal line indicates the median (50th percentile).  Targeting error is given in millimetres and measures the spatial discrepancy between the obtained and target activation locations in a spherical region of interest with 6 mm radius  cocentric with the target. The results obtained with L1L1(A) and L1L1(B) demonstrate a comparable sensitivity under lead field noise, suggesting that the limited dynamic range of L1L1(B) does not downgrade the L1L1 method's performance in situations involving uncertainty. Furthermore, both L1L1(A) and L1L1(B) find an overall elevated field ratio compared to RP and TLS with the greater PSNR level 50 dB, which allows for an improved +10 dB dynamic range for the optimization algorithm as compared to PSNR 40 dB. With respect to localization error, both L1L1 variants generally achieve accuracy comparable to or slightly better than RP and TLS across the considered targets and electrode configurations. In particular, L1L1(B) often exhibits a slightly more compact error distribution and fewer extreme outliers than L1L1(A), indicating that the reduced dynamic range constraint does not impair localization performance and may provide a modest stabilization effect under noisy lead field perturbations.}}
    \label{Fig:boxplot_noise_sensitivity}
\end{figure*}

Figure~\ref{Fig:boxplot_noise_sensitivity} further examines the noise sensitivity of the applied optimization methods using two specific target positions, illustrated in Figure~\ref{fig:target_positions_I_and_II}, covering a sample of twenty different noise realizations for both PSNR 40 and 50~dB. The results demonstrate that in a noisy setting, L1L1(A) and L1L1(B) produce similar distributions for \( \Gamma \), \( \Theta \), and localization accuracy \textcolor{black}{(defined here as the distance between the current density maximizer and the actual target position)}. Compared to RP and TLS, the L1L1 approach generally yields an elevated field ratio \( \Theta \), with the difference being marginal or absent at PSNR 40~dB but becoming more apparent at PSNR 50~dB. Accordingly, the potential benefit of using L1L1 depends critically on the level of modeling uncertainty and diminishes when the uncertainty becomes excessive. \textcolor{black}{With respect to localization performance, both L1L1 variants achieve accuracy that is comparable to or slightly better than RP and TLS across the considered targets and electrode configurations. The median targeting errors (in millimeters) remain similar between L1L1(A) and L1L1(B), suggesting that the reduced dynamic range constraint does not impair localization capability and may provide modest stabilization under noisy lead field perturbations.} Of the four optimization approaches tested, L1L1(B) provides the closest agreement between noiseless and noisy outcomes and proves particularly effective in maximizing the field ratio \( \Theta \) in the presence of noise.

\begin{table}[ht!]
\begin{footnotesize}
\centering
\caption{Maximum focused current density and field ratio described as $\Gamma_{max}$ and $\Theta_{max}$, respectively, for a 8- and 40-contacts multipolar DBS lead using the L1L1 method over a low-resolution lead field (LR-LF), with current dipole targeting aligned parallel and perpendicular to the DBS lead, for both noiseless (No Noise) and noisy (PSNR 40 dB) settings.}
\label{Tab:Together}
\renewcommand{\arraystretch}{1.25}
\resizebox{\columnwidth}{!}{%
\begin{tabular}{|ll|cccc|}
\hline
\multicolumn{2}{|c|}{\multirow{2}{*}{\textbf{8-contacts}}} & \multicolumn{4}{c|}{\textbf{Multipolar}} \\ \cline{3-6} 
\multicolumn{2}{|c|}{}                                     & \multicolumn{2}{c|}{L1L1(A), $\varepsilon \in [-160, 0]$} & \multicolumn{2}{c|}{L1L1(B), $\varepsilon \in  [-10, 0]$} \\ \hline
\thickhline
\multicolumn{1}{|l|}{\textbf{DV}} & \textbf{Orientation} & \multicolumn{1}{c|}{\textbf{No Noise}} & \multicolumn{1}{c|}{\textbf{PSNR 40 dB}} & \multicolumn{1}{c|}{\textbf{No Noise}} & \textbf{PSNR 40 dB} \\ \hline
\multicolumn{1}{|l|}{\multirow{2}{*}{$\Gamma_{max}$}} & Parallel      & \multicolumn{1}{c|}{2.42} & \multicolumn{1}{c|}{2.27} & \multicolumn{1}{c|}{2.27} & 2.27 \\ \cline{2-6} 
\multicolumn{1}{|l|}{}                          & Perpendicular & \multicolumn{1}{c|}{3.04} & \multicolumn{1}{c|}{2.35} & \multicolumn{1}{c|}{2.28} & 2.34 \\ \hline
\multicolumn{1}{|l|}{\multirow{2}{*}{$\Theta_{max}$}} & Parallel      & \multicolumn{1}{c|}{21.53} & \multicolumn{1}{c|}{5.56} & \multicolumn{1}{c|}{13.08} & 5.48 \\ \cline{2-6} 
\multicolumn{1}{|l|}{}                          & Perpendicular & \multicolumn{1}{c|}{27.07} & \multicolumn{1}{c|}{7.49} & \multicolumn{1}{c|}{13.18} & 7.22 \\ \hline
\end{tabular}
}

\vskip0.5cm

\resizebox{\columnwidth}{!}{%
\begin{tabular}{|ll|cccc|}
\hline
\multicolumn{2}{|c|}{\multirow{2}{*}{\textbf{40-contacts}}} & \multicolumn{4}{c|}{\textbf{Multipolar}} \\ \cline{3-6} 
\multicolumn{2}{|c|}{}                                     & \multicolumn{2}{c|}{L1L1(A), $\varepsilon \in  [-160, 0]$} & \multicolumn{2}{c|}{L1L1(B), $\varepsilon \in [-10, 0]$} \\ \hline
\thickhline
\multicolumn{1}{|l|}{\textbf{DV}} & \textbf{Orientation} & \multicolumn{1}{c|}{\textbf{No Noise}} & \multicolumn{1}{c|}{\textbf{PSNR 40 dB}} & \multicolumn{1}{c|}{\textbf{No Noise}} & \textbf{PSNR 40 dB} \\ \hline
\multicolumn{1}{|l|}{\multirow{2}{*}{$\Gamma_{max}$}} & Parallel      & \multicolumn{1}{c|}{2.20} & \multicolumn{1}{c|}{1.74} & \multicolumn{1}{c|}{2.02} & 2.24 \\ \cline{2-6} 
\multicolumn{1}{|l|}{}                          & Perpendicular & \multicolumn{1}{c|}{2.24} & \multicolumn{1}{c|}{2.29} & \multicolumn{1}{c|}{3.37} & 2.23 \\ \hline
\multicolumn{1}{|l|}{\multirow{2}{*}{$\Theta_{max}$}} & Parallel      & \multicolumn{1}{c|}{526.44} & \multicolumn{1}{c|}{7.56} & \multicolumn{1}{c|}{9.17} & 8.49 \\ \cline{2-6} 
\multicolumn{1}{|l|}{}                          & Perpendicular & \multicolumn{1}{c|}{536.05} & \multicolumn{1}{c|}{9.96} & \multicolumn{1}{c|}{10.97} & 8.47 \\ \hline
\end{tabular}%
}

\end{footnotesize}
\end{table}

\section{Discussion}
\label{sec:discussion}
In this study, we explore further application of the recently introduced metaheuristic L1L1 method \cite{galazprieto_2022_L1vsL2} to optimize current patterns and electrode configurations for directional DBS leads, similar to our previous study concerning optimization of multi-channel transcranial electric stimulation (MC-tES). For comparison, we also include optimization performance from applying {\em Reciprocity Principle} (RP) ~\cite{FERNANDEZCORAZZA2020116403} and {\em Ti\-kho\-nov Regularized Least Squares} (TLS)~\cite{dmochowski2011optimized} methods for retrieving bipolar or multipolar electrode configurations. \textcolor{black}{The former approach derive stimulation patterns directly from the forward electric field, whereas the latter least-squares formulation determine electrode currents by minimizing the mismatch between desired and predicted field distributions. The L1L1 method differs from these formulations in that it incorporates {\em a priori} uncertainty through an explicit detectability constraint on the lead-field components.} The L1L1 method was applied in two variants, L1L1(A) constrains \( \varepsilon\) to the interval \( [-160, 0] \) dB and L1L1(B) to \( [-10, 0] \) dB reflecting different levels of {\em a priori} modeling uncertainty \cite{schmidt2016uncertainty, athawale2019statistical}.

While the reciprocity principle implies that the maximal focused current density $\Gamma_{\mathrm{max}}$ in a linear quasi-static system is achieved by a two-point stimulation pattern, this result pertains to the unconstrained optimization of the stimulation vector and does not account for hardware-imposed spatial constraints. In MC-tES, where electrode locations can be adjusted, such extremal solutions can often be physically realized. In contrast, DBS employs fixed, implanted lead geometries, restricting the feasible stimulation space to current distributions across predefined contacts. Consequently, the localization and directionality of stimulation in DBS are governed jointly by (i) the spatial configuration of the lead, (ii) the orientation and position of the target dipolar source within the volumetric head model, and (iii) the number and configuration of active channels. Under these constraints, multi-contact current steering becomes necessary to approximate the ideal dipolar field, enabling directional control and improved focality. In this study, this effect is explicitly examined by comparing conventional 8-contact leads with a high-density 40-contact configuration, demonstrating how increased channel count expands the attainable subspace of stimulation patterns and improves the trade-off between intensity and focality. The three metrics used in our optimization framework to characterize the spatial current distributions obtained from the forward model are: {\em focused current density} \( \Gamma \), {\em nuisance current density} \( \Xi \), and {\em field ratio} \( \Theta = \Gamma / \Xi \). Within this constrained setting, the L1L1 method does not seek a purely extremal two-contact solution, but instead identifies sparse multi-contact configurations that optimally approximate the desired dipolar field while satisfying current balance and safety constraints.

\textcolor{black}{Uncertainty is introduced through perturbations of the lead-field matrix representing aggregate forward-model variability arising from multiple factors including numerical modeling discrepancies. Table~\ref{tab:uncertainty} summarizes representative magnitudes of modeling variability reported in DBS studies and provides quantitative context for the uncertainty assumptions adopted in this work. Variability arising from neural activation modeling has been reported in the range of 24--47\% ($-12.4$ to $-6.6$\,dB) for axonal activation thresholds and about 7--12\% ($-23$ to $-18$\,dB) for population activation metrics \cite{FrankemolleGilbert2021}. Differences associated with electrode geometry and stimulation configuration are reported around 10--20\% ($-20$ to $-14$\,dB) \cite{butson_2007_patient}. Tissue conductivity modeling typically produces smaller differences of about 3.8--5.6\% ($-28$ to $-25$\,dB) \cite{zhang2020steering}. Comparisons between simulations and intraoperative measurements show discrepancies of approximately 8--15\% ($-22$ to $-16$\,dB) \cite{Alonso2018}. Electrode localization uncertainty derived from imaging reconstruction can introduce positional errors corresponding to about 10--40\% variability in predicted stimulation fields \cite{Lasica2025LocalizationPrecision, athawale2019statistical}. In this work, these heterogeneous variability sources are interpreted as aggregate perturbations affecting the lead-field matrix, motivating the stochastic uncertainty model introduced in Section~\ref{sec:uncertain}; the detectability threshold is expressed in decibels (dB) and quantifies the relative attenuation of the lead-field response with respect to the target stimulation location..
}

\textcolor{black}{When the detectability constraint is relaxed (i.e., L1L1(A) case), the optimization procedure results in unrealistically high field‑ratio values, approximately fifty times larger than those obtained under a stricter threshold (i.e., L1L1(B) case). By introducing this detectability constraint, we are able to restrict the feasible set of DOFs to spatial locations whose lead‑field magnitude remains distinguishable under forward‑model perturbations, preventing the optimizer from exploiting weak lead‑field components that are highly sensitive to modeling uncertainty. In this framework, the feasible set \( \mathcal{V}_\delta \) corresponds to spatial locations where the lead‑field amplitude attenuates by a factor \( \delta \) for \( \delta \in [-40,-10] \) dB. These attenuation levels are consistent with spatial decay patterns reported in DBS modeling studies. In particular, attenuation around \( \delta \approx -40 \) dB corresponds to spatial regions comparable to a VTA radius of approximately 4.0~mm, as reported by Butson and McIntyre \cite{butson_2007_patient}, Miocinovic et~al.\ \cite{miocinovic_2006_subthalamic}, and M{\"a}dler and Coenen \cite{madler2012activated}. Consequently, the PSNR values used (40~dB and 50~dB) provide reasonable bounds for forward‑model uncertainty when defining the detectability threshold in the optimization framework. At PSNR 40~dB the dynamic range of \( \Theta \) is limited by the maximum lead‑field attenuation within the VTA, whereas PSNR 50~dB supports the assumption that each position inside the VTA can be independently and coherently targeted.
}

When comparing electrode configurations, increasing the number of contact points generally improves adaptability and control over the resulting electric field. Parallel configurations, i.e., alignment of the current dipole with the DBS lead, produced localized activation near the lead tip, whereas perpendicular configurations resulted in a broader and medially shifted activation pattern. This behavior can be attributed to shorter distances between electrodes and target locations, improved alignment between the induced field and the target current direction, and increased flexibility arising from a larger number of controllable degrees of freedom.

Overall, the findings provide a proof-of-concept for integrating {\em a priori} uncertainty modeling into numerically simulated DBS optimization. \textcolor{black}{Sources of uncertainty such as anatomical variability, tissue conductivity changes, and electrode placement errors can be viewed as perturbations of the lead-field matrix used in the optimization, consistent with prior DBS uncertainty analyses and localization studies \cite{schmidt2016uncertainty, athawale2019statistical, Lasica2025LocalizationPrecision}. Within this framework, the L1L1 method can improve focality and directional control over the stimulation, allowing more selective targeting of desired regions while reducing off-target activation in the computational model. These results therefore provide a computational proof-of-concept for incorporating uncertainty-aware constraints into DBS optimization beyond the deterministic forward-model assumptions often used in electric field modeling studies.}

\subsection{Limitations and Future Work}
\label{sec:Future}

While the current uncertainty modeling strategy provides a principled framework for constraining the optimization based on lead field attenuation, it also has several limitations. First, the threshold values used to define acceptable nuisance current levels are heuristic and may not fully capture the complexity of individual patient anatomies or tissue conductivities. This simplification assumes spatial uniformity in attenuation, which may not hold, for example, in the presence of anisotropic or heterogeneous tissue properties \cite{dannhauer2011modeling}. Additionally, the approach relies on a static surrogate model of the VTA and does not directly incorporate a probabilistic model of uncertainty \cite{schmidt2016uncertainty, athawale2019statistical}, which would allow for a more comprehensive treatment of variability. Potential alternative surrogate approaches exist, such as, threshold-based curve fitting \cite{Butson_2006_ElectrodeDesign}, and statistical emulators such as Gaussian process classifiers \cite{DeLaPava_2017_accelerating, Orozco_2019_KernelDBS}, which approach VTA in different ways.

The Gaussian noise model is adopted because it furnishes the maximally uninformative representation of uncertainty under a finite-variance constraint \cite{Jaynes2003, CoverThomas2006}. Given that the present optimization framework does not specify a detailed probabilistic law for perturbations, this choice avoids introducing unsupported structure and provides a neutral, information-theoretic baseline. More sophisticated uncertainty quantification techniques, such as Markov Chain Monte Carlo (MCMC)-based Bayesian sampling \cite{kaipio2006statistical} or deep-learning approaches to predictive uncertainty \cite{schmidt2016uncertainty, athawale2019statistical}, could be employed when richer data are available.

Compared to the reference methods, the L1L1 method is slightly more demanding in terms of computing time. However, the increase is not significant and does not limit its practical use. Other mathematical optimization techniques, such as TLS, may offer promising alternatives alongside the L1L1 method for further refining the {\em a priori} uncertainty modeling adopted in this study. A logical direction for future research is the development of {\em a posteriori} models, where the current approach could be validated in more advanced experimental or clinical settings. This would enable a deeper understanding of the potential hyperparameter realizations in real-world applications. Translating these computational findings into clinical practice will require rigorous validation through empirical testing and clinical trials. 

This manuscript does not address clinical implications; however, for the configuration studied, certain anatomical structures, such as the \textit{lamina medullaris} near the modeled lead trajectory, may require special consideration in clinical practice to avoid stimulation-induced side effects. Translating these methodological results to clinical applications will necessitate examining the relationship between the L1L1 method and uncertainties in individual factors such as lead placement, tissue conductivity and anisotropy, and patient-specific anatomical variations including thalamic segmentation. To support this, we employ finite element (FE) discretizations incorporating thalamic subdivisions through a boundary-fitted multi-compartment FE mesh together with a Complete Electrode Model (CEM) \cite{pursiainen2017forward}, as in the present study.

\section{Conclusions}
\label{sec:Conclusions}
This study represents an advancement in optimizing the volumetric current fields for DBS by leveraging the metaheuristic L1L1 method. Previous contact selection studies, such as \cite{Anderson_2018} and later \cite{Janson_2020_Activation}, have used convex optimization to determine optimal stimulation settings and this represents at addressing varying levels of {\em a priori} lead field uncertainty, which can be attributed as patient anatomy, misleading lead placement, and tissue properties, which often vary in clinical practice. 

Our findings demonstrate that the method can achieve relatively high field ratios with a carefully chosen activation threshold $\varepsilon$ value. This, in turn, has the potential to improve individualized stimulation outcomes. An ongoing work is to integrate patient-specific data, such as thalamic segmentation, to improve stimulation accuracy. This involves using boundary-fitted multi-compartment FE meshes \cite{galazprieto_2023_mesh} and a Complete Electrode Model \cite{pursiainen2017forward}, as in this study.

\section{Conflict of Interest}
\label{sec:conflict}
The authors certify that this study is a result of purely academic, open, and independent research. They have no affiliations with or involvement in any organization or entity with a financial interest or non-financial interest, such as personal or professional relationships, affiliations, knowledge, or beliefs in the subject matter or materials discussed in this manuscript.

\section{Author Contributions}
FGP: Conceptualization, Data curation, Formal analysis, Investigation, Methodology, Software, Visualization, Writing---original draft. AL: Conceptualization, Data curation, Investigation, Methodology, Project administration, Writing---original draft. MS: Formal analysis, Supervision, Validation, Writing---review \& editing. SP: Conceptualization, Formal analysis, Funding acquisition, Investigation, Methodology, Project administration, Resources, Software, Supervision, Validation, Visualization, Writing---review \& editing.

\section{Funding}
FGP, MS, and SP were supported by the Research Council of Finland (RCF) through the Center of Excellence in Inverse Modelling and Imaging 2018--2025 (RCF 353089), Exploratory Study for Radar Tomography of Dimorphos (RCF 359198), and the Flagship of Advanced Mathematics for Sensing, Imaging and Modelling (FAME) (RCF 359185); DAAD (German Academic Exchange Service) project "Non-invasively reconstructing and inhibiting activity in focal epilepsy" (354976, 367453); and the ERA Personalised diagnosis and treatment for refractory focal paediatric and adult epilepsy (PerEpi) project (RCF 344712). MS holds Academy Research Fellowship 2025--2029 (RCF 371055).

\section{Ethical Statement}
This work is based exclusively on mathematical modeling and numerical simulations, using an openly available MRI dataset \cite{piastra_maria_carla_2020_3888381} as input. No new data was collected, and no sensitive or personally identifiable information was processed at any stage of the study.

\section{Acknowledgments}
\label{sec:ack}
The authors would like to thank Prof. Dr. rer. nat. Carsten H. Wolters, researchers and clinicians from the Institute for Biomagnetism and Biosignalanalysis (IBB), and the PerEpi consortium for their continuous support, discussions, and feedback regarding non-invasive brain stimulation topics, regression analysis, and the seminars prepared during the elaboration of this study.

\appendix
\label{sec:Appendix}

\section{Mathematical Optimization Scheme}

\subsection{Formulation of the Lead Field Matrix}
\label{App:sec:Math_Lyx}
The linear forward mapping between electrode currents and the resulting volumetric current density distribution is given by the real-valued $ (N+M) \times K $ \emph{lead field matrix} $\mathbf{L}$:
\begin{equation}
    \mathbf{L}\,\mathbf{y} = \mathbf{x},
    \label{eq:leadfield_mapping}
\end{equation}
where $\mathbf{y} \in \mathbb{R}^{K}$ denotes the applied current pattern across $K$ electrode contacts, and $\mathbf{x} \in \mathbb{R}^{N+M}$ is the discretized volumetric current density with each entry corresponding to a local source element parameterized by three Cartesian \emph{dipolar degrees of freedom} (DOFs), representing current density along the $x$-, $y$-, and $z$-axes. For our optimization framework, $\mathbf{L}$ is partitioned into target and nuisance components:
\begin{equation}
    \mathbf{L} = 
    \begin{pmatrix}
        \mathbf{L}_1 \\
        \mathbf{L}_2
    \end{pmatrix},
    \qquad
    \mathbf{x} =
    \begin{pmatrix}
        \mathbf{x}_1 \\
        \mathbf{0}
    \end{pmatrix},
\end{equation}
where $\mathbf{L}_1 \in \mathbb{R}^{N\times K}$ models the target region, $\mathbf{L}_2 \in \mathbb{R}^{M\times K}$ models the nuisance region, and $\mathbf{x}_1 \in \mathbb{R}^{N}$ defines the desired field in the target region. To ensure that only the \textit{direction} of $\mathbf{x}_1$ is prescribed, we introduce the projection
\begin{equation}
    \mathbf{P} = \frac{\mathbf{x}_1 \mathbf{x}_1^{\mathsf{T}}}{\|\mathbf{x}_1\|_2^2}
    \in \mathbb{R}^{N\times N},
\end{equation}
which enforces alignment of the induced field with $\mathbf{x}_1$. The resulting reduced formulation is:
\begin{equation}
    \begin{pmatrix}
        \mathbf{P}\mathbf{L}_1 \\
        \mathbf{L}_2
    \end{pmatrix}
    \mathbf{y}
    =
    \begin{pmatrix}
        \mathbf{P}\mathbf{x}_1 \\
        \mathbf{0}
    \end{pmatrix}.
    \label{eq:projection_system_final}
\end{equation}
The first block enforces field orientation within the target region, while the second enforces suppression outside it.

%%%
\subsection{L1-norm Regularized L1-norm Fitting (L1L1)}
\label{App:sec:L1L1}
\textcolor{black}{The general form of the \emph{L1-norm regularized L1-norm fitting} (L1L1) method described by Galaz Prieto et al.~\citep{galazprieto_2022_L1vsL2} is a linear programming problem that determines the stimulation vector ${\bf y}$ so that the induced electric field ${\bf L}{\bf y}$ approximates a prescribed target field ${\bf x}$ while suppressing off-target stimulation. Using the lead-field formulation introduced in Sec.~\ref{App:sec:Math_Lyx}, the method minimizes a sparsity-regularized $L_1$ fitting functional defined on the focused and nuisance field components}. The task is
\begin{alignat}{2}
&& \underset{{\bf y}}{\textrm{minimize}} 
&\quad 
\left\|
\begin{pmatrix}
{\bf L}_1{\bf y}-{\bf x}_1 \\
\Psi_\varepsilon\!\left(\nu^{-1}{\bf L}_2{\bf y}\right)
\end{pmatrix}
\right\|_1
+ \alpha\,\zeta\,\|{\bf y}\|_1
\label{eq:L1L1Method_1}
\\
&& \textrm{subject to} 
&\quad {\bf y}\preceq \gamma{\bf 1},
\label{eq:L1L1Method_2}
\\
&& &\quad \|{\bf y}\|_1 \le \mu,
\label{eq:L1L1Method_3}
\\
&& &\quad \sum_{\ell=1}^{K} y_\ell = 0 .
\label{eq:L1L1Method_4}
\end{alignat}
\textcolor{black}{Constraint~\eqref{eq:L1L1Method_2}} limits the current amplitude at each electrode contact to $\gamma$. \textcolor{black}{Constraint~\eqref{eq:L1L1Method_3}} restricts the total injected current to the safety limit $\mu$. \textcolor{black}{Constraint~\eqref{eq:L1L1Method_4} enforces charge balance by requiring the net injected current across all electrodes to be zero}.
The regularization parameter $\alpha$ \textcolor{black}{controls} the $L_1$ regularization level relative to the scaling value $\zeta=\|{\bf L}\|_1$. The \textcolor{black}{operator}
\[
\Psi_\varepsilon[{\bf w}]_m=\max\{\,|w_m|,\varepsilon\,\}, 
\qquad m=1,\ldots,M ,
\]
introduces the \textcolor{black}{\emph{nuisance-field threshold} $0\le\varepsilon\le1$ with respect to} the scaling value $\nu=\|{\bf x}_1\|_\infty$. Consequently, nuisance-field entries satisfying $|({\bf L}_2{\bf y})_m|<\varepsilon\nu$ do not actively contribute to the objective function. The corresponding \emph{constraint support} is
\[
\{\,m:\,|({\bf L}_2{\bf y})_m|\ge\varepsilon\nu\,\},
\]
\textcolor{black}{which defines the subset of nuisance nodes that influence the optimization}. To evaluate stimulation performance and to quantify the trade-off between target and undesired current spread, we define three decision variables: the \emph{focused current density} $\Gamma$, the \emph{nuisance current density} $\Xi$, and the resulting \emph{field ratio} $\Theta$.

%%% Focus
The \emph{focused current density} $\Gamma$ (A/m\textsuperscript{2}) quantifies the degree of field alignment between the induced current distribution and the desired stimulation target. It is defined as
\begin{equation}
    \Gamma = \frac{{\bf x}_1^T {\bf L}_1 {\bf y}}{\| {\bf x}_1 \|_2},
    \label{eq:Gamma}
\end{equation}
where ${\bf x}_1 \in \mathbb{R}^N$ represents the target volumetric current density template specifying the spatial and directional characteristics of the desired field within the stimulation region. The inner product ${\bf x}_1^{T}{\bf L}_1{\bf y}$ measures the alignment between the induced field ${\bf L}_1{\bf y}$ and the prescribed target field pattern ${\bf x}_1$.

%%% Nuisance
The \emph{nuisance current density} $\Xi$ (A/m\textsuperscript{2}) characterizes the average magnitude of the induced field in off-target regions and is defined by
\begin{equation}
    \Xi = \frac{\| {\bf L}_2 {\bf y} \|_2}{\sqrt{M}}.
    \label{eq:Xi}
\end{equation}
Dividing the Euclidean norm $\|{\bf L}_2{\bf y}\|_2$ by $\sqrt{M}$ yields the root-mean-square (RMS) field amplitude per node,
\[
\textcolor{black}{\Xi=\sqrt{\frac{1}{M}\sum_{m=1}^{M}\big(({\bf L}_2{\bf y})_m\big)^2}}.
\]
This normalization makes $\Xi$ invariant to the number of discretization points in the nuisance domain. In cases where only a subset of nodes exhibits appreciable current density, the normalization may alternatively be restricted to the active \textit{constraint support}, i.e., replacing $M$ with the number of nodes satisfying $|({\bf L}_2{\bf y})_m|\ge\varepsilon\nu$.

%%% Ratio
The \emph{field ratio} $\Theta$ is introduced to characterize the overall selectivity of stimulation:
\begin{equation}
    \Theta = \frac{\Gamma}{\Xi}.
    \label{eq:Theta}
\end{equation}
A high $\Theta$ value reflects both strong current focusing within the target and minimal activation of non-target regions.

%%% Tikhonov-regularized
\subsection{Ti\-kho\-nov-Regularized Least Squares}
\label{App:sec:TLS}
\textcolor{black}{Tikhonov-regularized least-squares} (TLS) estimation \cite{dmochowski2011optimized} is formulated as the following quadratic optimization problem:
\begin{equation}
\min_{\mathbf{y}}
\left\{
\|\mathbf{L}_1\mathbf{y}-\mathbf{x}_1\|_2^2
+ \gamma^2\beta^2\|\mathbf{L}_2\mathbf{y}\|_2^2
+ \gamma^2\varsigma^2\|\mathbf{y}\|_2^2
\right\}.
\label{eq:TLS_minimization}
\end{equation}
where $\varsigma=\|\mathbf{L}\|_2$ \textcolor{black}{is the scaling factor used for model-dependent normalization}. To \textcolor{black}{control stimulation} focality, the nuisance-field weight $\beta \ge 0$ is treated as a variable parameter. The solution to (\ref{eq:TLS_minimization}) satisfies the linear system
\begin{equation}
\left(
\mathbf{L}_1^T\mathbf{L}_1
+ \gamma^2\beta^2\mathbf{L}_2^T\mathbf{L}_2
+ \gamma^2\varsigma^2\mathbf{I}_K
\right)\mathbf{y}
=
\mathbf{L}_1^T\mathbf{x}_1.
\label{eq:TLS_formula}
\end{equation}
\textcolor{black}{The regularization parameter $\gamma$ controls the strength of Tikhonov regularization, while the nuisance-field weight $\beta$ balances suppression of off-target stimulation relative to target-field fitting}.

%%%%%
% \section{two}
% The forward modeling uncertainty ranges summarized in Table~\ref{tab:forward_uncertainty_db} provide a physical interpretation of the attenuation feasibility constraints and noise levels employed in the Galaz Prieto uncertainty modeling framework. In the proposed method, uncertainty is incorporated by restricting the admissible attenuation range of the predicted electric field, expressed in decibels, thereby preventing the optimization from exploiting unrealistically large field ratios that may arise from nominal lead-field solutions. Specifically, the constrained formulation limits the feasible attenuation interval to values consistent with physiologically plausible and numerically stable field propagation, while the unconstrained variant permits attenuation over an extended dynamic range. The literature summarized in Table~\ref{tab:forward_uncertainty_db} shows that realistic forward modeling variability produces magnitude deviations on the order of approximately $0.1$--$2$~dB, corresponding to bounded multiplicative deviations in predicted field intensity. These deviations fall well within the attenuation ranges explicitly considered in the Galaz Prieto formulation and are further reflected in the PSNR levels used to evaluate robustness. Consequently, the uncertainty levels assumed in the optimization framework are consistent with empirically observed forward modeling variability, and the attenuation-based constraint ensures that the resulting stimulation configurations remain stable under realistic modeling uncertainty.

\bibliographystyle{ieeetr}
\bibliography{references.bib}
%\input{references}

% Directional Leads for Deep Brain Stimulation: Technical Notes and Experiences
% https://karger.com/sfn/article/99/4/305/295161/Directional-Leads-for-Deep-Brain-Stimulation

% Patient-Specific Simulations of Deep Brain Stimulation Electric Field with Aid of In-house Software ELMA
% https://ieeexplore.ieee.org/document/8856307

% Patient-specific analysis of the relationship between the volume of tissue activated during DBS and verbal fluency
% https://www.sciencedirect.com/science/article/pii/S1053811910003551

% Method for patient-specific finite element modeling and simulation of deep brain stimulation
% https://link.springer.com/article/10.1007/s11517-008-0411-2

% A Study on the Feasibility of the Deep Brain Stimulation (DBS) Electrode Localization Based on Scalp Electric Potential Recordings
% https://www.frontiersin.org/articles/10.3389/fphys.2018.01788/full

% Optimized programming algorithm for cylindrical and directional deep brain stimulation electrodes
% https://iopscience.iop.org/article/10.1088/1741-2552/aaa14b/meta
\end{document}